\documentclass[reqno,12pt,twoside]{amsbook}
\usepackage{graphicx}
\usepackage{latexsym}
\usepackage[centertags]{amsmath}
\usepackage{amsfonts}
\usepackage{amssymb}
\usepackage{amsthm}
\usepackage{newlfont}
\usepackage{ametsoc}

\usepackage{setspace}             
\usepackage[nolot,nolof,nocopyright]{brandeis}        

\def\mytopsep{3mm}

\newtheoremstyle{myplain}{\mytopsep}{\mytopsep}{\itshape}{0pt}{\bfseries}{.}{3mm}{}
\newtheoremstyle{mydefinition}{\mytopsep}{\mytopsep}{\normalfont}{0pt}{\bfseries}{.}{3mm}{}

\newtheoremstyle{myremark}{\mytopsep}{\mytopsep}{\normalfont}{0pt}{\bfseries}{.}{3mm}{}

\theoremstyle{myplain}

\newtheorem{thm}{Theorem}[section]
\newtheorem{cor}[thm]{Corollary}
\newtheorem{lem}[thm]{Lemma}
\newtheorem{prop}[thm]{Proposition}

\theoremstyle{mydefinition}
\newtheorem{dfn}[thm]{Definition}
\theoremstyle{myremark}
\newtheorem{rem}[thm]{Remark}
\newtheorem{exa}[thm]{Example}

\allowdisplaybreaks[1]

\renewcommand{\thesection}{\thechapter-\arabic{section}}

\makeatletter\@addtoreset{equation}{section}\makeatother

\def\rank{\mathop{\mathrm{rank}}}
\def\Oge{\mathop{\Omega}_\ge}
\def\Oeq{\mathop{\Omega}_=}

\def\bcup{\bigcup\limits}
\def\eql{\simeq}
\def\equiv{\sim}
\def\tomonoid{TO-monoid}
\def\sy{\mathfrak{S}}

\def\supp{\mathop{\mbox{supp}}}
\def\ord{\mathrm{ord}}
\def\mb{\mathbf}
\def\diag{\mathrm{diag }}
\def\alg{{\mathrm{alg}}}
\def\fra{{\mathrm{fra}}}

\def\NN{\mathbb{N}}
\def\tt{t^{-1}}
\def\QQ{\mathbb{Q}}
\def\PP{\mathbb{P}}
\def\xx{x^{-1}}
\def\yy{y^{-1}}
\def\CC{\mathbb{C}}
\def\ZZ{\mathbb{Z}}
\def\RR{\mathbb{R}}

\def\res{\mathop{\mathrm{Res}}}
\def\ct{\mathop{\mathrm{CT}}}
\def\PT{\mathop{\mathrm{PT}}}
\def\NT{\mathop{\mathrm{NT}}}
\def\pt{\mathop{\mathrm{PT}}}
\def\nt{\mathop{\mathrm{NT}}}
\def\CT{\mathop{\mathrm{CT}}}

\def\uprho{\mbox{}^\rho \,}

\def\poly{\mbox{\sf Poly}}
\def\frr{\mbox{\sf Frac}}
\def\rmd{\mathrm{rmd }}

\def\bous{Bousquet-M\'{e}lou}
\newcommand{\pad}[2]{\displaystyle\frac{\partial #2}{\partial #1}}
\newcommand{\ceiling}[1]{\,\lceil #1 \rceil\,}
\renewcommand{\ll}{\langle\!\langle}
\renewcommand{\gg}{\rangle\!\rangle}

\title{The Ring of Malcev-Neumann Series and the Residue Theorem
}

\thesisauthor{Guoce Xin}

\thesisadvisor{Ira M. Gessel}

\acknowledgments{\vspace{1cm}I am very grateful to my advisor, Ira
Gessel, for his guidance, suggestions, constant encouragements and
help. I thank Mireille Bousquet-M\`{e}lou, who expresses her
interests in my work and supplies me with reprints of her recent
work, which help me a lot on my research. I thank Richard Stanley
for providing me with useful references and friendly suggestions.
I thank Susan Parker for helping me on improving my teaching. I
thank Michael Cleber for giving me practical advices. I thank
Harry Tamvakis for introducing me the book, Integral
Representation of Combinatorial Sums, which speeds up my research
a lot. I thank my wife and colleague, Ji Li, for her patience,
constant support and encouragements.

}     

\thesisabstract{We develop a theory of the field of double Laurent
series, iterated Laurent series, and Malcev-Neumann series that
applies to most constant term evaluation problems. These include
(i) MacMahon's partition analysis, counting solutions of systems
of linear Diophantine equations or inequalities, counting the
number of lattice points in convex polytopes, (ii) evaluating
combinatorial sums and their generating functions, and proving
combinatorial identities, and (iii) lattice path enumeration such
as walks on the slit plane and walks on the quarter plane.

In the general setting of this new theory, the natural definition
of ``taking the constant term" of a formal series works well and
thus  the operators of taking constant terms commute with each
other. The proof of Bousquet-M\'{e}lou and Schaeffer's conjecture
about walks on the slit plane is included. In addition, the
counting problem of walks on the half plane avoiding the half line
is solved. Jacobi's multivariate residue theorem is generalized to
a field of Malcev-Neumann series, which gives a new interpretation
and a better understanding of the residue theorem. One application
of the residue theorem is a concise proof of Dyson's conjecture.

A new algorithm for partial fraction decompositions is developed.
This new algorithm is fast and uses little storage space. It also
results in an efficient algorithm for MacMahon's partition
analysis and related constant term evaluations.
}            

\date{2004}

\begin{document}

\thesisfront[May]{Department of Mathematics}{Ira M.
Gessel}{Michael Kleber, Department of Mathematics}{Richard P.
Stanley, MIT}{to}

\vfill\eject \setcounter{chapter}{-1} \setcounter{page}{1}
\chapter{Introduction}

This thesis is about combinatorial applications of formal Laurent
series. Our central topic is constant term evaluations, or
equivalently, residue evaluations. We will develop a general
setting for working with constant term evaluations that arose in
many areas. These include three major ones: (i) MacMahon's
partition analysis, counting solutions of systems of linear
Diophantine equations or inequalities, counting the number of
lattice points in convex polytopes, (ii) evaluating combinatorial
sums and their generating functions, and proving combinatorial
identities, and (iii) lattice path enumeration such as walks on
the slit plane and walks on the quarter plane.

Simply speaking, we mainly deal with formal Laurent expansion of
multivariate rational functions.

Let $K$ be a field. Starting from the field $K((x))$ of Laurent
series, we study the field $K((x))((t))$ of double Laurent series,
which is the field of Laurent series in $t$ with coefficients in
$K((x))$. Then we generalize to the multivariate case, the field
$K\ll x_1,\dots, x_n\gg$ of iterated Laurent series. Finally we
generalize to the ring of Malcev-Neumann series. The latter three
fields and rings have been little studied by combinatorists but
have many applications. Our general setting is in them.

\section{Connection to Complex Analysis}
To understand formal Laurent expansions of rational functions in a
simple fashion, we connect them with complex analysis. Note that
the arguments in this section are not rigorous.

First let $A$ and $B$ be two complex numbers. Then $A-B$ has a
reciprocal if $A\ne B$. We have the following geometric series
expansion
\begin{align}\label{e-0-AB}
\frac{1}{A-B}=\left\{%
\begin{array}{ll}
    \frac{1}{B} \frac{-1}{1-A/B}=-\sum_{n\ge 0} A^n/B^{n+1}, & {\text{if } A<B ,} \\
    \frac{1}{A} \frac{1}{1-B/A}=\sum_{n\ge 0} B^n/A^{n+1}, & {\text{if } A>B.} \\
\end{array}%
\right.
\end{align}

The observation is that in order to get a series expansion of
$1/(A-B)$, we need to know what is $A$ and $B$ is greater. Note
that the above expansions makes no sense when $\CC$ is replaced
with an arbitrary field $K$.

Now let $K$ be a field. By introducing a new variable $x$, and
treating $x$ as $o(1)$, or equivalently $x<c$ for all $0\ne c\in
K$, we informally get the field $K((x))$ of Laurent series.

Now let $A$ and $B$ be two series in $K((x))$. How can we expand
$1/(A-B)$ in $K((x))$? Informally, we have the expansions in
\eqref{e-0-AB}, except that the relation $A<B$ is replaced with
$A=o(B)$. When $A=O(B)$ we cannot expand $1/(A-B)$ in terms of $A$
and $B$. This argument can be made rigorous by the composition law
of $K[[x]]$.

How to generalize this idea to the two variable case? The obvious
way of letting $t=o(x)$ does not work, because we will have
trouble in expanding $1/(x^2-t)$. The solution is letting $t<c$
for all $0\ne c\in K((x))$. This is our field $K((x))((t))$ of
double Laurent series, i.e., the field of Laurent series in $t$
with coefficients in $K((x))$.

This idea naturally generalizes to the multivariate case, the
field of iterated Laurent series, and we always have the
expansions in \eqref{e-0-AB} depending on $A=o(B)$ or $B=o(A)$.

Let us recall the well-known result about residues in complex
analysis:
\begin{thm}
Let $\gamma$ be a simple curve in $\CC$. If $f$ is meromorphic
function that has no singularity on $\gamma$, then
$$\frac{1}{2\pi i}\int _\gamma f dz = \sum_{a\in E} \res_{z=a} f,$$
where $E$ is the set of singularities of $f$ that lie inside
$\gamma$.
\end{thm}

\begin{figure}[h]
 $$ \includegraphics[width=8cm]{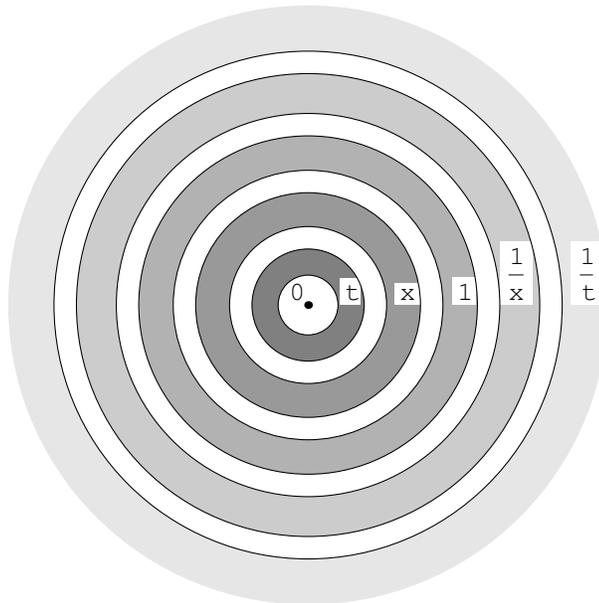}$$
  \caption{The plane of $K((x))((t))$}\label{f-kxt}
\end{figure}

In our general setting, e.g., for $F=F(x,t) $ in $K((x))((t))$,
$\res_{x} F $ is defined to be $[\xx] F$. It can be thought of as
$$\res_{x} F = \frac{1}{2\pi i} \int_ \gamma F(z,t) dz,$$
where $\gamma$ is the curve $|z|=x$, and the plane of complex
numbers should be replaced with the plane of $K((x))((t))$ as
shown in Figure \ref{f-kxt}.

The shaded regions are described as follows. Note that there are
gaps between those shaded regions.
\begin{multline*}
 \{\,0\,\}\hookrightarrow
tK((x))[[t]]\hookrightarrow
 tK((x))[[t]]\cup xK[[x,t]] \\
 \hookrightarrow
 tK((x))[[t]]\cup K[[x,t]]\hookrightarrow K((x))[[t]]\hookrightarrow K((x))((t)).
\end{multline*}

When integrating along the curve $|z|=x$, we need only consider
those singularities inside the curve. Since the singularities must
be independent of $x$, they belongs to $K((t))$. Using this
understanding, Theorems \ref{t-lagrange1} and \ref{t-1-mdividg}
can be thought as variations of the complex residue theorem.

The plane for $K\ll x_1,\dots ,x_n \gg$ can be drawn similarly,
but the shaded regions would be hard to describe.

\section{Structure of This Thesis}
This thesis consists of five chapters. Chapter 0 is this
introduction. We connect our theory to complex analysis in the
previous section. This connection will give a guide on how to
expand rational functions.

In Chapter 1 we rigorously develop the theory of $K((x))((t))$,
the field of double Laurent series.  The study of this field was
inspired by the application of $K[[x,t/x]]$ by \citep{ira}.
Results by \citep{bous} about walks on the slit plane stimulated
our research.

We proceed by introducing the basic concepts and operators, in
which three decompositions of double Laurent series are formally
given. The first one is used frequently. It says that we can
uniquely separate a given $f$ into two parts: one with only
nonnegative powers in $x$, and the other with only positive powers
in $x$. The second decomposition gives the concept of initial term
that evolves from that of $K((x))$. In terms of complex analysis,
a monomial $f$ is the initial term of $F$ if and only if
$F=f+o(f)$. This gives a guide on expanding $1/F$ into a double
Laurent series. The third decomposition comes from the unique
factorization lemma in \citep{ira} and \citep{bous}. It says that
if the initial term of $f$ is $1$ then $f$ can be uniquely
factored into three parts: one with only positive powers in $x$,
one being independent of $x$, and one with only negative powers in
$x$. It is obtained from the first decomposition by taking a
logarithm. It has many applications to lattice path enumeration.

An element $f(x,t)\in K((x))((t))$ has the form
$$f(x,t)=\sum_{n\ge n_0} a_n(x) t^n ,$$
where $a_n(x)$ is a Laurent series in $x$  for all $n$. It can
also be written in the form
$$f(x,t)=\sum_{m\in \ZZ} \sum_{n\ge n_0} a_{mn}x^mt^n,$$
where $a_{mn}$ is in $K$.

We define the constant term of $f$ in $x$ to be
$$\CT_x f(x,t)=\sum_{n\ge n_0 } \left(\CT_x a_n(x)\right) t^n=\sum_{n\ge n_0} a_{0n }t^n.
$$

Now let $F$ and $G$ in $K((x))((t))$ contain only nonnegative
powers in $x$. One basic problem is to evaluate the constant term
$\ct_x xF/G$, or equivalently $\res_x F/G$.

The most useful result in Chapter \ref{s-kxt} is Theorem
\ref{t-lagrange1}, which gives a formula for $\res F/G$ in terms
of $Y$ when $G$ has only one simple root $Y$ that is $o(x)$. It is
a generalization of the Lagrange inversion formula \citep[Theorem
5.4.2]{EC2}, and can be used to derive the multivariate Lagrange
inversion formula as described in Section \ref{ss-lag}. Other
applications of Theorem \ref{t-lagrange1} can be found in Chapter
\ref{s-lattice} on lattice path enumerations. It plays an
important role in proving a conjecture of \bous\ about walks on
the slit plane. See Section \ref{ss-conj} or \citep{xin}.

\vspace{3mm} Another useful result is Theorem \ref{t-mdivid},
which gives a formula for $\res F/G$ in terms of all the roots
that are $o(x)$.

As an application, we give a short proof of the well-known result:
{\em the diagonal of a rational power series in two variables is
algebraic.} See, e.g., \cite[Theorem 6.33]{EC2}. Note that we also
use Puiseux's Theorem.

In evaluating $\ct_x xF/G$, knowing the properties of the roots of
$G(x,t)$ will be helpful. This is the motivation of Section
\ref{s-puiseux}. We generalize Puiseux's Theorem a little bit and
use it to evaluate $\ct_x xF/G$.

\vspace{3mm} In Chapter \ref{s-multi} we study the field of
iterated Laurent series $K\ll x_1,\dots ,x_n\gg$. The fundamental
result is Proposition \ref{p-wellorder}, which says that a formal
series is an iterated Laurent series if and only if it has a
well-ordered support.

This result validates the application of the constant term
operator $\ct_{x_j}$ naturally defined by
$$\ct_{x_j} \sum_{(i_1,\dots ,i_n)\in \ZZ^n}
a_{i_1,\dots ,i_n} x_1^{i_1} \cdots x_n^{i_n} = \sum_{(i_1,\dots
,i_n)\in \ZZ^n, i_j=0} a_{i_1,\dots ,i_n} x_1^{i_1} \cdots
x_n^{i_n},$$ where $a_{i_1,\dots ,i_n}$ belongs to $K$.

This natural definition has some nice properties, such as
$\ct_{x_i}$ commutes with $\ct_{x_j}$, and $\ct_{x_i}$ commutes
with $\sum$. It is also consistent with the previous definitions.

For example, \citet{zeil} proved a conjecture of \cite{chan} by
showing an identity equivalent to
\begin{equation}
\ct_{x_1}\cdots \ct_{x_n} \frac{1}{\prod_{i=1}^n (1-x_i)}
\frac{1}{\prod_{i<j} (x_i-x_j)} =C_1\cdots C_{n-1},
\end{equation}
where $C_n$ is the Catalan number. As pointed out in
\citep{welleda}, this identity should be interpreted as taking
iterated constant terms; i.e., while applying $\ct_{x_n}$ to a
rational function, we expand it as a Laurent series in $x_n$. The
result is still a rational function and we can apply
$\ct_{x_{n-1}}$, \dots, $\ct_{x_1}$ iteratively. Note that in this
definition, $\ct_{x_i}$ does not commute with $\ct_{x_j}$.

Our approach is to treat a rational function as an element in
$K\ll x_1,\dots ,x_n \gg$, and then take the constant terms. So
after specifying the working field, the iterated constant term
operator is simply $\ct_{x_1,\dots ,x_n}$.

Once knowing this general setting,  the basic computational rules
are easy to use. In fact, the residue computation can be done
similarly as in complex analysis. The difference is that we shall
use the plane of $K\ll x_1,\dots ,x_n\gg$ instead of the plane of
complex numbers.

Section \ref{s-comsum} is the application to combinatorial sums.
We did not include many examples because much of this work has
been done in \citep{ego}. But we believe that our approach is
simpler.

Section \ref{s-parfrac} gives a new algorithm for partial fraction
decompositions of rational functions. This new algorithm is fast
and uses little storage space. We give a natural proof of the nice
reciprocity law for higher dimensional Dedekind sums in
\citep{zagier}.

Section \ref{s-Mac} is the application to MacMahon's partition
analysis, which has been given a new life by \citet{george6} in a
series of papers. The problem is reduced to evaluating the
constant term of a special type of rational function, which we
call the Elliott-rational functions. The denominators of these
rational functions have simple factors of one or two terms. The
constant terms of Elliott-rational functions are still
Elliott-rational. Thus taking constant term in several variables
can be reduced to iteratively taking constant term in one
variable.

Our approach is to embed the rational functions in a field of
iterated Laurent series, so that its series expansion is separated
from its rationality. More precisely, we first use partial
fraction decomposition and then apply its series expansion. This
approach results in an efficient algorithm as given in Section
\ref{s-Mac}.

Chapter 3 develops the most general setting: the ring of
Malcev-Neumann series (or MN-series for short). Let $K$ be a field
and let $G$ be a {\em totally ordered group}; i.e., $G$ has a
total ordering that is compatible with its group structure. Let
$K_w[G]$ be the set of all formal series in $G$ that have a
well-ordered support. \citet{malcev} and \citet{neumann} showed
that $K_w[G]$ is a division algebra that includes the group ring
$K[G]$ as a subalgebra. The importance of this result was to solve
an algebraic problem: $K[G]$ has no zero divisors if $G$ can be
made a totally ordered group.

Using the construction of MN-series in \citep{passmann}, we show
that the field $K$ can be replace with a commutative ring $R$ with
a unit, and that $G$ can be
 a monoid. So $R_w[G]$ is the ring of MN-series.

In Section \ref{s-MNresidue}, we show that under the reverse
lexicographical ordering, the field $K_w[\ZZ^n]$ is isomorphic to
$K\ll x_1,\dots ,x_n \gg$. So the field of iterated Laurent series
is a special case of the field of MN-series. Next we give the
residue theorem for the ring of MN-series. This is a twisted
multivariate residue theorem. This is the most significant result
in the thesis from several aspects. First of all, there is no
known analogous explanation of this twisting in complex analysis.
Next, as we will discuss in Section \ref{s-MNresidue}, this result
includes other (combinatorial) residue theorems as special cases,
and it has fewer conditions. Finally, our result asserts that the
residue theorem in fact involves two rings (or fields), which has
been overlooked by combinatorists.

Using our residue theorem, we give another view of the Lagrange
inversion formula in Section \ref{ss-lag}, and we give two proofs
of Dyson's conjecture in Section \ref{s-dyson}.

Section \ref{s-morris} simplifies the proof of the Morris identity
by \citep{welleda}. Section \ref{s-Macr} talks about the
theoretical aspects of MacMahon's partition analysis. We give a
new proof of the reciprocity theorem for a system of homogeneous
linear Diophantine equations by \cite{stanley-rec}.

Chapter 4 is the applications to lattice path enumeration. We use
the bridge lemma in \citep{bous} and the concept of Gessel pair,
that results from \citep{ira}, to work with lattice path
enumeration problems.

We simplify the previous works about walks slit plane in
\citep{bous,bouso}. Section \ref{s-4-exa} solves a problem
proposed by \citep{bouso} and solves some new types of lattice
path problems. Section \ref{ss-conj} solves a conjecture in
\citep{bous}. In all of this work, Theorem \ref{t-lagrange1} is a
basic tool.

Section \ref{s-quart} is about walks on the quarter plane, studied
by \citep{bousquart,bousquartf,bousquartD} in several papers.  We
give a simple description for the functional equation we need to
solve, and described the solution for a simple case.

\renewcommand{\theequation}{\thesection.\arabic{equation}}

\renewcommand{\theequation}{\thesection.\arabic{equation}}

\vfill\eject \setcounter{chapter}{0} 

\chapter{The Field of Double Laurent Series\label{s-kxt}}

\section{Notations and Background}
In this thesis,  $R$ is always a commutative ring with unit, and
$K$ is always a field of characteristic $0$. Let $t$ be a formal
variable. We review some conventional notation.
\begin{enumerate}
\item $R[t]$:  the ring
of polynomials in $t$ with coefficients in $R$.

\item $R(t)$: the ring of rational functions in $t$.

\item $R[[t]]$: the ring of formal power series in $t$.

\item $R[t,\tt]$ the ring of Laurent polynomials in $t$.

\item $R((t))$ the ring of Laurent power series in $t$.
\end{enumerate}

The ring $K[[x]]$ has been studied by many authors. Generating
functions of most combinatorial objects are in $K[[x]]$, or in
$K[[x_1,x_2,\ldots ,x_n]]$ for the multivariate case.

The ring $K[[x]]$ is a local ring. Its unique maximal ideal is
generated by $x$. Thus by adjoining $x^{-1}$, we can get its
quotient field $K((x))$, which is called the field of Laurent
series in $x$. We
 can identify $K((x))$ with
$K[[x]][\xx],$ the polynomial ring in $\xx$ with coefficients in
$K[[x]]$. We will see that in many situations, it is more
advantageous to work in $K((x))$, because of its field structure,
than in $K[[x]]$.

An element $\eta \in K((x))$ has the form
$$\eta=\sum_{n\ge n_0} a_n x^n,$$
 where $a_n\in K$ for all $n$. If $a_{n_0}\ne 0$, then we say that $\eta $ has order $n_0$,
 and $\eta$ can be written as $x^{n_0} \rho$, with $\rho $ an ordinary power series with
nonzero constant
 term. Moreover, $\eta^{-1}$ has order $-n_0$ since $\eta^{-1}=x^{-n_0}\rho^{-1}$.
We have the following three situations.
\begin{enumerate}
\item
If $\eta $ has positive order, then $\eta \in xK[[x]]$.
\item
If $\eta $ has order zero, then $\eta$  is a unit in $K[[x]]$.
\item
If $\eta $ has negative order, then $\eta^{-1} \in xK[[x]]$.
\end{enumerate}

Obviously, if $\eta(x)\in K((x))$, then $\eta(\xx)\in K((\xx))$.
The field $K((\xx))$ will turn out to be useful later. Now let us
look at some basic facts.

Clearly we have $K((x))\cap K((\xx))=K[x,\xx]$, the ring of
Laurent polynomials in $x$. Now let $\eta \in K[x,\xx]$. The
expansion of $\eta^{-1}$ in $K((x))$ is usually different from the
expansion of $\eta^{-1}$ in $K((\xx))$. For example, consider
$\eta=1-x \in K[x]$. The expansion of $\eta^{-1}$ in $K((x))$ is
$$\frac{1}{1-x}=\sum_{n\ge 0} x^n,$$
while the expansion of $\eta^{-1}$ in $K((\xx))$ is
$$\frac{1}{1-x}=\frac{\xx}{\xx-1} =-\frac{\xx}{1-\xx} =\sum_{n\ge 0} -x^{-n-1}.$$

So we shall specify the working field whenever the reciprocal
$\eta^{-1}$ comes into account.

\begin{rem}
The set of all elements of the form $\sum_{n\in \ZZ} a_n x^n$ is
not a ring under the usual multiplication.
\end{rem}

The field of double Laurent series $K((x))((t))$ is a field
extension of $K((x))$. It contains all the Laurent series in $t$
with coefficients in $K((x))$. The study of $K((x))((t))$ is
inspired by Gessel's work on the ring $K[[x,y/x]]$ in \citep{ira},
and stimulated by the work in \citep{bous}.


\vspace{3mm} Algebraic power series form a frequently used class
of generating functions in combinatorial theory. We include the
definition of \citep{EC2} as follows.

\begin{dfn}
Let $K$ be a field. A formal power series $\eta\in K[[x]]$ is said
to be {\em algebraic} if there exist polynomials
$P_0(x),P_1(x),\ldots ,P_d(x)\in K[x]$, not all $0$, such that
\begin{align}\label{e-algbraicdfn}
P_0(x)+P_1(x)\eta +\cdots +P_d(x) \eta^d &=0.
\end{align}
The smallest positive integer $d$ for which \eqref{e-algbraicdfn}
holds is called the degree of $\eta$.
\end{dfn}
Note that an algebraic series $\eta$ has degree one if and only if
$\eta$ is rational. We denote $K_{\alg}[[x]]$ the set of all
algebraic power series over $K$.

\begin{exa}
Let
$$\eta =(1-4x)^{-{1\over 2}}$$
Then we have $(1-4x)\eta ^2-1=0$. Hence $\eta$ is algebraic of
degree one or two. It is easy to check that $\eta$ is algebraic of
degree one if and only if the characteristic of $K$ is $2$.
\end{exa}


\section{Basic Concepts and Operators}
\subsection{Concepts}
 An element $f(x,t)\in K((x))((t))$ could be written in the form
$$f(x,t)=\sum_{n\ge n_0} a_n(x) t^n ,$$
where $a_n(x)$ is a Laurent series in $x$  for all $n$. It can
also be written as
$$f(x,t)=\sum_{m\in \ZZ} \sum_{n\ge n_0} a_{mn}x^mt^n,$$
where all the $a_{mn}$ are in $K$ and some of them are restricted
to be zero.

We denote by $[x^mt^n] f(x,t)$ the coefficient of $x^mt^n$ in
$f(x,t)$. Two elements in $K((x))((t))$ are equal if and only if
all of the corresponding coefficients are equal.

Note that $K((x))((t))$ is different from $K((t))((x))$. It is not
hard to see that the intersection of the two fields is $K((x,t))$,
which is called the ring of Laurent series in $x$ and $t$. It can
be identified with $K[[x,t]][\xx,\tt]$.

Similarly, we can consider the field $K((\xx))((t))$, which is
isomorphic to $K((x))((t))$. It is easy to see that $f(x,t)\in
K((x))((t))$ if and only if $f(\xx,t)\in K((\xx))((t))$. The map
induced by $x\to \xx$ connects these two fields. The intersection
$$K[x,\xx]((t))= K((x))((t))\cap K((\xx))((t))$$
is the ring of Laurent series in $t$ with coefficients that are
Laurent polynomials of $x$.

\vspace{-1mm} Polynomials in $x,t$ are clearly in $K((x))((t))$.
Several other basic series that we are going to use are listed as
follows.
\begin{align*}
\frac{1}{1-t} &= \sum_{n\ge 0} t^n ,\\
(1+t)^\alpha &= \sum_{n\ge 0} {\alpha \choose n} t^n, \quad \text{for all } \alpha \in \CC,\\
\log \frac{1}{1-t} &=\sum_{n\ge 1} \frac{1}{n} t^n,\\
e^t &= \sum_{n\ge 0} \frac{1}{n!} t^n
\end{align*}

We call the following the {\em composition law} of $K((t))$.
\begin{lem}\label{l-composition}
If $u=u(t)\in K[[t]]$ has constant term $0$, and $f(t)\in K((t))$,
then $f(u(t))\in K((t))$.
\end{lem}
\begin{proof}
Write $f(t)=t^{n_0} g(t)$ where $g(t)\in K[[t]]$ with nonzero
constant term. Then $g(u(t))\in K[[t]]$ by the composition law of
formal power series (see, e.g., \citep{EC1}). Hence
$f(u(t))=(u(t))^{n_0} g(u(t))$ belongs to $K((t))$ by its field
structure.
\end{proof}

In $K((x))((t))$, the composition law is just the application of
Lemma \ref{l-composition} on $K((x))$, and on $K((x))((t))$ by
passing the base field to $K((x))$.
\begin{prop}
If $f(x,t)\in K((x))((t))$, $g(x,t)\in tK((x))[[t]]$, and $h(x)\in
xK[[x]]$, then
 $f(h(x),g(x,t))\in K((x))((t))$.
\end{prop}

The field structure of $K((x))((t))$ and the composition law make
it possible for us to work with a large class of series. For
example, we can work with rational functions.

Every rational function $P(x,t)/Q(x,t)$ with $P(x,t),Q(x,t)\in
K[x,t]$ has a unique expansion in $K((x))((t))$. To expand it, we
write
$$Q(x,t)=\sum_{i=d}^m a_i(x) t^i,$$
where $a_d(x)\ne 0$. Then
$$
\frac{P(x,t)}{Q(x,t)}=\frac{P(x,t)}{a_d(x)t^d}\frac{1}{1+t\sum_{i=d+1}^m
a_i(x) t^{i-d-1} /a_d(x)}.
$$

By symmetry, every rational function also has an expansion in
$K((t))((x))$. But these two expansions are usually different. For
example, in $K((x))((t))$ we have
\begin{align*}
\frac{1}{t-x} = \frac{1}{-x}\frac{1}{1-t/x} =\sum_{n\ge 0}
-x^{1-n} t^n,
\end{align*}
but in $K((t))((x))$ we have
$$\frac{1}{t-x}=\frac{1}{t}\frac{1}{1-x/t}=\sum_{n\ge 0} t^{1-n} x^n.$$

For any element $f(x,t)\in K((x))((t))$, we have the following
three decompositions, which help us to understand the behavior of
double Laurent series. The first decomposition is straightforward
but frequently used.

\begin{lem}[First Decomposition] \label{l-1decom}
Any $f(x,t)$ in $K((x))((t))$ can be uniquely written as
$f_1+f_2$, where $f_1$ contains only nonnegative powers in $x$,
and $f_2$ contains only negative powers in $x$.
\end{lem}

We call $f_1$ the nonnegative part of $f$ in $x$, denoted by
$\pt_x f(x,t)$, and $f_2$ the negative part, denoted by $\nt_x
f(x,t)$. Thus $f(x,t)=\pt_x f(x,t)+\nt_x f(x,t)$ and
$$\pt_x f(x,t)= \sum_{n\ge n_0}\sum_{m\ge 0}  a_{mn}x^mt^n,$$
$$\nt_x f(x,t)=\sum_{n\ge n_0}\sum_{m<0}  a_{mn}x^mt^n.$$
Note that $f_1(x,t)$ is in $K[[x]]((t))$, but $f_2(x,t)$ is
actually in $K[\xx]((t))$.

If $f$ contains only nonnegative powers in $x$, then $f=\pt_x f$,
and we say that
 $f$ is $\pt$ in $x$. Similarly we can define $f$ to be $\nt$ in $x$
if $f$ contains only negative powers in $x$. Of course we can
define $\pt_t f(x,t)$ and $\nt_t f(x,t)$, but they are not as
useful.

On the issue of finding the right expansion for a reciprocal, the
second decomposition, as shown below, is going to be helpful.

\begin{lem}[Second Decomposition] \label{l-2decom}
Any $f(x,t)$ in $K((x))((t))$ can be uniquely written as
\begin{align}\label{e-seconddec}
f(x,t)=c x^{m_0}t^{n_0} g(x) h(x,t),
\end{align} where $c\in K$ is a constant,
$g(x)\in K[[x]]$ with constant term $g(0)=1$, and $h(x,t)\in
K((x))[[t]]$ with $h(x,0)=1$.
\end{lem}

In the second decomposition of $f(x,t)$, as given by
\eqref{e-seconddec}, we call $c x^{m_0}t^{n_0}$ the {\em initial
term} of $f(x,t)$, and $(m_0,n_0)$ the {\em order} of $f(x,t)$.
Moreover the second decomposition of $1/f(x,t)$ is given as
follows:
$$\frac{1}{f(x,t)}=\frac{1}{c} x^{-m_0}t^{-n_0} \frac{1}{g(x)} \frac{1}{h(x,t)},$$
where the meanings of $1/g(x)$ and $1/h(x,t)$ are clear.

In the following Lemma \ref{l-unifactor}, $h(x,t)$, as in Lemma
\ref{l-2decom}, can be decomposed further. This decomposition is
called the Unique Factorization Lemma by \citep{ira}, and by
\citep{bous}. It follows from the first decomposition through
taking a logarithm, and has some nice applications in lattice path
enumeration, as we shall discuss later in chapter \ref{s-lattice}.

\begin{lem}[Third Decomposition]\label{l-unifactor}
Let $h(x,t)$ be an element in $K((x))[[t]]$, in which the constant
term in $t$ is $1$, i.e., $h(x,0)=1$. Then $h$ has a unique
factorization in $K((x))[[t]]$ such that  $h=h_-h_0h_+$, where
$h_-\in K[\xx][[t]]$, $h_0\in K[[t]]$, and $h_+\in K[[x,t]]$.
Moreover, all the constant terms of $h_-,h_0, $ and $h_+$  are
$1$.
\end{lem}
\begin{proof}
Let $\log h=\sum_{i,j}b_{ij}x^it^j$. Then
\begin{align*}
h_- &=\exp\Big(\sum_{i<0,j>0} b_{ij} x^it^j \Big), \\
h_0 &=\exp\Big(\sum_{j\ge 1} b_{0j} t^j \Big), \\
h_+ &=\exp\Big(\sum_{i> 0,j>0} b_{ij} x^it^j \Big).
\end{align*}
The uniqueness follows from the first decomposition.
\end{proof}

The importance of the third decomposition is due to \citet{ira},
who connected it with the factorization of lattice paths. It is
also an important tool in the work of \citep{bous}.

\begin{rem}
Note that we shall still get a unique factorization if we group
$h_0$ and $h_+$ together. More precisely, if $h=h_1h_2$ with
$h_1\in K[\xx][[t]]$ and $h_2\in K[[x,t]]$, and both have initial
term $1$, then $h_1=h_-$ and $h_2=h_0h_+$. Similarly we can group
$h_-$ and $h_0$ together.
\end{rem}

\subsection{Operators}

One of the basic operators on $K((x))((t))$ is $\CT_x$, which
takes the constant term in $x$ of a series.
\begin{dfn}
For any $a(x)\in K((x))$, we denote by $\CT_x a(x)$  the constant
term $[x^0] a(x)$ of $a(x)$. Also for $f(x,t)\in K((x))((t))$, we
define
$$\CT_x f(x,t)=\sum_{n\ge n_0 } \left(\CT_x a_n(x)\right) t^n=\sum_{n\ge n_0} a_{0n }t^n.
$$
\end{dfn}

Clearly, we have $\ct_x f(x,t)=\left. \pt_x f(x,t)\right|_{x=0}.$
On the other hand, we will give a formula for $\pt_x f(x,t)$ in
terms of $\ct_x$.

For any Laurent series $h(x)$, the residue of $h(x)$ in $x$ is
defined to be
$$\res_x  h(x) = [x^{-1}] h(x) =\CT_x x h(x).$$
So essentially, the operator $\res_x$ plays the same role as the
operator $\CT_x$. Mathematicians are familiar with residue
computations, because in complex analysis, the residue can be
represented as an integral. For example, see \citep{ego}. The
operator $\CT_x$ is also frequently used, for it is more
convenient in many situations. For instance, if $f(x,t)$ is $\pt$
in $x$, then $\ct_x f(x,t)= f(0,t)$. This fact is seldom noticed
in residue computation, but easy to use in constant term
evaluation. We will use both operators. Note that $\res_x f=\CT_x
x f$, and $\CT_x f =\res_x x^{-1} f$.

\begin{dfn}
The {\em Hadamard product} of two series $f(t)=\sum_{n\ge 0} a_n
t^n$ and $g(t)=\sum_{n\ge 0} b_n t^n$ is defined to be
$$f(t)\odot g (t) =\sum_{n\ge 0} a_n b_n t^n .$$
\end{dfn}

The computation of Hadamard product can be converted into constant
term evaluation. We have
$$f(t)\odot g (t) = \ct_x f(t/x) g(x).$$
We can prove the following well-known result. See, e.g., \citep{EC1}.
\begin{thm}\label{t-1-hadamard}
If $f,g\in K[[t]]$ are also rational, then $f\odot g$ is
rational.
\end{thm}

We will give a more general form of this result in terms of
constant terms later. Now let us do an example.

If we want to compute $f\odot g$ for rational $f$ and $g$, we can
use the partial fraction method. For example, if $f$ and $g$ are
both
 the generating function of Fibonacci
numbers, i.e.
$$f(t)=g(t)=\frac1{1-t-t^2},$$
 then we can compute $f\odot g$ as follows.

Using Maple we can convert $f(t/x) g(x)$ into partial fraction in
$x$,
\begin{align*}
\begin{split}\left (1-{\frac {t}{x}}-{\frac {{t}^{2}}{{x}^{2}}}\right )^{-1}\left (
1-x-{x}^{2}\right )^{-1}\rule{9cm}{0pt}\\
={\frac {t \left (x+t-{t}^{2}\right )}{\left
(1-2t-2{t}^{2}+{t}^{3}\right ) x^2 \left( 1-t/x-{t/x}^{2}\right
)}}+{\frac {1-t+tx}{\left (1-2t-2{t}^{2}+{t}^{3}\right ) \left
(1-x-{x}^{2}\right )}}.
\end{split}
\end{align*}
We see that on the right hand side of the above equation, the
first term contains only negative powers in $x$ and the second
term contains only nonnegative powers in $x$. So by setting $x=0$
in the second term, we get
$$f(t)\odot f(t)= \ct_x f(t)f(t/x) =\frac {1-t}{1-2t-2{t}^{2}+{t}^{3}}
 .$$
The above argument can be used to compute the Hadamard product of
several rational functions.

\begin{dfn}
The diagonal of an element
$$f(x,t)=\sum_{m=-\infty}^{\infty}\sum_{n=n_0}^\infty a_{mn} x^m t^n$$
in $K((x))((t))$ is defined to be
$$\diag (f) (t)= \sum_{n=n_0}^\infty a_{nn} t^n.$$
\end{dfn}
The diagonal can be converted into constant term evaluation. We
have
$$\diag(f)(t)= \CT_x f(x,t/x).$$

If $f(x,t)$ is also in $K((t))((x))$, then $\diag(f)(x)=\CT_t
f(x/t,t)$ by symmetry.  We will see that some results on diagonal
are more suitably reformulated in terms of constant terms.

The partial differential operators $\pad{x}{}$ and $\pad{t}{}$ are
useful. One important fact is that
 $\res_x \pad{x}{f(x,t)}=0$ for any $f(x,t)\in K((x))((t))$, and similarly for $t$.
 As a direct consequence, we have that for any $f,g\in K((x))((t))$,
 $$\res_x \pad{x}f(x,t) \cdot g(x,t) = -\res_x f(x,t)\cdot \pad{x}{g(x,t)},$$
 since $\res_x \pad{x}{} [f(x,t) g(x,t)] =0$.
 This shows why it is sometimes more convenient to use residues than to use constant terms.

The integration operator $\int\; \cdot \; dx$ can only be applied
to elements
$f(x,t)$ with \\
$\res_x f(x,t)=0.$

The last useful operator in this section is the {\em divided
difference operator}.
\begin{dfn}
The divided difference operator $\partial_a$ with respect to $x$
defined on functions or series $f(x)$ is given by
$$\partial_a f(x)= \frac{f(x)-f(a)}{x-a}.$$
\end{dfn}
In this thesis, we are only going to use divided difference
operators with respect to one particular variable. When this
variable is clear, we will omit it.

Let $R$ be a commutative ring with unit. If $u$ is a new variable,
then $\partial_u$ is a linear operator from $R[[x]]$ to
$R[[x,u]]$. We have
\begin{align*}
\partial_u \sum_{n\ge 0} a_n x^n =\frac{\sum_{n\ge 0} a_n x^n
-\sum_{n\ge 0} a_n u^n}{x-u} = \sum_{n\ge 0} a_n
\left(x^{n-1}+x^{n-2}u+\cdots +u^{n-1}\right).
\end{align*}
If $u$ belongs to $xR[[x]]$, then $\partial_u f(x)\in R[[x]]$. If
$u=x$, then $\partial_x $ reduces to the derivative.

It is easy to see that $\partial_u f(x)$ is symmetric in $u$ and
$x$, and if $f(x)\in R[x]$ is a polynomial, then the degree of
$\partial_u f(x)$ in $x$ is one less than that of $f(x)$.
Moreover, the result of $\partial_u$ acting on a rational function
is still a rational function.

Let $$\Delta (z_1,\ldots ,z_n)=\prod_{i<j}(z_i-z_j)=
\left|\begin{array}{cccc}
   u_0^{n} & u_1^{n}  & \cdots &
   u_n^{n}\\
   u_0^{n-1} & u_1^{n-1}  & \cdots & u_n^{n-1}\\
   \vdots &\vdots   & \vdots & \vdots \\
   1& 1 & 1 & 1\end{array}\right|$$
be the Vandermonde determinant in $z_1,\ldots ,z_n$. Then we have
the following result, which will be used in the next section.

\begin{lem}
Taking divided difference with respect to $u_0$, we have
\begin{align}
\partial_{u_1}\partial_{u_2}\cdots \partial_{u_n} f(u_0)
&=
    \left|\begin{array}{cccc}
   f(u_0) &  f(u_1) & \cdots &
   f(u_n)\\
   u_0^{n-1} & u_1^{n-1}  & \cdots & u_n^{n-1}\\
   \vdots &\vdots   & \vdots & \vdots \\
   1& 1 & 1 & 1
   \end{array}\right|
  \Delta(z_1,\ldots ,z_n)   ^{-1} \label{e-1-7}\\
&=\sum_{i=0}^n \frac{f(u_i)}{\prod_{j\ne i}
(u_i-u_j)}.\label{e-1-8}
\end{align}
In particular,
\begin{align}\label{e-divid-h}
\partial_{u_1}\partial_{u_2}\cdots \partial_{u_n} u_0^{m+n}=h_m(u_0,\dots ,u_n)
=\sum_{0\le i_1\le \cdots \le i_m\le n} u_{i_1}\cdots u_{i_m},
\end{align}
which is a complete symmetric function of $u_0,\dots ,u_n$.
\end{lem}
\begin{proof}
The equivalence of \eqref{e-1-7} and \eqref{e-1-8} follows by
expanding the determinant by the first row. It is a well-known
result in the theory of symmetric function that equation
\eqref{e-1-7} reduces to \eqref{e-divid-h} when
$f(u_0)=u_0^{m+n}$.

Now by linearity, it suffices to show that the lemma is true for
all integers $m$ and $f(u_0)=u_0^m$. For $m\ge 0$, we can look at
the result of applying $\partial_{u_1}\partial_{u_2}\cdots
\partial_{u_n}$ to the generating function $1/(1-u_0z)$. By
induction on $n$, we have
\begin{align}
\partial_{u_1}\partial_{u_2}\cdots \partial_{u_n} \frac{1}{1-u_0z}= \frac{1}{(1-u_0z)\cdots (1-u_nz)}.
\end{align}
Equation \eqref{e-divid-h} hence follows by equating coefficients
of $z$. The case $m<0$ is similar.
\end{proof}

\section{Computational Rules in $K((x))((t))$}
In this section we shall establish the computation rules and the
residue theorem in the field of double Laurent series. These rules
will be generalized in the next two chapters. In all situations,
we shall see that the right expansion of a reciprocal is important
to our computations.

We start from the following easy fact: If $u$ is independent of
$x$, then
$$\ct_x \Big(\sum_{n\ge 0} u^n/x^n\Big) \cdot x^k = \left\{ \begin{array}{ll}
u^k & \text{if } k\ge 0, \\
0 & \text{if } k<0.
\end{array}\right. $$

Thus by linearity, we have the following.
\begin{lem}\label{l-basic}
If $f(x)=\sum_{n\ge 0} a_n x^n$, i.e., $f$ is $\pt$ in $x$, then
$$\ct_x \Big(\sum_{n\ge 0} \frac{u^n}{x^n} \Big)\cdot
f(x)= f(u).$$
\end{lem}

Note that in the field $K((x))((u))$ we can say that
$$\sum_{n\ge 0} \frac{u^n}{x^n}=\frac{1}{1-u/x}=\frac{x}{x-u}.$$
But the above equation is not true in $K((u))((x))$.

In order to use rational functions, we shall specify the ring in
which we are working. In this section, we  are working in
$K((x))((t))$ or $K((\xx))((t))$, i.e., series in $t$.

We have three situations for the expansion of $1/(x-u)$, as stated
in the following lemma.
\begin{lem}\label{l-3abc}
Let $R$ be a commutative ring with unit. Suppose $u=u(t)$ is a
Laurent series in $t$, and $v(t)=1/u(t)$. Then we have
\begin{enumerate}
\item[a)]
If $u(t)\in tR[[t]]$, then $(x-u(t))^{-1}$ belongs to $\xx
R[\xx][[t]]$. The following expansion is valid in both
$R((x))[[t]]$ and $R((\xx))[[t]]$.
\begin{align}\label{e-a}
\frac{1}{x-u(t)}=\frac{1}{x} \frac{1}{1-u(t)\xx} =\xx \sum_{n\ge
0}  x^{-n}u(t)^n .
\end{align}

\item[b)]
If $v(t) \in t R[[t]]$, then $(x-u(t))^{-1}$ belongs to
$R[x][[t]]$. The following expansion is valid in both
$R((x))[[t]]$ and $R((\xx))[[t]]$.
\begin{align}\label{e-b}
\frac{1}{x-u(t)} &= \frac{v(t)}{v(t)x-1} = -\frac{v(t)}{1-x v(t)}=
 -\sum_{n\ge 0} x^n u(t)^{-n-1}.
\end{align}

\item[c)]
If $u(t) \in R[[t]]$ and $u(0)\ne 0$, then the expansion of
$(x-u(t))^{-1}$ in $R((x))[[t]]$ is \eqref{e-b}, and the expansion
in $R((\xx))[[t]]$ is \eqref{e-a}.
\end{enumerate}
\end{lem}

\begin{rem}
Note that only in case $c)$ do we have different expansions of
$1/(x-u)$ in the two rings $R((x))[[t]]$ and $R((\xx))[[t]]$. This
small difference results in parallel theories.
\end{rem}

\noindent {\bf Convention}: If we write $\ct_x$, then we are
working in $R((x))((t))$. If we write $\ct_{\xx} f(x,t)$, then we
are working in $R((\xx))((t))$. So $f$ is $\pt$ in $\xx$ means
that $f$ belongs to the ring $R((\xx))((t))$, and that $f$
contains only nonnegative powers in $x$.

\vspace{3mm} Now applying Lemma \ref{l-basic} to the above lemma,
we have
\begin{cor}\label{c-base}
If $Q(x,t)$ is $\pt$ in $x$ and $u=u(t)\in tR[[t]]$, then
\begin{equation}\label{e-ctqq}
\ct_x \frac{x}{x-u} Q(x,t)=Q(u,t).
\end{equation}
If $Q(x,t)$ is $\pt$ in $\xx$ and $u=u(t)\in R[[t]]$, then
$$\ct_{\xx} \frac{x}{x-u} Q(x,t)=Q(u,t).$$
\end{cor}

\begin{rem}
Note that if $u$ has constant term nonzero, the composition law
will not guarantee the existence of $Q(u,t)$, but the condition of
$Q(x,t)$ being $\pt_{\xx}$ is sufficient. The condition of
$Q(x,t)$ being $\pt_x$ is equivalent to $t^kQ(x,t)\in R[[x,t]]$
for some integer $k$. Since $t$ can be factored out when taking
the constant term in $x$, we can simply say that $Q(x,t) \in
R[[x,t]]$ instead of $Q(x,t)$ is $\pt_x$ in the two variable case.
A similar situation does not happen in the  multivariate case.
\end{rem}

Now we show that the  $\pt$ operator can also be realized by
$\ct$. Let $Q(x,t)=Q_1(x,t)+Q_2(x,t)$ be the first decomposition
of $Q(x,t)$ in $x$. Then $Q_2(x,t)$ is $\nt$ in $x$, which implies
that  $Q_2(x,t)\frac{x}{x-u}$ is also $\nt$ in $x$, and thus has a
zero constant term in $x$. Now we have
$$ \ct_{x} Q(x,t)\frac{x}{x-y}=\ct_x Q_1(x,t)\frac{x}{x-y}=Q_1(y,t)= \pt_y Q(y,t).$$

We can generalize Corollary \ref{c-base} in two directions. One
is the following result.
\begin{thm}\label{t-mdivid}
If $u_1,\ldots ,u_n\in tR[[t]]$ for all $n$, then for any $Q(x,t)$
that is $\pt_x$ we have
$$\pt_x Q(x,t)\prod_{i=1}^n\frac1{ x-u_i}=
\partial_{u_n}\partial_{u_{n-1}}\cdots \partial_{u_1} Q(x,t).$$
In particular,
$$\ct_x xQ(x,t)\prod_{i=1}^n\frac1{ x-u_i}=\sum_{i=1}^n Q(u_i,t)
\prod_{1\le j\le n \atop j\ne i}\frac{1}{u_i-u_j}.$$
\end{thm}

We need the following lemma.
\begin{lem}\label{l-l11}
If $Q(x,t)$ is $\pt$ in $x$, and if $u=u(t)$ is a formal power
series in $t$ with constant term $0$, then
\begin{align}\label{e-l11a}
\frac1{x-u}Q(x,t)=\xx \frac{ Q(u,t)}{1-u\xx} + \partial_u Q(x,t),
\end{align}
in which the first part is $\nt$ in $x$ and the second part is
$\pt$ in $x$.
\end{lem}
\begin{proof}
Equation \eqref{e-l11a} follows from direct computation. We have
\begin{align*}
\xx \frac{Q(u,t)}{1-u\xx} +\frac{Q(x,t)-Q(u,t)}{x-u} &=
\frac{Q(u,t)}{x-u}+\frac{Q(x,t)-Q(u,t)}{x-u}=\frac{Q(x,t)}{x-u}.
\end{align*}
Now clearly on the left most sides of the above equation, the
first term is $\nt$ in $x$ and the second term is $\pt$ in $x$.
This completes the proof.
\end{proof}

\begin{rem} Note that equation \eqref{e-ctqq} can be obtained by a specialization in
\eqref{e-l11a}. We have
\begin{align*}
\ct_x \frac{x}{x-u} Q(x,t) & = \left. \pt_x \frac{1}{x-u} xQ(x,t)\right|_{x=0} \\
&= \left. \frac{xQ(x,t)-uQ(u,t)}{x-u} \right|_{x=0} =Q(u,t)
\end{align*}
\end{rem}

\begin{proof}[Proof of Theorem \ref{t-mdivid}]
The proof is by repeatedly using equation \eqref{e-l11a}. In the
following computation, ``other terms" refers to
 terms with only negative powers in $x$.
Using the fact that $1/(x-u_i)$ belongs to $\xx R[\xx][[t]]$, we
have
\begin{align*}
&\frac1{(x-u_1)(x-u_2)\cdots (x-u_n)}Q(x,t) \\
&\qquad\qquad\qquad=
 \frac1{(x-u_n)(x-u_{n-1})\cdots (x-u_2)}(\partial_{u_1} Q(x,t)
+\mbox{other terms} )\\
&\qquad\qquad\qquad=\frac1{(x-u_n)(x-u_{n-1})\cdots
(x-u_2)}\partial_{u_1} Q(x,t) +\mbox{other terms} .
\end{align*}
Since $\partial_{u_1}Q(x,t)$ is still in $R[[x,t]]$, we can repeat
the above computation, and get
$$\frac1{(x-u_n)(x-u_{n-1})\cdots (x-u_1)} Q(x,t)
=\partial_{u_n}\partial_{u_{n-1}}\cdots \partial_{u_1} Q(x,t)+
\mbox{other terms}.
$$
This completes the proof of the first part of the theorem. Now \\
 $\partial_{u_n}\partial_{u_{n-1}}\cdots \partial_{u_1} x Q(x,t)\in
R[[x,t]]$, and
$$\partial_{u_n}\partial_{u_{n-1}}\cdots \partial_{u_1} x Q(x,t)
=\sum_{i=0}^n u_i Q(x_i) \prod_{1\le j\le n \atop j\ne i}
\frac{1}{u_i-u_j},$$ where we identify $t$ as $u_0$. Then by
setting $u_0=t=0$, we get the constant term
$$\sum_{i=1}^n Q(x_i) \prod_{1\le j\le n \atop j\ne i} \frac{1}{u_i-u_j}. $$
\end{proof}

\begin{rem}
We shall mention two things about this theorem. First, $u_i$ is
allowed to be $0$. Second, $u_i$ is allowed to be equal to $u_j$
for some $j\ne i$. For example, to deal with the case
 $u_i=u_j$,  we replace $u_j$ with $u_i+a$, let
$a$ approach $0$, and apply L'H\^opital's rule.
\end{rem}

As an application, we give a short proof of the well-known theorem
about the diagonal. See \citep[Theorem 6.33]{EC2}. We also need
the Puiseux's Theorem, which will be discussed later in the last
section of this chapter.

\begin{thm}
The diagonal of a rational power series in two variables is
algebraic.
\end{thm}

This result is clearly a consequence of the following theorem.

\begin{thm}
If $f(x,t)\in K((x))((t))$ is rational, then $\ct_x f(x,t)$ is
algebraic.
\end{thm}
\begin{proof}[Sketch of the proof]
Write $f(x,t)=N(x,t)/D(x,t)$ as a quotient of two polynomials. By
Puiseux's Theorem, there exists a positive integer $M$ such that
we can factor $D(x,t)$ as $A(t)(x-u_1)\cdots (x-u_m) (x-v_1)\cdots
(x-v_n)$, where $A(t)$ is a rational function in $t$, and
$u_i,v_j$ lie in $K((t^{1/M}))$ with $u_i$ having positive order
and $v_i$ having nonpositive order.

Thus Theorem \ref{t-mdivid} can be applied after multiplying $f$
by $x/(x-u_0)$, where $u_0=0$. The result will be a rational
function in the $u_i$'s, $v_i$'s and $t$, possibly obtained after
some derivatives and specializations, and hence is algebraic.
\end{proof}

Gessel observed a more practical method when dealing with the
situation in Theorem \ref{t-mdivid}. The idea is to use partial
fraction decomposition, together with the following lemma.

\begin{lem}\label{l-c1-m}
If $u\in tR[[t]] $ then for any nonnegative integer $k$ and
$Q(x,t)$ that is $\pt$ in $x$, we have
$$
\ct_x xQ(x,t)\frac{1}{(x-u)^{k+1}} = \frac{1}{k!}\left.
\frac{\partial^{k}}{\partial x^k} Q(x,t) \right|_{x=u}.$$
\end{lem}
The proof of this Lemma is trivial by linearity.

Since we have a partial fraction decomposition of the product
$\prod_{i=1}^n{ (x-u_i)}^{-1}$, we can apply Lemma \ref{l-c1-m} to
the evaluation of $\ct_x xQ(x,t)\prod_{i=1}^n{ (x-u_i)}^{-1}$.

Let $G(x,t)$ belong to $K[[x,t]]$. If $X=X(t)$ satisfies
$G(X,t)=0$, then we say that $X$ is a {\em root} of $G(x,t)$ for
$x$. Such $X$ is usually a fractional Laurent series, which we
will discuss later, so the order of $X$ is well defined. If $X$ is
a root of $G(x,t)$ for $x$ and if $X$ has positive order, then we
say $X$ is a {\em positive root} of $G(x,t)$.

Using the following well-known result (see e.g., \citep[Theorem
4.2]{ira}), we can generalize Corollary \ref{c-base} in another
direction. See Theorem \ref{t-lagrange1} below.
\begin{lem}
If $G(x,t)\in R[[x,t]]$ and $G(x,0)$ can be written as
$ax+\text{\rm higher terms}$ with $a\ne 0$, then $G(x,t)$ has a
unique positive root $X(t)$ for $x$, and this $X=X(t)$ belongs to
$tR[[t]]$.
\end{lem}

Theorem \ref{t-lagrange1} below is the most useful result in
chapter \ref{s-kxt}. It is a generalization of the Lagrange
inversion formula. In the case that $F$ is independent of $t$ and
$G(x,t)=x-tH(x)$, where $H(x)$ is a power series, we can easily
derive Lagrange's inversion formula. See Stanley \citep[Theorem
5.4.2]{EC2}. Moreover, the multivariate Lagrange inversion formula
can be deduced from it, as discussed further in section
\ref{ss-lag}.

Other applications of Theorem \ref{t-lagrange1} can be found in
chapter \ref{s-lattice} on lattice path enumeration. In the proof
of a conjecture in \citep{bous} about walks on the slit plane,
Theorem \ref{t-lagrange1} plays an important role.

\begin{thm}\label{t-lagrange1}
Let $G(x,t), F(x,t) \in K[[x,t]]$. If  $G(x,0)$ can be written as \\
$ax+\text{higher terms}$ with $a\ne 0$, then
\begin{equation}
\CT_{x} \frac{x}{G(x,t)} F(x,t)= \left.
\frac{F(x,t)}{\displaystyle{\partial\over
\partial x} G(x,t)} \right|_{x=X},
\end{equation}
where $X=X(t)$ is the unique element in $tK[[t]]$ such that
$G(X,t)=0$.
\end{thm}
\begin{proof}
Since $X(t)$ is the unique root of $G(x,t)=0$ that lies in
$tK[[t]]$, we have
$$\frac{G(x,t)}{x-X}=\frac{G(x,t)-G(X,t)}{x-X}=\partial_X G(x,t). $$
This is an element in $K[[x,t]]$ with nonzero constant term. Thus
applying Corollary \ref{c-base}, we get
\begin{align*}
\CT_{x} \frac{x}{G(x,t)} F(x,t) &= \CT_x \frac{x}{x-X}
\left(\frac{G(x,t)}{x-X}\right)^{-1}
F(x,t) \\
&=\left. \left(\frac{G(x,t)}{x-X}\right)^{-1}
F(x,t) \right|_{x=X} \\
&= \left. \frac{F(X,t)}{1-t {\partial\over
\partial x} G(x,t)} \right|_{x=X}
\end{align*}
\end{proof}

\begin{rem}
In Theorem \ref{t-lagrange1}, if we are working in the ring
$K((x^{-1}))$, we shall require that $X$ is the unique root of
$G(x,t)$ that lies in $K[[t]]$.
\end{rem}

The well-known rule for change variables in the computation of
residues is the following.
\begin{thm}
Let $K$ be a field, and $h(x)\in K((x))$. Suppose $n\in \ZZ$ be
such that $h(x)/x^n \in K[[x]]$ has nonzero constant term. Then
for any $\Phi(x)\in K((x))$, we have
\begin{equation}
\res_x \Phi(h) h'(x) =n \res_y \Phi(y),
\end{equation}
provided $\Phi(h(x))\in K((x))$.
\end{thm}
There is also a similar result for residues in $K((x))((t))$.
\begin{thm}
Let $K$ be a field, and $f(x,t)\in K((x))((t))$. Suppose the
initial term of $h(x)$ is $cx^bt^a$. Then for any $\Phi(x,t)\in
K((x))((t))$, if $\Phi(h,t)\in K((x))((t))$, then we have
\begin{equation}\label{e-jacobi1}
\res_x \Phi(h,t) \pad{x} h(x,t) =b \res_y \Phi(y,t).
\end{equation}
\end{thm}
\begin{proof}
By linearity, it suffices to show that this is correct for
$\Phi(x,t)=x^it^j$. Since $t$ can be factored out, we can assume
$\Phi(x,t)=x^i$. Then the right-hand side of \eqref{e-jacobi1}
becomes $0$ for $i\ne -1$ and $b$ for $i=-1$. Now let us compute
the left-hand side. If $i\ne -1$, then
$$\res_x h^i \pad{x} h =\res_x \frac1{i+1} \pad{x} h^{i+1} =0.$$
If $i=-1$, then $h^{-1}\pad{x}{h} = \pad{x}{\log  h}$. But $\log
h$ is not in $K((x))((t))$. We can overcome this by using the
formula
$$(fg)^{-1}\pad{x}{fg} = f^{-1}\pad{x}f+g^{-1} \pad{x}g,$$
which can be easily checked.

By the second decomposition, $h(x,t)$ can be uniquely factored as
$$h(x,t)=c t^a x^b h_1(x)h_2(x,t),$$
where $h_1(x)\in K[[x]]$ with constant term $1$, $h_2(x,t)\in
K((x))[[t]]$ with constant term $1$, and $c\in K$ is a constant.
Hence $\log h_1$ and $\log h_2$ belong to $K((x))[[t]]$, and we
have
\begin{align*}
\res_x h^{-1}(x,t) \pad{x} h(x,t) &= \res_x (ct^a x^b)^{-1}
\left(\pad{x}{ct^ax^b}\right) +  h_1^{-1}(x) \pad{x} h_1(x)
+h_2^{-1}(x,t) \pad{x} h_2(x,t) \\
&= \res_x b x^{-1}+\pad{x} \log  (h_1(x))+\pad{x}\log  (h_2(x,t))
=b
\end{align*}
\end{proof}

\section{Binomial Coefficients and Combinatorial Sums}

Binomial coefficients $\binom{n}{k}$ are the most frequently used
numbers in combinatorics. They are defined by
$$\binom{n}{k}= \frac{n(n-1)\cdots (n-k+1)}{k!}.$$
This holds for all nonnegative integers $k$ and complex numbers
$n$. From the well-known binomial
 theorem, we see that
 \begin{align}
\binom{n}{k}=\ct_\alpha
\frac{(1+\alpha)^n}{\alpha^k}.\label{e-binom}
\end{align}
 Starting from this identity, we can prove many identities
 involving binomial coefficients.

\begin{exa}
Compute $f(n)=\sum_{k=0}^n \binom{n}{k}.$
\end{exa}
The clever way is to use the formula
$$(1+x)^n =\sum_{k=0}^n \binom{n}{k}x^k.$$
By setting $x=1$, we get $f(n)=2^n$. This is a specialization of a
more general formula. But we are not always so lucky to find the
corresponding general formula. Here we give two methods to apply
Theorem \ref{t-lagrange1} on this trivial example. We shall see
that working in the field $K((\xx))((t))$ might be better than in
$K((x))((t))$.

Method $1$ is to show that $\sum_{n\ge 0} f(n)x^n=1/(1-2x)$. The
working field is $K((\alpha))((x)).$
\begin{align*}
\sum_{n\ge 0} f(n)x^n & = \sum_{n\ge 0} \sum_{0\le k\le n} \binom{n}{k} x^n\\
&= \sum_{n\ge 0} \sum_{0\le k\le n} \ct_\alpha \frac{(1+\alpha)^n}{\alpha^k} x^n\\
&= \ct_\alpha  \sum_{k\ge 0} \alpha^{-k} \sum_{ n\ge k} (1+\alpha)^n x^n\\
&= \ct_\alpha  \sum_{k\ge 0} \alpha^{-k} \frac{(1+\alpha)^kx^k }{1-(1+\alpha )x}\\
&= \ct_\alpha \frac{\alpha}{\alpha-(1+\alpha)x}\cdot
\frac{1}{1-(1+\alpha)x}.
\end{align*}
Now the term after the ``$\cdot$" contains only positive powers in
$\alpha$. Solving the denominator $\alpha-(1+\alpha)x$ for
$\alpha$, we get $\alpha=x/(1-x)$, which is in $x\CC[[x]]$. Thus
we can apply Theorem \ref{t-lagrange1} and get
\begin{align*}
\sum_{n\ge 0} f(n)x^n = \frac{1}{1-x}\cdot \left.
\frac{1}{1-(1+\alpha)x} \right|_{\alpha=x/(1-x)} =\frac{1}{1-2x}.
\end{align*}

Method 2. The working field is $K((\alpha^{-1}))$. Since
$\binom{n}{k}=0$ for $k>n$, we have
\begin{align*}
f(n) & = \sum_{k\ge 0} \binom{n}{k} = \sum_{k\ge 0} \ct_\alpha
(1+\alpha)^n \alpha^{-k}
= \ct_{\alpha^{-1}}\frac{\alpha}{\alpha-1} \cdot (1+\alpha)^n \\
&= \frac{1}{1} \left. (1+\alpha)^n \right|_{\alpha=1} =2^n.
\end{align*}
Comparing the above two method, we see that in some cases, it is
much simpler to work in $K((\alpha^{-1}))$ than to work in
$K((\alpha))$.

The next example shows that the residue theorem might simplify the
computation a lot. The drawback is that in general we might not
know how to change the variables.
\begin{exa}
Compute $f(n)=\sum_{k=0}^{n-1}\binom{n+k-1}{k}2^{-k}.$ \citep[p.
98, Exer 5.53]{EC1}
\end{exa}

Method 1. We compute the generating function of $f(n).$
\begin{align*}
\sum_{n\ge 0}f(n)x^n & =\sum_{n\ge 0}x^n
\sum_{k=0}^{n-1}\ct_\alpha (1+\alpha)^{n+k-1}\alpha^{-k}2^{-k}\\
&= \ct_\alpha \sum_{k\ge 0} (1+\alpha)^{k-1}\alpha^{-k}2^{-k}\sum_{n>k}x^n(1+\alpha)^n\\
&=\ct_\alpha \sum_{k\ge 0} (1+\alpha)^{k-1}\alpha^{-k}2^{-k}
\frac{(1+\alpha)^{k+1}x^{k+1}}{1-(1+\alpha)x}\\
&= \ct_\alpha \frac{2\alpha}{2\alpha -(1+\alpha)^2x} \cdot
\frac{x}{1-(1+\alpha)x}.
\end{align*}
Now compute the positive root (root with positive order) of
$2\alpha -(1+\alpha)^2x$ for $\alpha$, and denote it by $A$. Then
$A=\frac{1-x-\sqrt{1-2x}}{x}$. Apply Theorem \ref{t-lagrange1},
and simplify. We get
\begin{align*}
\sum_{n\ge 0}f(n)x^n & = \frac{x}{1-2x}
\end{align*}

Method 2. We use the residue theorem. Observe that
$$f(n)=[x^{n-1}](1-x)^{-1}(1-x/2)^{-n}.$$
Then we have
$$
f(n) =\res_x \frac{1}{(x-x^2/2)^n} \frac{1}{1-x} .
$$

Change variables by $y=x-x^2/2$. Then ${dy}/{dx}=1-x$, and
$(1-x)^2= 1-2y$. Hence
$$f(n)=\res_y y^{-n} (1-2y)^{-1}=2^{n-1}.$$

\begin{exa}
Show that $[x^{n-1}] (1+x)^{2n-1} (2+x)^{-n}=\frac{1}{2}$ is an
identity in $K((x))$. \citep[p. 98, Exer 5.57]{EC2}
\end{exa}

We use the residue theorem.
\begin{align*}
[x^{n-1}] (1+x)^{2n-1} (2+x)^{-n} & = \res_x (1+x)^{2n-1}
(2x+x^2)^{-n}
\end{align*}
Change variables by $y=2x+x^2$. Then ${dy\over dx}=2+2x$ and
$(1+x)^2=1+2y$. The above becomes
\begin{align*}
 \res_y (1+x)^{2n-1} y^{-n} (2+2x)^{-1} = {1\over 2} \res_y (1+y)^{n-1} y^{-n} ={1\over 2}.
\end{align*}

\vspace{4mm} We can see that the advantage of working in
$K((\xx))((t))$ happens when we can extend
 a finite sum to an infinite sum, in which case the formula for the sum of  a geometric series
 has a simple form. But the formula \eqref{e-binom} does not work in many situations.
For example, let us investigate the following summation.
\begin{align}\label{ex-fib}
\sum_{k=0}^{\lfloor n/2 \rfloor} \binom{n-k}{k}.
\end{align}

This summation can be extended to $0\le k \le n$, but not to all
nonnegative integer $k$,
 because $\binom{n-k}{k}$ is nonzero for $k>n$.

In order to extend this sum to all integer $k$, we need to
interpret  $\binom{n}{k}$ as zero when $n$ is a negative integer.
This can be done in the field $K((\alpha^{-1}))$, for we have
\begin{lem}
\begin{align}
\ct_{\alpha^{-1}} \frac{(1+\alpha)^n}{\alpha^k} =
\left\{\begin{array}{cc}
\binom{n}{k} & \text{ if } n\ge 0 \\
0 & \text{ if } n<0
\end{array}\right..
\end{align}
\end{lem}
\begin{proof}
The lemma is clearly true for nonnegative integer $n$. When $n$ is
a negative integer, $m=-n$ is a positive integer, and we have the
following expansion in $K((\alpha^{-1}))$.
\begin{align*}
 \frac{(1+\alpha)^n}{\alpha^k}= \alpha^{-k}\frac{1}{(1+\alpha)^m}=
  \alpha^{-k-m} \frac{1}{(1+1/\alpha)^m}= \alpha^{-k-m} \sum_{i\ge 0 } \binom{m+i-1}{i} \alpha^{-i}.
\end{align*}
Hence the lemma follows.
\end{proof}

Now the evaluation of \eqref{ex-fib} can be done as follows.
\begin{align*}
\sum_{k\ge 0} \ct_{\alpha^{-1}} \alpha^{-k} (1+\alpha)^{n-k}& =
\ct_{\alpha^{-1}} \sum_{k\ge 0} (1+\alpha)^n
\alpha^{-k}(1+\alpha)^{-k}\\
&= \ct_{\alpha^{-1}} \frac{\alpha}{\alpha^2+\alpha-1} \cdot
(1+\alpha)^{n+1}.
\end{align*}
Using the quadratic formula, we can solve for $x$ in the
denominator and get two roots $\frac{-1\pm \sqrt{5}}{2}$, denoted
by $A$ and $B$. Then applying Theorem \ref{t-mdivid}, we get
$$\frac{(1+A)^{n+1}-(1+B)^{n+1}}{A-B},$$
which turns out to be a Fibonacci number.

\section{Fractional Laurent Series and Puiseux's Theorem\label{s-puiseux}}
\subsection{Motivation and Background}
In both Theorem \ref{t-mdivid} and Theorem \ref{t-lagrange1}, the
third
 decomposition is obtained first. The evaluation of the constant term of
 $F(x,t)/G(x,t)$, where both $F$ and $G$ belong to $K[[x,t]]$, seems easier if we know its third
decomposition. It turns out that the positive roots of $G(x,t)$,
i.e., the $X$ satisfying
 $G(X,t)=0$ and having positive order, play a central role.
Theorem \ref{t-lagrange1} deals with a special case of such
evaluations. Our purpose in this section is to deal with a more
general case.

The root of a polynomial or power series $G(x,y)$ can be expressed
as a fractional Laurent series. Puiseux's Theorem \ref{t-puiseux}
deals with the case that $G(x,y)$ belongs to $K[[x]][y]$. In most
cases, this is sufficient. But we would like to consider a larger
set of $G(x,y)$. This results in a more general form of  Puiseux's
Theorem.

Before going further, let us review some basic concepts.

A {\em fractional Laurent series} (or Puiseux series) $\eta$ has
the form
$$\eta =\sum_{n\ge n_0} a_n x^{n/N}$$
for some $N\in \NN$. Let $K^{\fra}((x))$ (respectively,
$K^{\fra}[[x]]$)
 denote the set of all fractional Laurent series (respectively, fractional power
 series) over $K$. More precisely,
$$K^{\fra}((x))=\bcup_{N\ge 1} K((x))[x^{1/2},x^{1/3},\dots ,x^{1/N}].$$
Or in modern terminology, $K^\fra((x))$ is a direct limit.


Similarly we can define $K^\fra[[x]]$. It is clear that
$K^{\fra}((x))$ is the quotient field of the ring $K^{\fra}[[x]]$,
which contains only nonnegative powers in $x$. Note that
$\sum_{N\ge 1}x^{1/N}$ is not a fractional series in our sense of
the term.

For completeness, we include the following result. See
\citep[]{EC2}
 \begin{prop}
 The field $K^{\fra}((x))$ is an algebraic extension of $K((x))$;
 i.e., every $\eta \in K^{\fra}((x))$ satisfies an equation
 $$P_0(x)+P_1(x)\eta +\cdots +P_d(x) \eta^d =0,$$
 where $P_i(x)\in K((x))$ and not all $P_i(x)=0$.
 \end{prop}
\begin{proof}
Let $\eta =\sum_{n\ge n_0} a_n x^{n/N} \in K^{\fra}((x))$. There
are then unique series $\eta_0,\eta_1,\ldots \eta_{n-1} \in
K((x))$ such that
$$\eta=\eta_0+x^{1/N}\eta_1 +x^{2/N}\eta_2 +\cdots +x^{(N-1)/N}\eta_{N-1}.$$
Since $x^{i/N}$ are clearly algebraic over $K((x))$ for
$i=1,2,\ldots ,N-1$, the theorem follows from a general result in
field theory: For any extension field $E$ of any field $F$, the
elements of $E$ that are algebraic over $F$ form a subfield of $E$
containing $F$.
\end{proof}

\begin{thm}[Puiseux's Theorem]\label{t-puiseux}
Let $K$ be an algebraically closed field of characteristic zero
$($e.g., $K=\CC)$. Then the field $K^{\fra}((x))$ is algebraically
closed.
\end{thm}
There are many proofs of this theorem. One uses Newton polygons.
Here we will give another approach to Puiseux's Theorem. This new
approach handles a more general case. It bypasses the Newton
polygon argument. Of course Newton polygons will give us more
details about the roots.

Let $G(x,y)\in K[[x,y]]$. We say that $Y=Y(x)\in K^\fra ((x))$ is
a root of $G(x,y)$ if $G(x,Y(x))=0$. Puiseux's Theorem
characterizes all the roots of $G(x,y)\in K[[x]][y]$. Here we want
to characterize the roots of $G(x,y)\in K[[x,y]]$.

A problem arises in this consideration: the substitution of $Y\in
K^\fra ((x))$ for $y$ in $G(x,y)$ in general is not well defined.
This problem exists even if we only consider $Y\in K^\fra [[x]]$.

To avoid this situation, we define $Y$ to be a {\em positive root}
of $G(x,y)\in K[[x,y]]$ if $Y\in K^\fra[[x]]$ with $Y(0)=0$, for
such a root has positive order. By the composition law, the
substitution of such $Y$ for $y$ always results in an element in
$K^\fra [[x]]$. So one of our tasks  is to characterize all the
positive roots of $G(x,y)\in K[[x,y]]$. We can also consider
 roots $Y\in K^\fra[[x]]$ of $G(x,y)\in K[y][[x]]$, with which restriction,
 the substitution of $Y$ for $y$ is always valid. Note that the ring
 $K[y][[x]]$ contains the ring $K[[x]][y]$.
Finally, we will show that our theory implies Puiseux's Theorem.

Now let us study the roots of $G(x,y)\in K[[x,y]]$. If $Y\in
K[[x^{1/N}]]$ is a root of $G(x,y)$, then $G(x,Y)=0$. By setting
$x=0$, we get $G(0,Y(0))=0$. This is to say that $Y(0)$ is a root
of $G(0,y)$, which belongs to $K[[y]]$. In combinatorics, infinite
sum of nonzero elements in $K$ does not make sense. Thus it makes
no sense to say, for instance, that $\pi$ is a root of $\sin y$ in
the ring of formal power series. In fact, we have the following:

\begin{lem}
Let $Y=Y(x)$ be a fractional power series with nonzero constant
term, then $G(x,Y)$ makes sense if and only if $G(x,y)$ belongs to
$K[y][[x]]$.
\end{lem}
\begin{rem}
The lemma is trivial when $Y$ is a nonzero constant. The general
case follows from a more general result. At this moment, let us
take this lemma as a fact.
\end{rem}

So when characterizing those roots of $G(x,y)$ with nonzero
constant term, we require that $G(x,y)$ belong to $K[y][[x]]$, but
when characterizing positive roots of $G(x,y)$, there is no
restriction.

\subsection{Main Results}
To state our main results, we need some concepts. Let us first
establish some basic properties of the roots of $G(x,y)$ that lie
in $K^\fra[[x]]$.

In what follows, we always assume that $G(0,y)$ is a polynomial
unless specified otherwise. We will see that this assumption will
not lose any roots from reduction 1 in the next subsection.

\begin{lem}\label{l-pu-pre}
If for some positive integer $N$, $Y\in K[[x^{1/N}]]$ is a root of
$G(x,y)$, then $F(x,y):=G(x,y)/(y-Y)$ belongs to $K[[x^{1/N},y]]$,
and $\deg F(0,y)=\deg G(0,y)-1$.
\end{lem}
\begin{proof}
Write $G(x,y)=A(y)+x H(x,y)$ with $A(y)\in K[y]$ and $H(x,y)\in
K[[x,y]]$. In the following, the divided difference is taken with
respect to $y$. Since $G(x,Y)=0$, we have
\begin{align*}
F(x,y):=\frac{G(x,y)}{y-Y}&=\frac{G(x,y)-G(x,Y)}{y-Y}
=\partial_{Y} G(x,y),
\end{align*}
which belongs to $K[[x,y,Y]]\subset
K[[x,y,x^{1/N}]]=K[[x^{1/N},y]]$.

By setting $x=0$, we see that $F(0,y)=G(0,y)/(y-Y(0))$ belongs to
$K[y]$. Thus $\deg F(0,y)=\deg G(0,y)-1$.
\end{proof}

The following is the main lemma for us in studying the positive
roots of $G(x,y)$. We will give five reduction procedures to prove
this result in the next subsection.

\begin{lem}[Main Lemma]\label{l-pu-1}
If for some positive integer $d$, we can write $G(0,y)=y^dA(y)$,
where $A(y)$ belongs to $K[[y]]$, then $G(x,y)$ has at least one
positive root.
\end{lem}

Using Lemma \ref{l-pu-1}, we can show that:
\begin{lem}\label{l-pu-multi}
If $G(x,y)=y^d+x H(x,y)\in K[[x,y]]$ for some positive integer
$d$, then there are positive integers $N$, $k_1,\dots, k_r$, with
$k_1+\cdots +k_r=d$, and distinct $Y_i\in K[[x^{1/N}]]$ with
constant term $0$ for $1\le i \le r$, and $F(x,y)\in
K[[x^{1/N},y]]$ with constant term $1$, such that $G(x,y)$ can be
written as:
\begin{align}\label{e-pu-1}
G(x,y)=(y-Y_1)^{k_1}\cdots (y-Y_r)^{k_r} F(x,y),
\end{align}
and this form is unique up to the order of the factors.
\end{lem}
If \eqref{e-pu-1} is a factorization of $G(x,y)$ as above, then
 we say the {\em multiplicity}
of $Y_i$ is $k_i$.
\begin{proof}[Proof of Lemma \ref{l-pu-multi}]
By Lemma \ref{l-pu-1}, we can suppose that $Z_1 \in
K[[x^{1/N_1}]]$ is a root of $G(x,y)$ with constant term $0$.

Let
$$G_1(x,y):=\frac{G(x,y)}{y-Z_1}.$$
 Then by
Lemma \ref{l-pu-pre}, $G_1(x,y) \in K[[x^{1/N_1},y]]$. By setting
$x=0$, we get $G_1(0,y)=y^{d-1}.$

Thus we can repeat the above argument and get $N_{i}$ for
$i=2,3,\ldots ,d$ such that
$$ G_i(x,y):=\frac{G_{i-1}(x,y)}{y-Z_i},$$
with $N_{i-1}$ dividing $N_i$, $Z_i\in K[[x^{1/N_i}]]$,
$G_i(x,y)\in K[[x^{1/N_i},y]]$, and $G_i(0,y)=y^{d-i}$.

Now let $N=N_d$, and let $F(x,y)=G_d(x,y)$. Then
\begin{equation}
G(x,y)= F(x,y){\prod_{i=1}^{d} (y-Z_i)}.
\end{equation}
Clearly every $Z_i$ is a root of $G(x,y)$. Equation \eqref{e-pu-1}
is hence obtained by collecting equal terms of  $Z_i$'s and then
renaming.

The uniqueness follows from the following two facts. First,
$K[[x^{1/N},y]]$ is an integral domain. Second, If $Y\in
K^\fra[[x]]$ is a positive root of $G(x,y)$, then $Y=Y_i$ for some
$i$. For otherwise, by setting $y=Y$ in \eqref{e-pu-1}, we get
$F(x,Y)=0$, which contradicts the fact that $F(0,y)=1$.
\end{proof}

Now we can give our main results.
\begin{thm}[Generalized Puiseux Theorem]\label{t-pu-const}
For all $G(x,y)\in K[[x,y]]$, the number of positive roots
$($counted with multiplicity$)$ of $G(x,y)$ that lie in
$K^\fra[[x]]$ equals the order of $G(0,y)$.

For all $G(x,y)\in K[y][[x]]$, the number of roots $($counted with
multiplicity$)$ of $G(x,y)$ that lie in $K^\fra[[x]]$ and have
constant term $\alpha $ equals the multiplicity of $\alpha$ as a
root of $G(0,y)$.
\end{thm}
\begin{proof}
The first part follows from Lemma \ref{l-pu-multi}. The second is
obtained by a linear transformation. For details, see reduction 3
in the next subsection.
\end{proof}

To see that Theorem \ref{t-pu-const} implies Puiseux's Theorem, we
prove it as follows.
\begin{proof}[Proof of Puiseux's Theorem]
Suppose $G(x,y)\in K[[x]][y]$ is of degree $d$ in $y$. Let $s$ be
the degree of $G(0,y)$, and let
$$ G(0,y)=a(y-\alpha _1)^{k_1}\cdots (y-\alpha_r)^{k_r},$$
be the factorization of $G(0,y)$ in $K[y]$. Since
$G(x,y-\alpha)\in K[[x]][y]$ for all $\alpha \in K$, by Theorem
\ref{t-pu-const}, it has $k_i$ roots with constant term $\alpha_i$
for $i=1,\dots ,r$. Thus, we get $s$ roots of $G(x,y)$ that lie in
$K^\fra[[x]]$, among which, the number of positive roots equals
the order of $G(0,y)$.

Now consider $H(x,y):=y^dG(x,1/y)\in K[[x]][y]$. Then the order of
$H(0,y)=y^d G(0,1/y)$ is $d-s$. By theorem \ref{t-pu-const},
$H(x,y)$ has $d-s$ positive roots that lie in $K^\fra [[x]]$,
 the reciprocals of which are
clearly roots of $G(x,y)$ that have negative order.
\end{proof}

A direct consequence of the above argument is the following.
\begin{cor}
Suppose that $G(x,y)\in K[[x]][y]$. If $d$ is the degree of
$G(x,y)$ in $y$, then among all of the $d$ roots of $G(x,y)$,
$\ord(G(0,y))$ roots have positive order,
$\deg(G(0,y))-\ord(G(0,y))$ roots have zero order, and
$d-\deg(G(0,y))$ roots have negative order.
\end{cor}

Let us see some consequence of the generalized Puiseux Theorem
\ref{t-pu-const}.
\begin{cor}\label{c-pu-fac}
Any $G(x,y)\in K[[x,y]]$ can be uniquely factored as
$$G(x,y)=(y-Y_1)\cdots (y-Y_d)F(x,y),$$
where $d$ equals the order of $G(0,y)$, $Y_i\in K^\fra[[x]]$ with
$Y_i(0)=0$ for all $i=1,\ldots ,d$, and $1/F(x,y)\in K[[x,y]]$.
Moreover, $(y-Y_1)\cdots (y-Y_d)$ belongs to $ K[[x]][y]$.
\end{cor}
\begin{proof}
We only show that $(y-Y_1)\cdots (y-Y_d)\in K[[x]][y]$, for the
other part is easy.

Let $f=F(0,0)$ which is not zero. In the field $K((y))((x))$, the
constant term of $G(x,y)/(fy^d)$  in $x$ is $1$, and
$$\frac{G(x,y)}{fy^d}=\left(1-\frac{Y_1}{y}\right)\cdots \left(1-\frac{Y_d}{y}\right)\cdot \frac{F(x,y)}f.$$
By the uniqueness of the third decomposition (Lemma
\ref{l-unifactor}), we must have $(1-\frac{Y_1}{y})\cdots
(1-\frac{Y_d}{y}) \in K[\yy][[x]]$. Therefor $(y-Y_1)\cdots
(y-Y_d)\in K[[x]][y]$.
\end{proof}

\begin{cor}\label{t-pu-factorable}
Let $G(x,y)\in K[[x,y]]$. If the order of $G(0,y)$ is positive and
less than the degree of $G(x,y)$ in $y$, then $G(x,y)$ is not
irreducible in $K[[x,y]]$.
\end{cor}

\begin{prop}\label{p-fac-1}
Let $H(x,y) \in K[[x,y]]$. Working in $K((y))((x))$, we have an
expression of  $(1-xH(x,y)/y^d)_-$ in $K((y))^\fra ((x))$ as
follows.
\begin{align}\label{e-p-fac-1}
\left(1-x\frac{H(x,y)}{y^d}\right)_- =
\left(1-\frac{Y_1}{y}\right)\cdots \left(1-\frac{Y_d}{y}\right),
\end{align}
where $Y_1,\cdots Y_d \in K^\fra [[x]]$ are the $d$ positive roots
of $y^d-xH(x,y)$.
\end{prop}
\begin{proof}
Let $G(x,y)=y^d-x H(x,y)\in K[[x,y]]$. Then from Lemma
\ref{l-pu-multi}, we have a decomposition
$$G(x,y)=(y-Y_1)\cdots (y-Y_r) F(x,y),$$
where $Y_1,\ldots ,Y_d$ are the $d$ roots of $G(x,y)$ and
$F(x,y)\in K[[x^{1/N},y]]$ has constant term $1$. Thus
$$\frac{G(x,y)}{y^d}=1-x\frac{H(x,y)}{y^d} = (1-\frac{Y_1}{y})\cdots
(1-\frac{Y_d}{y}) F(x,y)$$ is the third decomposition in
$K((y))((x^{1/N}))$. Thus equation \eqref{e-p-fac-1} follows from
the uniqueness (Lemma \ref{l-unifactor}).
\end{proof}

\begin{thm}\label{t-lagmm}
Suppose $H(x,y)\in K[[x,y]]$ and $Y_1, Y_2,\ldots , Y_r$ in
${K}[[x^{1/N}]]$ be the $r$ positive roots of $y^r-xH(x,y)$ . Then
for $k>0$, we have the following identity.
\begin{equation}
\sum_{i=1}^r Y_i^k = k [y^{-k}] \log  \frac{1}{1-xH(x,y)/y^r}
\end{equation}
\end{thm}
\begin{proof}
By Theorem \ref{l-pu-multi}
 we have a decomposition
 $$y^r-xH(x,y)=(y-Y_1)(y-Y_2)\cdots (y-Y_r) F(x,y),$$
 where $F(x,y)\in {K}[[x^{1/N}]]$ has constant term $1$.
 Then we have
\begin{equation}
\log  \frac{1}{1-xH(x,y)/y^r}= \log  F(x,y)+\sum_{i=1}^r \log
\frac{1}{1-Y_i/y}\label{e-lagmm}
\end{equation}
The first term on the right hand side of \eqref{e-lagmm} contains
only positive powers in $y$, and the other terms contains only
negative powers in $y$. The theorem then follows by equating
coefficients of $y^{-k}$ on both sides of \eqref{e-lagmm}.
\end{proof}

In the special case of  $H(x,y)$ being a polynomial, we can say
something more.

\begin{thm}
\label{t-1-3rd0} Suppose that $H(x,y)$ is a polynomial. Let
$Y_1,\dots,Y_r$ be the $r$ positive roots of $G(x,y)=y^r-xH(x,y)$,
and let $Y_{r+1},\dots, Y_d$ be the other roots of $G(x,y)$, where
$d$ is the degree of $G(x,y)$ in $y$. Then
\begin{align}
\left(1-xH(x,y)/y^r\right)_0&=(-1)^{d-r} Y_{r+1}\cdots Y_{d}\,
[y^d]\,
\left(y^r-xH(x,y)\right)\label{e-1-3rd1}\\
&=(-1)^r Y_1^{-1}\dots Y_r^{-1} \,[y^0] \, \left(
y^r-xH(x,y)\right).\label{e-1-3rd2}
\end{align}
\end{thm}
\begin{proof}
As a polynomial in $y$, we have the factorization
\begin{align}\label{e-1-3rd0}
y^r-xH(x,y)=A(x)(y-Y_1)\cdots(y-Y_r)(1-y/Y_{r+1})\cdots (1-y/Y_d),
\end{align}
 where
$A(x)$ is fractional power series and we can check that the
initial term of $A(x)$ is $1$. Thus we have the following third
decomposition with respect to $y$:
$$1-\frac{xH(x,y)}{y^r}= A(x)(1-Y_1/y)\cdots(1-Y_r/y)(1-y/Y_{r+1})\cdots (1-y/Y_d). $$
Thus $A(x)=\left(1-xH(x,y)/y^r\right)_0$.

By comparing coefficient in $y^d$ on both sides of equation
\eqref{e-1-3rd0}, we get
$$[y^d]\, \left(y^r-xH(x,y)\right)= (-1)^{d-r}A(x)Y_{r+1}^{-1} \cdots Y_{d}^{-1}.$$
Equation \eqref{e-1-3rd1} thus follows. Equation \eqref{e-1-3rd2}
follows from the fact that
$$(-1)^d Y_1\cdots Y_d \,[y^d]\, \left(y^r-xH(x,y)\right) =[y^0] \, \left(y^r-xH(x,y)\right).$$
\end{proof}

When $r=1$, Theorem \ref{t-lagmm} reduces to the following, which
was first obtained in \citep{ira} and was shown to be equivalent
to Lagrange inversion formula.

\begin{prop}
Let $H(x,y)$ belong to $K[[x,y]]$, and let $Y$ be the unique
positive root of $y-xH(x,y)$. Then
$$ Y^k =k [y^{-k}] \log \frac{1}{1-xH(x,y)/y}.$$
\end{prop}

Now let us see an application of Theorem \ref{t-1-3rd0}.
\begin{exa}
In a complete solution to the so-called tennis ball problem
\citep[Theorem 1]{anna-marc}, the final generating function is
given by
$$Q(z)= \frac{-1}{z}(1-w_1)\cdots (1-w_l),$$
where $w_1,\dots ,w_l$ are fractional power series that satisfying
the equation
$$(w-1)^l-zw^{k+l}=0.$$

By changing variables $w=1+u$ and $w_i=1+u_i$, we have
$$Q(z)= \frac{-1}{z} (-1)^l u_1\cdots u_l,$$
with $u_i$ being the positive roots of $u^l-z(1+u)^{k+l}$. Thus
applying Theorem \ref{t-1-3rd0}, we get
$$\log Q(z) =\ct_u \log \frac{1}{1-z(1+u)^{k+l}/u^l}=\sum_{n\ge 0} \frac{1}{n} \binom{nk+nl}{nl}.$$
Therefore, we get the following concise formula:
\begin{align}
Q(z)= \exp \Big(\sum_{n\ge 0} \frac{1}{n} \binom{nk+nl}{nl} z^n
\Big)
\end{align}

This result is very similar to the generating function of paths
with steps $(1,k)$ or $(1,-l)$ that start at $0,0$, end on the
horizontal axis, and never goes below the horizontal axis. Let
$R(t)$ be the generating function, then we have \citep{bizley}
$$ R(z) =\exp \Big(\sum_{n\ge 0} \frac{1}{nk+nl} \binom{nk+nl}{nl} z^n \Big). $$
\end{exa}

\begin{thm}\label{t-1-mdividg}
Suppose that $H(x,y)\in K[[x,y]]$. If $Y_1, Y_2,\ldots , Y_r$ are
the $r$ distinct positive roots of $y^r-xH(x,y)$ in
${K}[[x^{1/N}]]$. Then for any $\Phi(x,y)\in K[[x,y]]$, we have
\begin{equation}\label{e-16-divid}
\ct_y \frac{y}{y^r-xH(x,y)}\Phi(x,y) =\sum_{i=1}^r
\left.\frac{\Phi(x,y)}{ry^{r-1}-x{\partial \over \partial y}
H(x,y)}\right|_{y=Y_i}
\end{equation}
\end{thm}
\begin{proof}
Clearly, we have a factorization of $y^r-xH(x,y)$
$$y^r-xH(x,y)=(y-Y_1)\cdots (y-Y_r) F(x,y),$$
where $F(x,y)$ has constant term $1$ and belongs to $K[[x,y]]$. So
we have
$$\ct_y \frac{y}{y^r-xH(x,y)}\Phi(x,y)= \ct_y \frac{y}{(y-Y_1)\cdots (y-Y_r)}\cdot
\frac{1}{F(x,y)}\Phi(x,y).$$ Now $1/F(x,y)$ contains only
nonnegative powers in $y$, and so does $\Phi(x,y)$, and $Y_i$ has
positive order. In the field $K((y))^\fra[[x]]$, we can apply
Theorem \ref{t-mdivid}. Equation \eqref{e-16-divid} then follows
by putting $F(x,y)=G(x,y)/\prod_{i=1}^r(y-Y_i)$, and applying the
L'H\^opital's rule.
\end{proof}

The application of this theorem will appear in the chapter on
lattice path enumeration.

\subsection{Proof of the Main Lemma}

The main lemma says that if the order of $G(0,y)$ is positive,
then $G(x,y)$ has at least one positive root. The basic idea of
proving this lemma is: Find the initial term of the assumed
positive root and then inductively find the next one. This idea
also works when the characteristic of $K$ is nonzero.

We shall give five reduction procedures to prove the main lemma.
Reduction 1 reduces the general $G(x,y)$ to the case $G(0,y)=y^d$
for some positive integer $d$. Reduction 2 deals with the base
case when $G(0,y)=y$. Reduction 3 says that using a linear
transformation, the computation of roots of $G(x,y)$ with constant
term $\alpha$ can be converted into the computation of positive
roots of some $H(x,y)\in K[[x,y]]$. Reduction 4 handles the case
$G(x,y)=y^d+xH(x,y)$ for some $H(x,y)\in K[[x,y]]$ with $H(0,0)\ne
0$. Reduction 5 will be used to deal with the case
$G(x,y)=y^d+xH(x,y)$ for some $H(x,y)\in K[[x,y]]$ with
$H(0,0)=0$. The first 4 reduction procedures are routine. The
fifth
 is complicated. Note that reduction 5 covers all the cases of $G(0,y)=y^d$.

Now we begin to give these reductions.

Reduction 1. $G(x,y)\longrightarrow \frac{G(x,y)}{A(y)}, \text{
where } G(0,y)=y^d A(y)$.

For any $G(x,y)\in K[[x,y]]$, $G(0,y)$ clearly belongs to $
K[[y]]$. Thus it can be written as $y^d A(y)$ for some $d\ge 0$
and $A(y)\in K[[y]]$ with constant term nonzero. Since $A(y)$ is
invertible, i.e., $1/A(y)\in K[[y]]$, replacing $G(x,y)$ by
$G(x,y)/A(y)$ will not change any positive roots.

Now we can assume $G(0,y)=y^d$. Note that if $d=0$, i.e.
$G(0,y)=1$, $G(x,y)$ does not have any positive roots since it is
invertible.

Reduction 2. $ \text{If } G(0,y)=y, \text{ then  apply Lemma
\ref{l-pu-d1} to get the unique root}.$

This is the base case of $d=1$. We have the following well-known
result, e.g., \citep[Theorem 4.2]{ira}.
\begin{lem}\label{l-pu-d1}
For any $H(x,y)\in K[[x,y]]$, there is a unique $Y\in K^\fra
[[x]]$ with $Y(0)=0$ such that $Y+x H(x,Y)=0$. Moreover, this $Y$
belongs to $K[[x]]$.
\end{lem}
The proof of this lemma is by assuming $Y(x)=a_1x+a_2x^2+\cdots $,
and solve for $a_1$, $a_2,$ and so on subsequently. In fact, every
$a_i$ is obtained after finitely many additions and
multiplications. In other words, $a_i$ lies in the ring generated
by the coefficients of $H(x,y)$.

Reduction 3. $ G(x,y)\longrightarrow H(x,y):=G(x,y-\alpha),
Y:=Y'+\alpha. $

Let $G(x,y)\in K[[x,y]]$ with $G(0,y)\in K[y]$. If $\alpha\ne 0$
is a root of $G(0,y)$, and $H(x,y)=G(x,y-\alpha)$ belongs to
$K[[x,y]]$, (i.e. $G(x,y)\in K[y][[x]]$),
 then $Y'$ is a root of $H(x,y)$ with
constant term $0$ if and only if $Y=Y'+\alpha $ is a root of
$G(x,y)$.

Reduction 4. If $d\ge 2$, and $G(x,y)=y^d+xH(x,y)$ with
$H(0,0)\ne 0$ then apply Lemma \ref{l-pu-nonzero}.

This case is covered by the following lemma.
\begin{lem}\label{l-pu-nonzero}
If $G(x,y)=y^d+xH(x,y)\in K[[x,y]]$ with $d\ge 2$ and $H(0,0)\ne
0$, then $G(x,y)$ has $d$ distinct positive roots.
\end{lem}
\begin{proof}
Let $c=H(0,0)$. Since $c\ne 0$, $H(x,y)/c$ has constant term $1$
and thus $(H(x,y)/c)^{1/d}$ (with constant $1$) is well defined
 in $K[[x,y]]$. Thus we have
 $$y^d +x H(x,y)= \prod _{i=1}^d \left(y+\zeta _i x^{1/d}
 \left(\frac{H(x,y)}{c}\right)^{1/d}\right)$$
 where
 $y^d+c=\prod _{i=1}^d (y+\zeta _i)$. Then solving the above $d$
 factors will give us
$d$ distinct elements,
 $Y_i=-\zeta_i x^{1/d}+\cdots$ for $1\le i \le d$, that lie in
 $ K[[x^{1/d}]]$. All these $Y_i$  are clearly
roots of $G(x,y)$. Thus by Lemma \ref{l-pu-multi}, $G(x,y)$ has no
other roots.
\end{proof}

\begin{rem}
Lemma \ref{l-pu-nonzero} does not apply when the characteristic of
$K$ is $p\ne 0$. The reason is that if $p$ divides $d$ then we can
not take the $d$-th root. We shall see a counterexample later.
\end{rem}

Reduction 5 needs more explanations. Let $G(x,y)=y^r+xH(x,y)\in
K[[x,y]]$ with $r\ge 2$ and $H(0,0)=0$. The basic idea is to
factor out a power of $x$ in the assumed positive root $Y$ of
$G(x,y)$.

Rewrite $G(x,y)$ in the following form
$$G(x,y)=y^r-x\sum_{i=0}^\infty x^{n_i} A_i(x)y^i,$$
where $A_i(x)\in K[[x]]$ has nonzero constant term and $n_i\ge 0$
for all $i$. Let
$$s=\min_{0\le j\le r-1} \frac{n_j+1}{r-j}
\text{ and } y=x^s \tilde{y}.$$ Then $s>0$ and
\begin{align*}
G(x,y) &=x^{sr} \tilde{y}^r+x\sum_{i\ge 0} x^{n_i}A_i(x)
x^{si}\tilde{y}^i = x^{sr} \left( \tilde{y}^r+\sum_{i\ge 0}
x^{n_i+1+si-sr}A_i(x) \tilde{y}^i\right)
\end{align*}
Clearly, $n_i+1+si-sr\ge 0$ for all $i\ge r$. By the choice of
$s$, $n_i+1+si-sr\ge 0$ for $0\le i \le r-1$, and the equality
holds for at least one $i$ with $0\le i\le r-1$. Let $j$ be the
smallest such that $n_j+1-s(r-j)=0$, which yields
$s=\frac{n_j+1}{r-j}$.

So if we let
$$\tilde{G}(x,\tilde{y})
=\tilde{y}^r+\sum_{i\ge 0} x^{n_i+1+si-sr}A_i(x) \tilde{y}^i,$$
then $\tilde{G}(x,\tilde{y})\in K[[x^{1/(r-j)},y]]$, and we have
the relation
$$x^{sr}\tilde{G}(x,\tilde{y})=G(x,x^{s}\tilde{y}).$$
This relation guarantees that $\tilde{G}(x,\tilde{y}+\alpha) \in
K[[x^{1/(r-j)},\tilde{y}]]$ for any $\alpha \in K$, since
$$x^{sr}\tilde{G}(x,\tilde{y}+\alpha )=
G(x,x^{s}(\tilde{y}+\alpha))=G(x,y+\alpha x^{s})$$ belongs to $
K[[x^{1/(r-j)},y]]$. Thus the condition in reduction $3$ is
satisfied.

Moreover, if we denote $\tilde{G}(0,\tilde{y})$ by $B(\tilde{y})$,
then it is a polynomial and
$$B(\tilde{y})=\tilde{y}^r +\sum_{i} A_i(0) \tilde{y}^i, $$
where the sum runs over all $0\le i\le r-1$ such that
$n_i+1+si-sr=0$.  Clearly, $B(\tilde{y})$ has  the highest term
$\tilde{y}^r$ and lowest term $A_j(0)\tilde{y}^j$.

Pick a root $\alpha$ of $B(y)$ and apply reduction $3$ and then
reduction $1$. Denote the result of the above procedure by
$G'(x,y)$, and the assumed root of $G'(x,y)$ by $Y'$. Then
$G'(0,y)=y^{k}$ for some $k$ with $0\le k\le r$, and
$Y=x^sY'+\alpha x^s$.

{\bf Observation:} The only chance for $k$ to be $r$ is when
$B(y)=(y-\alpha)^r$ for some $\alpha \neq 0$. In this case, $j=0$
and $A_{r-1}(0)=r(-\alpha)^{r-1}\ne 0$ (which
 is not true when the characteristic of $K$ is not $0$),
 and hence $n_{r-1}+1+s(r-1)-sr=0$,
which implies that $s=n_{r-1}+1$ is a positive integer.

$\textbf{Reduction $5$: } G(x,y)\longrightarrow G'(x,y),  \text{
and } Y:=x^sY'+\alpha x^s.$

Now we can prove our main lemma.
\begin{proof}[Proof of the Main Lemma \ref{l-pu-1}]
We prove this lemma by induction on $d$. Lemma \ref{l-pu-d1} shows
that the lemma is true for the base case $d=1$. Now suppose it is
true for $1,2,\ldots ,d-1$. Then we need to show it is true for
$d$.

We apply the following reduction procedure.
\begin{enumerate}

\item Apply reduction $1$ to make $G(0,y)=y^d$ for some $d\ge 0.$

\item If $G(0,y)=y$, then apply reduction $2$ to get the unique
root.

\item If $G(x,y)=y^d+x H(x,y)$ with $d\ge 2$ and $H(0,0)\ne 0$, then
apply reduction $4$ to get $d$ distinct roots.

\item If $G(x,y)=y^d+x H(x,y)$ with $H(0,0)=0$ and $d\ge 2$, then
apply reduction $5$.
\end{enumerate}
Steps 1, 2, and 3 will give us the result immediately. In step 4,
if applying reduction 5 gives us $G'(0,y)=y^k$ for some $k<d$,
then we can get a root by the induction hypothesis.

So the above reduction procedure will stop unless beginning at
some point, every application of reduction $5$ results in some
$G'(0,y)=y^d$. Therefore, we can assume that at this point,
$G_0(x,y)=y^d+xH_0(x,y)$, and that $G_{i+1}(x,y)$ is obtained from
$G_i(x,y)$ by applying reduction 5 for $i=0,1,2,\ldots $.

When applying reduction $5$ on $G_i(x,y)$, we get a positive
rational number $s_{i+1}$, and $G_{i+1}(x,y)=y^d+xH_{i+1}(x,y)$,
and $\alpha_{i+1}\ne 0$. The relation between the assumed root
$Y_i$ of $G_i(x,y)$ and $Y_{i+1}$ of $G_{i+1}(x,y)$ is given by
$$Y_{i}=\alpha _{i+1}x^{s_{i+1}}+Y_{i+1}.$$

Now let
$$Y=\sum_{i\ge 1} \alpha_i x^{s_1+s_2+\cdots +s_i}.$$
From the construction of $Y$, we see that
 $G_0(x,Y)=0$, provided that $Y$ is a fractional power series.
In fact $Y$ is a power series, since from the observation, every
$s_{i+1}$ is a positive integer.
\end{proof}

\begin{exa}
We consider the positive roots of
$G(x,y)=y^3+x^3+x^4+3x^2y+3xy^2-y^4e^{xy}.$
\end{exa}
First, we shall apply reduction 5. Using the notation in reduction
5, we have $r=3$, $n_0=2, n_1=1, n_2=0$, and hence
$\frac{n_j+1}{r-j}$ equals $1$ for $j=0,1,2$. Therefore $s=1$, and
we shall let $y=xy_1$. Now
\begin{align*}
G(x,y)&=x^3y_1^3+x^3+x^4+3x^3y_1+3x^3y_1^2-x^4y_1^4e^{x^2 y_1} \\
&=x^3\left( y_1^3+1+x+3y_1+3y_1^2 -xy_1^4e^{x^2y_1}\right).
\end{align*}

Let $G_1(x,y_1)= y_1^3+1+3y_1+3y_1^2+x -xy_1^4e^{x^2y_1}. $ Then
$G_1(x,Y_1)=0$ if and only if $G(x,Y)=0$, where $Y=xY_1$.

To find the roots of $G_1(x,y_1)$ for $y_1$, we use reduction 3.
Let $y_1=y_2-1$. Then \begin{align*} G_1(x,y_1)=G_2(x,y_2)&=y_2^3
+x(1-(y_2-1)^4 e^{x^2(y_2-1)})\\
&=y_2^3+x\Big(1-e^{-x^2}+e^{-x^2}(4-x^2)y_2-e^{-x^2}(6-4x^2+\frac{x^4}{2})y_2^2\\
&\qquad\qquad\qquad\qquad\qquad\qquad\qquad +\text{higher order
terms}\Big).
\end{align*}
To find the roots of $G_2(x,y_2)$, we need to use reduction 5
again. This time $r=3, n_0=2,n_1=0,n_2=0$, and $\frac{n_j+1}{r-j}$
equals $1,1/2,1$ for $j=0,1,2$ respectively. Therefore $s=1/2$,
and we shall let $y_2=x^{1/2}y_3$. After some algebraic
manipulations, we get
\begin{multline*}
   G_3(x,y_3)= y_3^3+\Big(x^{-1/2}(1-e^{-x^2})+e^{-x^2}
(4-x^2)y_3 \\
 -e^{-x^2}(6-4x^2+x^4/2)x^{1/2}y_2^2+x^{1/2}\cdot
     \text{higher ordered terms} \Big).
\end{multline*}

Now $G_3(0,y_3)=y_3^3+4y_3=y_3(y_3+2\sqrt{-1})(y_3-2\sqrt{-1}).$
One positive root of $G_3(x,y_3)$ for $y_3$ can be found
immediately by reduction 2. The other two zero order roots have
constant terms $2\sqrt{-1}$ and $-2\sqrt{-1}$.

Recalling that $y=xy_1,$ $ y_1=y_2-1$, and $y_2=x^{1/2}y_3$, we
find that the three positive roots of $G(x,y)$ are
$-x+2\sqrt{-1}x^{3/2}+\text{higher order terms}$,
$-x-2\sqrt{-1}x^{3/2}+\text{higher order terms}$, and $-x+x^2\cdot
T(x)$ with $T(x)\in K[[x^{1/2}]]$. All of the three positive roots
are in $K[[x^{1/2}]]$.

\begin{exa}
Let $K$ be an algebraically closed field with characteristic $p\ne
0$. Consider the roots of $G(x,y)=y^p-y-\xx$.
\end{exa}

\citep{chevalley} proved that $y^p-y-\xx$ does not have a root for
$y$ in $K^\fra((x))$. In addition, he gave the following
factorization in a certain field:
$$y^p-y-\xx= \prod_{i=0}^{p-1} \Big(i+ \sum_{k=1}^\infty x^{-1/p^k}\Big).$$

We will describe how to obtain this result. At this moment let us
see that reduction 4 fails in this situation.

Obviously, $Y=Y(x)$ is a root of $y^p-y-\xx$ if and only if it is
a root of $x(y^p-y)-1$, if and only if $1/Y$ is a root of
$x(1-y^{p-1})-y^p$. So $y^p-x(1-y^{p-1})$ has no roots in
$K^\fra((x))$. Reduction 4 fails because $(1-y^{p-1})^{1/p}$ is
not a power series.

Now let us see how to obtain the roots of $G(x,y)$. Assume that
$Y=Y(x)$ is a root of $G(x,y)$ and that $Y=ax^s+\text{higher order
terms}$. It is easy to see that $s=-1/p$ and $a=1$. Thus
substituting $Y=x^{-1/p}+Y_1$ for $y$ in $G(x,y)$, we see that
$Y_1$ is a root of $y^p-y-x^{-1/p}$. This above argument applies
repeatedly and we can assume that $Y=Y_0+\sum_{k\ge 1}
x^{-1/p^k}$, where the order of $Y_0$ is great than $-1/p^k$ for
any positive integer $k$. Now substituting this $Y$ for $y$ in
$G(x,y)$, we see that $Y_0$ is a root of $y^p-y$. Thus $Y_0$
equals $0,1,\dots, p-1$, and we have the desired factorization.

\renewcommand{\theequation}{\thesection.\arabic{equation}}
\vfill\eject \setcounter{chapter}{1} 
\chapter{The Field of Iterated Laurent Series\label{s-multi}}
After studying the field of double Laurent series, it is natural
to study the multivariate theory. The proofs of many combinatorial
identities involve more than two variables. The theory we are
going to develop in this chapter has three major applications. The
first application is to the evaluation of combinatorial sums
\citep{ego}; the second is to MacMahon's partition analysis, which
has been restudied by \citep{george6} in a series papers; and the
last is to lattice path enumeration, which will be carried out in
Chapter \ref{s-lattice}.

\section{The Fundamental Structure of $K\ll x_1,\dots,x_m\gg$}
In what follows, we denote $m$-vectors by bold face letters. Thus
$\mb{n}$ denotes the vector $(n_1,n_2,\ldots , n_m)$. Then
$\mb{x^n}:= x_1^{n_1}x_2^{n_2}\cdots x_m^{n_m}$ and
$\mb{n!}:=n_1!\, n_2!\cdots n_m!$. We also identify $\mb{n}_i$
with $n_i$.

By a formal Laurent series in $\mb{x}$, we mean a series that can
be written in the form
$$\sum_{n_1=-\infty}^{\infty}\cdots \sum_{n_m=-\infty}^{\infty} a_{n_1\ldots n_m}x_1^{n_1}
\cdots x_m^{n_m},$$ where $a_{n_1\ldots n_m}$ are elements in $K$.
Obviously, the set of all formal Laurent series in $\mb{x}$ does
not form a ring. However, some of its subsets do. In fact, one
well-known ring is the ring of Laurent series in $\mb{x}$, denoted
by $K((\mb{x}))=K((x_1,\ldots ,x_m))$, which is a subset of the
set of formal Laurent series. A formal Laurent series belongs to
$K((\mb{x}))$ if and only it has a lower bound for the power of
each $x_i$. Indeed, we have the following identification
$$K((x_1,\ldots ,x_n))=K[[x_1,\ldots ,x_n]][x_1^{-1},\ldots x_n^{-1}].$$

Suppose that $\mb{x}=(x_1,x_2,\ldots ,x_m )$ is an ordered set of
formal variables. We define $K_m\ll \mb{x}\gg =K\ll x_1,x_2,\ldots
,x_m\gg $ inductively by $K_m\ll \mb{x}\gg =K_{m-1}\ll \mb{x}\gg
((x_m))$, with $K_0\ll \mb{x}\gg =K$. So $K_1\ll \mb{x}\gg
=K((x_1))$ is the field of Laurent series in $x_1$, and $K_2\ll
\mb{x}\gg =K((x_1))((x_2))$ is the field of double Laurent series
in $x_1,x_2$, which has
 been studied in
chapter \ref{s-kxt}.

Clearly, $K_m\ll \mb{x}\gg $ is a field. We call $K_m\ll \mb{x}\gg
$
 the \emph{field of iterated Laurent series}. We shall see that
 many rings, such as the ring of polynomials $K[\mb{x}]$, the field of rational functions
$ K(\mb{x})$, the ring of Laurent series $K((\mb{x}))$, and the
ring $K[\lambda^{-1},\lambda][[\mb{x}]]$,
  can be embedded into the field of iterated Laurent series.
Thus the results on the field of iterated Laurent series apply to
many situations. Right now we are going to focus on the field of
iterated Laurent series to develop the general theory of this
field. This does not seem to have been done before.

Now let us look at some simple properties of iterated Laurent
series. An element $f(\mb{x})$ belongs to $K_m\ll \mb{x}\gg $ if
and only if it can be written in the form
$$f(\mb{x})=\sum_{n_m\ge N_{m}} a_{n_m} x_m^{n_m},$$
where $a_{n_m}\in K_{m-1}\ll \mb{x} \gg $. So $f(\mb{x})$ is
firstly regarded as a Laurent series in $x_m$, then a Laurent
series in $x_{m-1}$, and so on.

Similar to the two variable case, we have the composition law.

\begin{prop}
If $f(\mb{x})\in K_m\ll \mb{x}\gg $, and $g_i\in x_iK_{i-1}\ll
\mb{x}\gg [[x_i]]$ for
 $i=1,2,\ldots ,m$, then
$f(g_1,g_2,\ldots,g_m)\in K_m\ll \mb{x}\gg $.
\end{prop}

This law is in fact the application of the composition law of the
one variable Laurent series. It is not so useful since it does not
implies the composition law of the ring $K[[x_1,\dots ,x_m]]$. A
general composition law will be given in the next Chapter.

Clearly we can write $f$ as a formal Laurent series
$$ f(\mb{x})=\sum_{(i_1,\dots,i_m)\in \ZZ^m} a_{i_1,\dots,i_m} x_1^{i_1}\cdots x_m^{i_m}.$$
But it is not clear what the restrictions on these coefficients
is.  The structure of $K((x))$ is clear, and the structure of
$K((x))((t))$ is simple enough for our purpose. But for the three
variable case, it is not obvious whether the obvious definition of
the operator $\ct_{x_i}$ works or not. In fact, the obvious
definition works. To see this, we need to describe the structure
of $K\ll \mb{x}\gg$ more clearly.

Recall that a totally ordered set $S$ is {\em well-ordered} if
each nonempty subset of $S$ contains a minimal element. Applying
the basic theory of well-ordered sets, we get the fundamental
structure (Proposition \ref{p-wellorder}) for the field of
iterated Laurent series, which is going to play an important role
in our further development.

Let $M$ be the group of monomials in $x_1,\dots ,x_m$ with usual
multiplication, and let $\ZZ^m$ be the group written additively.
Clearly $M$ is isomorphic to $\ZZ^m$. The {\em reverse
lexicographic ordering} $<$ on $\ZZ^m$ is defined by
$(n_1,n_2,\ldots ,n_m) < (k_1,k_2,\ldots, k_m)$ if and only if
there is an $i$ such that $n_i<k_i$ and $n_j=k_j$ for all $j>i$.
This ordering is clearly a total ordering on $\ZZ^m$ that is
compatible with its group structure. Transferring this total
ordering to $M$, we get a total ordering ``$\preceq$", which plays
a central rule when expanding $1/f$ into an iterated Laurent
series. Thus if $1\le i<j\le m$, then for any positive integer
$s$, we have $x_i^s \prec x_{j}$, and the expansion of
$1/(x_i^s-x_j)$ is given by $x_i^{-s}\sum_{n\ge 0} \left(x_j/x_i^s
\right)^n$.
 The analogous situation for complex
variables would be informally written as $1>\!\!> x_1>\!\!> \cdots
>\!\!> x_n$ when expanding rational functions into Laurent series,
where $>\!\!> $ means ``much greater". See \citep{wilson} and
\citep[p. 231]{stanley-rec}.

The {\em order} of a monomial $c\mb{x^n}$ (where $c\in K$) is
defined to be $\mb{n}$. We say that the order of $\mb{x^n}$ is
smaller than the order of $\mb{x^k}$ if $\mb{n}<\mb{k}$ in the
reverse lexicographic ordering, or equivalently, $\mb{x^n}\prec
\mb{x^k}$.

Suppose that $b_\mb{n}\in K$ and that
$$f(\mb{x})=\sum_{\mb{n} \in \ZZ^m} b_\mb{n} \mb{x^n}$$
is a formal series. Then the {\em support} of $f$
 is defined to be the set $\{\, \mb{n}\mid b_{\mb{n}} \ne 0 \,\}$.

Now we can give the fundamental structure of the field $K\ll
\mb{x}\gg $, the proof of which will be provided later.
\begin{prop}[Fundamental structure]\label{p-wellorder}
A formal Laurent series in $\mb{x}$ belongs to $K\ll \mb{x} \gg$
if and only if it has a well-ordered support.
\end{prop}

This new result not only gives an overall view of iterated Laurent
series, but also validates the following natural definition.

\begin{dfn}\label{dfn-natural}
The operator $\ct_{x_j}$ acts on a formal series in $x_1,\dots
,x_m$ by
$$\ct_{x_j} \sum_{(i_1,\dots ,i_m)\in \ZZ^m}
a_{i_1,\dots ,i_m} x_1^{i_1} \cdots x_m^{i_m} = \sum_{(i_1,\dots
,i_m)\in \ZZ^m, i_j=0} a_{i_1,\dots ,i_m} x_1^{i_1} \cdots
x_m^{i_m},$$ where $a_{i_1,\dots ,i_m}$ belongs to $K$.
\end{dfn}

This natural definition has some obvious commutativity properties.
(See $P2$ and $P3$ below.) But the set of all formal series in
$x_1,\dots ,x_n$ does not form a ring, which means that we cannot
apply multiplication.

From the fundamental structure and the simple and useful fact that
{\em any subset of a well-ordered set is well-ordered}, it is easy
to see the following three properties hold.

\begin{enumerate}
\item[$P1$.] The operator $\ct_{x_i}$ results in an
iterated Laurent series when acting on an iterated Laurent series.

\item[$P2$.] The operator $\ct_{x_i}$ commutes with
$\sum$.

\item[$P3$.] The operator $\ct_{x_i}$ commutes with
$\ct_{x_j}$.
\end{enumerate}

Property $P1$ is necessary to make our definition applicable, and
it is nontrivial for $n\ge 3$ without the fundamental structure.
The commutativity property $P2$ is the key to converting many
problems into simple algebraic computations. The commutativity
property $P3$ may significantly simplify the constant term
evaluations.

Let us compare with another definition of $\ct_{x_i}$ by an
example.

\begin{exa}
\citet{zeil} proved a Conjecture of Chan et al. by showing an
identity that is equivalent to the following
\begin{equation}\label{e-zeil0}
\ct_{x_1}\cdots \ct_{x_n} \frac{1}{\prod_{i=1}^n (1-x_i)}
\frac{1}{\prod_{i<j} (x_i-x_j)} =C_1\cdots C_{n-1},
\end{equation}
where $C_n$ is the Catalan number. As pointed out in
\citep{welleda}, this identity should be interpreted as taking
iterated constant terms; i.e., in applying $\ct_{x_n}$ to the
displayed rational function, we expand it as a Laurent series in
$x_n$; the result is still a rational function and we can apply
$\ct_{x_{n-1}}$, \dots, $\ct_{x_1}$ iteratively. The $\ct_{x_i}$
does not commute with $\ct_{x_j}$.

Our approach is to expand rational functions in $K\ll x_1,\dots
,x_n \gg$ and then take the constant term in $x_1,\dots ,x_n$. So
after specifying the working field, the iterated constant term
operator is simply $\ct_{x_1,\dots ,x_n}$.
\end{exa}

\begin{proof}[Proof of Proposition \ref{p-wellorder}]
We proceed by induction on $m$. If $m=1$, then we are considering
$K((x_1))$. The proposition is clearly true. Suppose it is true
for $m-1$. Now we prove that it is also true for $m$.

Suppose that $b_\mb{n}\in K$, and that
$$f(\mb{x})=\sum_{\mb{n} \in \ZZ^m} b_\mb{n} \mb{x^n}.$$

On the one hand, if $f(\mb{x})\in K\ll \mb{x} \gg $, then we have
$$f(x)=\sum_{n_m\ge N_{m}} a_{n_m} x_m^{n_m}=\sum_{\mb{n} \in \ZZ^m} b_\mb{n} \mb{x^n},$$
with $a_{n_m}\in K_{m-1}\ll x\gg $ and $N_m$ an integer. Let $S$
be any subset of $P:=\{\, \mb{n}\mid b_{\mb{n}} \ne 0 \,\}$. Thus
for all $\mb{n}\in P$, $\mb{n}_m$ is greater than or equal to $
N_m$. Thus $\min_{\mb{n}\in S} \mb{n}_m$, denoted by $B_m$, exists
and is $\ge N_m$.
 Let $S'=\{\, \mb{n}\in S\mid  \mb{n}_m=B_m\,\}$. Since $S'$ is a subset of the set of  powers of
the nonzero terms in $a_{B_m}\in K_{m-1}\ll \mb{x}\gg $, it is a
well-ordered set. Using induction on $m$, we see that $S'$ has a
minimum, written as $(B_1,\ldots ,B_{m-1})$. Then $(B_1,\dots,
B_m)$ is the minimum of $S$.

On the other hand, if $P:=\{\, \mb{n}\mid b_{\mb{n}} \ne 0 \,\}$
is a well-ordered set, then it has a minimum, say $(N_1,\ldots
,N_m)$. Therefor, $N_m$ is the minimum of $\{\, \mb{n}_m\mid
\mb{n}\in P\,\}$, and
$$a_{n_m}x_m^{n_m}=\sum_{\mb{k}\in \ZZ^m, \mb{k}_m=\mb{n}_m} b_\mb{k} x_1^{k_1}\cdots
x_{m-1}^{k_{m-1}} x_m^{n_m}.
$$
Since the set of powers of the nonzero terms in $a_{n_m}x_m^{n_m}$
is a subset of $P$, it is well-ordered. By induction on $m$, we
get that $a_{n_m}\in K\ll x_1,\ldots ,x_{m-1}\gg $, and that $f$
can be written as
$$f(\mb{x})=\sum_{\mb{n}\in \ZZ^m, n_m\ge N_m} a_{n_m} x_m^{n_m} .$$
\end{proof}

Now for any $f(\mb{x})\in K\ll \mb{x}\gg $, we define $\ord(f)$ to
be the minimum of the support of $f$. If $\ord(f)=\mb{k}$, then we
call $\mb{b_kx^k} $ the initial term of $f$. It is clear that the
initial term of $fg$ equals the initial term of $f$ times the
initial term of $g$.

Similar to the case of double Laurent series, we have the
following three decompositions for iterated Laurent series. The
first decomposition follows from the fundamental structure
(Proposition \ref{p-wellorder}).

\begin{lem}[First Decomposition in $K\ll \mb{x}\gg $]
For each $1\le i\le m$, and $f\in K\ll \mb{x}\gg $, $f$ can be
uniquely written as $f=f_1+f_2,$ where $f_1$ contains only
nonnegative powers in $x_i$ and $f_2$ contains only negative
powers in $x_i$.
\end{lem}
Thus we can define $\pt_{x_i} f(\mb{x})=f_1$ and $\nt_{x_i}
f(\mb{x})=f_2$. If $f=\pt_{x_i}f$ (or $f=\nt_{x_i}f$), then we say
$f$ is $\pt$ (or $\nt$) in $x_i$, just the same as in the two
variable case. We shall mention that without using our fundamental
structure, it is not obvious that $\pt_{x_i} f(\mb{x}) \in K\ll
\mb{x} \gg$ for $m \ge 3$.

Now $\ct_{x_i} f(\mb{x})$ is in $K\ll x_1,\ldots ,x_{i-1}
,x_{i+1}, \ldots ,x_m\gg $, and $\ct_{x_{i_1},\ldots ,x_{i_r}}
f(\mbox{x})$ is independent of $x_{i_1},\ldots x_{i_r}$. The
residue is defined by
$$\res_{ x_{i_1},\ldots , x_{i_r}} f(\mb{x})
= \ct_{x_{i_1},\ldots ,x_{i_r}} x_{i_1}\cdots x_{i_r} f(\mb{x}).$$

In the following, we will see that the second decomposition is
useful in expanding $1/f$.

\begin{lem}[Second Decomposition in $K\ll \mb{x}\gg $]
\label{l-m-sec} If $f(x)\in K\ll \mb{x}\gg $, then $f$ can be
uniquely factored into the form
\begin{equation}
\label{e-mdec2} f(\mb x)=a\mb{x^k} f_1(x_1)f_2(x_1,x_2)\cdots
f_m(x_1,x_2,\ldots ,x_m),
\end{equation}
so that $a\in K$, and $f_i(x_1,\ldots ,x_i)\in K_{i-1}\ll
\mb{x}\gg [[x_i]]$ with constant term $1$ for all $i$. Moreover,
$a\mb{x^k}$ is the initial term of $f(\mb{x})$, and the second
decomposition of  $1/f(\mb{x})$ is given by
$$\frac{1}{f(\mb{x})}=\frac{1}{a}\mb{x}^{-\mb{k}} \frac{1}{f_1(x_1)}
\frac{1}{f_2(x_1,x_2)}\cdots \frac{1}{f_m(x_1,x_2,\ldots x_m)} .$$
\end{lem}
\begin{proof}
We prove this lemma by induction on $m$. It is trivial for $m=1$.
When $m=2$, we have shown it in chapter \ref{s-kxt}. Now suppose
it is true for $m-1$. We want to show it is true for $m$.

Since $K\ll \mb{x}\gg$ can be written as $K_{m-1}\ll \mb{x}\gg
((x_m))$, the Laurent series in $x_m$ with coefficients in
$K_{m-1}\ll \mb{x}\gg $, and $K_{m-1}\ll \mb{x}\gg $ is also a
field, we can write $f(\mb{x})=bx_m^{N_m} f_m(\mb{x})$, where
$b\in K_{m-1}\ll \mb{x}\gg $ and $f_m(x)\in K_{m-1}\ll \mb{x}\gg
[[x_m]]$ with constant term $1$. By induction, we can write
$$b=ax_1^{k_1}x_2^{k_2}\cdots x_{m-1}^{k_{m-1}} f_1f_2\cdots f_{m-1}$$
with $a\in K$ and $f_i\in K_{i-1}\ll \mb{x}\gg [[x_i]]$. So we
have the decomposition \eqref{e-mdec2}.

Now let $f(\mb{x})=a'\mb{x^{r}} g_1\cdots g_m$ be another
decomposition, and let $c=f(\mb{x})/ (x_{m}^{r_{m}} g_m)$. Then
$c\in K_{m-1}\ll \mb{x}\gg $, and $f=c x_m^{r_m}g_m$. By the
uniqueness of the decomposition of $K_{m-1}\ll \mb{x}\gg ((x_m))$,
we must have $b=c$ and $x_m^{k_m}f_m=x_m^{r_m}g_m$. Since $f_m$
and $g_m$ are both in $K_{m-1}\ll \mb{x} \gg [[x_m]]$ with
constant term (in $x_m$) $1$, we have $k_m=r_m$ and $f_m=g_m$. By
induction, $g_i=f_i$ and $r_i=k_i$ for all $i$ and $a=a'$. This
shows uniqueness.

It is clear that the initial term of $f_i$ is $1$ for all $i$, so
$a\mb{x^k}$ is the initial term of $f(\mb{x})$. The remaining
assertions are obvious.
\end{proof}

\begin{rem}\label{r-secdec}
With the above notation, we can see that $\log f_i(x_1,\ldots
,x_i)$ belongs to the ring $x_i K\ll x_1,\ldots ,x_{i-1}\gg
[[x_i]]$. Therefore $\log (f_1\cdots f_m)$ equals the sum of $\log
f_i$, and hence belongs to $K\ll \mb{x}\gg$. This fact will be
used later for proving a generalized residue theorem.
\end{rem}

\begin{lem}[Third Decomposition in $K\ll \mb{x}\gg $]
\label{l-mfactorization} If $f(\mb{x})\in K\ll \mb{x}\gg $ has
initial term  $1$, then for each $i$ with $1\le i\le m$, we have a
unique decomposition in $K\ll \mb{x}\gg $ $f=f_{i+}\,  f_{i0} \,
f_{i-} $, where $f_{i+}$ contains only positive powers in $x_i$,
$f_{i0}$ does not contain $x_i$, $f_{i-}$ contains only negative
powers in $x_i$,  and each of them has initial term $1$.
\end{lem}
\begin{proof}
Similar to the two variable case, this lemma follows from the
first decomposition through taking a logarithm. The difference is
that we need to show that $\log f$ belongs to $K\ll \mb{x}\gg $,
which follows from Remark \ref{r-secdec}.
\end{proof}

Applications of the third decompositions will be give in the
chapter about lattice path enumeration.

We conclude this section by giving some properties of the operator
``$\ord$". We have the following properties of the operation
``ord":
\begin{enumerate}
\item $\ord (fg)=\ord(f) +\ord(g)$.
\item $\ord(f+g)\ge \min(\ord(f),\ord(g))$, the greater only happens when the sum
of the initial terms of $f$ and $g$ equals $0$.
\item For any $ N\in \ZZ$ we have $\ord(f^N)=N\ord(f)$. In particular, $\ord(f^{-1})=-\ord(f)$.
\end{enumerate}
The first two properties are obvious. The third property is
trivial when $N$ is nonnegative. So it suffices to show that
$\ord(f^{-1})=-\ord(f)$, which follows from the second
decomposition (Lemma \ref{l-m-sec}).

\section{Basic Computational Rules}
Depending on the working field, rational functions
$Q(x_1,x_2,\ldots,x_m)$ may have as many as $m!$ different
expansions. More precisely, if $\sigma$ is a permutation of $[m]$,
then $Q(\mb{x})$ will have a unique expansion in $K\ll
x_{\sigma_1},x_{\sigma_2},\ldots ,x_{\sigma_m}\gg $. The
expansions of $Q(\mb{x})$ for different $\sigma$ are usually
different. So we need to specify the working field whenever a
reciprocal comes into account. Note that the intersection of all
these $m!$ sets $K\ll x_{\sigma_1},x_{\sigma_2},\ldots
,x_{\sigma_m}\gg $ is the ring of Laurent series $K((x_1,x_2,\dots
,x_m))$.

The computational rules in the working field $K\ll x_1,\dots
,x_m\gg$ are listed as follows, where $F$ and $G$ are in $K\ll
x_1,\dots ,x_m\gg$. These rules are similar to those for the field
of double Laurent series.
\begin{enumerate}
\item[Rule 1:] (linearity) For any $a,b$ that are independent of $x_i$,
$$\ct_{x_i} \left( a F(\mb{x})+bG(\mb{x})\right)= a \ct_{x_i} F(\mb{x}) + b \ct_{x_i} G(\mb{x}).$$

\item[Rule 2:] If $F$ can be written as $\sum_{k\ge 0}
a_k x_i^k$, then
$$\displaystyle \ct_{x_i} F = \left. F
\right|_{x_i=0}.$$

\item[Rule 3:]
$$\res_{x_i} \frac{\partial F(\mb{x})}{\partial x_i} G(\mb{x})= -\res_{x_i}
F(\mb{x}) \frac{\partial G(\mb{x})}{\partial x_i}.$$

\item[Rule 4:]
Suppose $F$ is $\pt$ in $x_i$. If $G$ can be factored in $K\ll
x_1,\dots ,x_m\gg$ as $(x_i-u)H$ such that $u$ is independent of
$x_i$ and $\ord (u) > \ord(x_i)$, and $1/H$ is $\pt$ in $x_i$,
then
$$\ct_{x_i} F(\mb{x}) \frac{x_i}{G(\mb{x})}=\left. \frac{F(\mb{x})}{\displaystyle
\frac{\partial G(\mb{x})}{\partial x_i}} \right|_{x_i=u}$$

\end{enumerate}

Rule 3 follows from the well-known property of the residues
$$\res_{x_i} \frac{\partial H(\mb{x})}{\partial x_i} =0.$$
Rule 4 is a reformulation of Theorem \ref{t-lagrange1} in the
multivariate case.

\section{Application to the Evaluation of Combinatorial Sums \label{s-comsum}}
One major application of our theory is on the evaluation of
combinatorial sums. To apply our theory, we first use the binomial
theorem and the formula for geometric series to convert the sums
into constant terms, and then Theorem \ref{t-lagrange1} (rule 4).

Let $\alpha$ be short for $\alpha_1,\dots,\alpha_r$. The working
field in this section is always $K\ll \alpha, \mb{x}\gg$.

\begin{exa}
Saalsch\"{u}tz's Theorem is equivalent to the following identity.
\begin{align*}
\sum_{k\ge 0} (-1)^k \binom{a+k-1}{k} \binom{a+e}{n-k}
\binom{d+e+k-1}{e} = \binom{d-a+n-1}{n} \binom{d+n-1}{e-n},
\end{align*}
where the sum ranges from $0$ to $n$.
\end{exa}
\begin{proof}
We prove this identity by showing that both sides have the same
generating function. The generating function for the left side can
be evaluated as follows:
\begin{align*}
& \ct_{\alpha_1,\alpha_2,\alpha_3} \sum_{a,d,e,n\ge 0}\sum_{k=0}^n
{\frac {\left (-1\right )^{e}\left (1+\alpha_{1}\right )^{-a}\left
( 1+\alpha_{2}\right )^{a+e}\left (1+\alpha_{3}\right
)^{-d-k}{x_{{1
}}}^{a}{x_{2}}^{d}{x_{3}}^{e}{x_{{4}}}^{n}}{{\alpha_{1}}^{k}{
\alpha_{2}}^{n-k}{\alpha_{3}}^{e}}}
 \\
&=  \ct_{\alpha_1,\alpha_2,\alpha_3} {\frac
{\alpha_{1}\alpha_{2}\alpha_{3} \left (1+\alpha_{1}\right )\left
(1+\alpha_3 \right )^{2}}{\left (\alpha_{2}-x_{{4}} \right )\left
(\alpha_{1}+\alpha_{1}\alpha_{3}-x_{{4}}\right )}}\cdot \\
& \qquad \qquad\qquad \quad \frac{1}{\left
(1+\alpha_{1}-(1+\alpha_{2})x_{1}\right )\left (
1+\alpha_{3}-x_{{2 }}\right )\left (\alpha_{{3
}}+(1+\alpha_{2})x_{3}\right )} .
\end{align*}
Now we take the constant term in $\alpha_1$ first. Only the second
factor in the denominator can result in negative powers in
$\alpha_1$, and it has a unique root $x_4/(1+\alpha_3)$, whose
order is higher than that of $\alpha_1$. Thus we can apply Theorem
\ref{t-lagrange1} and get
$$ \ct_{\alpha_2,\alpha_3}{\frac{\alpha_{2}\alpha_{3}\left (1+\alpha_{3}\right )\left (1+
\alpha_{3}+x_{{4}}\right )}{\left
(\alpha_{3}+x_{3}+\alpha_{2}x_{3} \right ) \left
(1+\alpha_{3}-x_{2}\right ) \left
(1+\alpha_{3}-(1+\alpha_2)(1+\alpha_3)x_1+x_{{4}} \right )\left
(\alpha_{2}-x_{{4}}\right  )}}.
$$
Similarly, we take the constant term in $\alpha_2$. Only the last
factor in the denominator can result in negative powers in
$\alpha_2$, and it has a unique root $x_4$. Applying Theorem
\ref{t-lagrange1} we get
$$\ct_{\alpha_3} {\frac {\alpha_{3}\left (1+\alpha_{3}\right )\left (1+\alpha_{{3}
}+x_{{4}}\right )}{ \left (\alpha_{3}+x_{3}+x_{3}x_{{4}}\right )
\left (1+\alpha_{3}-x_{2}\right )\left
(1+\alpha_{3}-x_{1}+x_{{4}}-\alpha_{3}x_{1}
-x_{1}x_{{4}}-\alpha_{3}x_{1}x_{ {4}} \right )}} .
$$
Now only the first factor in the denominator can result in
negative powers in $\alpha_3$, and it has a unique root
$x_3+x_3x_4$. Applying Theorem \ref{t-lagrange1} we get the final
generating function
$${\frac {\left (1-x_{3}\right )\left (1-x_{3}-x_{3}x_{{4}}
\right )}{\left (1-x_{1}-x_{{3
}}+x_{1}x_{3}+x_{1}x_{3}x_{{4}}\right ) \left
(1-x_{2}-x_{3}-x_{3}x_{{4}}\right )}} .$$

For the right side, we can evaluate the generating function as
follows.
\begin{align*}
&\ct_{\alpha_1,\alpha_2} {\frac {\left (-1\right )^{e}\left
(1+\alpha_{1}\right )^{-d+a} \left (1+\alpha_{2}\right
)^{-d-n}{x_{1}}^{a}{x_{2}}^{d}{x_{3}
}^{e}{x_{{4}}}^{n}}{{\alpha_{1}}^{n}{\alpha_{2}}^{e-n}}}\\
&\qquad = \ct_{\alpha_2}\ct_{\alpha_1} {\frac
{\alpha_{1}\alpha_{2}\left (1+\alpha_{ {1}}\right )\left
(1+\alpha_{2}\right )^{2}}{\left (\alpha_{1}+
\alpha_{1}\alpha_{2}-\alpha_{2}x_{{4}}\right )\left
(1-x_{1}-\alpha_{1}x_{1} \right )\left
(1+\alpha_{1}+\alpha_{2}+\alpha_{1}\alpha _{2}-x_{2}\right )\left
(\alpha_{2}+x_{3}\right )}} .
\end{align*}

Only the first factor in the denominator will result in negative
powers in $\alpha_1$, and it has a unique root
$\alpha_2x_4/(1+\alpha_2)$. Applying Theorem \ref{t-lagrange1} we
get
$$\ct_{\alpha_2}{\frac {\alpha_{2}\left (1+\alpha_{2}\right )\left (1+\alpha_{2}+\alpha_{{2
}}x_{{4}}\right )}{\left (\alpha_{2}+x_{3}\right ) \left
(1+\alpha_{2}-x_{2}+\alpha_{2}x_{{4}}\right )\left
(1+\alpha_{2}-x_{{1}
}-\alpha_{2}x_{1}-\alpha_{2}x_{1}x_{{4}}\right )}}.$$ Only the
first factor in the denominator will result in negative powers in
$\alpha_2$, and it has a unique root $-x_3$. Applying Theorem
\ref{t-lagrange1} we get the final generating function
$${\frac {\left (1-x_{3}\right )\left (1-x_{3}-x_{3}x_{{4}}
\right )}{\left (1-x_{1}-x_{{3
}}+x_{1}x_{3}+x_{1}x_{3}x_{{4}}\right ) \left
(1-x_{2}-x_{3}-x_{3}x_{{4}}\right )}} .$$ Saalsch\"{u}tz's Theorem
thus follows.
\end{proof}

\begin{exa}
Evaluate the generating function
$$\sum_{m,n=0}^\infty \sum_{a=0}^{m-2}\sum_{b=0}^{n-2}
\binom{n+a-1}{a}\binom{m+b-1}{b}\binom{m+n-a-b-4}{n-b-2} x^my^n.$$
This evaluation arose in counting directed convex polyominoes with
certain parameters.
\end{exa}

First we convert the sum into a constant term evaluation. We get
$$\ct_{\alpha_1,\alpha_2,\alpha_3} \sum_{m,n=0}^\infty \sum_{a=0}^{m-2}\sum_{b=0}^{n-2}
{\frac {\left (1+\alpha_{1}\right )^{n+a-1}\left (1+\alpha_{2}
\right )^{m+b-1}\left (1+\alpha_{3}\right
)^{m+n-a-b-4}{x}^{m}{y}^{n
}}{{\alpha_{1}}^{a}{\alpha_{2}}^{b}{\alpha_{3}}^{n-b-2}}}.
$$
The summation can be computed by first summing on $m\ge a+2$,
 $n\ge b+2$, and then summing on $a\ge 0$ and $b\ge 0$. We get

\begin{multline*}
\ct_{\alpha_1,\alpha_2,\alpha_3}{\frac
{\alpha_{1}\alpha_{2}\alpha_{3} \left (1+\alpha_{1}\right )\left
(1+\alpha_{2}\right ){x}^ {2}{y}^{2}}{ \left (1-x(1+\alpha_{2}
+\alpha_{3}+\alpha_{2}\alpha_{3})\right ) \left
(\alpha_{3}-y(1+\alpha_{{1
}}+\alpha_{3}+\alpha_{1}\alpha_{3})\right )}}\cdot \\
 \frac{1}{
\left (\alpha_{1} -x(1+\alpha_{1}+\alpha_{2}+\alpha_{1}\alpha_{2})
\right )\left
(\alpha_{2}-y(1+\alpha_{1}+\alpha_{2}+\alpha_{1}\alpha_{2} )\right
)}.
\end{multline*}
We first take the constant term in $\alpha_1$. Only the third
factor in the denominator will result in negative powers in
$\alpha_1$, which is linear and has a unique root $x\left
(1+\alpha_{2}\right)/\left(1-x(1+\alpha_{2})\right) $. Thus
applying Theorem \ref{t-lagrange1} we get
\begin{multline*}
\ct_{\alpha_2,\alpha_3}{\frac
{\alpha_{2}\alpha_{3}(1+\alpha_2){y}^{2}{x}^{2}}{ \left
(\alpha_{2}-y(1+\alpha_{2})-x(\alpha_{2}+{\alpha_{2}}^{2}) \right
)}}\cdot
\\
\frac{1}{\left (\alpha_{3}-y(1+\alpha_{3})-x(\alpha_{3}+\alpha_{2}
\alpha_{3})\right )\left (1-x(1+\alpha_{2}+\alpha_{3}+
\alpha_{2}\alpha_{3})\right )}.
\end{multline*}

Now taking the constant term in $\alpha_3$ is better than
$\alpha_2$. Only the second factor in the denominator will result
in negative powers in $\alpha_3$, which has a unique root
$y/(1-x-\alpha_{2}x-y) $. Applying Theorem \ref{t-lagrange1}, we
get
$$\ct_{\alpha_2} {\frac {\alpha_{2}\left (1+\alpha_{2}\right ){x}^{2}{y}^{2}}{
\left
(1-2\,x-y-2\alpha_{2}x+{x}^{2}+2\alpha_{2}{x}^{2}+{\alpha_{2}}^{2}{x}^{2}\right
)\left (\alpha_{2}-x(\alpha_{2}+{
\alpha_{2}}^{2})-y(1+\alpha_{2}\right )}} .$$ To take the constant
term in $\alpha_2$, we need to solve for $\alpha_2$ in the first
factor of the denominator. Only one root has order great than
$\alpha_2$, which may be found by the quadratic formula,
$$A={\frac {1-x-y-\sqrt {(1-x-y)^2-4xy}}{2x}}
.$$ Applying Theorem \ref{t-lagrange1}, and simplifying, we get
the generating function
$$\frac {xy \left(1-x-y-\sqrt {(1-x-y)^2-4xy}\right)}{2\left((1-x-y)^2-4xy\right)}
.$$ \qed

\begin{exa}
Super Catalan numbers $S(m,n)$ are defined
\begin{equation}
S(m,n)=\frac{(2m)!\,(2n)!}{m!\,n!\,(m+n)!}\label{e-scatalan1}.
\end{equation}
They were first stated to be integers by \citep{catalan}.
\end{exa}

We compute the generating function of $S(m,n)$ as follows. It is
easy to check that
$$S(m,n)= (-1)^n 4^{m+n} \binom{m-1/2}{m+n} .$$
Thus we have
\begin{align*}
\sum_{m,n\ge 0} S(m,n)x^my^n &=\sum_{m,n\ge 0} (-1)^n 4^{m+n} \binom{m-1/2}{m+n} x^my^n\\
&= \ct_{\alpha} \sum_{m,n\ge 0} (-1)^n 4^{m+n} (1+\alpha)^{m-1/2} \alpha^{-m-n} x^my^n \\
&=\ct_{\alpha} (1+\alpha)^{-1/2} \frac{1}{1-4(1+\alpha)x/\alpha}\frac{1}{1+4y/\alpha} \\
&= \ct_{\alpha} \frac{\alpha}{(1-4x)\sqrt{1+\alpha}}\cdot
\frac{\alpha}{\left(\alpha-4x/(1-4x)\right)(\alpha +4y)}
\end{align*}
Thus we have two roots $A_1=4x/(1-4x),A_2=-4y$ for $\alpha$ in the
denominator that will result in negative powers in $\alpha$. We
can use partial fractions or apply Theorem \ref{t-mdivid} to get
\begin{align*}
\frac{1}{1-4x}
\frac{\frac{A_1}{\sqrt{1+A_1}}-\frac{A_2}{\sqrt{1+A_2}}}{A_1-A_2}
=
\frac{1}{x+y-4xy}\left(\frac{x}{\sqrt{1-4x}}+\frac{y}{\sqrt{1-4y}}\right).
\end{align*}

\section{A New Algorithm for Partial Fraction
Decompositions\label{s-parfrac}}

The original purpose of this section is for the application of our
theory to MacMahon's partition analysis. But these results are of
independent interest.

The partial fraction decomposition (or expansion) of a one
variable rational function is very useful in mathematics. For
example, it is crucial to get the partial fraction decomposition
of a rational function when integrating it. Kovacic's algorithm
\citep{kovacic} for solving the differential equation
$y''(x)+r(x)y(x)=0$ requires the full partial fraction expansion
of $r(x)$ over the complex numbers.

The classical algorithm for partial fraction expansion relies on
the following theorem. To make it simple, we consider rational
functions in $\CC(t)$.

\begin{thm}
If $a_1,\dots ,a_n$ are $n$ distinct numbers, $m_1,\dots ,m_n$ are
positive integers, and the degree of $p(t)$ is less than
$m_1+\cdots +m_n$, then there are unique numbers $A_{i,j}$, where
$1\le i\le n$ and $1\le j\le m_i$, such that
\begin{equation}
\label{e-parfrac1} \frac{p(t)}{(t-a_1)^{m_1}\cdots (t-a_n)^{m_n}}
= \sum_{i=1}^n  \sum_{j=1}^{m_i} \frac{A_{i,j}}{(t-a_i)^j}.
\end{equation}
\end{thm}

The classical algorithm multiplies both sides by the denominator,
and then equates coefficients to solve a large system of linear
equations for the $A_{i,j}$'s.


It is the key observation of our new algorithm that linear
transformations will keep the structure of the partial fraction
expansion. We illustrate this idea by an example and will give a
precise argument later. See Lemma \ref{l-ppfraction-tr}.

Example: The partial fraction expansion of $f(t)$ is
$A/(t-a)+B/(t-b)$ if and only if the partial fraction expansion of
$f(t+c)$ is $A/(t+c-a)+B/(t+c-b)$. So if $a_i$ is not $0$ for all
$i$, then we can compute the partial fraction expansion of
$f(t+a_1)$, after that, replacing $t$ with $t-a_1$.

For example, let $f(t)=(t-a)^{-10}(t-b)^{-20}$. Maple will get
stuck when converting $f(t)$ into partial fractions, in which it
needs to solve a system of linear equations of $30$ unknowns. But
Maple can convert $f(t+b)$ into partial fractions quickly, and the
replacing of $t$ by $t-b$ costs little time. This is because after
that transformation most of coefficients in those $30$ linear
equations become $0$.

\subsection{The Theorems and the Algorithm}
In this section we develop a completely new algorithm for
computing partial fraction decompositions of rational functions.
This new algorithm not only has theoretical applications, but also
is very fast. When the base field is algebraically closed, our
algorithm is surprisingly simple. When the base field is not
algebraically closed, we also have a fast algorithm, and we will
explain how to compute the full partial fraction decompositions of
rational functions.

Denote by $F(t)$ the left hand side of equation
\eqref{e-parfrac1}. Let $M$ be the degree of the denominator of
$F(t)$, which is $m_1+m_2+\cdots +m_k$. Compared with the
classical algorithm for obtaining the partial fraction
decomposition of $F(t)$, our new algorithm has three improvements.
This comparison is under the assumption of fast multiplications of
(usually rational) numbers. In the following, when we say that an
algorithm takes $O(m)$ time, we mean that the algorithm will do
  $O(m)$ multiplications.

\begin{enumerate}
\item The new algorithm is fast. The classical algorithm needs to
solve $M$ linear equations of $M$ unknowns, which takes $O(M^3)$
time by using the Gaussian elimination algorithm. See
\citep[Property 37.1]{sedgewick}. But our algorithm only takes
about $O(M^2)$ time.

\item The new algorithm needs little storage space.
The classical algorithm needs to record all of the $M^2$
coefficients in those $M$ linear equations. So the storage space
is about $O(M^2)$. But our new algorithm needs only to record two
polynomials of degree $m$, where $m$ is the maximum of the
$m_i$'s. So the storage space is only $O(m)$.

\item The new algorithm computes the partial fraction
expansion at different $a_i$'s separately, so it is more suitable
for parallel programming.
\end{enumerate}

\vspace{3mm} Let $K$ be any field, and $t$ be a variable. It is
well-known that $K[t]$ has many nice properties. Here we use the
fact that $K[t]$ is a unique factorization domain.

In what follows, the degree of an element $r\in K[t]$, denoted by
$\deg(r)$, is the degree of $r$ as a polynomial in $t$. The degree
of the $0$ polynomial is treated as $-\infty$. We start with the
division theorem in $K[t]$.

\begin{thm}\label{t-division}
Let $D ,N \in K[t]$ and suppose $D\ne 0$. There is a unique pair
$(p,r)$ such that $p,r\in K[t]$, $N=Dp+r$, and $\deg r< \deg D$.
\end{thm}

\begin{rem}
In the above theorem, $r$ is called the {\rm remainder}. The
well-known division algorithm computes both $p$ and $r$ for given
$N$ and $D$. It is easy to see that this will take $O(\deg(p)
\deg(D))$ time. If we only care about $p$, or only care about $r$,
there exist faster algorithms, especially in some special cases.
We will discuss this later.
\end{rem}

A rational function $N/D$ with $N,D\in K[t]$ is said to be {\em
proper}
 if $\deg N <\deg D$. A proper rational function is
simply called a {\em proper fraction}. The unit $1$ is not proper,
but $0$ is considered to be proper. It is clear that the sum of
proper fractions is  a proper fraction, and the product of proper
fractions is  a proper fraction. But the set of all proper
fractions does not form a ring, for $1$ does not belong to it.

By Theorem \ref{t-division}, any rational function $N/D$ can be
uniquely written as the sum of a polynomial and a proper fraction.
Such a decomposition is called a {\em ppfraction} (short for
polynomial and proper fraction) of $N/D$. If $N=Dp+r$ with
$\deg(r)<\deg(D)$, then $N/D=p+r/D$ is a ppfraction. We denote by
$\poly(N/D)$ the polynomial part of $N/D$, and by $\frr(N/D)$ the
fractional part of $N/D$.

Recall the following well-known result in algebra.
\begin{lem}\label{l-ppfraction-k}
Let $N,D\in K[t]$ with $D\ne 0$. If $D=D_1\cdots D_k$ is a
factorization of $D$ in $K[t]$, and all the $D_i$ are pairwise
relatively prime, then $N/D$ can be uniquely written as
\begin{equation}\label{e-ppfraction-k}
\frac{N}{D}= p+\frac{r_1}{D_1}+\cdots +\frac{r_k}{D_k},
\end{equation}
where $r_i$ is a polynomial of degree smaller than $\deg (D_i)$
for all $i$, and $p$ equals the polynomial part of $N/D$. We call
such decomposition the {\rm ppfraction expansion of $N/D$ with
respect to $(D_1,\dots ,D_k)$}.
\end{lem}

Suppose that $D=D_1D'$ and that $D_1$ and $D'$ are relatively
prime. Then we have a ppfraction of $N/D$ with respect to
$(D_1,D')$: $N/D=\poly(N/D)+r_1/D_1+r'/D'.$ In such a
decomposition, we  call $r_1/D_1$  the {\em fractional part of
$N/D$ with respect to $D_1$}, and denote it by $\frr(N/D,D_1)$. If
$D_1=(t-a)^m$ for some $a\in K$, then we simply denote it by
$\frr(D/D,t=a)$.

Clearly $\frr(N/D,1)$ is always $0$, and $\frr(N/D,D_1)$ is always
a proper fraction with denominator $D_1$. We have the following
simple property.
\begin{lem}\label{l-ppfraction-M}
Let $M,N,D, D_1 \in K[t]$ with $D\ne 0$ and $D_1$ dividing $D$. If
$D_1$ and $D/D_1$ are relatively prime, then $\frr(MN/D,D_1)=
\frr(M \cdot \frr(N/D,D_1)).$
\end{lem}
\begin{proof}
Let $D'=D/D_1$, and let $N/D=p+r_1/D_1+r'/D'$ be the ppfraction
expansion of $N/D$ with respect to $(D_1,D')$. Then
$\frr(N/D,D_1)=r_1/D_1$. Now
\begin{align*}
    \frac{MN}{D}&=Mp+\frac{Mr_1}{D_1}+\frac{Mr'}{D'} \\
    &= Mp+\poly\left(\frac{Mr_1}{D_1}\right)
+\poly\left(\frac{Mr'}{D'}\right) +
\frr\left(\frac{Mr_1}{D_1}\right)
+\frr\left(\frac{Mr'}{D'}\right),
\end{align*} in which the sum of the first
three terms is a polynomial, the fourth term is a proper fraction
with denominator $D_1$, and the fifth term is a proper fraction
with denominator $D'$. Hence $\frr(MN/D,D_1)=\frr(M r_1/D_1)$ as
desired.
\end{proof}

\begin{thm}\label{t-ppfraction-main1}
For any $N,D\in K[t]$ with $D\ne 0$, if $D_1,\dots D_k \in K[t]$
are pairwise relatively prime, and $D=D_1\cdots D_k$, then
$$\frac{N}{D}= \poly\left(\frac{N}{D}\right) +\frr \left(\frac{N}{D}, D_1 \right)+\cdots +
\frr \left(\frac{N}{D}, D_k \right)$$ is the ppfraction expansion
of $N/D$ with respect to $(D_1,\dots ,D_k)$.
\end{thm}
\begin{proof}
Suppose that
$$\frac{N}{D}= p+\frac{r_1}{D_1}+\cdots +\frac{r_k}{D_k}$$
is the ppfraction expansion of $N/D$ with respect to $(D_1,\dots
,D_k)$. Let $D'=D_2\cdots D_k$. Then $D_1$ and $D'$ are relatively
prime and $r_2/D_2+\cdots +r_k/D_k$ is a proper fraction with
denominator $D'$. Denote it by $r'/D'$. By the uniqueness of
ppfraction of $N/D$ with respect to $(D_1,D')$, we have $r_1/D_1=
\frr(N/D,D_1)$. Similarly $r_i/D_i= \frr(N/D,D_i)$ for all $i$.
\end{proof}

Thus to find the ppfraction expansion of $N/D$ with respect to
$(D_1,\dots ,D_k)$, it suffices to find $\poly(N/D)$, which can be
easily done by the division algorithm, and $\frr(N/D,D_i)$ for
every $i$. From this idea, we can give a fast algorithm for
computing the ppfraction expansion with respect to $D_1,\dots
,D_k$. For this problem, the classical way is to assume that $N/D$
is a proper fraction, assume also that $r_i(t)= a_{i,0}+a_{i,1} t+
\cdots +a_{i,d_i-1} t^{d_i-1}$ for every $i$, where
$d_i=\deg(D_i)$, then solve a system of linear equations in
$d_1+\cdots +d_k$ indeterminates by equating coefficients of the
equation $N=r_1D/D_1+\cdots +r_kD/D_k$.

\begin{thm}\label{t-ppfraction-main2}
Suppose that $N\in K[t]$ and that $D=D_1\cdots D_k$ is a
factorization of $D$ in $K[t]$ such that $D_1$ is relatively prime
to $D_i$ for $i=2,\dots ,k$. Suppose also that $1/(D_1D_i)=
s_i/D_1+r_i/D_i$, which is not required to be a ppfraction
expansion. Then $\frr(N/D,D_1) =\frr(Ns_2s_3\cdots s_k/D_1)$.
\end{thm}
\begin{proof}
We have
\begin{align*}
\frac{1}{D} &= \frac{1}{D_1D_2} \cdot \frac{1}{D_3\cdots D_k} =
\frac{s_2}{D_1D_3\cdots D_k} +\frac{r_2}{D_2D_3\cdots D_k}
\end{align*}
Applying a similar procedure successively to $D_1D_3, D_1D_4,
\dots , D_1D_k$ in the first term, we get
$$\frac{1}{D}=\frac{s_2s_3\cdots s_k}{D_1}+\frac{s_2\cdots s_{k-1}r_k}{D_k}
+\frac{s_2\cdots s_{k-2}r_{k-1}}{D_{k-1}D_k} +\cdots
+\frac{r_2}{D_2\cdots D_k}.$$ Now denoted by $T_i$ the $i$th term
on the right hand side of the above equation. Then $T_1$ is a
rational function with denominator $D_1$, and for $i\ge 2$, $T_i$
is a rational function with denominator $D_kD_{k-1}\cdots
D_{k-i+2}$, which divides $D'$, where $D'=D_2\cdots D_k$. Thus
$\frr(T_2)+\cdots +\frr(T_k)$ is a proper fraction with
denominator $D'$. Denote it by $r'/D'$. Now
$$1/D= \poly(T_1)+\cdots +\poly(T_k) +\frr(T_1)+\frr(T_2)+\cdots +\frr(T_k).$$
The sum of the polynomial part of all the $T_i$'s has to be equal
to the polynomial part of $1/D$, which is $0$. Thus we get $1/D=
\frr(T_1)+r'/D'$. This is a ppfraction expansion of $1/D$ with
respect to $(D_1,D')$. So
$$\frr(1/D,D_1)= \frr(T_1) =\frr(s_2s_3\cdots s_k/D_1).$$
Thus by Lemma \ref{l-ppfraction-M} $\frr(N/D,D_1)=
\frr(Ns_2s_3\cdots s_k/D_1)$.
\end{proof}

Given relatively prime polynomials $D_1$ and $D_2$, we can use the
classical method to find $s_2,r_2$ such that $1/(D_1D_2)=
s_2/D_1+r_2/D_2$ with $\deg(s_2)<\deg(D_1)$ and
$\deg(r_2)<\deg(D_2)$. Alternatively, we can write the equation in
the form $1=s_2 D_2+r_2 D_1$ and find $s_2$ and $r_2$ by the
Euclidean algorithm.

From Theorem \ref{t-ppfraction-main2}, after solving $k-1$ linear
equations, with the $i$th having $\deg(D_1)+\deg(D_i)$
indeterminates for $i=2,\dots ,k$,  we can compute
$\frr(N/D,D_1)$, which is equal to the fractional part of
$Ns_2\cdots s_k/D_1$. This  algorithm is much more efficient than
the classical method for large $k$.

\vspace{3mm} If $D=a p_1^{m_1} \cdots p_k^{m_k}$, where $a\in K$,
 is a factorization of $D$ into monic primes in $K[t]$,
then $p_1^{m_1}, \dots , p_k^{m_k}$ are pairwise relatively prime.
Let $D_i=p_i^{m_i}$, and let $r_i$ be a polynomial with
$\deg(r_i)<\deg(D_i)$. Then every $r_i/D_i$ can be uniquely
written in the form $\sum_{j=1}^{m_i} A_j/p_i^j$ with $\deg
(A_j)<\deg(p_i)$ for all $j$. The {\em partial fraction expansion}
of $N/D$ is the result of applying the above decomposition to the
ppfraction of $N/D$ with respect to $(D_1,\dots ,D_k)$. In this
case, we can use the following lemma to reduce the problem to
computing only the partial fraction expansion of $1/(p_ip_j)$ for
all $i\ne j$.

\begin{lem}\label{l-frac-power}
Let $p,q\in K[t]$ be relatively prime polynomials. If $r$ and $s$
are two polynomials such that $1/(pq)=r/p+s/q$, then for any
positive integers $m,n$,
\begin{equation}\label{e-ppfraction-pq}
\frac{1}{p^mq^n} =\frac{1}{p^m} \sum_{i=0}^{m-1} \binom{m+i}{i}
r^ns^{i}p^{i} +\frac{1}{q^n} \sum_{j=0}^{n-1} \binom{n+j}{j}
r^{j}s^m q^{j}.
\end{equation}
\end{lem}
\begin{proof}
Using the formula $1/(pq)=r/p+s/q$, we have
$$\frac1{p^mq^n} =\frac{1}{pq}\cdot \frac{1}{p^{m-1}q^{n-1}}
= \frac{r}{p^mq^{n-1}}+\frac{s}{p^{m-1}q^n}.$$ If we let
$A(m,n)=1/(p^mq^n)$, then the above equation is equivalent to
$$A(m,n)=rA(m,n-1)+sA(m-1,n).$$
Using this recursive relation, we can express $A(m,n)$ in terms of
$A(0,j)$ and $A(i,0)$, where $1\le j\le n$ and $1\le i\le m$.

Either using induction or a combinatorial argument, we can easily
get
$$A(m,n)= \sum_{i=0}^{m-1} \binom{m+i}{i} r^{n}s^{i}A(m-i,0)
+\sum_{j=0}^{n-1} \binom{n+j}{j} r^{j}s^m A(0,n-j).$$ Equation
\eqref{e-ppfraction-pq} is just a restatment of the above
equation.
\end{proof}

\vspace{3mm} Let $b\in K$ and $\tau_b$ by the transformation
defined by $\tau_b\  f(t)=f(t+b)$ for any $f(t)\in K[t]$ or
$f(t)\in K(t)$. Then $\tau_b$ is clearly an automorphism on $K[t]$
and on $K(t)$, and its inverse is $\tau_{-b}$. The following
properties can be easily checked for any $p,q\in K[t]$ and $b\in
K$.
\begin{enumerate}
\item  $p$ is prime in $K[t]$ if and only if
$\tau_b\  p$ is.
\item $\tau_b\  \gcd(p,q) =\gcd(\tau_b\  p,\tau_b\  q)$.
\item $\deg (\tau_b\ p)=\deg(p)$.
\item $p/q$ is a proper fraction if and only if $\tau_b\ p/q$ is.
\end{enumerate}
Thus for any $N,D\in K[t]$ with $D\ne 0$, $N/D= p+r_1/D_1+\cdots
+r_k/D_k$ is the ppfraction expansion of $N/D$ if and only if
$\tau_b\  N/D= (\tau_b\  p)+(\tau_b\  r_1/D_1)+\cdots +(\tau_b\
r_k/D_k)$ is a ppfraction expansion of $\tau_b\  N/D$. The partial
fraction expansion can be obtained by first computing the partial
fraction expansion of $\tau_b\  N/D$, then applying $\tau_{-b}$ to
the result. Choosing $b$ appropriately can simplify the
computation. The above argument gives us the following lemma.
\begin{lem}\label{l-ppfraction-tr}
For any $N,D,D_1\in K[t]$ with $D\ne 0$, $D/D_1\in K[t]$, and
\\
$\gcd (D_1,D/D_1)=1$, we have
$$\frr(N/D,D_1)=\tau_{-b}\  \frr(\tau_b\  N/D, \tau_b\  D_1).$$
\end{lem}

\vspace{3mm} Now consider the case when $K$ is algebraically
closed. This is the simplest case, since every prime in $K[t]$ is
linear and can be written as $t-a$ for some $a\in K$.


Let $\ceiling{t^m}$ be the map from $K[[t]]$ to $K[t]$ given by
setting  $t^n=0$ for all $n\ge m$. More precisely,
$$\ceiling{t^m} \sum_{n\ge 0} a_n t^n =\sum_{n=0}^{m-1} a_n t^n.$$
where $a_i\in K$ for all $i$. The following properties can be
easily checked for all $f,g\in K[[t]]$.
\begin{enumerate}
\item  $\ceiling{t^m} (f+g) =\ceiling{t^m} f +\ceiling{t^m} g $.
\item  $\ceiling{t^m} (fg) =\ceiling{t^m} (\ceiling{t^m} f \ceiling{t^m} g)$.
\item  If $0<k< m$ then $\ceiling{t^m} t^k f =t^k \ceiling{t^{m-k}} f$.
\item  If $g(0)\ne 0$, then $\ceiling{t^m} f/g =\ceiling{t^m} (\ceiling{t^m} f/ \ceiling{t^m} g)$.
\end{enumerate}

\begin{prop}\label{p-0-multiply}
Let $P(t)$ and $Q(t)$ be two polynomials of degree $m-1$, then
$\ceiling{t^m}  P(t)Q(t)$ can be computed in $O(m^{1.58})$ time.
\end{prop}

It is known that the multiplication of two polynomial of degree
$m-1$ can be done in $O(m^{1.58})$ time. See, e.g.,
\citep[Property 36.1]{sedgewick}.

Fact: Let $P(t)$ and $Q(t)$ be two polynomials in $t$ of degree
$m-1$. To obtain $\ceiling{t^m} P(t)Q(t) $ needs only half of the
time to obtain $P(t)Q(t)$. This fact follows from the following
observation.

Bisect $P(t)$ into $P(t)= P_1(t)+t^{d}P_2(t)$, where $d=\lfloor
m/2\rfloor$ and $P_1(t)$ and $P_2(t)$ are both polynomials of
degree no more than $d$. Bisect $Q(t)$ into $Q(t)=Q_1(t)+t^d
Q_2(t)$ similarly. Then we have
$$\ceiling{t^{m_0}} P(t)Q(t) =( \ceiling{t^m} P_1(t)Q_1(t))+t^d \left(\ceiling{t^{m-d}}
(P_1(t)Q_2(t)+P_2(t)Q_1(t))\right).$$ Now it is easy to analyze
this to see the fact is true.

\begin{thm}\label{t-ppfraction-main-0}
Let $N,D\in K[t]$ and $D=t^{m}E$ with $E\in K[t]$ and $E(0)\ne 0$.
Then
$$ t^{m} \frr(N/D,t^{m})=\ceiling{t^{m}} \frac{N(t)}{E(t)}.$$
\end{thm}

\begin{proof}
Since $E(0)\ne 0$, $t^{m}$ and $E$ are relatively prime. Let
\begin{align}\label{e-ppfraction-th0}
\frac{N(t)}{D(t)}= p(t) +\frac{r(t)}{t^m} +\frac{s(t)}{E(t)}
\end{align}
be the ppfraction of $N/D$ with respect to $(t^m,E)$. Thus $\deg
(r(t))<m$, and $r(t)=t^m \frr(N/D,t^m)$.

Because $K(t)$ can be embedded into the field of Laurent series
$K((t))$, equation \eqref{e-ppfraction-th0} is also true as an
identity of $K((t))$.
 On the right hand side of
equation \eqref{e-ppfraction-th0}, when expanded as Laurent series
in $K((t))$, the the second term contains only negative powers in
$t$, and the other term contains only nonnegative powers in $t$.
Therefore, $r(t)/t^m$ equals the negative part of $N/D$ when
expanded as a Laurent series. More precisely, for $i=1,\dots , m$,
we have
$$[t^{-i}] \frac{N(t)}{D(t)} = [t^{-i}] \frac{ r(t)}{ t^m}.$$
This is equivalent to $[t^{m-i}] N(t)/E(t)= [t^{m-i}] r(t)$ for
$i=1,\dots, m$. Now $r(t)$ is a polynomial of degree at most
$m-1$, and $N(t)/E(t) \in K[[t]]$, so
$$r(t)= \ceiling{t^m} \frac{N(t)}{E(t)}.$$
\end{proof}

\begin{rem}
The idea of using Laurent expansion to obtain part of the partial
fraction expansion of rational functions appeared in the proof of
\citep[Theorem 4.4]{gessel2}.
\end{rem}

Gessel observed that this same idea can also be used to compute
the polynomial part of a rational function. And it is fast when
the polynomial part has small degree.
\begin{prop}\label{p-ppfraction-ira}
If $R(t)$ is a rational function in $K(t)$, then the polynomial
part $P(t)$ could be computed by the following equation.
$$t^{-1}P(t^{-1})= \frr( t^{-1} R(t^{-1}),t=0).$$
\end{prop}
\begin{proof}
Let $R(t)=P(t)+N(t)/D(t)$ be the ppfraction of $R(t)$, and let
$p=\deg (P)$, $d=\deg(D)$, and $n=\deg(N)$. Then $n<d$. Now we
have
$$t^{-1}R(t^{-1})=t^{-1}P(t^{-1})+t^{-1} N(t^{-1})/D(t^{-1}) =t^{-1}P(t^{-1})+
t^{d-n-1} \tilde{N}(t)/(\tilde{D}(t)),$$ where
$\tilde{D}(t)=t^dD(t^{-1})$, and similarly for $\tilde{N}(t)$.

Apply ppfraction expansion to the second term. Since
$\tilde{D}(t)$ has nonzero constant term, it is relatively prime
to $t^{p+1}$. Now it is clear that $t^{-1}P(t^{-1})$ is the
fractional part of $t^{-1}R(t^{-1})$ with respect to $t^{p+1}$.
\end{proof}

\begin{exa}
It is easy to check that
 $$R(t)={\frac {{t}^{3}+2\,{t}^{2}-3\,t+4}{{t}^{2}-4\,t+2}}
=t+6+{\frac {-8+19\,t}{{t}^{2}-4\,t+2}}.$$ Now we compute it by
Proposition \ref{p-ppfraction-ira}.
\end{exa}
\begin{align*}
t^{-1}R(t^{-1})& ={\frac {1+2\,t-3\,{t}^{2}+4\,{t}^{3}}{t^2\left
(1-4\,t+2\,{t}^{2}\right
)}} \\
t^2 \frr( t^{-1} R(t^{-1}),t^2) &= \ceiling{t^2} {\frac
{1+2\,t-3\,{t}^{2}+4\,{t}^{3}}{\left (1-4\,t+2\,{t}^{2}\right
)}}\\
&=\ceiling{t^2} \frac{1+2t}{1-4t} =1+6t.
\end{align*}
So we obtain that the polynomial part of $R(t)$ is $t+6$.

Note that when expanded as Laurent series in $t$, we have
$$\ceiling{t^{m_0}} \frac1{(t-a_i)^{m_i}}=\sum_{j=0}^{m_0-1} (-1)^{m_i}
\binom{m_i-1+j}{j} \frac{t^j}{a_i^{m_i+j}}.$$ Hence by Theorem
\ref{t-ppfraction-main2}, we get
\begin{cor}\label{c-ppfraction-0}
Let $N\in K[t]$ and $D=t^{m_0}(t-a_1)^{m_1}\cdots (t-a_k)^{m_k}$
with all the $a_i$'s being distinct and not equal to $0$. Then
$$t^{m_0} \frr\left(\frac{N}{D},t^{m_0}\right)=\ceiling{t^{m_0}}
Ns_1\cdots s_k,$$ where
$$s_i=\sum_{j=0}^{m_0-1} (-1)^{m_i}
\binom{m_i-1+j}{j} \frac{t^j}{a_i^{m_i+j}}.$$
\end{cor}
We following the notation in Corollary \ref{c-ppfraction-0}.
Because the ratios of the consecutive terms in the $s_i$ above are
simple rational numbers multiplied by $t$, the construction of
$s_i$ needs only $O(m_0)$ time. Thus from Proposition
\ref{p-0-multiply}, the computation of $\frr(N/D,t^{m_0})$ can be
done in $O(km_0^{1.58})$ time.

Therefore, combining Theorem \ref{t-ppfraction-main1}, Lemma
\ref{l-ppfraction-tr} and Corollary \ref{c-ppfraction-0}, we
obtain an algorithm for computing the partial fraction
decomposition of a proper rational function of the general form
$$F(t)=\frac{N(t)}{(t-a_1)^{m_1}\cdots (t-a_k)^{m_k}}.$$

\begin{enumerate}
\item Let $S:=0$
\item For $i$ from $1$ to $k$ do $G(t) :=F(t+a_i) ,$
 $S:=S+\frr (G(t),t^{m_i} )$ next $i$.
 \item Return $S$.
\end{enumerate}

The computation of $\frr(G(t),t^{m_i})$ will take $O(km_i^{1.58} )
$ time. Sum on all $i$ this gives us
$k/2(m_1^{1.58}+m_2^{1.58}+\cdots +m_k^{1.58})$. Now the only part
left is the computation of $F(t+a_i) $ for all $i$, which can be
easily checked to be no more than $O(M^2)$. So in any case, our
new algorithm will take no more than $O(M^2)$ time.

This new algorithm also enables us to work with some
 difficult rational functions by hand.

\begin{exa}
Compute the partial fraction expansion of $f(t)$, where
$$f(t)= \frac{t}{(t+1)^2(t-1)^3(t-2)^5}.$$
\end{exa}

Solution. Clearly, the polynomial part of $f(t)$ is $0$. Although
applying Corollary \ref{c-ppfraction-0} is faster, we compute the
fractional part of $f(t)$ at $t=-1$ and $t=1$ differently. For the
fractional part of $f(t)$ at $t=-1$, we apply $\tau_{-1}$, and
compute $\frr (f(t-1),t^2)$ by Theorem \ref{t-ppfraction-main-0}.
We have
\begin{align*}
t^2 \frr(f(t-1),t^2) &= \ceiling{t^2} \frac{t-1}{(t-2)^3(t-3)^5}\\
  &= \ceiling{t^2} \frac{ t-1}{ (-8+12t)((-3)^5+3^4\cdot 5t)} \\
  &= \ceiling{t^2} \frac{t-1}{8\cdot 3^5(1-19/6t)} \\
  &= \ceiling{t^2} \frac{(t-1)(1+19/6t)}{8\cdot 3^5} = -\frac{1}{8\cdot 3^5}(1+\frac{13t}{6}).
\end{align*}
Thus
$$\frr(f(t), (t+1)^2)= -\frac{1}{2^3\cdot 3^5 (t+1)^2}-\frac{13}{2^4 \cdot 3^6 (t+1)}.$$

Similarly, we can compute the fractional part of $f(t)$ at $t=1$.
We have
\begin{align*}
t^3 \frr(f(t+1),t^3) &= \ceiling{t^3} \frac{t+1}{(t+2)^2(t-1)^5}\\
  &= \ceiling{t^3} \frac{ t+1}{ (t^2+t+4)(-10t^2+5t-1)} \\
  &= \ceiling{t^3} \frac{t+1}{-4+16t-21t^2} \\
  &= -\frac{1}{4}\ceiling{t^3} (t+1)(1+4t-\frac{21}{4}t^2 +16t^2) \\
  &=-\frac{1}{4}(1+5t+\frac{59}{4} t^2).
\end{align*}
Thus
$$\frr(f(t), (t-1)^3)=  -\frac{1}{4(t-1)^3}-\frac{5}{4(t-1)^2}-
\frac{59}{16(t-1)}.$$ The fractional part of $f(t)$ at $t=2$ can
be obtained similarly, but it is better to use Corollary
\ref{c-ppfraction-0}. In fact, this computation becomes quite
complicated. Although it is still possible to work by hand, we did
use Maple.
\begin{align*}
&t^5\frr(f(t+2),t^5) \\
\qquad\qquad &=\ceiling{t^5} (t+2) \left(\frac{1}{9}-{\frac
{2t}{27}}+\frac{{t}^{2}}{27}-{\frac {4{t}^{3}}{243}}+{\frac
{5{t}^{4}}{729}} \right) \left(1-3t+6{t}^{2}-10{t}^{3}+15{t}^{4}
 \right) \\
\qquad\qquad &= \frac{2}{9}-{\frac {19}{27}}t+{\frac
{13}{9}}{t}^{2}-{\frac {593}{243}}{ t}^{3}+{\frac
{2689}{729}}{t}^{4}.
\end{align*}
Apply theorem \ref{t-ppfraction-main1}, we get the partial
fraction expansion of $f(t)$, which is too lengthy to be worth
giving here.

\vspace{3mm} Now we come back to the general case. In Maple, the
full partial fraction expansion of a rational function will
involve a form like
$$\sum_{\alpha= \text{\rm root of } p(t)}\sum_{j=1}^m \frac{h_j(\alpha)}{(t-\alpha)^j},$$
where $p(t)$ is a prime polynomial, and $h_j(t)$ will be a
polynomial of degree no more than $\deg(p(t))$. This expansion is
useful in some situations. We can also get this kind of expansion
by applying Theorem \ref{t-ppfraction-main-0}. This is best
illustrated by an example.

\begin{exa}
Compute the partial fraction expansion of $f(t)$, where
$$f(t)=\frac{t}{(t^2-t-1)^2(t^2-t+2)}.$$
\end{exa}

Solution. Suppose $\alpha$ is a root of the prime polynomial
$p(t):=t^2-t-1$. Since $K(\alpha)$ is a field, and
$\alpha^2=\alpha+1$, we can use this relation to get rid of all
terms containing $\alpha^n$ for $n\ge 2$. Because $p(t)$ is a
prime polynomial, $\alpha$ can only be a simple root of $p(t)$.
Then $t$ divides $p(t+\alpha)$ and $p(t+\alpha)/t$ has nonzero
constant term. In the present example,
$$p(t+\alpha)= (t+\alpha)^2-(t+\alpha)-1= t(t+2\alpha-1).$$
Note that the constant term of $p(t+\alpha)$ is always $0$.

Clearly, $\tau_\alpha\  (t^2-t+2)$ has constant term nonzero, for
otherwise it will not be relatively prime to $p(t)$. In the
present situation,
$$ (t+\alpha)^2-(t+\alpha)+2=t^2+(2\alpha -1)t+3.$$
By Lemma \ref{l-ppfraction-tr} and Theorem
\ref{t-ppfraction-main-0}, we can work in $K(\alpha)[[t]]$.
\begin{align*}
\ceiling{t^2} t^2 f(t+\alpha) &= \ceiling{t^2}
\frac{t+\alpha}{(t+2\alpha -1)^2(t^2+(2\alpha-1)t+3)}
\\
&= \frac{1}{15} \ceiling{t^2} \frac{t+\alpha}{1+(2\alpha-1)11t/15}\\
&= \frac{1}{15} \ceiling{t^2} (t+\alpha)(1-11(2\alpha-1)t/15) \\
&= \frac{1}{15}\alpha +\frac{(-11\alpha+7)}{15^2}t.
\end{align*}
Thus the fractional part of $f(t)$ at $\alpha$ that satisfies
$p(\alpha)=0$ can be written as
$$  \frac{\alpha}{15(t-\alpha)^2}
+\frac{(7-11\alpha)}{225(t-\alpha)}.$$ Similarly, the fractional
part of $f(t)$ at $\beta$ that satisfies $\beta^2-\beta+2=0$ can
be written as
$$ \left({\frac {4}{63}}-{\frac {1}{63}}\beta\right)\left (t-\beta
\right )^{-1}
 .$$
Together with the fact that the polynomial part of $f(t)$ is
clearly $0$, the full partial fraction expansion of $f(t) $ is
hence
$$f(t)= \sum_{\alpha^2-\alpha-1=0} \left[\frac{\alpha}{15(t-\alpha)^2}
+\frac{(7-11\alpha)}{225(t-\alpha)}\right]
+\sum_{\beta^2-\beta+2=0} \frac{4-\beta}{63(t-\beta) }.$$

\subsection{Applications to Generalized Dedekind Sums\label{ss-dedekind}}
\begin{prop}\label{p-dedekind1}
If the denominator of $R(t)$ is relatively prime to $p(t)$, and
$p(t)$ has only nonzero simple roots, then
\begin{equation}
\label{e-dedekind1} \sum_{p(\alpha)=0} \frac{R(\alpha)}{\alpha
p'(\alpha)} = -\frr(p^{-1}(t)R(t),p(t)) |_{t=0}.
\end{equation}
\end{prop}
\begin{proof}
Let $F(t)= R(t)/p(t)$, and let $\alpha$ be a root of $p(t)$. Then
$t$ divides $p(t+\alpha)$ and the constant term of $p(t+\alpha)/t$
is
$$\lim_{t\to 0} \frac{p(t+\alpha)}{t} =\lim_{t\to 0} p'(t+\alpha)= p'(\alpha),$$
where $p'(t)$ is the first derivative of $p(t)$ with respect to
$t$. Since $\alpha $ is a simple root, $p'(\alpha)\ne 0$.

By Theorem \ref{t-ppfraction-main-0}, $\frr(F(t+\alpha),t)$ is
then equal to
$$\ceiling{t^1} \frac{1}{p'(\alpha)+\text{higher terms}} R(t+\alpha) = \frac{R(\alpha)}{p'(\alpha)}.$$
Thus the ppfraction expansion of $F(t)$ at $p(t)$ can be written
as
$$\frr(F(t), p(t))=\sum_{p(\alpha)=0} \frac{R(\alpha)}{(t-\alpha)p'(\alpha)}.$$
If $p(0)\ne 0$, then by setting $t=0$, we get
$$\frr(F(t),p(t))|_{t=0}=-\sum_{p(\alpha)=0} \frac{R(\alpha)}{\alpha p'(\alpha)},$$
which is equivalent to \eqref{e-dedekind1}.
\end{proof}

\begin{cor}\label{c-dedekind}
If $R(\alpha)$ has no poles at $\alpha$ with $\alpha^n=1$, then
\begin{equation}
\label{e-dedekind2} \sum_{\alpha^n=1} R(\alpha)= -n
\frr(R(t)/(t^n-1),t^n-1) |_{t=0}.
\end{equation}
If $R(\alpha)$ has no poles at $\alpha$ with $\alpha^n=1$ except
$\alpha=1$, then
\begin{equation}
\label{e-dedekind3} \sum_{\alpha^n=1,\alpha\ne 1} R(\alpha)=
-n\frr(R(t)/(t^n-1), t^{n-1}+\cdots +t+1)|_{t=0}.
\end{equation}
\end{cor}
\begin{proof}
For the first part, let $p(t)=t^n-1$. Then $tp'(t)=nt^n$. When
$\alpha^n=1$, $\alpha p'(\alpha)=n$. Therefore equation
\eqref{e-dedekind2} follows from Proposition \ref{p-dedekind1}.

For the second part, let $p(t)=t^{n-1}+\cdots +t+1=(t^n-1)/(t-1).$
Then $tp'(t)= n t^n/(t-1)- t(t^n-1)/(t-1)^2.$ With the condition
that $\alpha^n=1$ and $\alpha\ne 1$, we have $\alpha
p'(\alpha)=n/(\alpha-1)$. Hence equation \eqref{e-dedekind3}
follows from Proposition \ref{p-dedekind1}.
\end{proof}

{\em Generalized Dedekind sums} are sums of the following form:
$$\sum_{\alpha^n=1,\alpha\ne 1} R(\alpha),$$
where $R(t)$ is a rational function. Sometimes $\alpha$ is allowed
to be $1$. This kind of sums has been studied by many authors.

One important class of generalized Dedekind sums is the class of
higher dimensional Dedekind sums, which are defined by
\begin{equation}
\label{e-dedekind-def} d(n;a_1,\dots ,a_m) = \sum_{\alpha^n=1,
\alpha\ne 1} \prod_{i=1}^m \frac{\alpha^{a_i}+1}{\alpha^{a_i}-1},
\end{equation}
where $n$ and $a_i$'s are positive integers, and $n$ is relatively
prime to $a_i$ for all $i$. For other equivalent definitions, see
\citep{zagier}.

Don Zagier gave a nice reciprocity law for higher dimensional
Dedekind sums in \citep{zagier}. The proof used a kind of residue
theorem.
\begin{thm}
If $a_0,\dots ,a_m$ are pairwise coprime positive integers, then
\begin{equation}
\sum_{j=0}^m \frac{1}{a_j} d(a_j;a_0,\dots ,\hat{a}_j,\dots, a_m)
=\phi_n(a_0,\dots ,a_n), \label{e-dedekind-rec}
\end{equation}
where the hat over $a_j$ denotes its omission from the list, and
$\phi_n$ is a certain rational function in $n+1$ variables.
\end{thm}

This theorem seems more naturally to be discovered by using
partial fraction expansion. Let
$$F(t)=\prod_{i=0}^{m} \frac{t^{a_i}+1}{t^{a_i}-1}.$$
Then by Corollary \ref{c-dedekind}, it is easy to see that
$$d(a_0;a_1,\dots ,a_m) =-\frac{a_0}2 \frr(F(t),t^{a_0}+\cdots +t+1) |_{t=0}.$$
Note that $F(t)$ is symmetric in $a_0,\dots ,a_n$. Thus we have
$$\sum_{j=0}^m \frac{1}{a_j} d(a_j;a_0,\dots ,\hat{a}_j,\dots, a_m)
=-\frac{1}{2} \sum_{i=0}^m \frr(F(t),t^{a_i}+\cdots +t+1)
|_{t=1}.$$ Now $F(t)$ has a ppfraction expansion of the form
$$F(t)= \poly(F(t))+ \frr(F(t),(t-1)^{m+1})+ \sum_{i=0}^m \frr(F(t),t^{a_i}+\cdots +t+1).$$
It is easy to see that $\poly(F(t))=1$ and $F(0)=1$. Thus by
setting $t=0$, we obtain that
$$\sum_{j=0}^m \frac{1}{a_j} d(a_j;a_0,\dots ,\hat{a_j},\dots, a_m)
=\frac{1}{2} \frr(F(t),(t-1)^{m+1})|_{t=0}.$$

Note that Zagier used a residue theorem to express this in terms
of Bernoulli numbers.

\section{Applications to MacMahon's Partition Analysis\label{s-Mac}}

\subsection{Background}
\begin{dfn}
 An Elliott-rational function is a rational
function that can be written in such a way that its denominator
can be factored into products of one monomial minus another, with
the $0$ monomial allowed.
\end{dfn}
In the one-variable case, this concept reduces to the generating
function of a quasi-polynomial. There is much interest in problems
of counting solutions to systems of linear Diophantine equations
and inequalities, and counting lattice points in convex polytopes.
Such counting problems can be converted into evaluating the
constant term of certain Elliott-rational functions. This
conversion has been known as MacMahon's partition analysis, and
has been given a new life by \citeauthor{george6} in a series of
papers
\citep{george6,MR2003d:05007,MR2003e:11110,MR2003d:05017,MR2002h:11100,MR2003d:05008,MR2002g:05014,MR2001j:05009,MR99j:05012}.

MacMahon's idea was to introduce new variables
$\lambda_1,\lambda_2,\dots $ to replace linear constraints. For
example, suppose we want to count the nonnegative integral
solutions to the linear equation $2a_1-3a_2+a_3+2=0$. We can
compute the generating function of such solutions as the
following:

$$\sum_{a_1,a_2,a_3\ge 0 \atop 2a_1-3a_2+a_3+2=0} x_1^{a_1}x_2^{a_2}x_3^{a_3}
=\sum_{a_1,a_2,a_3 \ge 0} \ct_\lambda \lambda^{2a_1-3a_2+a_3+2}
x_1^{a_1}x_2^{a_2}x_3^{a_3}.$$ Now apply the formula for the sum
of a geometric series. It becomes
$$\ct_\lambda \frac{\lambda^2}{(1-\lambda^2 x_1)(1-\lambda^{-3}x_2)(1-\lambda x_3)}.$$
The above expression is a power series in $x_i$ but not in
$\lambda$.

It is clear that if there are $r$ linear equations, we can resolve
them by introducing $r$ variables $\lambda_1,\dots ,\lambda_r$.
Thus counting solutions of a system of linear Diophantine
equations can be converted into evaluating the constant term of an
Elliott-rational function.

So the central problem in this section is to evaluate the constant
terms of Elliott-rational functions. One important result to this
problem is the following.

\begin{thm}\label{t-elliott}
If $F$ is Elliott-rational, then the constant terms of $F$ are
still Elliott-rational.
\end{thm}

This result follows from  ``The method of Elliott" (see \cite[p.
111--114]{mac}) developed from the following identity. Note that
we have not specified the working field yet.

\begin{lem}[Elliott Reduction Identity]
For positive integers $j$ and $k$,
$$\frac{1}{(1-x\lambda^j)(1-y\lambda^{-k})}=\frac{1}{1-xy\lambda^{j-k}}
\left(\frac{1}{1-x\lambda^j}+\frac{1}{1-y\lambda^{-k}}-1\right).$$
\end{lem}

Elliott's argument is that after finitely many applications of the
above identity to an Elliott-rational function, we will get a
summation of rational functions, in which the denominators
contains either all factors of the form $1-x\lambda^i$, or all
factors of the form $1-y/\lambda^i$. Now taking the constant term
of each summand is easy.

Theorem \ref{t-elliott} reduces the evaluation of $\ct_{\Lambda}
F$ to the univariate case $\ct_{\lambda}F$ by iteration.
Unfortunately, the Elliott reduction algorithm is not efficient in
practice. Other algorithms have been developed, and computer
programs have been set up, such as the ``Omega" package
\citep{george6}. But we can do much better by the partial fraction
method and working in a field of iterated Laurent series.

Before going further, let us review some of the work in
\citep{george6}.
 The key ingredient in
their argument is MacMahon's Omega operator $\Omega_\ge$.
\begin{dfn}
The operator $\Omega_\ge$ is defined by
$$\Oge \sum_{s_1=-\infty}^\infty \cdots \sum_{s_r=-\infty}^{\infty} A_{s_1,\dots ,s_r}
\lambda_1^{s_1} \cdots \lambda_r^{s_r} := \sum_{s_1=0}^\infty
\cdots \sum_{s_r=0}^{\infty} A_{s_1,\dots ,s_r},$$ where the
domain of the $A_{s_1,\dots ,s_r}$ is the field of rational
functions over $\CC$ in several complex variables and $\lambda_i$
are restricted to a neighborhood of the circle $|\lambda_i|=1.$ In
addition, the $A_{s_1,\dots ,s_r}$ are required to be such that
any of the $2^{r}-1$ sums
$$\sum_{s_{i_1}=0}^\infty \cdots \sum_{s_{i_j}=0}^{\infty} A_{s_{i_1},\dots
,s_{i_j}}$$
 is absolute convergent within the domain of the definition of $A_{s_1,\dots ,s_r}$.
\end{dfn}
Another operator $\Oeq$ is given by
$$\Oeq \sum_{s_1=-\infty}^\infty \cdots \sum_{s_r=-\infty}^{\infty} A_{s_1,\dots ,s_r}
\lambda_1^{s_1} \cdots \lambda_r^{s_r} := A_{0,\dots ,0}.$$

It was emphasized in \citep{george6} that it is essential to treat
everything analytically rather than formally because the method
relies on unique Laurent series representations of rational
functions.

It is not hard to see their definition always works if we are
working in a ring such as the ring of formal power series in
$\mb{x}$ with coefficients Laurent polynomials in $\Lambda$, where
$\mb{x}$ is short for $x_1,\dots, x_n$ and $\Lambda$ is short for
$\lambda_1,\dots ,\lambda_r$. In fact, this approach was used in
\citep{han}.

By Theorem \ref{t-elliott}, it suffices to consider the case of
$r=1$, since the general case can be done by iteration. In the
previous work by Andrews et al. or by Han, the problem was reduced
to evaluating the constant term (with respect to $\lambda$) of a
rational function of the form
\begin{align}\label{e-2-mac-g}
\frac{\lambda^k}{\prod_{1\le i\le m}
(1-\lambda^{j_i}x_i)\prod_{1\le i\le n} (1-y_i/\lambda^{k_i})}.
\end{align}
This treatment has assumed the obvious geometric expansion. In our
terms, $1$ is the initial term of each factor in the denominator.

\vspace{3mm} We find it better to do this kind of work in a
certain field of iterated Laurent series, because in such a field,
we can use the theory of partial fraction decompositions in
$K(\lambda)$ for any field $K$ and any variable $\lambda$.

We illustrate this idea by solving a problem in \citep[p.
2]{george6} with the partial fraction method.

\noindent {\bf Problem} Find all nonnegative integer solutions
$a,b$ to the inequality $2a\ge 3b$.

First of all, using geometric series summations  we translate the
problem into a form which MacMahon calls the {\em crude generating
function}, namely
$$f(x,y):=\sum_{a,b\ge 0, 2a-3b\ge 0} x^a y^b =\Oge \sum_{a,b\ge 0} \lambda^{2a-3b} x^ay^b
=\Oge \frac{1}{(1-\lambda^2 x)(1-\lambda^{-3}y)},$$ where
everything is regarded as a power series in $x$ and $y$ but not in
$\lambda$.

Now by converting into partial fractions in $\lambda$, we have
$$\frac{1}{(1-\lambda^2 x)(1-\lambda^{-3}y)}=
\frac{y(1+\lambda x^2y+ \lambda^2 x)}{(1-x^3y^2)(\lambda^3-y)}+
\frac{1+\lambda x^2y}{ (1-x^3y^2)(1-\lambda^2 x)}.$$ When the
right-hand side of the above equation is expanded as a power
series in $x$ and $y$, the second term contains only nonnegative
powers in $\lambda$, and the first term,
$$\frac{y(1+\lambda x^2y+ \lambda^2 x)}{(1-x^3y^2)(a^3-y)}=\frac{y}{1-x^3y^2}
\frac{\lambda^{-3}+\lambda^{-2}x^2y +\lambda^{-1} x}{
1-\lambda^{-3} y}$$ contains only negative powers in $\lambda$.
Thus by setting $\lambda=1$ in the second term, we obtain
$$f(x,y)=\frac{1+x^2y}{
(1-x^3y^2)(1- x)}.$$ By a geometric series expansion, it is easy
to deduce that
$$\{\, (a,b)\in \NN^2: 2a\ge 3b \,\} = \{\, (m+n+\lceil n/2 \rceil ,n): (m,n)\in \NN^2\,\}.$$

\subsection{Algorithm by Partial Fraction Decomposition}
Working in the field of iterated Laurent series has two
advantages. First, the expansion of a rational function into
Laurent series is determined by the total ordering ``$\preceq$ "
on monomials, so we can temporarily forget its expansion as long
as we work in this field. Second, the fact that $F$ is a rational
function in $\lambda$ with  coefficients in a certain \emph{field}
permits us to apply the theory of partial fraction decompositions.

Note that the idea of using partial fraction decompositions in
this context was first adopted in \citep[p.
229--231]{stanley-rec}, but without the use of computers,  this
idea was thought to be impractical.

MacMahon's partition analysis always works in a ring like
$K[\Lambda, \Lambda^{-1}][[\mb{x}]]$, where $\Lambda^{-1}$ is
short for $\lambda^{-1}_1,\dots ,\lambda^{-1}_r$. This ring can be
embedded into a field of iterated Laurent series, such as $K\ll
\Lambda,\mb{x} \gg$.

While working in the field of iterated Laurent series, MacMahon's
operators can be realized as the following.
\begin{align}
\Oge F(\Lambda, \mb{x}) &=\left. \PT_{\lambda} F(\Lambda, \mb{x})
\right|_{\Lambda= (1,\dots
,1)},\\
\Oeq F(\Lambda, \mb{x}) &=\ct_\Lambda F(\Lambda,\mb{x})=\left.
\PT_{\lambda} F(\Lambda, \mb{x}) \right|_{\Lambda= (0,\dots ,0)}.
\end{align} So it suffices to find $\pt_\Lambda F$.

In fact, it is well-known that $\Oge$ can be realized by $\Oeq$.
This is just like the fact that $\pt$ can be realized by $\ct$ as
we described in chapter 1. So either an algorithm for $\pt_\Lambda
F$ or an algorithm for $\ct_\Lambda F$ will be sufficient for our
purpose. Generally speaking, $\pt$ is more suitable for the
algorithm, and $\ct$ is more suitable for theoretical analysis.
This will be seen from our further discussion.

Now we need an algorithm to evaluate $\pt_\lambda F(\lambda)$ with
$$F(\lambda)= \frac{P(\lambda)}{\prod_{1\le i\le n} (\lambda^{j_i}-z_i)} $$
where $P(\lambda)$ is a polynomial in $\lambda$, $j_i$ are
nonnegative integers, and $z_i$ are independent of $\lambda$. Note
that we allow $z_i$ to be zero, so that the case of $P(\lambda)$
being Laurent polynomial is covered. Also note that our approach
is different from the previous algorithms, which deal with
rational functions expressed as in \eqref{e-2-mac-g}.

We have the following result.
\begin{thm}\label{t-ct-F}
Suppose that the factors in the denominator of $F$ are pairwise
relatively prime, and that the partial fraction decomposition of
$F$ is
$$F=f(\lambda)+\sum_{1\le i\le n} \frac{p_i(\lambda)}{\lambda^{j_i}-z_i},$$
where $f(\lambda)$ is a polynomial in $\lambda$, and
$p_i(\lambda)$ is a polynomial of degree less than $j_i$ for each
$i$. Then
$$\pt_\lambda F=f(\lambda) +\sum_{i} \frac{p_i(\lambda)}{\lambda^{j_i}-z_i},$$
where the sum ranges over all $i$ such that $z_i \prec
\lambda^{j_i}$.
\end{thm}
\begin{proof}
The condition that $z_i$ is independent of $\lambda$ implies that
either $ \lambda^{j_i} \prec z_i$ or $z_i \prec \lambda^{j_i}$. In
the former case, we observe that the expansion of
$p_i(\lambda)/(\lambda^{j_i}-z_i)$ into Laurent series contains
only negative powers in $\lambda$, hence has no contribution when
applying $\pt_\lambda$. In the latter case, the expansion contains
only nonnegative powers in $\lambda$. Thus the the theorem
follows.
\end{proof}

Now we need an efficient algorithm for the partial fraction
decompositions. The classical algorithm does not seem to work
efficiently. This is the motivation of our new algorithm for
partial fraction decomposition in last section.

By Theorem \ref{t-ppfraction-main2}, we need two formulas to
develop our algorithm. One is  for the fractional part of
$p(\lambda)/(\lambda^j-a)$, and the other for  the partial
fraction decomposition of $(\lambda^j-a)^{-1}(\lambda^k-b)^{-1}$.
These are given as Propositions \ref{p-2-frac} and \ref{p-2-pfrac}
respectively.

Let $\rmd(n,k)$ be the remainder of $n$ when divided by $k$. We
have

\begin{prop}\label{p-2-frac}
The fractional part of $p(\lambda)/(\lambda^j-a)$ can be obtained
by replacing $\lambda^d$ with $\lambda^{\rmd(d,j)} a^{\lfloor d/j
\rfloor}$ in $p(\lambda)$ for all $d$.
\end{prop}
\begin{proof}
By linearity, it suffice to show that the remainder of $\lambda^d$
when divided by $\lambda^j-a$
 equals
$\lambda^{\rmd(d,j)} a^{\lfloor d/j \rfloor}$, which is trivial.
\end{proof}
It is easy to see that this operation takes time linear in the
number of nonzero terms of $p(\lambda)$, where we assumed fast
arithmetic operations.
\begin{rem}\label{r-2-fracm}
Observe that the numerator of the fractional part of
$p(\lambda)/(\lambda^j-a)$ is always a Laurent polynomial in all
variables.
\end{rem}

\begin{lem}\label{l-2-pfrac}
For positive integers $j$ and $k$, if $a^k\ne b^j$, then the
following is a partial fraction expansion.
\begin{equation}
\frac{1}{(\lambda^j-a)(\lambda^k-b)} = \frac1{b^j-a^k}
\frr\left(\frac{\sum_{i=0}^{k-1} \lambda^{ij}a^{k-1-i}}
{\lambda^{k}-b}\right) -\frac1{b^j-a^k}
\frr\left(\frac{\sum_{i=0}^{j-1} \lambda^{ik}b^{j-1-i}}
{\lambda^{j}-a}\right)
\end{equation}
\end{lem}
\begin{proof}
First we show that if $a^k\ne b^j$, then $\lambda^j-a$ and
$\lambda^k-b$ are relatively prime. If not, say $\xi$ is their
common root in a field extension, then $\xi^j=a$ and $\xi^k=b$.
Thus we have $a^k=(\xi^j)^k=\xi^{jk}=(\xi^k)^j=b^j$, a
contradiction.

We have \begin{align*}
\frac{b^j-a^k}{(\lambda^j-a)(\lambda^k-b)}&=\frac{\lambda^{jk}-a^k}{(\lambda^j-a)(\lambda^k-b)}-
\frac{\lambda^{jk}-b^j}{(\lambda^j-a)(\lambda^k-b)}\\
&=\frac{\sum_{i=0}^{t-1} \lambda^{ij}a^{k-1-i}} {\lambda^{k}-b}
-\frac{\sum_{i=0}^{s-1} \lambda^{ik}b^{j-1-i}}{\lambda^j-a}.
\end{align*}
Now the polynomial part of
$\frac{b^j-a^k}{(\lambda^j-a)(\lambda^k-b)}$ is clearly $0$. Thus
the sum of the polynomial parts of the two terms on the right side
of the above equation also equals $0$. So taking the fractional
part of both sides and then dividing both sides by $b^j-a^k$ gives
the desired result.
\end{proof}

Now if $\gcd(j,k)$ is not $1$, then we can replace
$\lambda^{\gcd(j,k)}$ with $\mu$ and apply the above lemma. This
gives us the following result.

Let
$$\mathcal{F}(\lambda^j-a,\lambda^k-b)=\frac{\sum_{i=0}^{j'-1}
\lambda^{ik'}b^{j'-1-i}}{a^{k'}-b^{j'}},$$ where $j'=j/\gcd(j,k)$
and $k'=k/\gcd(j,k)$.

\begin{prop}\label{p-2-pfrac}
For positive integers $j$ and $k$, if $a^k\ne b^j$, then we have
\begin{align}\label{e-2-pfracm}
\frr\left( \frac{1}{(\lambda^j-a)(\lambda^k-b)}, \lambda^j-a
\right)&= \frr \left(\frac{\mathcal{F}(\lambda^j-a,\lambda^k-b)
}{\lambda^{j}-a} \right),
\end{align}
\end{prop}
\begin{rem}
Note that a similar result appeared in \citep{george6}, but their
proof was lengthy.
\end{rem}

Now by Theorem \ref{t-ppfraction-main2}, we have the following:
\begin{thm}\label{t-parfrac-F1}
With the notation of Theorem \ref{t-ct-F}, the polynomial
$p_s(\lambda)$ equals the remainder of
$$P(\lambda) \prod_{i=1,i\ne s}^n \mathcal{F}(\lambda^{j_s}-a_s, \lambda^{j_i}-a_i),$$ when divided by
$\lambda^{j_i}-z_i$ as a polynomial in $\lambda$.
\end{thm}

In Theorem \ref{t-ct-F}, we assumed that $\lambda^{j_i}-z_i$ and
$\lambda^{j_k}-z_k$ are relatively prime. Now let us consider the
case that $\lambda^{j_i}-z_i$ and $\lambda^{j_k}-z_k$ have a
nontrivial common factor. This happens if and only if
$z_i^{j_k}=z_k^{j_i}$, which can be easily checked. If they are
identical, then we can combine them together and apply Lemma
\ref{l-frac-power}. Otherwise, we can temporarily regard $z_i$ and
$z_j$ as two different variables. After the computation, we
replace them.

Thus the above argument, Theorem \ref{t-ct-F}, and
\ref{t-parfrac-F1} together will give us an efficient algorithm
for evaluating $\ct_\lambda F$.

\begin{rem}
From Remark \ref{r-2-fracm}, Theorem \ref{t-ct-F}, and Theorem
\ref{t-parfrac-F1}, we see that $\pt_{\lambda}F$ is
Elliott-rational when $F$ is. This is another way to prove Theorem
\ref{t-elliott}.
\end{rem}

\begin{exa}
Evaluate the constant term of $F(\Lambda)$, where
$$F(\Lambda)=\frac{1}{(1-\frac{\lambda_2x}{\lambda_1^2})(1-\frac{\lambda_3x}{\lambda_1^2})
(1-\frac{\lambda_1x}{\lambda_2^2})(1-\frac{\lambda_3x}{\lambda_2^2})
(1-\frac{\lambda_1x}{\lambda_3^2})(1-\frac{\lambda_2x}{\lambda_3^2})}.$$
\end{exa}
Although is in $K[\Lambda,\Lambda^{-1}][[x]]$, we shall work in
$K\ll \Lambda,x\gg$.

First, we take the constant term in $\lambda_1$. Applying Theorems
\ref{t-ct-F} and \ref{t-parfrac-F1} to the factors of $F(\Lambda)$
containing $\lambda_1$, we get
\begin{multline*}
  \ct_{\lambda_1}  \frac{1}{(1-\frac{\lambda_2x}{\lambda_1^2})(1-\frac{\lambda_3x}{\lambda_1^2})
(1-\frac{\lambda_1x}{\lambda_2^2}) (1-\frac{\lambda_1x}{\lambda_3^2})}\\
    =-{\frac {{\lambda_{{3}}}^{2}{\lambda_{{2}}}^{7}}{ \left( {\lambda_{{2}
}}^{3}-{x}^{3} \right)  \left(
{\lambda_{{2}}}^{4}-\lambda_{{3}}{x}^{3 } \right)  \left(
{\lambda_{{2}}}^{2}-{\lambda_{{3}}}^{2} \right) }}+{ \frac
{{\lambda_{{2}}}^{2}{\lambda_{{3}}}^{7}}{ \left( {\lambda_{{3}}
}^{4}-\lambda_{{2}}{x}^{3} \right)  \left(
{\lambda_{{3}}}^{3}-{x}^{3 } \right)  \left(
{\lambda_{{2}}}^{2}-{\lambda_{{3}}}^{2} \right) }}.
\end{multline*}
Denote by $F_1$ and $F_2$ the above two summands. At this stage,
we note that the expansion of $(\lambda_2^2-\lambda^3)^{-1}$ dones
not exist in $K[\Lambda,\Lambda^{-1}][[x]]$, and there is no
advantage in getting rid of the factor $\lambda_2^2-\lambda_3^3$
in the denominator by combining the above two summands into one
rational function.

Now we have \begin{align}\label{e-exa-mac11}
    \ct_\Lambda F(\Lambda)
    =\ct_{\lambda_2,\lambda_3}
\frac{\lambda_2^2\lambda_3^2F_1}{(\lambda_2^2-\lambda_3x)(\lambda_3^2-\lambda_2x)}
+\ct_{\lambda_2,\lambda_3}
\frac{\lambda_2^2\lambda_3^2F_2}{(\lambda_2^2-\lambda_3x)(\lambda_3^2-\lambda_2x)}.
\end{align}

We shall take the constant term in $\lambda_2$ first, since in the
expansion of
\begin{multline*}
\frac{\lambda_2^2\lambda_3^2F_1}{(\lambda_2^2-\lambda_3x)(\lambda_3^2-\lambda_2x)}\\
=-{\frac {{\lambda_{{3}}}^{4}{\lambda_{{2}}}^{9}}{ \left(
{\lambda_{{2} }}^{3}-{x}^{3} \right)  \left(
{\lambda_{{2}}}^{4}-\lambda_{{3}}{x}^{3 } \right)  \left(
{\lambda_{{2}}}^{2}-{\lambda_{{3}}}^{2}
\right)(\lambda_2^2-\lambda_3x)(\lambda_3^2-\lambda_2x) }},
\end{multline*}
only one factor, $\lambda_3^2-\lambda_2x$, in the denominator will
produce nonnegative powers in $\lambda_2$. Our procedure gives the
first term in \eqref{e-exa-mac11} as
$$\ct_{\lambda_3}-{\frac {{\lambda_{{3}}}^{16}{x}^{2}}{ \left( {x}^{2}-{\lambda_{{3}}}^
{2} \right)  \left( {x}^{6}-{\lambda_{{3}}}^{6} \right)  \left( -{
\lambda_{{3}}}^{7}+{x}^{7} \right)  \left(
-{\lambda_{{3}}}^{3}+{x}^{3 } \right) }}=0,
$$
and the second term in \eqref{e-exa-mac11} as
$$\ct_{\lambda_3}\left[{\frac {{\lambda_{{3}}}^{10}}{ \left(x^3 -{\lambda_{{3}}}^{3}
 \right) ^{2} \left( {x}^{2}-{\lambda_{{3}}}^{2} \right) ^{2}}}
-{\frac {{\lambda_{{3}}}^{16}{x}^{2}}{ \left( {x}
^{3}-{\lambda_{{3}}}^{3} \right)  \left(
{x}^{6}-{\lambda_{{3}}}^{3} \right)  \left( {x}^
{2}-{\lambda_{{3}}}^{2} \right)  \left({x}^{7}
-{\lambda_{{3}}}^{7}
 \right) }}\right]
 =1.
$$
Note that in evaluating the constant terms in the above two
Elliott-rational functions, we need only their polynomial parts.
Thus $\ct_{\Lambda} F=1.$

To see this in another way, we solve the corresponding linear
equations
\begin{align} \left[ \begin {array}{cccccc}
-2&-2&1&0&1&0\\\noalign{\medskip}1&0&-2
&-2&0&1\\\noalign{\medskip}0&1&0&1&-2&-2\end {array} \right] \cdot
 [ x_{{1}}\ \ x_{{2}}\ \ x_{{3}}\ \ x_{{4}}\ \ x_{{5}}
\ \ x_{{6}} ]^T =0.\label{ex-mac-lin}
\end{align}

The solution is
$$\left\{ {x_4}=\frac32{x_1}+7/2{ x_6}+3{ x_5},{
x_2}=-\frac32{ x_1}-\frac32{ x_6}-{x_5},{x_3}=-{x_1}-3{x_6}-3{
x_5}
 \right\}$$ with free parameters $x_1,x_5,x_6$. Now it is easy to
 see that $x_i=0$ for all $i$ is the only nonnegative integral solution of
\eqref{ex-mac-lin}.

\section{About the Residue Theorem}

As for the field of double Laurent series, we need a residue
theorem for the field of iterated Laurent series.

\citep{jac} proved the following theorem for $m= 3$:
\begin{thm}
Let $f_1(\mb{x}), \ldots , f_m(\mb{x})$ be Laurent series and let
$\mb{n}^{(i)}\in \ZZ^m$ be such that
$f_i(\mb{x})/\mb{x}^{\mb{n}^{(i)}}$ is a formal power series with
nonzero constant term. Then for any Laurent series $\Phi(\mb{y})$
such that $\Phi(\mb{f})$ belongs to $K(( x_1,\dots ,x_m))$,
\begin{equation}
\res_{\mb{x}} \left| \pad{x_{j}} {f_i}\right| \Phi(\mb{f}) =\left|
n_j^{(i)}\right| \res_{\mb{y}} \Phi(\mb{y}).
\end{equation}
\end{thm}
This is a theorem on the ring of multivariate Laurent series. The
diagonal (Good's) Lagrange inversion formula can be easily derived
from it. (There is a good summary for this in \citep{gessel-res}.)
We shall discuss this later.

The term {\em homogeneous Laurent series} was introduced in
\citep{reversion}. They used ``homogeneous expansion" to explain
the residue theorem in the ring of homogeneous Laurent series and
derived a simple formula for the non-diagonal Lagrange inversion
formula. A homogeneous Laurent series is better understood by
adding a redundant variable $t$. It is defined to be a Laurent
series in $t$, with coefficients in $K((x_1,\dots ,x_m))$, such
that in each nonzero term, the sum of the powers of the $x$'s
equals the power of $t$. The set of homogeneous Laurent series
form a ring, and we denote it by $K_h((x_1,\dots ,x_m))$. This $t$
plays an important role in expanding reciprocals. Note that in
\citep{reversion}, the redundant variable $t$ was replaced by $1$,
and the ring of homogeneous Laurent series was denoted by
$K(((x_1,\dots ,x_m)))$.

Because the residue theorems are developed over rings
($K((x_1,\dots ,x_m))$ \\
and $K_h((x_1,\dots ,x_m))$), they can be applied only if every
$f_i$ has a reciprocal in the corresponding ring. Now we are going
to give a residue theorem for the field of iterated Laurent
series, in which this restriction no longer exists since we are
working in a field.

Now let us see the residue theorem for the field $K\ll x_1,\ldots
,x_m\gg $ relative to $x_1,\dots ,x_m$.

\begin{prop}\label{p-residue}
Let $F_1,\dots ,F_m$ be iterated Laurent series. Suppose that the
initial term of $F_i$ is $f_i=a_ix_1^{n_{i,1}}\cdots
x_m^{n_{i,m}}$, where $a_i$ is independent of $x_1,\dots ,x_m$.
Then for any formal Laurent series $\Phi(y_1,\dots ,y_m)$ such
that $\Phi(f_1,\dots ,f_m)$ converges, we have
\begin{equation}\label{e-residue}
\res_{x_1,\dots ,x_m} \left| \frac{\partial F_i}{\partial x_j}
\right|_{1\le i,j \le m} \Phi(F_1,\dots ,F_m) =\left|
n_{i,j}\right|_{1\le i,j \le m} \res_{F_1,\dots ,F_m}
\Phi(F_1,\dots ,F_m).
\end{equation}
\end{prop}
\noindent On the right hand side of \eqref{e-residue}, every $F_i$
is temporarily regarded as a new variable.

The proof of this proposition will not be given here, because we
are going to give a more general result in the next chapter. At
this moment, we only give some remarks on this proposition. Note
that in Proposition \ref{p-residue}, $\Phi(x_1,\dots ,x_m)$ need
not belong to $K\ll x_1,\dots ,x_m \gg$.

There are several deficiencies of Proposition \ref{p-residue}.

First: With respect to $x_1,\dots ,x_m$, we can give a residue
theorem for the field of $K\ll x_1,\dots ,x_{m+n}\gg$. But to
state it clearly will be lengthy. This same situation persists if
we want to give a residue theorem for $K\ll x_{\sigma(1)}, \dots,
x_{\sigma(n+r)}\gg$, where $\sigma \in \sy_{n+r}$.

Second: The condition that $\Phi(F_1,\dots, F_m)$ belongs to $K\ll
x_1,\dots ,x_m \gg$ is not desirable. At least we should have a
simple criterion.

Now let us see the following phenomenon, which need an
explanation.

I will describe the basic idea of our residue theorem by a simple
example in the field of double Laurent series $K((x))((t))$.

\begin{exa}\label{ex-residue}
Let $F=x^2t$ and $G=xt^2$. Our residue theorem gives us the
identity
\begin{equation}\label{e-2-exa}
\ct_{x,t} \Phi (F,G) =\ct_{F,G} \Phi(F,G)
\end{equation}
for any rational $\Phi$, where on the left hand side, we are
taking the constant term inside $K((x))((t))$. We claim that on
the right hand side, the constant term cannot always be taken in
$K((F))((G))$ or in $K((G))((F)).$ This can be seen from the
following two examples.
\end{exa}

 First example: let $\Phi(F,G) :=\frac{F}{F+G}$. Direct computation in
$K((x))((t))$ shows that
$$\ct_{x,t} \Phi(F,G) =\ct_{x,t}
\frac{x^2t}{x^2t+xt^2}=\ct_{x,t} \frac{1}{1+t/x}=\ct_{x,t}
\sum_{n\ge 0} \left(-t/x\right)^n =1.$$ Equation \eqref{e-2-exa}
is true in $K((F))((G))$ but false in $K((G))((F))$. The correct
expansion is
$$ \ct_{F,G} \Phi(F,G)=\ct_{F,G} \frac{1}{1+G/F}=\ct_{F,G} \sum_{n\ge 0}
\left(-G/F\right)^n =1.$$

Second example: let $\Phi(F,G):=\frac{F^2}{F^2+G}$. Direct
computation shows that
$$\ct_{x,t}\Phi(F,G)=\ct_{x,t}\frac{x^4t^2}{x^4t^2+xt^2}=\ct_{x,t}\frac{x^3}{x^3+1}=0.$$
Equation \eqref{e-2-exa} is false in $K((F))((G))$ but true in
$K((G))((F))$. The correct expansion is
$$ \ct_{F,G} \Phi(F,G)=\ct_{F,G} \frac{F^2}{G}\frac{1}{F^2/G+1}=\ct_{F,G} \sum_{n\ge 0}
(-1)^{n} \left(F^2/G\right)^{n+1} =0.$$

In these two examples, only one expansion of $\Phi(F,G)$ into a
series in $F$ and $G$ works.  Writing such expansions as series in
$x$ and $t$ gives elements in $K((x))((t))$. In fact, the correct
expansions we used for these two examples have a consistency. In
the first example, $F=x^2t \prec xt^2=G$, and in the second
example $G=xt^2 \prec x^4t^2=F^2$. The conclusion is that the
expansion of $\Phi (F,G)$ on the right hand side of
\eqref{e-2-exa} is determined in $K((x))((t))$.

To improve the above situation and give a nice residue theorem is
the motivation of the next chapter.

\renewcommand{\theequation}{\thesection.\arabic{equation}}
\vfill\eject \setcounter{chapter}{2} 
\chapter{The Ring of Malcev-Neumann Series and the Residue Theorem \label{s-main}}

In the last chapter, we developed the theory of the field of
iterated Laurent series. It has many applications, as we have
already seen, but at the same time, there is something missing in
it. First, the field $K((\xx))((t))$ turned out to be useful, but
in the multivariate case, we have $2^{n}$ fields: $K\ll
x^{e_1}_1,\dots ,x_n^{e_n}\gg$ with $e_i$ being $\pm 1$. This
makes it hard to describe the general theory.
Second, the residue theorem 
needs to be further developed.

In searching for a satisfactory solution for the above two
problems, the ring of Malcev-Neumann series (or MN-series for
short) arises naturally. With this tool, fields like $K\ll
x_1,x_2^{-1},x_3^{-1}\gg$ can be easily described through an
endomorphism of $\ZZ^3$. As for the residue theorem, we will see
that it indeed involves two fields, which had been commonly
overlooked by combinatorists.

The algebra of MN-series was first constructed by \citep{malcev}
and \citep{neumann}. See \citep{passmann} for further references.
It is defined in the following fashion.

Let $R$ be a commutative ring with unit, and let $G$ be a group.
Define $R[G]$ to be the set of all elements of the form
$$\sum_{g\in G} a_g g,$$ where $a_g$ belongs to $R$ and $g$ is regarded as
a symbol, such that only finitely many $a_g$'s are nonzero. Then
under the natural addition and multiplication (by linearly
extending the multiplication of $G$), $R[G]$ is a $R$-algebra,
called the group algebra of $G$.

Now suppose $R$ is a field, denoted by $K$. The MN-series was
developed to answer a problem in algebra: Can we embed $K[G]$ into
a $K$-division algebra? The answer is yes for a special class of
groups as we shall explain later.

The structure of this chapter goes roughly like this: the
construction of the ring of MN-series is included in the first
section; in the second section, we will give the residue theorem
for MN-series; then, we will take a different point of view about
the Lagrange inversion formula; finally, we will discuss the
theoretical aspects of MacMahon's partition analysis.

\section{The Construction of the Ring of MN-series}

The construction of MN-series that we are going to give comes from
some similar ideas in \citep{passmann}. The new points are that we
construct the ring from a totally ordered monoid instead of a
totally ordered group, and that we use the finite decreasing chain
condition for well-ordered sets, which makes the proof clearer
than using the definition directly as in \citep{passmann}.

Recall that
a {\em partial ordering} $\le $ on a set $S$ is a relation on $S$ that is reflexive ($x\le x$ 
for all $x\in S$), antisymmetric ($x\le y$ and $y\le x$ implies
$x=y$) and transitive ($x\le y$ and $y\le z$ implies $x\le z$).

A {\em poset} $(S,\le)$ (short for partially ordered set) is a set $S$ together with a partial 
ordering $\le $ on $S$. If for all $x,y\in S$, either $x\le y$ or
$y\le x$ holds, then we call $(S,\le)$ a totally ordered set, and
$\le $ a total ordering of $S$.

Let $\le$ be a partial ordering on $S$. We say that $S$ is a
well-ordered set if every nonempty subset $A$ of $S$ has a
smallest element $a$. (Thus $a\le b$ for every $b\in A$.) In this
case, $\le$ is also called a well-ordering of $S$.

We are going to study the properties of well-ordered sets. The
following property is trivial but important.
\begin{itemize}
\item Any subset of a well-ordered set is well-ordered.
\end{itemize}

The following equivalent definition of a well-ordered set is
useful.
\begin{prop}
Let $\le $ be a total ordering on $S$. Then $S$ is well-ordered if
and only if $S$ does not contain an infinite decreasing sequence.
\end{prop}
\begin{proof}
If $S$ has an infinite decreasing sequence, say $a_1>a_2>\cdots $,
then $\{\, a_i \,\}_{i\ge 1}$ is a subset of $S$ without smallest
element. Hence $S$ is not well-ordered.

On the other hand, if $S$ is not well-ordered, then $S$ has a
nonempty subset $A$ which has no smallest element. Pick an element
from $A$, say $a_1$. Since $a_1$ is not the smallest element in
$A$, we can pick $a_2<a_1$ from $A$. This procedure will continue,
and we will get an infinite decreasing sequence $a_1>a_2>\cdots $.
\end{proof}

Examples.
\begin{enumerate}
\item Totally ordered finite sets are  well-ordered.

\item Under the natural order, $\NN$ is the simplest infinite well-ordered set.

\item Under the natural order, $\ZZ, \QQ,$ and $\RR$ are not well-ordered sets.
\end{enumerate}

Now let $\le $ be a total ordering on $S$, but not necessarily a
well-ordering. Consider the set $W_S$ of all well-ordered subsets
of $S$.

\begin{lem}\label{l-weltop}
If $w_\alpha \in W_S$ for all $\alpha$, then $\cap_\alpha
w_\alpha$ is also in $W_S$; if $w_1,w_2\in W_S$, then $w_1\cup
w_2$ belongs to $W_S$.
\end{lem}
\begin{proof}
The first part is obvious. We prove the second part by
contradiction. If $w_1\cup w_2$ is not well-ordered, then there is
an infinite decreasing chain $a_1>a_2>\cdots $ in $w_1\cup w_2$.

Picking out all elements in $w_1$, we get a sequence $a_{i_1}>
a_{i_2}>\cdots $ in $w_1$. Since $w_1 $ is well-ordered, this
decreasing sequence has to be finite. Similarly, picking out all
elements in $w_2$, we get a finite decreasing sequence $a_{j_1}>
a_{j_2}>\cdots $. But every element of the infinite set $\{ \,a_n|
n\ge 1 \,\}$ decreasing sequences as two sets is $\{\,
a_n\,\}_{n\ge 1}$ is in one of these two finite sequences, a
contradiction.
\end{proof}

From the above lemma, we see that $W_S$ is closed under infinite
intersection and finite union. Therefore, we have the following:

\begin{prop}
For any totally ordered set $S$, $W_S \cup \{\, S\,\}$ is the set
of all closed sets of a topology on $S$.
\end{prop}

We call this topology the {\em well-ordered topology}, denoted by
$T_w(S)$. Note that the closure of any well ordered subset is
itself and that the closure of any other subset is $S$.

Examples.
\begin{enumerate}
\item If $S$ itself is well-ordered, then all subsets of $S$ are closed,
and $T_w(S)$ is the discrete topology.

\item For the set of integers $\ZZ$ under the natural order, the closed sets in
$T_w(\ZZ)$ are all subsets of $\ZZ$ that have a least element.

\item For the set of rational numbers $\QQ$, the elements of
$W_\QQ$ do not have a simple description. Any subset of $\QQ$ that
has both a minimal element and a maximal denominator is
well-ordered. For if we let $A$ be a subset of $\QQ$, with minimal
element $a$ and maximal denominator $d$, then $d!\, A$ is a subset
of $\ZZ$ with minimal element $d!\, a$. So it is well-ordered. But
the converse is not true. For example, $\{\,1-2^{-n}\,\}_{n\ge 1}$
is an increasing sequence and hence well-ordered. But it has no
maximal denominator.
\end{enumerate}

The study of the well-ordered topology might be  interesting. It
would be good to give a simple description of the well-ordered
subsets of $\QQ$, or even $\RR$.

\vspace{3mm} A {\em monoid} is a semigroup with a unit. A {totally
ordered monoid} or TO-monoid is a monoid $G$ equipped with a total
ordering $\le$ that is compatible with the multiplication of $G$;
i.e., for all $x,y,z\in G$, $x<y$ implies that $zx<zy$ and that
$xz<yz$. An immediate consequence is that if $x<y$ and $x'<y'$,
then $xx'<yy'$. For $x<y$ implies $xx'<yx'$, $x'<y'$ implies
$yx'<yy'$, and  the transitivity of $<$ yields $xx'<yy'$.

If a TO-monoid $G$ is also abelian and written additively, then
$<$ is said to be {\em translation invariant}; i.e., $x<y$ implies
$x+z<y+z$. Similarly we can define a TO-group. The abelian groups
$\ZZ,\QQ,$ and $\RR$ are all totally ordered abelian groups.

Given two subsets $X$ and $Y$ of $G$, we define the product
$X\cdot Y$ to be the set $\{\, xy: x\in X,y\in Y\,\}$. We also
define $X^{\cdot n}$ to be the product $X \cdots  X$ ($n$ times).

Now we consider $W_G$, the set of all well-ordered subsets of $G$.
The following lemma will be useful.

\begin{lem}\label{l-infsubseq}
If $S$ is a totally ordered set, then any infinite sequence
$a_1,a_2,\dots$ in $S$ contains at least one of the following.
\begin{enumerate}
\item An infinite increasing subsequence.
\item An infinite constant subsequence.
\item An infinite decreasing subsequence.
\end{enumerate}
\end{lem}
\begin{proof}
Suppose $a_1,a_2,\dots$ has neither an infinite decreasing
subsequence nor an infinite constant subsequence. We want to show
that it has an infinite increasing subsequence.

Since it contains no infinite decreasing subsequence, it has a
smallest element, say $a_{i_1}$. For otherwise we can construct an
infinite decreasing subsequence. Deleting the first $i_1$ elements
from $\{\, a_{n} \,\}_{n\ge 1}$ leaves an infinite sequence. Since
there are only finitely many $a_n$'s that equal $a_{i_1}$,
deleting all of them still results in an infinite sequence. In
this new sequence, every element is greater than $a_{i_1}$, and
there is no infinite decreasing or constant subsequence. Thus we
can repeat the above procedure, and get an infinite increasing
subsequence $a_{i_1}<a_{i_2}<\cdots $.
\end{proof}

\begin{prop}
\label{p-welplus} If $G$ is a \tomonoid\ and $w_1,w_2\in W_G$,
then $w_1\cdot w_2\in W_G$.
\end{prop}
\begin{proof}
If not, we can assume that
$$a_1b_1>a_2b_2>\cdots $$
is an infinite decreasing sequence with $a_i\in w_1$ and $b_i\in
w_2$.

Since $w_1 $ is well-ordered, the infinite sequence $\{\, a_n
\,\}_{n\ge 1}$ has no infinite decreasing sequence. By Lemma
\ref{l-infsubseq}, it has an infinite weakly increasing
subsequence, say $a_{i_1}\le a_{i_2}\le \cdots$. Together with the
condition that $a_{i_1}b_{i_1}> a_{i_2}b_{i_2} >\cdots $, we get
an infinite decreasing sequence $b_{i_1}> b_{i_2}> \cdots$ in
$w_2$.  This contradicts the fact that $w_2$ is well-ordered.
\end{proof}

\vspace{3mm} Now  we can construct the ring of MN-series. Let $G$
be a totally ordered monoid, and let $R$ be a commutative ring
with unit.

A formal series $\eta$ on $G$ has the form
$$\eta =\sum_{g\in G} a_g g, $$
where $a_g\in R$ and $g$ is regarded as a symbol. The support of
$\eta$ is defined to be
$$\supp (\eta)= \{\, g\in G: a_g\ne 0\, \}.$$
A {\em Malcev-Neumann series} is a formal series on $G$ that has a
well-ordered support. We define $R_w[G]$ to be the set of all such
MN-series.

If $\eta\in R_w[G]$, then we can define the {\em order} of $\eta$
to be $\ord(\eta)=\min \left(\supp(\eta)\right)$. The {\em initial
term} of $\eta$ refers to the term with the smallest order. It is
clear that $\ord(\eta \tau)=\ord(\eta)\ord(\tau)$. We denote by
$[g] \eta$ the coefficient of $g$ in $\eta$.

\begin{thm}
If $G$ is a \tomonoid, then under the natural addition and
multiplication, $R_w[G]$ is a ring.
\end{thm}
\begin{proof}
By linearity, it suffices to show that $R_w[G]$ is closed under
addition and multiplication. Let $\eta,\tau \in R_w[G]$, and let
$A= \supp (\eta)$ and $B=\supp (\tau)$. Then $A$ and $B$ are
well-ordered. Since $[g] \eta+\tau = [g] \eta +[g] \tau$, $\supp
(\eta+\tau)$ is contained in $A\cup B$, which is well-ordered by
Lemma \ref{l-weltop}. So $\eta+\tau$ belongs to $R_w[G]$.

For the multiplication, we have
\begin{equation}
[g] \eta \tau= \sum_{ab=g} [a] \eta \cdot [b] \tau
,\label{e-welprod}
\end{equation}
where the sum can be restricted to $a\in A$ and $b\in B$, for
otherwise the summand is zero. So the support of $\eta \tau$ is
contained in $A\cdot B$, which is well-ordered by Proposition
\ref{p-welplus}.

Now we show that the summation on the right hand side of
\eqref{e-welprod} is a finite sum; i.e., for any $g\in G$ there
are only finitely many $(a,b)\in A\times B$ such that $ab=g$. If
not, suppose $g=a_nb_n$ for $n=1,2\dots $. Then $\{\, a_n\,
\}_{n\ge 1}$ is an infinite sequence of distinct elements of $A$,
which is well-ordered.
 By Lemma \ref{l-infsubseq}, it contains an infinite increasing subsequence, say
 $a_{i_1}<a_{i_2}<\cdots.$ But then $b_{i_1}>b_{i_2}>\cdots $
 forms an infinite decreasing sequence of $B$. This contradicts the fact that
 $B$ is well-ordered.
\end{proof}

The ring $R_w[G]$ has some nice properties. For example, it
contains $R[G]$ as a subring, because elements in $R[G]$ have
finite support.

If $G$ is also a group and $K$ is a field, then $K_w[G]$ is
maximal in the sense that if $\eta =\sum_{g\in G} a_g t^g$ is not
in $K_w[G]$, then adding $\eta $ into $K_w[G]$ cannot form a ring.
For if $\supp (\eta)$ is not well-ordered, we can assume that
$g_1>g_2>\cdots $ is an infinite decreasing sequence in
$\supp(\eta)$. Let $\tau= \sum_{n\ge 1} a_g^{-1} g^{-1}_n$. Note
that $\tau \in R_w[G]$, since $g^{-1}_1<g^{-1}_2<\cdots $ is well
ordered. But the constant term of $\eta \tau$ equals an infinite
sum of $1$'s, which diverges.

Let $\eta_1,\eta_2,\dots $ be a series of elements in $R_w[G]$.
Then we say that $\eta_1+\eta_2+\dots $ exists or {\em strictly
converges} to $\eta \in R_w[G]$, if for every $g\in G$, there are
only finitely many $i$ such that $[g] \eta_i\ne 0$, and
$\sum_{i\ge 1} [g] \eta_i = [g] \eta$. Note that $\sum_{n\ge 1}
2^{-n}$ does not strictly converge to $1$.

Let $f(t)=\sum_{n\ge 0} b_n t^n$ be a formal power series in $t$,
and let $\eta\in R_w[G]$. Then we define the composition  $f\circ
\eta$ to be
$$f\circ \eta := f(\eta)= \sum_{n\ge 0} b_n \eta^n $$
if it exists.

We have the following composition law for $R_w[G]$.
\begin{thm}
\label{t-welcomposition} If $f\in R[[t]]$ and $\eta \in R_w[G]$
with $\ord (\eta)>1$, then $f\circ \eta$ strictly converges in
$R_w[[G]]$.
\end{thm}

The proof of this theorem consists of two parts: one is to show
that the support of $f\circ \eta$ is well-ordered; the other is to
show that for any $g\in G$, $[g] f\circ \eta$ is a finite sum of
elements in $R$.

\begin{prop}\label{p-welpowerunion}
If $A\in W_G$, and $A>1$, i.e., for all $a\in A,$ $a>1$, then $\bcup_{n\ge 0} A^{\cdot n} \in 
W_G$.
\end{prop}

In order to prove this proposition, we introduce a new concept.
Let $S$ be a subset of $G$. If $S>1$, then we say that $S$ is {\em
Archimedean} if for all $x,y\in S$, there is a positive  integer
$n$ such that $x<_G y^n$.

\begin{lem}\label{l-welarki}
If $A\in W_G$, $A>1$ and $A$ is Archimedean, then $\bcup_{n\ge 0}
A^{\cdot n}\in W_G$.
\end{lem}
\begin{proof}
If not, we shall have an infinite decreasing sequence in
$\bcup_{n\ge 0} A^{\cdot n}$, say
$$a^{(1)} > a^{(2)} > \cdots $$
with $a^{(i)}\in A^{\cdot n_i}$ for some positive integer $n_i$
for all $i$. Since $A$ is well-ordered and $A>1$, we can assume
that $1< a\in A$ is the smallest element of $A$.

 Write $a^{(1)}=d_1d_2\cdots d_{n_1}$, where $d_i$ in $A$. Then by the assumption that
 $A$ is Archimedean, there are positive integers $m_j$ for $j=1,2,\cdots ,n_1$
 such that $d_j< a^{m_j}$. Let $m=m_1+m_2+\cdots +m_{n_1}$. Then
 $a^{(1)}<a^m $.

 Clearly $a^{(i)}\ge a^{n_i}$, so $n_i<m$ for all $i$,
 and the infinite decreasing sequence $a^{(1)} > a^{(2)} > \cdots $ in fact
 belongs to $\bcup_{n=0}^m A^{\cdot n}$, which is a finite union of well-ordered sets
 (by Proposition \ref{p-welplus}), and is hence well-ordered (by Lemma \ref{l-weltop}),
 a contradiction.
\end{proof}

For $x,y\ge 1$, we define $x<<y$ to mean that for all $n\in \NN$,
$x^n<_G y$. If there are positive integers $m$ and $n$ such that
$x<_G y^n$ and $y<_G x^m$, then we say that $x$ and $y$ are {\em
Archimedean equivalent}, denoted by $x\equiv y$. The following
properties are clear for any $x,y,z\ge 1$.

\begin{enumerate}
\item Exactly one of the three conditions holds: $x>>y$, or $x\equiv y$ or $x<<y$.
\item $\equiv $ is an equivalence relation.
\item If $x\equiv z$, and $y\equiv z$ or $y<<z$, then $xy\equiv yx \equiv z$.
\item If $x<<y$ and $y<<z$, then $x<<z$.
\item If $x\equiv y$, and $y<< z$, then $x<<z$.
\end{enumerate}

One consequence of $(3)$ is the following. Suppose $x_1x_2\cdots
x_n$ is a product of elements in $G$. Let $n_0$ be such that
$x_{n_0}$ is the largest among all the $x_i$'s. Then for $1\le
i\le n$, either $x_i\equiv x_{n_0}$ or $x_i<< x_{n_0}$. Using
$(3)$ inductively, we see that $x_{n_0}\equiv x_1x_2\cdots x_n$.

\begin{proof}[Proof of Proposition \ref{p-welpowerunion}]
We give a proof by contradiction. Suppose that $\bcup_{n\ge 0}
A^{\cdot n}$ is not well-ordered. Then we shall have an infinite
decreasing sequence in $\bcup_{n\ge 0} A^{\cdot n}$, say
$$a^{(1)} > a^{(2)} > \cdots $$
with $a^{(i)}$ a finite product of terms in $A$.

In every product $a^{(i)}$, there is at least one factor
 that is equivalent to $a^{(i)}$.
Let $a_i$ be the rightmost one. Then we can write $a^{(i)}
=x^{(i)}\cdot a_i\cdot y^{(i)}$, with
$a_i\in A$, $x^{(i)}$ and $y^{(i)}$ being finite products, and  $a_i \equiv x^{(i)}\cdot a_i 
\equiv a^{(i)}>>y^{(i)}$.

Case $1$. If there is an $N$ such that $a_1>>a_N$, then we also
have $a_1>> a_n$ for all $n\ge N$. Therefore we get an infinite
decreasing sequence $a^{(N)}>a^{(N+1)}> \cdots $ in $\cup_{n\ge 0}
A^{\cdot n}$, whose terms are all $<<a_1$. Record the above
sequence as
$$a_1 >>b^{(1)}> b^{(2)}>\cdots .$$

Case $2$. If $a_1\equiv a_i$ for all $i$, we let $A_1$ be the set
of all elements in $A$ that are Archimedean equivalent to $a_1$.
Then $A_1$ is well-ordered, greater than $1$, and Archimedean. By
Lemma \ref{l-welarki}, $\bcup_{n\ge 0} A_1^{\cdot n}$ is
well-ordered. So the infinite sequence $x^{(1)}\cdot
a_1,x^{(2)}\cdot a_2,\dots $ in $\bcup_{n\ge 0} A_1^{+n}$ has an
infinite weakly increasing subsequence, say $x^{(j_1)}\cdot
a_{j_1}\le x^{(j_2)}\cdot a_{j_2} \le \cdots $. Then $y^{(j_1)} >
y^{(j_2)} >\cdots $, whose terms are all $<<a_1$,
 is an infinite decreasing sequence in $\bcup_{n\ge 0} A^{+n}$. We can still record
 it as
 $$a_1 >>b^{(1)}> b^{(2)}>\cdots .$$

 Thus in either case, we can repeat the argument to get
 $a_1>>b_1>>c_1>>\cdots $, which is an infinite decreasing sequence in $A$.
 This contradicts the fact that $A$ is well-ordered.
\end{proof}

\begin{proof}[Proof of Theorem \ref{t-welcomposition}]
Suppose $A=\supp(\eta)$. Then $A$ is well-ordered and greater than
$1$. Let $a=\min A=\ord (\eta)$, and let $W=\bcup_{n\ge 0}
A^{\cdot n}$. By Proposition \ref{p-welpowerunion}, we see that
$W$ is well-ordered. Let $f=\sum_{n\ge 0} b_nt^n$. Then
$$[g] f\circ \eta =\sum_{a_1a_2\cdots a_m=g} b_m [{a_1}] \eta \cdots [{a_m}]\eta ,$$
where the sum is over all $m\in \NN$ and $a_1,a_2,\dots ,a_m\in
A$, we see that $\supp(f\circ \eta)$ is a subset of $W$, and hence
is well-ordered.

Now we let $P(g)$ be the claim that there are only finitely many tuples \\
$(m,a_1,a_2,\dots ,a_m)$, where $m\in \NN$ and $a_i\in A$, such
that $a_1a_2\cdots a_m=g$. In order to show that $f\circ \eta$ is
well defined, it suffices to show that $P(g)$ is true for all
$g\in W$.

It is clear that the least element of $W$ is $1$. In this case
$m=0$, so $P(1)$ is true. The second least element of $W$ is $a$.
In this case $m=1$ and $a_1=a$, so $P(a)$ is true. If $P(g)$ is
false for some $g\in W$ and $g>1$, then we can assume that $g$ is
the smallest such, for $W$ is well-ordered. So there are
infinitely many tuples, say for $i=1,2,\dots$,
$$a^{(i)}_1a^{(i)}_2\cdots a^{(i)}_{m_i} =g ,$$
where $m_i$ is a positive integer and $a^{(i)}_j\in A$  for all
$1\le j \le m_i$.

Consider $a^{(1)}_1, a^{(2)}_1,\dots $ as an infinite sequence in $A$. By Lemma 
\ref{l-infsubseq}, it contains an infinite increasing subsequence
or an infinite constant subsequence or both.

If $a^{(i_1)}_1< a^{(i_2)}_1<\cdots $ is an infinite increasing
subsequence, then $a^{(i_1)}_2\cdots a^{(i_1)}_{m_{i_1}}>
a^{(i_2)}_2\cdots a^{(i_2)}_{m_{i_2}}> \cdots$ is an infinite
decreasing sequence in $W$. This contradicts the fact that $W$ is
well-ordered.

If $a^{(i_1)}_1= a^{(i_2)}_1=\cdots $ is an infinite constant
sequence, then let $h=a^{(1)}_2\cdots a^{(1)}_{m_1}<g$. We have
$h\in W$ and $a^{(i_j)}_2a^{(i_j)}_3\cdots a^{(i_j)}_{m_{i_j}} =h$
for all $j$. Thus $P(h)$ is false. This contradicts the assumption
that $g$ is the smallest for $P(g)$ to be false. So $P(g)$ is true
for all $g\in W$.
\end{proof}

\begin{cor}
For any $\eta\in R_w[G]$ with initial term $1$,  $\eta^{-1}\in
R_w[G]$.
\end{cor}
\begin{proof}
Write $\eta=1-\tau$. Then $\tau\in R_w[G]$ and $\ord(\tau)>0$. By Theorem 
\ref{t-welcomposition}, $\sum_{n\ge 0} \tau^n$ strictly converges
in $R_w[G]$. One can check that $(1-\tau) \cdot \sum_{n\ge 0}
\tau^n =1.$
\end{proof}

So for any $\eta\in R_w[G]$ with initial term $f$, $\eta$ is
invertible if and only if $f$ is invertible. Write $f=a_g g$ with
$a_g\in R$. Then $f$ is invertible if and only if $a_g$ is
invertible in $R$ and $g$ is invertible in $G$. Thus if $G$ is a
group, then $\eta$ is invertible if and only if the coefficient
$a_g$ is invertible. Hence we have the following.

\begin{cor}
If $K$ is a field, and $G$ is a totally ordered group, then
$K_w[G]$ is a $K$-division algebra. Moreover, if $G$ is also
abelian, then $K_w[G]$ is a field.
\end{cor}

\begin{dfn}
If $G$ and $H$ are two TO-monoids, then the \emph{Cartesian
product} $G\times H$ is defined to be the set $G\times H$ equipped
with the usual multiplication and the reverse lexicographic order,
i.e., $(x_1,y_1)\le (x_2,y_2)$ if and only if $y_1<_H y_2$ or
$y_1=y_2$ and $x_1\le x_2$.
\end{dfn}

We define $G^n$ to be the Cartesian product of $n$ copies of $G$.
It is an easy exercise to show the following.

\begin{prop}
The Cartesian product of finitely many TO-monoids is a TO-monoid.
\end{prop}

One important example is that $\ZZ^n$ is a totally ordered abelian
group.

When considering the ring $R_w(G\times H)$, it is natural to treat
${(g,h)}$ as $gh$, where $g$ is identified with $(g,1)$ and $h$ is
identified with $(1,h)$. With this identification, we have the
following.

\begin{prop}\label{p-malcev-prod}
The ring $R_w[G\times H]$ is the same as the ring
$\left(R_w[G]\right)_w[H]$ of Malcev-Neumann series on $H$ with
coefficients in $R_w[G]$.
\end{prop}
\begin{proof}
Let $\eta \in R_w[G\times H]$, and let $A=\supp(\eta)$. Let $\rho$
be the second projection of $G\times H$, i.e., $\rho(g,h)=h.$

We first show that $\rho(A)$ is well-ordered. If not, then we have
an infinite decreasing sequence in $H$, say
$\rho(g_1,h_1)>\rho(g_2,h_2)>\cdots $, which by definition becomes
$h_1>h_2>\cdots  $. Then in the reverse lexicographic order, this
implies that $(g_1,h_1)>(g_2,h_2)>\cdots $ is an infinite
decreasing sequence of $A$, a contradiction. So $\rho(A)$ is
well-ordered.

Now $\eta$ can be written as
$$\eta =\sum_{h\in \rho(A)} \biggl(\sum_{g\in G, (g,h)\in A}  a_{g,h}
g\biggr) h.$$ Since for each $h\in \rho(A)$, the set $\{\, g\in G:
(g,h)\in A \,\}$ is a clearly a well-ordered subset of $G$,
$\sum_{g\in G, (g,h)\in A}  a_{g,h} g$ belongs to $R_w[G]$ for
every $h$, and hence $\eta\in R_w[G]_w[H]$.

Now let $\tau=\sum_{h\in D} b_h h \in R_w[G]_w[H]$, where
$D=\supp(\tau)$ is a well ordered subset of $H$, and $b_h\in
R_w[G]$. Let $B_h$ denote the support of $b_h$.
 We need to show that $\bcup_{h\in D} B_h\times \{\,h \,\})$ is well-ordered in $G\times H$.
Let $A$ be any subset of $\bcup_{h\in D} B_h\times \{\,h \,\})$.
We show that $A$ has a smallest element. Since $\rho(A)$ is a
subset of the well-ordered set $D$, we can take $h_0$ to be the
smallest element of $\rho(A)$. The set $A\cap B_{h_0}\times
\{\,h_0 \,\})$ is well-ordered for it is a subset of the
well-ordered set $B_{h_0}\times \{\,h_0 \,\})$. Let $(g_0,h_0)$ be
the smallest element of $A\cap B_{h_0}\times \{\,h_0 \,\})$. Then
$(g_0,h_0)$ is also the smallest element of $A$.
\end{proof}

\section{The Residue Theorem for MN-series\label{s-MNresidue}}

The main topic of this section is to describe and prove the
residue theorem for MN-series. We describe the theorem in the
first subsection and give the proof in the second.

\subsection{The Residue Theorem}
From now on, we let $R$ be a commutative ring with unit, and by a
monoid, we mean an abelian monoid written additively. Thus if
$\mathcal{G}$ is a monoid, then $R[\mathcal{G}]$ is a commutative
monoid ring.

Now let $\mathcal{G}$ be a totally ordered monoid. Then the ring
of MN-series $R_w[\mathcal{G}]$ is commutative. In order to
distinguish between the addition in $R$ and the addition in
$\mathcal{G}$, we replace $g$ by $t^g$. Thus $t^{g_1}
t^{g_2}=t^{g_1+g_2}$. Correspondingly, a formal series $\eta$ on
$\mathcal{G}$ has the form
$$ \sum_{g\in \mathcal{G}} a_g t^g,$$
where $a_g \in R$ and $t^g$ is regarded as a symbol. We also call
$t^g$ a monomial. Other terminologies are defined correspondingly.

Examples.
\begin{enumerate}
    \item $K_w[\ZZ]\simeq K((x))$ is the field of Laurent series.
    \item $K_w[\QQ]$ strictly contains the field $K^\fra ((x))$ of fractional Laurent
series, and is more complicated.
    \item By Proposition
\ref{p-malcev-prod}, $K_w[\ZZ^2]\simeq K((x_1))_w[\ZZ]\simeq
K((x_1))((x_2))$. Using induction, $K_w[\ZZ^n] \simeq K\ll
x_1,x_2,\dots,x_n\gg,$ is a field of iterated Laurent series,
which turns out to be the most useful special case.
\end{enumerate}

\vspace{3mm}

Observe that any submonoid of a \tomonoid\ is still a \tomonoid
under the induced total ordering. Let $\mathcal{G}$ be a TO-monoid
and let $\mathcal{H}$ be a monoid. If $\rho:\mathcal{H}\to
\mathcal{G}$ is an injective homomorphism, then
$\rho(\mathcal{H})\eql \mathcal{H}$ is a submonoid of
$\mathcal{G}$. We can thus regard $\mathcal{H}$ as a submonoid of
$\mathcal{G}$ through $\rho$. The induced ordering $\le^\rho$ on
$\mathcal{H}$ is given by $h_1\le^\rho h_2 \Leftrightarrow
\rho(h_1)\le_\mathcal{G} \rho(h_2)$. Thus $\mathcal{H}$ is a
TO-monoid under $\le^\rho$. Clearly a subset $A$ of
$(\mathcal{H},\le^\rho)$ is well-ordered if and only if $\rho(A)$
is well-ordered in $(\mathcal{G},\le_\mathcal{G})$.

Let $\mathcal{G}$ be a \tomonoid . We can give $\mathcal{G}$ a
different ordering so that under this new ordering $\mathcal{G}$
is still a \tomonoid . For instance, the total ordering
$\hat{\le}$ defined by $g_1\le g_2 \Leftrightarrow g_2\hat{\le }
g_1$ is clearly such an ordering. One special class of total
orderings is interesting for our purpose. If $\rho:\mathcal{G}\to
\mathcal{G}$ is an injective endomorphism, then the induced
ordering $\le^\rho$ is also a total ordering on $\mathcal{G}$. We
denote the corresponding ring of MN-series by
$R_w^\rho[\mathcal{G}]$.

For example, if $\mathcal{G}=\ZZ^n$, then any nonsingular matrix
$M\in GL(\ZZ^n)$ induces an injective endomorphism. In particular,
$K_w[\ZZ^2]\eql K\ll x,t \gg$ is the field of double Laurent
series, and $K_w^\rho[\ZZ^2]\eql K\ll \xx,t\gg$, where the matrix
corresponding to $\rho$ is the diagonal matrix $\diag (-1,1)$. It
is easy to see that $K\ll x_1^{e_1},\dots ,x_n^{e_n}\gg$ with
$e_i=\pm 1$ are special fields of MN-series, and the corresponding
matrices are diagonal matrices with entries $\pm 1$.

In order to state the residue theorem, we need more concepts.
Consider the following situation. Let $\mathcal{G}$ and
$\mathcal{H}$ be monoids with $\mathcal{H}\eql \ZZ^n$, and suppose
that we have a total ordering $\le $ on the direct sum
$\mathcal{G}\oplus \mathcal{H}$ such that $\mathcal{G}\oplus
\mathcal{H}$ is a \tomonoid . We identify $\mathcal{G}$ with
$\mathcal{G}\oplus 1$ and $\mathcal{H}$ with $1\oplus
\mathcal{H}$. Let $e_1,e_2,\dots ,e_n$ be a basis of
$\mathcal{H}$. Let $\rho$ be the endomorphism on
$\mathcal{G}\oplus \mathcal{H}$ that is generated by $\rho(e_i)=
g_i+\sum_j m_{ij} e_j$ for all $i$, where $g_i\in \mathcal{G}$,
and $\rho(g)=g$ for all $g\in \mathcal{G}$. Then $\rho$ is
injective if the matrix $M=\left(m_{ij}\right)_{1\le i,j\le n}$
belongs to $GL(\ZZ^n)$.

It is natural to use new variables $x_i$ to denote $t^{e_i}$ for
all $i$. Thus monomials in $R_w[\mathcal{G}\oplus \mathcal{H}]$
can be represented as $ t^gx_1^{k_1}\cdots x_n^{k_n}$.
Correspondingly, $\rho$ acts on monomials by $\rho( t^g)=t^g$ for
all $g\in \mathcal{G}$, and $\rho(x_i)= t^{g_i} x_1^{m_{i1}}\cdots
x_n^{m_{in}}.$

\noindent {\bf Notation}: If $f_i$ are monomials, we use $\mb{f}$
to denote the homomorphism $\rho$ generated by $\rho(x_i)=f_i$.

\vspace{3mm} An element $\eta$ of $R_w[\mathcal{G}\oplus
\mathcal{H}]$ can be written as
$$\eta=\sum_{\mb{k}\in \ZZ^n} \sum_{g\in \mathcal{G}} a_{g,\mb{k}} t^g x_1^{k_1}\cdots x_n^{k_n}
= \sum_{\mb{k}\in \ZZ^n} b_{\mb{k}} \mb{x}^{\mb{k}},$$ where
$a_{g,\mb{k}}\in R$ and $b_{\mb{k}}\in R_w[\mathcal{G}]$. We call
the $b_{\mb{k}}\mb{x}^{\mb{k}} $  an \emph{$x$-term} of $\eta$.
Since the set $\{\, \ord(b_{\mb{k}}\mb{x}^{\mb{k}}): \mb{k}\in
\ZZ^n\}$ is a subset of $\supp(\eta)$, it is well-ordered and
hence has a least element. Because of the different powers in the
$x$'s, no two of $\ord(b_{\mb{k}})\mb{x^k}$ are equal. So we can
define the $x$-{\em initial} term of $\eta$ to be the $x$-term
that has the least order.

Now the operators $\frac{\partial}{\partial x_i}$, $\ct_{x_i}$,
$\res_{x_i}$, $\pt_{x_i}$, $\nt_{x_i}$, and $\ct_{\mb{x}}$ are all
well defined in the ring of MN-series $R_w[\mathcal{G}\oplus
\mathcal{H}]$. Similar to the field of iterated Laurent series, we
have the basic computational rules for the field of MN-series.

\begin{lem}[Computational Rules]\label{l-comput-rule}
In a field of MN-series $K_w[\mathcal{G}\oplus \mathcal{H}]$ with
$\mathcal{H} \simeq \ZZ$, we identify $t^{(0,1)}$ with $x$, where
$(0,1)\in \mathcal{G}\oplus \mathcal{H}$. Let $F$ and $G$ be two
elements in $K_w[\mathcal{G}\oplus \mathcal{H}]$.
\begin{enumerate}
\item[Rule 1:] (linearity) For any $a,b$ that are independent of $x$,
$$\ct_{x} \left( a F({x})+bG({x})\right)= a \ct_{x} F({x}) + b \ct_{x_i} G({x}).$$

\item[Rule 2:] If $F$ can be written as $\sum_{k\ge 0}
a_k x^k$, then
$$\displaystyle \ct_{x} F = \left. F
\right|_{x=0}.$$

\item[Rule 3:]
$$\res_{x} \frac{\partial F({x})}{\partial x} G({x})= -\res_{x}
F({x}) \frac{\partial G({x})}{\partial x}.$$

\item[Rule 4:]
Suppose that $F$ is $\pt$ in $x$. If $G$ can be factored as
$(x-u)H$ such that $u$ is independent of $x$ and $\ord (u) >
\ord(x)$, and $1/H$ is $\pt$ in $x$, then
$$\ct_{x} F({x}) \frac{x}{G({x})}=\left. \frac{F({x})}{\displaystyle
\frac{\partial G({x})}{\partial x}} \right|_{x=u}$$
\end{enumerate}
\end{lem}

\vspace{3mm} Now we come back to the multivariate case, and
suppose $F_i\in R_w[\mathcal{G}\oplus \mathcal{H}]$ for all $i$.

 \begin{dfn}
 The Jacobian determinant (or simply Jacobian) of $\mb{F}$ with respect to $\mb{x}$ is defined to be
$$J\left({\mb{F}}|{\mb{x}}\right):=
J\left(\frac{F_1,F_2,\ldots ,F_n}{x_1,x_2,\dots ,x_n}\right)
=\det\left(\pad{x_j}{F_i} \right)_{1\le i,j\le n}.$$
\end{dfn}
 When the $x$'s  are clear, we write $J(F_1,F_2,\dots
,F_n)$ for short.

\begin{dfn}
If the $x$-initial term of $F_i$ is $a_i x_1^{b_{i1}}\cdots
x_n^{b_{in}}$, then the Jacobian number of $\mb{F}$ with respect
to $\mathbf{x}$ is defined to be
$$j\left({\mb{F}}|{\mb{x}}\right):=
j\left(\frac{F_1,F_2,\ldots ,F_n}{x_1,x_2,\dots ,x_n}\right)
=\det\left(b_{ij} \right)_{1\le i,j\le n}.$$
\end{dfn}

\begin{dfn}
The log Jacobian of $F_1,\dots ,F_n$ is defined to be
$$LJ(F_1,\dots,F_n):=\frac{x_1\cdots x_n}{F_1 \cdots F_n} J(F_1,\dots ,F_n).$$
\end{dfn}
We call it the log Jacobian because formally it can be written as
\citep{wilson}
$$LJ(F_1,\dots,F_n)=J\left(\frac{\log F_1,\ldots ,\log F_n}{\log x_1,
\dots ,\log x_n}\right),$$ since
$$\frac{\partial \log F}{\partial \log x}= \frac{\partial \log F}{\partial F}
\frac{\partial F}{\partial \log x}=\frac{1}{F}\frac{\partial
F}{\partial  x} \frac{\partial x}{\partial \log
x}=\frac{x}{F}\frac{\partial F}{\partial x}.$$

\begin{rem}
Generally speaking, the Jacobian is convenient in residue
evaluations, while the log Jacobian is convenient in constant term
evaluations.
\end{rem}

The following lemma is devised for the proof of our residue
theorem. It is also a kind of composition law.
\begin{lem}
\label{l-malcev-comp} Let $\Phi$ be a formal series in $x_1,\dots,
x_n$ with coefficients in $R_w[\mathcal{G}]$. Then $\Phi(F_1,\dots
,F_n)\in R_w[\mathcal{G}\oplus \mathcal{H}]$ if and only if
$\Phi(f_1,\cdots ,f_n)\in R_w[\mathcal{G}\oplus \mathcal{H}]$,
where $f_i$ is the $x$-initial term of $F_i$ for all $i$. Moreover
if $j(F_1,\dots ,F_n)\ne 0$, then $\Phi(F_1,\dots F_n)\in
R_w[\mathcal{G}\oplus \mathcal{H}]$ if and only if $\Phi(x_1,\dots
,x_n)\in R_w^\mb{f} [\mathcal{G}\oplus \mathcal{H}]$.
\end{lem}

This lemma reduces the convergence of $\Phi(F_1,\dots, F_n)$ to
that of $\Phi(f_1,\dots ,f_n)$. For example, the ring of formal
power series $R[[x_1,\dots ,x_n]]$ is isomorphic to $R_w[\NN^n]$,
where $\NN^n$ itself is well-ordered under the reverse
lexicographic ordering. If $\Phi$ is a formal power series in
$\mb{x}$, then $\Phi(f_1,\dots ,f_n)$ is also a formal power
series when $f_i$ are monomials in $R[[x_1,\dots ,x_n]]$. Thus
Lemma \ref{l-malcev-comp} implies the composition law of
$R[[x_1,\dots ,x_n]]$.

\begin{proof}[Proof of Lemma \ref{l-malcev-comp}]
Write every $F_i$ as $f_i(1+\tau_i)$, where $f_i$ is the
$x$-initial term and $\ord(\tau_i)>1$ or $\tau_i=0$.

For the first part, we show that if $\Phi(F_1,\dots ,F_n) \in
R_w[\mathcal{G}\oplus \mathcal{H}]$, then replacing $F_i$ by
$F_i(1+\tau)$ with $\ord(\tau)>0$ results in an element of
$R_w[\mathcal{G}\oplus \mathcal{H}]$. Then the first part follows
by replacing $F_i$ by $f_i=F_i(1+\tau_i)^{-1}$, (or conversely,
$f_i$ by $F_i=f_i(1+\tau_i)$) one by one for $i$ from $1$ to $n$.

We deal with the case $i=1$ as follows. The case of arbitrary $i$
is similar. Let $A=\supp(\Phi(F_1,\dots, F_n))$ and
$T=\supp(\tau)$. Then by assumption, $A$ is well-ordered, and $T$
is positive and well-ordered. We can write
$$\Phi(F_1,\dots ,F_n) =\sum_{k\in \ZZ} d_k F_1^{k},$$
where $d_k$ is a formal series in $F_2,\dots ,F_n$ with
coefficients in $R_w[\mathcal{G}]$. Then
\begin{align}\label{e-3-pcriterion}
\Phi(F_1(1+\tau),\dots ,F_n) =\sum_{k\in \ZZ} d_k
F_1^{k}(1+\tau)^k =\sum_{k\in \ZZ} d_kF_1^k \sum_{l\ge 0}
\binom{k}{l} \tau^l
\end{align}
Now we see that the support of $\Phi(F_1(1+\tau),\dots ,F_n)$ is a
subset of
$$\bcup_{l\ge 0} (A+T^{+l})=A+\bcup_{l\ge 0} T^{+l},$$ which is well-ordered by
Proposition \ref{p-welplus} and Proposition \ref{p-welpowerunion}.

To see that the coefficient of $t^{g+h}$ is a finite sum for every
$g$ and $h$, we observe that replacing each $\binom{k}{l}$ by $1$
will not decrease the number of summands. The right side of
equation \eqref{e-3-pcriterion} then becomes
$$\Big(\sum_{k\in \ZZ} d_k F_1^{k}\Big) \Big(\sum_{l\ge 0} \tau^{l}
\Big),$$ in which the coefficient of $t^{g+h}$ is a finite sum,
because it is a product of two elements in $R_w[\mathcal{G}\oplus
\mathcal{H}]$.

For the second part, if $j(F_1,\dots ,F_n)\ne 0$, then $\rho:
x_i\to {f_i}$ induces an injective endomorphism on
$\mathcal{G}\oplus \mathcal{H}$. We see that
 $\supp(\Phi(f_1,\dots, f_n))$  is well-ordered in $\mathcal{G}\oplus \mathcal{H}$ if and only
if $\rho \left(\supp(\Phi(x_1,\dots, x_n))\right)$ is
well-ordered. This, by definition, is to say that $
\Phi(x_1,\dots, x_n)\in R_w^\mb{f} [\mathcal{G}\oplus
\mathcal{H}].$ The lemma now follows from the first part.
\end{proof}

\vspace{3mm} \noindent {\bf Notation.} Starting with a totally
ordered abelian monoid $\mathcal{G}\oplus \mathcal{H}$ as
described above, let $\Phi$ be a formal series on
$\mathcal{G}\oplus \mathcal{H}$. When we write $\ct_x^\rho
\Phi(x_1,\dots ,x_n)$, we mean both that $\Phi(x_1,\dots ,x_n)$
belongs to $ R_w^\rho[\mathcal{G}\oplus \mathcal{H}]$, and that
the constant term is taken in this ring. When $\rho$ is the
identity map, it is omitted. When we write $\ct_\mb{F}
\Phi(F_1,\dots ,F_n)$, it is assumed that $\Phi(x_1,\dots ,x_n)\in
R_w^\mb{f}[\mathcal{G}\oplus \mathcal{H}]$, and we are taking the
constant term of $\Phi(x_1,\dots ,x_n) $ in the ring
$R_w^\mb{f}[\mathcal{G}\oplus \mathcal{H}]$. Or equivalently, we
always have
$$\ct_\mb{F} \Phi(F_1,\dots ,F_n) = \ct_{\mb{x}}\mbox{}^\mb{f} \Phi(x_1,\dots ,x_n).$$
This treatment is particularly useful  when dealing with rational
functions, as we shall see soon.

\vspace{3mm} Now comes our residue theorem for
$R_w[\mathcal{G}\oplus \mathcal{H}]$, in which we will see how an
element in one ring is related to an element in another ring
through taking the constant terms.
\begin{thm}[Residue Theorem]\label{t-MNresidue}
Suppose for each $i$, $F_i\in R_w[\mathcal{G}\oplus \mathcal{H}]$
has $x$-initial term $f_i=a_i x_1^{b_{i1}}\cdots x_n^{b_{in}}$
with $a_i, a_i^{-1}\in R_w[\mathcal{G}]$. If $j(F_1,\dots ,F_n)\ne
0$, then
 for any
 $\Phi(\mb{x})\in R_w^\mb{f} [\mathcal{G}\oplus \mathcal{H}]$, we have
\begin{align}\label{e-wel-residue}
\ct_\mathbf{x} \Phi(F_1,\dots ,F_n) LJ(F_1,\dots ,F_n)
=j(F_1,\dots ,F_n)\ct_\mb{F} \Phi(F_1,\dots ,F_n) .
\end{align}
\end{thm}
\begin{proof}[Proof of Theorem \ref{t-MNresidue}]
With the hypothesis, both sides of equation \eqref{e-wel-residue}
converge. In fact, Lemma \ref{l-malcev-comp} is designed for this
convergence.

Now by multilinearity, it suffices to show the theorem is true for
monomials $\Phi$. This is the main topic of the next subsection.
See Lemmas \ref{l-3-residue-n} and \ref{l-3-residue-y}.
\end{proof}

\begin{rem}
If $j(F_1,\dots ,F_n)= 0$, then $\Phi(F_1,\dots ,F_n)$ is only
well defined in some special cases.
\end{rem}

\begin{rem}
If $\Phi(x_1,\dots ,x_n)$ is a Laurent polynomial, then
$\Phi(F_1,\dots ,F_n)$ always exists. In this case, it is not
necessary to consider the map $\mb{f}$.
\end{rem}

Now let $K$ be a field and let $\mathcal{G}$ be a group. We are
going to consider both fields $K_w[\mathcal{G}\oplus \mathcal{H}]$
and $K_w^\mb{f}[\mathcal{G}\oplus \mathcal{H}]$, where $\mb{f}$ is
an injective endomorphism. Both fields contain
$K_w[\mathcal{G}][\mathcal{H}]$ as a subring, and thus contain the
quotient field of $K_w[\mathcal{G}][\mathcal{H}]$, which is the
field of rational functions. The operator $\ct_{x_i}$ is always
well defined. But the results of $\ct_{x_i}$ acting on a rational
function $\Phi$ will be different when working in different
fields. More precisely, let $\eta$ be the denominator of $\Phi$.
 Because of the different orderings, $\ord (\eta)$ in
$K_w[\mathcal{G}\oplus \mathcal{H}]$ is usually different from
$\ord^{\mb{f}}(\eta)$ in $K_w^\mb{f}[\mathcal{G}\oplus
\mathcal{H}]$. Thus $\eta^{-1}$ has different expansions in the
two fields.

\begin{rem}
If $R=K$ is a field and $\mathcal{G}$ is a group, then
$K_w^{\mb{f}}[\mathcal{G}\oplus \mathcal{H}]$ is a field for all
injective $\mb{f}$. Thus $K_w^\mb{f}[\mathcal{G}\oplus
\mathcal{H}]$ contains all rational $\Phi$. In applications of
this theorem, we need only to expand $\Phi$ correctly in a
specified field.
\end{rem}

\subsection{The Proof of the Theorem}

The proof of our residue theorem and lemmas basically comes from
\citep{reversion}, except for the proof of Lemma
\ref{l-3-residue-y}.

In what follows, we suppose $F_i, F_i^{-1}\in
R_w[\mathcal{G}\oplus \mathcal{H}]$ for all $i$.

The following properties of Jacobians can be easily checked.
\begin{lem}\label{l-lres}
Let the Jacobian be defined as in the previous subsection. Then
\begin{enumerate}
\item $J(F_1,F_2,\ldots ,F_n)$ is $R_w[\mathcal{G}]$-multilinear.
\item $J(F_1,F_2,\ldots ,F_n)$ is alternating; i.e.,
$J(F_1,F_2,\ldots ,F_n)=0$ if $F_i=F_j$ for some $i\ne j$.
\item $J(F_1,F_2,\ldots ,F_n)$ is anticommutative; i.e.,
$$J(F_1,\ldots,F_i,\ldots ,F_j,\ldots,F_n)=-J(F_1,\ldots,F_j,\ldots ,F_i,\ldots,F_n).$$
\item $($Composition Rule$)$ If $g(z)\in K((z))$ is a series in one variable, then
$$J(g( F_1),F_2,\ldots ,F_n)=g'(F_1) J(F_1,F_2,\ldots ,F_n).$$
\item $($Product Rule$)$ $$J(F_1G_1,F_2,\ldots ,F_n)=F_1J(G_1,F_2,\ldots ,F_n)+
G_1J(F_1,F_2,\ldots ,F_n).$$
\item $J(F_2^{-1},F_2,\ldots ,F_n)=0$.
\end{enumerate}
\end{lem}

\begin{lem} \label{l-residue-f} If all $F_i$ are $x$-monomials, then
\begin{align}\label{e-3-jacobian}
LJ(F_1\dots ,F_n)=j(F_1,\dots ,F_n).
\end{align}
\end{lem}
\begin{proof}
Suppose that  for every $i$, $F_i=a_i x_1^{b_{i1}}\cdots
x_n^{b_{in}}$, where $a_i$ is in $R_w[\mathcal{G}]$. Factoring
$F_i=a_ix_1^{b_{i1}}\cdots x_n^{b_{in}}$  from the $i$th row of
the Jacobian matrix for all $i$ and then factoring $x_j^{-1}$ from
the $j$th column for all $j$, we get
$$J(F_1,F_2,\ldots ,F_n)=\frac{F_1\cdots F_n}{x_1\cdots x_n}\det(b_{ij}).$$
Equation \eqref{e-3-jacobian} is just a rewriting of the above
equation.
\end{proof}

\begin{lem}\label{l-residue}
$$\res_x J(F_1,\ldots , F_n)=0.$$
\end{lem}
\begin{proof}
By multilinearity, it suffices to check monomials $F_i$. Suppose
that they are given as in Lemma \ref{l-residue-f}. Then equation
\eqref{e-3-jacobian} gives us
$$J(F_1,\dots ,F_n)= j(F_1,\dots ,F_n)\frac{F_1\cdots F_n}{x_1\cdots x_n}.$$
More explicitly,
$$J(F_1,\dots ,F_n)= a_1\cdots a_n \det(b_{ij}) x_1^{-1+\sum b_{i1}} \cdots x_n^{-1+\sum
{b_{in}}}j(F_1,\dots ,F_n). $$ If $\sum {b_{i1}}=\sum
{b_{i2}}=\cdots =\sum {b_{in}}=0$, then the Jacobian number is
$0$, and therefore the residue is $0$. Otherwise, at least one of
the $x_i$'s has exponent $\ne -1$,
 so the residue is $0$ by definition.
\end{proof}

\begin{lem}\label{l-3-residue-n}
For all integers $e_i$ with at least one of $e_i\ne -1$, we have
\begin{align}\label{e-residueth}
\res_x F_1^{e_1}\cdots F_n^{e_n} J(F_1,\ldots ,F_n)= 0.
\end{align}
\end{lem}
\begin{proof}
The clever proof in \citep[Theorem 1.4]{reversion} also works
here.

Permuting the $F_i$ and using $(3)$ of Lemma \ref{l-lres}, we may
assume that $e_1\ne -1$,\dots, $e_j\ne -1$, but $e_{j+1}=\cdots
=e_{n}=-1$, for some $j$ with $1\le j\le n-1$. Setting
$G_i=\frac{1}{e_i+1}F_i^{e_i+1}$ for $i=1,\ldots , j$, we have
$$F_1^{e_1}F_2^{e_2}\cdots F_n^{e_n} J(F_1,F_2,\ldots ,F_n)=
F_{j+1}^{-1}\cdots F_n^{-1} J(G_1,\ldots, G_j,F_{j+1}, \ldots
,F_n).$$ Then applying the formula
$$F_{j+1}^{-1} J(G_1,\ldots, G_j,F_{j+1}, \ldots ,F_n)=
J(F_{j+1}^{-1}G_1,G_2\ldots, G_j,F_{j+1}, \ldots ,F_n)$$
repeatedly for $j+1,j+2,\dots ,n$, we get
$$ J(F_{j+1}^{-1}\cdots F_n^{-1}G_1,G_2, \ldots, G_j,F_{j+1}, \ldots ,F_n).$$
The result now follows from Lemma \ref{l-residue}.
\end{proof}

For the case $e_1=e_2=\cdots =e_n=-1$, we have
\begin{lem}\label{l-3-residue-y}
\begin{align}
\label{e-3-residue-y}
 \res_{\mb{x}}F_1^{-1}\cdots F_n^{-1}J(F_1,\dots ,F_n)=j(F_1,\dots, F_n).
\end{align}
\end{lem}
The simple proof for this case in \citep{reversion} does not apply
in our situation. The reason will be explained in Proposition
\ref{p-3-ljacobian}.

Note that Lemma \ref{l-3-residue-y} is equivalent to saying that
\begin{align}
\label{e-3-lresidue-y}
 \ct_{\mb{x}}LJ(F_1,\dots ,F_n)=j(F_1,\dots, F_n).
\end{align}

\begin{proof}
Let $f_i:=a_ix_1^{b_{i1}}\cdots x_n^{b_{in}}$ be the $x$-initial
term of $F_i$. Then $F_i=f_i B_i$, where $B_i\in
R_w[\mathcal{G}\oplus \mathcal{H}]$ has $x$-initial term $1$. By
the composition law,
 $\log ( B_i) \in
R_w[\mathcal{G}\oplus \mathcal{H}] $. Now applying the product
rule, we have
\begin{align*}
&F_1^{-1}\cdots F_n^{-1}J(F_1,F_2,\ldots ,F_n) \\
&\qquad \qquad= f_1^{-1}F_2^{-1}\cdots F_n^{-1} J(f_1,F_2,\ldots,
F_n)
+B_1^{-1}F_2^{-1}\cdots F_n^{-1} J(B_1,F_2,\ldots, F_n)\\
&\qquad\qquad=f_1^{-1}F_2^{-1}\cdots F_n^{-1} J(f_1,F_2,\ldots,
F_n)+ F_2^{-1}\cdots F_n^{-1} J( \log (B_1),F_2,\ldots, F_n).
\end{align*}
From Lemma \ref{l-3-residue-n}, the last term in the above
equations has no contribution to the residue in $x$, and hence can
be discarded.

The same procedure can be applied to $F_2,F_3,\ldots ,F_n$.
Finally we will get
$$\res_x F_1^{-1}\cdots F_n^{-1}J(F_1,F_2,\ldots ,F_n)=
\res_x f_1^{-1}\cdots f_n^{-1} J(f_1,f_2,\ldots ,f_n),$$ which is
equal to the Jacobian number by Lemma \ref{l-residue-f}.
\end{proof}

The following proposition gives a good reason for using the log
Jacobian.
\begin{prop}\label{p-3-ljacobian}
The $x$-initial term of the log Jacobian $LJ(F_1,\dots ,F_n)$
equals the Jacobian number $j(F_1,\dots ,F_n)$ when it is nonzero.
\end{prop}
\begin{proof}
From the definition,
$$LJ(F_1,\dots ,F_{n})= \frac{x_1\cdots x_n}{F_1\cdots F_n}
\det \left( \displaystyle\frac{\partial F_i}{\partial x_j}
\right).$$ To obtain the $x$-initial term, we replace every term
with its $x$-initial term. The result will be of the least order
unless it is zero. Therefore by Lemma \ref{l-residue-f}, we can
write
$$LJ(F_1,\ldots ,F_n )= j(F_1,\dots ,F_n)+ \text{higher order terms}.$$
To show that $j(F_1,\dots ,F_n)$ is the $x$-initial term, we need
to show that all the other term that are independent of $x$
cancel. (Note that we do not have this trouble when all the
coefficients belong to $R$.) This is equivalent to saying that
$$ \ct_{\mb{x}}LJ(F_1,\ldots ,F_n )= j(F_1,\dots ,F_n),$$
which follows from Lemma \ref{l-3-residue-y}.
\end{proof}

\begin{exa}
Consider the field $K\ll x,t \gg$. Let $F=x^2+xt+x^3t$. Then the
initial $x$-term of $F$ is $x^2$. Now let us see what happens to
the log Jacobian $LJ(F|x)$ of $F$ with respect to $x$.
\end{exa}
\begin{align*}
LJ(F|x)=\frac{x}{F}\frac{\partial F}{\partial x}
&=\frac{x(2x+t+3x^2t)}{x^2(1+t/x+xt)} \\
&=(2+t/x+3xt) \sum_{k\ge 0} (-1)^k (t/x+xt)^k
\end{align*}
It is not clear that $2$ is the unique term in the expansion, but
all the other terms cancel. We continue to check as the following.
\begin{align*}
\ct_x LJ(F|x) &= \ct_x (2+t/x+3xt) \sum_{k\ge 0} (-1)^k (t/x+xt)^k \\
&= 2 \sum_{k\ge 0 }\binom{2k}{k}t^{2k}- t\sum_{k\ge
0}\binom{2k+1}{k}t^{2k+1}-3t \sum_{k\ge 0} \binom{2k+1}{k+1} t^{2k+1} \\
&= 2+ \sum_{k\ge 1} \left( 2 \binom{2k}{k}-4
\binom{2k-1}{k}\right)t^{2k}.
\end{align*}
Now it is easy to see that the terms not containing $x$ in the
expansion of the log Jacobian really cancel.

From Theorem \ref{t-MNresidue} and Lemma \ref{l-residue-f}, we see
directly the following result.
\begin{cor}\label{c-3-residue-m}
If $F_i$ are all $x$-monomials in $K_w[\mathcal{G}\oplus
\mathcal{H}]$, and $\Phi\in K_w^\mb{F}[\mathcal{G}\oplus
\mathcal{H}]$, which indicates that $j(F_1,\dots ,F_n)$ is
nonzero, then
$$\ct_{\mb{x}} \Phi(F_1,\dots ,F_n)=\ct_{F_1,\dots ,F_n} \Phi(F_1,\dots ,F_n).$$
\end{cor}
This is saying that change of variables by monomials will not
change the constant terms. Now it is easy to understand the
phenomenon of Example \ref{ex-residue}.

In the case that all $F_i$ are monomials in $K[\mb{x},\mb{\xx}]$
with $j(\mb{F})\ne 0$, $\Phi$ is in $K[\mb{x},\mb{\xx}]$ if and
only $\Phi(F_1,\dots ,F_n)$ is. We always have
$$\ct_{F_1,\dots ,F_n} \Phi(F_1,\dots ,F_n)=\ct_{x_1,\dots ,x_n} \Phi(x_1,\dots ,x_n).$$
More generally, we have the following result, which will be used
later.

\begin{cor}\label{c-3-monomial}
Suppose $\mb{y}$ is another set of variables. If $\Phi \in
K[\mb{x},\mb{\xx}]\ll \mb{y}\gg$, and if $F_i$ are all monomials
in $\mb{x}$ with $j(\mb{F})\ne 0$, then
$$\ct_{\mb{x}} \Phi(F_1,\dots ,F_n)=\ct_{\mb{x}} \Phi(x_1,\dots ,x_n).$$
\end{cor}

\begin{exa}
Evaluate the following constant term in $\CC\ll x,y,t\gg$.
\begin{align}\label{exa-3-residue}
\ct_{x,y}-{x}^{3}{e^{{\frac {t}{xy}}}} \left( 3\,xy-2\,t \right)
\left( {x}^{3 }y{e^{{\frac {t}{xy}}}}-tx-ty \right) ^{-1} \left(
x-y \right) ^{-1}
 \left( -1+{x}^{3}{e^{{\frac {t}{xy}}}} \right) ^{-1}.
\end{align}
This is an example that is hard to evaluate without using our
residue theorem.
\end{exa}
Let $F=x^2ye^{\frac{t}{xy}}$, $G=xy^2 e^{\frac{t}{xy}}$. It is
easy to compute the log Jacobian and the Jacobian number. We have
$$LJ(F,G|x,y)=3-\frac{2t}{xy}, \text{ and }j(F,G|x,y)=3. $$
We can check that \eqref{exa-3-residue} can be written as
$$\ct_{x,y} \frac{F^3G}{(F^2-(F+G)t)(F-G)(G-F^2)}LJ(F,G|x,y).$$
Thus by the residue theorem, the above constant term equals
\begin{align}\label{exa-3-residue1}
\ct_{F,G} \frac{3F^3G}{(F^2-(F+G)t)(F-G)(G-F^2)}=\ct_{F,G}
\frac{3}{ (1-\frac{(F+G)t}{F^2})(1-\frac{G}{F})(1-\frac{F^2}{G})},
\end{align}
where on the right hand side of \eqref{exa-3-residue1}, we can
check that $1$ is the initial term of each factor in the
denominator.

At this stage, we can use the series expansion to obtain the
constant term. But we will evaluate it by the computational rule 4
in Lemma \ref{l-comput-rule}.

Starting from the left hand-side of \eqref{exa-3-residue1}, we
first take the constant term in $G$. We can solve for $G$ in the
denominator since all these three factors are linear in $G$. Only
one root, $F^2$, has higher order than $G$. Thus we can apply rule
4 and get
\begin{align*} \ct_{F,G}
\frac{3F^3G}{(F^2-(F+G)t)(F-G)(G-F^2)}&=\ct_F
\frac{3F^3}{(F^2-(F+F^2)t)(F-F^2)}\\
&=\ct_F \frac{3F}{(F-(1+F)t)(1-F)} \\
&=\frac{3}{(1-t)(1-\frac{t}{1-t})},
\end{align*}
where in the last step, we applied rule 4 again. One can check
that the two roots of the denominator for $F$ are $t/(1-t)$ and
$1$, and that only the former root has higher order than $F$.

After simplification, we finally get
$$
\ct_{x,y}-{x}^{3}{e^{{\frac {t}{xy}}}} \left( 3\,xy-2\,t \right)
\left( {x}^{3 }y{e^{{\frac {t}{xy}}}}-tx-ty \right) ^{-1} \left(
x-y \right) ^{-1}
 \left( -1+{x}^{3}{e^{{\frac {t}{xy}}}} \right)
 ^{-1}=\frac{3}{1-2t}.$$

\section{Another View of Lagrange's Inversion Formula\label{ss-lag}}
Let $F_1,\ldots ,F_n$ be power series in variables $x_1,\ldots
,x_n$ of the form $F_i=x_i+$ ``higher degree terms",  with
indeterminate coefficients for each $i$. It is known, e.g.,
\citep[Proposition 5, p. 219]{jac}, that $\mathbf{F}=(F_1,\ldots
,F_n)$ has a unique compositional inverse, i.e., there exists
$\mathbf{G}=(G_1,\ldots , G_n)$ where each $G_i$ is a power series
in $x_1,\ldots ,x_n$ such that $F_i(G_1,\ldots ,G_n)=x_i$
 and $G_i(F_1,\ldots ,F_n)=x_i$ for all $i$.

Lagrange inversion gives a formula of $G$'s in terms of $F$'s.

The above case is known as non-diagonal case. The diagonal case is
when $F_i$ divides $x_i$ for every $i$, or equivalently, $F_i=x_i
H_i$, where $H_i\in K[[x_1,\dots ,x_n]]$ with constant term $1$.

The formula of Good deals with the case when $x_i$ in fact divides
$F_i$. Such a formula is called diagonal (or Good's) Lagrange
inversion formula. This formula can be easily derived by the
ordinary residue theorem. We can illustrate this in our terms.

In the diagonal case, we can suppose that $F_i=x_iH_i$, where
$H_i$ is in $K[[x_1,\ldots ,x_n]]$ with constant term $1$.
Consider this in the field $K\ll x_1,x_2,\ldots ,x_n\gg $. Then
$x_i$ is the initial term of $F_i$, and the Jacobian number
$j(F_1,\dots ,F_n)=1$.

Change variables by $y_i=F_i(x)$, we will have $x_i=G_i(y)$. Then
\begin{align*}
[y_1^{k_1}\cdots y_n^{k_n}] G_i(y) &= \res_y y_1^{-1-k_1}\cdots y_n^{-1-k_n} G_i(y)\\
&= \res_x F_1^{-1-k_1}\cdots F_n^{-1-k_n} x_i J(\mathbf{F}),
\end{align*}
 where $J(\mathbf{F})$ is the Jacobian of $F_1,\ldots, F_n$.

Now let us consider the non-diagonal case. In this case, we cannot
apply the residue Theorem \ref{t-MNresidue} directly, because when
working in $K\ll x_1,\ldots, x_n\gg $,  we might meet the
situation that the Jacobian number equals $0$. For example, if
$x_n$ does not divide $F_n$, then it is easily seen that the power
of $x_n$ in the initial term of $F_i$ is zero for all $i$. So the
Jacobian number of $F_1,\dots ,F_n$ is $0$.

This difficulty can be overcome by introducing a new variable $t$.
After we get a suitable formula, replace $t$ by $1$. The result
obtained this way is equivalent to the homogeneous expansion
introduced in \citep{reversion}.

The working field is $K\ll x_1,x_2,\ldots ,x_n,t\gg $. In stead of
dealing with $F_1,\dots ,F_n$ directly, we consider the
compositional inverse of the system $y_i=F_i(x_1t,\ldots ,x_nt)$.
Clearly if there is a solution, we shall have
$x_it=G_i(\mathbf{y})$. Then by setting $t=1$, we will get the
desired result.

Since the initial term of $F_i(x_1t,\ldots ,x_nt)$ is $x_it$, the
Jacobian number is $1$. It is also easy to see that
$J(\mathbf{F}(t\mb{x}))=t^n J(\mathbf{F})|_{\mb{x}=t\mb{x}}$. So
we have the same formula, but interpreted differently. Setting
$t=1$ in the result is valid, since the power in $t$ equals the
sum of powers in the $x_i$'s. This is equivalent to the
homogeneous expansion.

Let $\Phi \in K[[y_1,\ldots ,y_n]]$. We get the formula
\begin{equation}
[y_1^{k_1}\cdots y_n^{k_n}] \Phi(\mathbf{G}) =\res_x
F_1^{-1-k_1}\cdots F_n^{-1-k_n} \Phi(\mathbf{x}) J(\mathbf{F}).
\end{equation}

Multiplying both sides of the above equation by $y_1^{k_1}\cdots
y_n^{k_n}$, and summing on all nonnegative integers
$k_1,k_2,\ldots ,k_n$, we get
\begin{align}\label{e-3-another}
\Phi(\mathbf{G}(\mathbf{y}))=\res_x \frac{1}{F_1-y_1}\cdots
\frac{1}{F_n-y_n}J(\mathbf{F}) \Phi(\mathbf{x}),
\end{align}
which is true as power series in the $y_i$'s.

It's natural to ask if we can get this formula directly from the
Residue Theorem. The answer is yes. The argument is given as
follows.

Working in $K\ll x_1,\ldots ,x_n,y_1,\ldots ,y_n\gg $. We make the
change of variables by $z_i=F_i-y_i$. Then
$x_i=G_i(\mb{y}+\mb{z})$, and the initial term of $F_i-y_i$ is
$x_i$, for $y_i$ has higher order. Thus the Jacobian number is
$1$. The Jacobian determinant still equals to $J(\mathbf{F})$.
Applying the residue theorem, we get
$$\res_x \frac{1}{F_1-y_1}\cdots \frac{1}{F_n-y_n}J(\mathbf{F})
\Phi(\mathbf{x}) =\res_z \frac{1}{z_1z_2\cdots z_n}
\Phi(\mathbf{G}(\mb{y}+\mb{z})).$$ Since $\Phi(G(y+z))$ is in
$K[[\mb{y},\mb{z}]]$. The final result is obtained by setting
$\mb{z}=\mb{0}$ in $\Phi(G(\mb{y}+\mb{z}))$.

Note that $J(\mathbf{F})\in K[[\mb{x}]]$ has constant term $1$.
Therefore $J(\mathbf{F})^{-1}\Phi(\mb{x})$ is also in
$K[[\mb{x}]]$. Hence we can reformulate \eqref{e-3-another} as
$$\res_\mb{x} \frac{1}{F_1-y_1}\cdots \frac{1}{F_n-y_n}\Phi(\mb{x})=
\Phi(\mb{x})J(\mathbf{F})^{-1}|_{\mb{x}=\mb{G}}.$$

\vspace{3mm} Here is another way to prove Lagrange's Inversion
formula. We only give the proof for the case $n=2$. The general
case is similar by induction.

Applying Theorem \ref{t-lagrange1} with respect to $x_1$, we get
$$\CT_{x_1,x_2} x_1x_2 \frac{1}{F_1-y_1} \frac{1}{F_2-y_2}\Phi(x_1,x_2)
=\CT_{x_2} x_2\frac{1}{\pad{x_1}{F_1(x_1,x_2)}} \frac{1}{F_2-y_2}
\Phi(x_1,x_2)|_{x_1=H_1},$$ where $H_1=H_1(x_2,y_1)\in
K[[x_2,y_1]]$, so that
$$F_1(H_1,x_2)-y_1=0.$$
Now let $H_2=H_2(y_1,y_2)\in K[[y_1,y_2]]$, so that
$$F_2(H_1,x_2)-y_2=0.$$
Applying Theorem \ref{t-lagrange1} with respect to $x_2$, we get
$$\CT_{x_1,x_2} x_1x_2 \frac{1}{F_1-y_1} \frac{1}{F_2-y_2}\Phi(x_1,x_2)
=\frac{1}{\pad{x_1}{F_1(H_1,H_2)}} \frac{1}{\frac{\partial
F_2(H_1,H_2)}{\partial x_2}} \Phi(H_1,H_2).$$ Now it is routine to
check that $\pad{x_1}{F_1(H_1,H_2)}\frac{\partial
F_2(H_1,H_2)}{\partial x_2}$ equals the Jacobian.

\section{About Dyson's Conjecture\label{s-dyson}}
We give an example of the application of the residue theorem. The
following is a conjecture of Dyson.
\begin{thm}\label{t-dyson}
Let $a_1,\ldots ,a_n$ be $n$ nonnegative integers. Then the
following equation holds as Laurent polynomials in $\mb{z}$.
\begin{equation}
\CT_{\mb{z}} \prod_{1\le i\ne j \le n}
\left(1-\frac{z_i}{z_j}\right)^{a_j} =
 \frac{(a_1+a_2+\cdots a_n)!}{a_1!\, a_2!\, \cdots a_n!}.\label{e-dyson}
\end{equation}
\end{thm}

For $n=3$ this assertion is equivalent to the familiar Dixon
identity:
\begin{equation}
\sum_{j}(-1)^j \binom{a+b}{a+j}\binom{b+c}{b+j}\binom{c+a}{c+j}
=\frac{(a+b+c)!}{a!\, b!\, c!}.
\end{equation}
Theorem \ref{t-dyson} was proved by \citet{wilson} and
\citet{gunson} independently. A similar proof was given in
\citep{ego}. Theses proofs use integrals of analytic functions. A
simple induction proof was found by \citet{good1}. We are going to
give a proof by using the Residue Theorem for Malcev-Neumann
series.

Let $\mb{z}$ be the vector $(z_1,z_2,\ldots, z_n)$. If $\mb{z}$
appears in the computation, we use $\mb{z}$ for the product
$\mb{z^1}=z_1z_2\cdots z_n$. We use similar notation for $\mb{u}$.

Let $\Delta(\mb{z})=\Delta (z_1,\ldots
,z_n)=\prod_{i<j}(z_i-z_j)=\det(z_i^{n-j})$ be the Vandermonde
determinant in $\mb{z}$, and let
$\Delta_j(\mb{z})=\Delta(z_1,\ldots ,\hat{z}_j,\ldots ,z_n)$,
where $\hat{z}_j$ means to omit $z_j$. We introduce new variables
$u_j=(-1)^{j-1}z_j^{n-1}\Delta_j(\mb{z})$. Then they satisfy the
equations
$$\Delta(\mb{z})=\sum_{j=1}^n (-1)^{j-1} z_j^{n-1}\Delta_j(\mb{z})=u_1+u_2+\cdots +u_n,$$
$$u_1\cdots u_n=\prod_{j=1}^n  (-1)^{j-1}z_j^{n-1}\Delta_j(\mb{z}) =(-1)^{\binom{n}{2}}
\mb{z}^{n-1} (\Delta(\mb{z}))^{n-2}.$$ We also have
$$\prod_{i=1,i\ne j}^n \left(1-\frac{z_i}{z_j}\right) =(-1)^{j-1}z_j^{n-1}
\frac{\Delta(\mb{z})}{\Delta_j(\mb{z})}=\frac{u_1+u_2+\cdots
+u_n}{u_j}.$$

Thus equation \eqref{e-dyson} is equivalent to
$$\CT_{\mb{z}} \frac{(u_1+u_2+\cdots +u_n)^{a_1+a_2+\cdots +a_n}}{u_1^{a_1}\cdots u_n^{a_n}}
=\frac{(a_1+a_2+\cdots a_n)!}{a_1!a_2!\cdots a_n!},$$ which is a
direct consequence of the multinomial theorem and the following
proposition.
\begin{prop}\label{p-dyson1}
For any series $\Phi(\mb{z})\in K^\mb{u}\ll \mb{z} \gg$, we have
$$\ct_{\mb{z}} \Phi(u_1,\dots ,u_n) =\ct_{\mb{u}} \Phi(u_1,\dots ,u_n).$$
\end{prop}

In fact, we can prove a more general formula. Let $r$ be an
integer and let $u^{(r)}_j = (-1)^{j-1}z_j^{r}\Delta_j(\mb{z})$.
Then $u^{(r)}_1+\cdots +u^{(r)}_n$ equals $ h_{r-n+1}
(z_1,z_2,\dots ,z_n)\Delta(\mb{z})$ for $r\ge n-1$ and equals $0$
for $0\le n\le n-2$. We have the following generalization.
\begin{thm}\label{t-dyson-g}
If $r$ is not equal to one of $0,1,\cdots ,n-2,$ or
$-\binom{n-1}{2}$, then for any series $\Phi(\mb{z})\in K^\rho\ll
\mb{z} \gg$, where $\rho(z_i)= u_i^{(r)}$, we have
$$\ct_{\mb{z}} \Phi(u^{(r)}_1,\dots ,u^{(r)}_n) =\ct_{\mb{u}^{(r)}} \Phi(u^{(r)}_1,\dots
,u^{(r)}_n).$$
\end{thm}
Note that Proposition \ref{p-dyson1} is the special case for
$r=n-1$ of Theorem \ref{t-dyson-g}. By Theorem \ref{t-MNresidue},
the above result is equivalent to saying that the log Jacobian is
a nonzero constant. To show this, we use the argument by
\citep{wilson}.

\begin{lem}\label{l-dyson-sym}
Let $G(x_1,\dots ,x_n)$ be a function of $n$ variables such that
\begin{enumerate}
\item $G$ is a symmetric function of $x_1,\dots ,x_n$.
\item $G$ is a ratio of two polynomials in the $x$'s.
\item $G$ is homogeneous of degree $0$ in the $x$'s.
\item The denominator of $G$ is $\Delta(x_1,\dots ,x_n)$.
\end{enumerate}
Then $G$ is a constant.
\end{lem}
\begin{proof}
Since the denominator of $G$ changes sign when the values of any
pairs $x_i,x_j$ are exchanged, the numerator must also change sign
under such an exchange. Thus, the numerator vanishes when
$x_i=x_j$. Hence the numerator has $x_i-x_j$ as a factor for any
$i$ and $j$, i.e., it has the entire denominator as a factor. So
$G$ is a polynomial. Together with the degree $0$ condition, $G$
must be a constant.
\end{proof}

\begin{proof}
In order to compute the log Jacobian, we let
$$J=\det(J_{ij})=\det \left( \frac{\partial \log u^{(r)}_i}{\partial \log z_j}\right).$$ Then
$J_{ii}=r$ and $J_{ij}=\sum_{k\ne i} \frac{z_i}{z_k-z_j}$ for
$i\ne j$. We first show that $J$ is a constant by Lemma
\ref{l-dyson-sym}. It is easy to see that $J$ satisfies the
conditions $1,2$ and $3$ in Lemma \ref{l-dyson-sym}. Now we show
that the denominator of $J$ is $\Delta(\mb{z})$, so that we can
claim that the Jacobian is a constant, and hence equals the
Jacobian number.

Evidently $J$ is the ratio of two polynomials in the $\mb{z}$'s,
whose denominator is a product of factors $z_i-z_j$ for some $i\ne
j$. From the expression of $J_{ij}$, we see that $z_i-z_j$ only
appears in the $i$th or the $j$th column. Every $2$ by $2$ minor
of the $i$th and $j$th  columns are of the following form, in
which we assume that $k$ and $l$ are not one of $i$ and $j$.
$$\left|\begin{array}{cc}
J_{ki} & J_{kj} \\
J_{li} & J_{lj}
\end{array} \right|
=\left|\begin{array}{cc} \frac{z_k}{z_j-z_i}+\sum_{s\ne i,j}
\frac{z_k}{z_s-z_i} & \frac{z_k}{z_i-z_j}+\sum_{s\ne i, j}
\frac{z_k}{z_s-z_j} \\
\frac{z_l}{z_j-z_i}+ \sum_{s\ne i,j} \frac{z_l}{z_s-z_i} &
\frac{z_l}{z_i-z_j} +\sum_{s\ne i, j} \frac{z_l}{z_s-z_j}
\end{array} \right|,
$$
in which the terms containing $(z_i-z_j)^2$ as the denominator
cancel. Therefore, expanding the determinant according to the
$i$th and  $j$th column, we see that $\Delta(\mb{z})$ is the
denominator of $J$.

Now the initial term of $z_i-z_j$ is $z_i$ if $i<j$. We see that
the initial term of $u^{(r)}_1$ is $z_1^rz_2^{n-2}z_3^{n-3}\cdots
z_{n-1}$. Similarly we can get the initial term for $u^{(r)}_j$.
The Jacobian number, denoted by $j(r)$, is thus the determinant
$$j(r)=\det \left(
\begin{array}{ccccc}
r & n-2 & n-3&  \cdots &0 \\
n-2 &r &  n-3&  \cdots & 0 \\
\vdots &\vdots &\vdots  & \vdots &\vdots \\
n-2 & n-3 & n-4 & \cdots & r
\end{array}
\right),
$$
where the displayed matrix has diagonal entries $r$, and other
entries in each row are $n-2,n-3,\dots , 0$, respectively from
left to right.

Since the row sum of each row is $r+\binom{n-1}{2}$,
$j(-\binom{n-1}{2})=0$. We claim that $j(r)=0$ when $r=0,1,\dots ,
n-2$. For in those cases, $u_1^{(r)}+\cdots +u_{n}^{(r)} =0$. This
implies that the Jacobian is $0$, and hence $j(r)=0$. We can
regard $j(r)$ as a polynomial in $r$ of degree $n$, and we have
already got $n$ zeros. So $j(r)=r(r-1)\cdots (r-n+2)
(r+\binom{n-1}{2}).$ up to a constant. This constant equals $1$
through comparing the leading coefficient of $r$.

In particular, $j(n-1)= \binom{n}{2} (n-1)!= \frac{n-1}{2} n!.$
Note that in \citep{ego}, the constant was said to be
$\frac{n-3}{2}n!$, which is wrong.
\end{proof}

Another proof of Dyson's conjecture by our residue theorem is to
use the change of variables by \citep{wilson}.

Let
$$v_j=\prod_{1\le i\le n, i\ne j} (1-z_j/z_i)^{-1}.$$
Then the initial term of $v_j$ is $z^{n-j}z_{j+1}\cdots z_n$ up to
a constant. Since the order of $v_n$ is $\mb{0}$, we have to
exclude $v_n$ from the change of variables, for otherwise, the
Jacobian number will be $0$. In fact, we have the relation
$v_1+v_2+\cdots +v_n=1$, which can be easily shown by Lemma
\ref{l-dyson-sym}.

Dyson's conjecture is equivalent to
\begin{align}\label{e-3-dyson-v}
\ct_{\mb{z}} \prod_{j=1}^n v_i^{-a_j} = \frac{(a_1+a_2+\cdots
a_n)!}{a_1!a_2!\cdots a_n!}
\end{align}

\begin{proof}[Another Proof of Dyson's Conjecture]
Using Lemma \ref{l-dyson-sym} and Wilson's argument, we can
evaluate the following log Jacobian. (Or see \citep{wilson} for
details.)
$$\frac{\partial (\log v_1,\log v_2 ,\dots ,\log v_{n-1})}{\partial(
\log z_1,\log z_2,\dots ,\log z_{n-1})}=(n-1)! v_n.$$

Then by the residue theorem
$$\ct_z \Phi(v_1,\dots ,v_{n-1},z_n) =\ct_{v_1,\dots,v_{n-1},z_n}
(1-v_1-\cdots -v_{n-1})^{-1} \Phi(v_1,\dots, v_{n-1},z_{n}).$$

In particular, (since the initial term of $1-v_1-\cdots -v_{n-1}$
is $1$), we have:
\begin{align*}
\ct_{\mb{z}} \prod_{j=1}^n v_i^{-a_j}
&=\ct_{v_1,\dots,v_{n-1},z_n}
(1-v_1-\cdots -v_{n-1})^{-a_n-1}\prod_{j=1}^{n-1} v_i^{-a_j}\\
&=[v_1^{a_1}\cdots v_{n-1}^{a_{n-1}}] \sum_{m\ge 0 }
\binom{a_n+m}{a_n}
(v_1+\dots +v_{n-1})^m\\
&=\binom{a_n+a_1+\cdots + a_{n-1}}{a_n} \binom{a_1+\cdots +
a_{n-1}}{a_1,\dots,a_{n-1}}.
\end{align*}
Equation \eqref{e-3-dyson-v} then follows.
\end{proof}

\section{About Morris's Identity\label{s-morris}}
We give a simplified proof of the following form of Morris's
identity \citep[Theorem 27, Corollary 28]{welleda}, which was
proved by using \emph{total residue}.

\begin{thm}If $k_1,k_2,k_3\in \NN$ and $k_1+k_2\ge 2$, then
\begin{multline}
 \ct_x \prod_{i=1}^r x_i^{-k_1+1}\prod_{i=1}^r (1-x_i)^{-k_2}
\prod_{i<j}(x_i-x_j)^{-k_3} \\
=\prod_{j=0}^{r-1}
\frac{\Gamma(1+\frac{k_3}{2})\Gamma(k_1+k_2-1+(r+j-1)\frac{k_3}{2})}{
\Gamma(1+(j+1)\frac{k_3}{2})\Gamma(k_1+j\frac{k_3}{2})\Gamma(k_2+j\frac{k_3}{2})}.
\end{multline}
\end{thm}

Let $P_{l,r}$ be the symmetric function defined by
$$P_{l,r}=\sum_{w\in \sy_r} w\cdot (x_1 x_2\cdots x_l)=l!\,(r-l)!\, e_l,$$
where $\sy_r$ is the symmetric group on $1,2,\dots ,r$ and $w$
acts by permuting the indexes of the $x$'s, and $e_l$ is the
elementary symmetric function. In particular,
$$P_{0,r}=r! \qquad P_{r,r}= r!\, x_1x_2\cdots x_r.$$
When $r$ is fixed, we write $P_l$ for $P_{l,r}$. Let
$$\phi_r(l,k_1,k_2,k_3)=
\frac{P_l}{\prod_{i=1}^r x_i^{k_1-1} \prod_{i=1}^r (1-x_i)^{k_2}
\prod_{i<j}(x_i-x_j)^{k_3}},$$ where $k_1,k_2,$ and $k_3$ are
nonnegative integers. If $k_3$ is odd, this function is
anti-symmetric in $x_1,\dots ,x_r$. If $k_3$ is even, this
function is symmetric.

Now let $C_r(l,k_1,k_2,k_3)$ be the constant term of
$\phi_r(l,k_1,k_2,k_3)$.

The following is in \citep[Theorem 27]{welleda}.
\begin{thm}
Let $k_1,k_2,k_3\ge 0,$ $0\le l\le r$. The constants
$C_r(l,k_1,k_2,k_3)$ are uniquely determined by the relations:
\begin{enumerate}
\item $C_r(r,k_1,k_2,k_3)=r! \, C_r(0,k_1-1,k_2,k_3).$
\item $C_r(r-1,1,k_2,k_3)=C_{r-1}(0,k_3,k_2,k_3).$
\item $C_r(0,1,k_2,0)=r!$.
\item $C_1(l,0,k_2,k_3)=0.$
\item For $1\le l \le r$,
\begin{multline*}
\left(k_1+k_2-2+{k_3\over 2}(2r-l-1)\right) C_r(l,k_1,k_2,k_3)\\
=\left(k_1-1+{k_3\over 2}(r-l)\right) C_r(l-1,k_1,k_2,k_3).
\end{multline*}
\end{enumerate}
\end{thm}
We only prove these five relations, which in fact give a recursive
formula for $C_r(l,k_1,k_2,k_3)$. Note that the proof of relation
2 by using total residue was lengthy in \citep{welleda}.
\begin{proof}
Relation $1$ follows directly from the definition.

Now $\prod_{i=1}^r (1-x_i)^{-k_2}  \prod_{i<j}(x_i-x_2)^{-k_3}$ is
always a power series in $x_r$. So if $k_1=0$ or $k_1=1$, then we
can get the constant term in $x_r$ by setting $x_r=0$. When
$k_1=0$, we get relation $4$. When $k_1=1$, we get
\begin{multline*}
\ct_x \frac{\sum_{w\in \sy_r}w \cdot ( x_1\cdots
x_{r-1})}{\prod_{i=1}^r (1-x_i)^{k_2}
\prod_{i<j}(x_i-x_2)^{k_3}}\\
=\ct_{x_1,\dots ,x_{r-1}} \frac{\sum_{w\in \sy_{r-1}}w\cdot(
x_1\cdots x_{r-1})}{\prod_{i=1}^{r-1} x_i^{k_3}\prod_{i=1}^{r-1}
(1-x_i)^{k_2} \prod_{i<j}(x_i-x_j)^{k_3}}.
\end{multline*}
This implies relation $2$.

Relation $3$ is equivalent to
$$\ct_x \prod_{i=1}^r \frac{1}{(1-x_i)^{k_2}}=1,$$
which is obvious.

Now we show relation $5$. Let $l>0$, and let
$$U= \frac{1}{\prod_{i=1}^r x_i^{k_1} \prod_{i=1}^r (1-x_i)^{k_2}
\prod_{i<j}(x_i-x_j)^{k_3}}.$$ Then $C_r(l,k_1,k_2,k_3)=\res_{x}
P_l U$. We have
\begin{align*}
&\frac{\partial}{\partial x_1} (1-x_1) x_1 x_2\cdots x_l U\\
=& (k_2-1) x_1\cdots x_l U+ (1-k_1)(1-x_1)x_2\cdots x_l
U-k_3(1-x_1)x_1\cdots x_l \sum_{j=2}^r
\frac{U}{x_1-x_j} \\
=& (k_1+k_2-2) x_1\cdots x_l U+(1-k_1) x_2\cdots x_l
U-k_3(1-x_1)x_1\cdots x_l \sum_{j=2}^r \frac{U}{x_1-x_j}
\end{align*}
If $k_3$ is odd, then $U$ is antisymmetric. Anti-symmetrizing over
$\sy_r$, we get
\begin{align*}
&\sum_{w\in \sy_r} (-1)^w w\cdot \left(\frac{\partial}{\partial
x_1} (1-x_1) x_1 x_2\cdots x_l U
\right) \\
&= (k_1+k_2-2)P_l U+(1-k_1) P_{l-1}U -k_3 \sum_{w\in \sy_r} w\cdot
(1-x_1)x_1\cdots x_l \sum_{j=2}^r \frac{U}{x_1-x_j} .
\end{align*}
To  compute
$$\sum_{w\in \sy_r} w\cdot (1-x_1)x_1\cdots x_l \sum_{j=2}^r
\frac{1}{x_1-x_j} U ,$$ we first sum over the transpositions
$(j,1)$. For $2\le j\le l$, we use the formula
$$\frac{(1-x_1) x_1x_j}{x_1-x_j} +\frac{(1-x_j)x_1x_j}{x_j-x_1} =-x_1x_j.$$
For $j>l$, we use the formula
$$\frac{(1-x_1)x_1}{x_1-x_j}+\frac{(1-x_j)x_j}{x_j-x_1} =1-x_1-x_j.$$
We obtain that
\begin{align*}
2 \sum_{w\in \sy_r} w\cdot (1-x_1)x_1\cdots x_l \sum_{j=2}^r
\frac{1}{x_1-x_j} = (-(l-1)-2(r-l)) P_l U+(r-l) P_{l-1}.
\end{align*}
Thus finally we get
 relation $5$ when $k_3$ is odd.

If $k_3$ is even, then $U$ is symmetric. Symmetrizing over
$\sy_r$, we get
\begin{align*}
&\sum_{w\in \sy_r} w\cdot \left(\frac{\partial}{\partial x_1}
(1-x_1) x_1 x_2\cdots x_l U
\right) \\
&= (k_1+k_2-2)P_l U+(1-k_1) P_{l-1}U -k_3 w\cdot (1-x_1)x_1\cdots
x_l \sum_{j=2}^r \frac{1}{x_1-x_j} U.
\end{align*}
The rest of the proof of relation 5 for $k_3$ even proceeds as in
the case of $k_3$ odd.
\end{proof}

\section{MacMahon's Partition Analysis Revisited\label{s-Macr}}

In section \ref{s-Mac}, we discussed the algorithmic aspect of
MacMahon's partition analysis. In this section, we shall discuss
the theoretical aspect. Some work was done in \citep[p.
229--231]{stanley-rec} by using residue computations and partial
fraction decompositions. We are going to work in a field of
MN-series. The foundation of this part is Theorem \ref{t-elliott},
which says that the constant term of an Elliott-rational function
is still Elliott-rational. This statement is true for any field of
MN-series. Our goal in this section is to give new proof of the
reciprocity theorem for a system of homogeneous linear Diophantine
equations. See Theorem \ref{t-3-recip} below.

First, we shall clarify the notation. Let $\rho$ be an injective
endomorphism of $\ZZ^{r+n}$, or more generally a total ordering on
the group of monomials that is compatible with its group
structure. We use $\Lambda$ to denote the vector $(\lambda_1,\dots
,\lambda_r)$ and $\mb{x}$ to denote the vector $(x_1,\dots ,x_n)$.
Then $\CC^\rho \ll \Lambda, \mb{x} \gg$ is a field of MN-series.
The field $\CC(\Lambda,\mb{x})$ of rational functions can be
embedded into $\CC^\rho \ll \Lambda, \mb{x} \gg$, and any rational
function $F(\Lambda, \mb{x})$ has a unique expansion in
$\CC^\rho\ll \Lambda, \mb{x} \gg$.

It is convenient for our purposes to denote by $K$ the field
$\CC(\mb{x})$. The field of rational functions
$\CC(\Lambda,\mb{x})$ can be identified with $K(\Lambda)$. Usually
we are taking constant terms in the ring
$\CC[\Lambda,\Lambda^{-1}][[\mb{x}]]$, where $\Lambda^{-1}$ refers
to $(\lambda_1^{-1},\dots ,\lambda_r^{-1})$. This ring can be
embedded into $\CC\ll \Lambda, \mb{x} \gg$, where $\rho$ is
omitted since it is the identity map.

\subsection{The Case of $r=1$}

In this case, we need not restrict ourselves to Elliott-rational
functions. Thus we need to consider the following problem.

\vspace{3mm} \noindent \emph{Problem:} Given
 a rational function $Q(\lambda)$ (short for $Q(\lambda,\mb{x})$)
 of $\lambda$ and $\mb{x}$,
compute $\PT^\rho _\lambda Q(\lambda,\mb{x})$. Recall that
$\PT^\rho_\lambda$ indicates that $Q(\lambda,\mb{x})$ is treated
as an element of $\CC^\rho \ll \lambda,\mb{x} \gg$.

To deal with this problem, we shall understand that $Q(\lambda)$
is not only an element of $K(\lambda)$, but also an element of
$\CC^\rho\ll \lambda,\mb{x}\gg$. As an element of $K(\lambda)$,
$Q(\lambda)$ can be written as $p(\lambda)/q(\lambda)$, where
$p(\lambda)$ and $q(\lambda)$ are both in $K[\lambda]$. As an
element of $\CC^\rho\ll \lambda,\mb{x}\gg$, the denominator
$q(\lambda)$ plays a central role.

Recall that $\CC^\rho\ll \lambda,\mb{x}\gg$ is equipped with an
operator $\ord^\rho$ and a total ordering on its monomials. Let us
write $q(\lambda)=\sum_{i=0}^{d} a_i \lambda^i$, with $a_i\in
\CC(\mb{x})$ and $a_d\ne 0$. To expand $Q(\lambda)$ into a series
in $\CC^\rho \ll \lambda,\mb{x} \gg$, we need to find the
$\lambda$-initial term $a_j\lambda^j$, or equivalently, the $j$
such that $\ord^\rho (a_j \lambda^j)$ is smaller than $\ord^\rho
(a_i\lambda^i)$ for all $i\ne j$. This can be achieved because of
the different powers in $\lambda$. Then
$$\frac{1}{q(\lambda)} =\frac1{a_j\lambda^j} \frac{1}{1+\sum_{i\ne j}
a_i/a_j\lambda^{i-j} } =\frac1{a_j\lambda^j} \sum_{k\ge 0} (-1)^k
\Big(\sum_{i\ne j} a_i/a_j\lambda^{i-j}\Big)^k.$$

It is now clear that we have the following three situations.
\begin{enumerate}
\item If $j$ equals $0$, then for any polynomial $p(\lambda)$,
$p(\lambda)/q(\lambda)$ contains only nonnegative powers in
$\lambda$. In this case, we say that $1/q(\lambda)$ is $\PT^\rho$
in $\lambda$.

\item If $j$ equals $d$, then for any polynomial $p(\lambda)$ of degree
in $\lambda$ less than $d$, $p(\lambda)/q(\lambda)$ contains only
negative powers in $\lambda$. In this case, we say that
$1/q(\lambda)$ is $\NT^\rho$ in $\lambda$.

\item If $j$ equals neither $0$, nor $d$, then
 $1/q(\lambda)$ contains both positive and
negative powers in $\lambda$. Thus $1/q(\lambda)$ is neither
$\PT^\rho$ nor  $\NT^\rho$ in $\lambda$.
\end{enumerate}

\begin{lem}
Let $q_1$ and $q_2$ be polynomials. Then for any fixed $\rho$
\begin{itemize}
    \item Both $1/q_1(\lambda)$ and
$1/q_2(\lambda)$ are $\pt^\rho$ in $\lambda$ if and only if $1/(q
_1q_2)$ is.
    \item Both $1/q_1(\lambda)$ and
$1/q_2(\lambda)$ are $\nt^\rho$ in $\lambda$ if and only if $1/(q
_1q_2)$ is.
    \item For all the other cases, $1/(q_1q_2)$ is neither $\pt^\rho $
    in $\lambda$ nor $\nt^\rho$ in $\lambda$.
\end{itemize}
\end{lem}
\begin{proof}
We prove the first case for $\pt$ as follows. The other cases are
similar. Write $$q_1=\sum_{i=0}^{d_1} a_i \lambda^i, \quad
q_2=\sum_{i=0}^{d_2} b_i \lambda^i, \quad \text{ and
}q_1q_2=\sum_{i=0}^{d_1+d_2} c_i \lambda^i.$$ Suppose that
$a_{j_1}\lambda^{j_1}$ and $b_{j_2}\lambda^{j_2}$ are the
$\lambda$-initial term of $q_1$ and $q_2$ respectively. Now if we
expand the product $q_1q_2$ but do not collect terms, then
$a_{j_1}b_{j_2}\lambda^{j_1+j_2}$ is the unique term with the
least order. So the order of $c_{j_1+j_2}\lambda^{j_1+j_2}$ has to
equal the order of $a_{j_1}b_{j_2}\lambda^{j_1+j_2}$. This implies
that the $\lambda$-initial term of $q_1q_2$ is
$c_{j_1+j_2}\lambda^{j_1+j_2}$. The assertion for $\pt$ in the
lemma hence follows from the fact that $j_1+j_2=0\Leftrightarrow
j_1=0 \text{ and } j_2=0$. (Remember that $j_1,j_2\ge 0$).
\end{proof}
 A direct consequence of the above lemma is the following corollary.
 \begin{cor}
If $1/q_1(\lambda)$ is $\PT^\rho$ in $\lambda$ and
$1/q_2(\lambda)$ is $\NT^\rho$ in $\lambda$, then $q_1(\lambda)$
and $q_2(\lambda)$ cannot have a nontrivial common divisor in
$K[\lambda]$, i.e., they are relatively prime.
\end{cor}

\begin{dfn}
If $q(\lambda)$ can be factored as $q_1(\lambda)q_2(\lambda)$ such
that $1/q_1(\lambda)$ is $\NT^\rho$ in $\lambda$ and
$1/q_2(\lambda)$ is $\PT^\rho$ in $\lambda$, then we say that
$q(\lambda)$ is $\rho$-{\rm factorable}, and
$q(\lambda)=q_1(\lambda)q_2(\lambda)$ is a $\rho$-{\rm
factorization}. Such factorization is unique (if it exists) up to
a constant in $K$.
\end{dfn}

\begin{thm}
Let $p(\lambda), q(\lambda)\in K[\lambda]$. If $q(\lambda)$ is
$\rho$-factorable, then $\CT^\rho_\lambda p(\lambda)/q(\lambda)$
is in $K$, i.e., is rational.
\end{thm}
\begin{proof}
Suppose $q(\lambda)=q_1(\lambda)q_2(\lambda)$ is such a
$\rho$-factorization. Since  $1/q_1(\lambda)$ is $\PT^\rho $ in
$\lambda$ and $1/q_2(\lambda)$ is $\NT^\rho$ in $\lambda$,
$q_1(\lambda)$ and $q_2(\lambda)$ are relatively prime in
$K[\lambda]$. Thus we have the unique partial fraction expansion
in $K(\lambda)$:
\begin{align} \label{e-3-frac-r}
\frac{p(\lambda)}{q(\lambda)}
=p_0(\lambda)+\frac{p_1(\lambda)}{q_1(\lambda)}
+\frac{p_2(\lambda)}{q_2(\lambda)},
\end{align}
where $p_i$ are polynomials in $\lambda$ for $i=0,1,2$ and
 $\deg p_i(\lambda)< \deg q_i(\lambda)$ for
$i=1,2$. Since when expanded as series in $\CC^\rho\ll
\lambda,\mb{x}\gg$, $p_1(\lambda)/q_1(\lambda)$ contains only
negative powers in $\lambda$, and $p_0$ and
$p_2(\lambda)/q_2(\lambda)$ contains only nonnegative powers in
$\lambda$, we have
$$\pt_\lambda \uprho \frac{p(\lambda)}{q(\lambda)} =
p_0(\lambda)+\frac{p_2(\lambda)}{q_2(\lambda)}.$$ Thus
$\ct^\rho_\lambda=p_0(0)+p_2(0)/q_2(0)$ is in $\CC(\mb{x})$.
\end{proof}

This result clearly implies Theorem \ref{t-1-hadamard}.

\begin{cor}\label{c-mac-rec}
Suppose the degree of $p(\lambda)$ is less than the degree of
$q_1(\lambda)q_2(\lambda)$, and $p(0)=0$. If $1/q_1(\lambda)$ is
$\NT^\rho$ in $\lambda$, but is $\PT^\sigma$ in $\lambda$, and
$1/q_2(\lambda)$ is $\PT^\rho$ in $\lambda$ but is $\NT^\sigma$ in
$\lambda$, then
$$\ct_\lambda \uprho \frac{p(\lambda)}{q_1(\lambda)q_2(\lambda)}
=-\ct_\lambda \mbox{}^\sigma \,
\frac{p(\lambda)}{q_1(\lambda)q_2(\lambda)},$$ where the equation
is regarded as an element of $K$.
\end{cor}
\begin{proof}
From the hypothesis, it is easy to see that for $i=1$ or $2$
$q_i(0)\ne 0$, and $q_i(\lambda)$ can not be of degree $0$. Thus
the corollary follows from equation \eqref{e-3-frac-r} by setting
$\lambda =0$.
\end{proof}

As an element of $K[\lambda]$, $q(\lambda)$ can be factored into
the product of irreducible polynomials. Let
$q(\lambda)=q_1(\lambda)\cdots q_k(\lambda)$ be such a
factorization. Then $q(\lambda)$ is $\rho$-factorable if and only
if every $1/q_i$ is either $\PT^\rho$ or $\NT^\rho$. When this is
true, the $\rho$-factorization can be obtained by collecting
similar terms.

All Elliott-rational functions are $\rho$-factorable for any
$\rho$. For in such a function, the denominator is a product of
the form $\lambda^j-a$, where $a\in K$ and $j$ is a positive
integer. Thus for any $\rho$, $1/(\lambda^j-a)$ is either
$\NT^\rho$ or $\PT^\rho$ in $\lambda$.

More precisely, any Elliott-rational function $F$ can be written
as follows:
\begin{equation}\label{e-2-Mac-F2}
F=\frac{p(\lambda)}{(\lambda^{j_1}-a_1)\cdots (\lambda^{j_n}-a_n)
(\lambda^{k_1}-b_1)\cdots (\lambda^{k_m}-b_m)},
\end{equation}
where $p(\lambda)$ is a polynomial of $\lambda$, $j_i$ and $k_i$
are positive integers, $m$ and $n$ are nonnegative integers, and
$a_i, b_l\in K$. For a particular $\rho$, we require that
$1/(\lambda^{j_i}-a_i)$ is $\NT^\rho$ in $\lambda$, and
$1/(\lambda^{k_i}-b_i)$ is $\pt^\rho$ in $\lambda$. Note that
$a_1$ can be $0$. The conclusion is that
 $\ct^\rho_\lambda F$ is always rational.

For any  total ordering $\rho$ on the monomials of $K(\lambda)$,
we let $\hat{\rho}$ be the total ordering such that
$\hat{\rho}(m_1)\le \hat{\rho} (m_2)$ if and only if $\rho(m_1)\ge
\rho(m_2)$ for all monomials $m_1,m_2$. Then we have a sort of
reciprocity formula.
\begin{cor}
Let $F(\lambda)$ be of the form \eqref{e-2-Mac-F2}. If $F(0)=0$,
and $F(\lambda)$ is a proper rational function in $\lambda$, then
for any $\rho$, we have the reciprocity
$$\ct_\lambda \uprho F(\lambda) =-\ct_\lambda \mbox{}^{\hat{\rho}} \, F(\lambda),$$
where both sides are regarded as elements in $K$.
\end{cor}

\subsection{The General Case}
MacMahon's partition analysis can be applied to solve a system of
linear Diophantine equations or inequalities. It is well-known
that inequalities can be replaced with equations by introducing
new variables.

Solving  linear Diophantine equations means finding all vectors
$\alpha\in \NN^n$ that satisfy $A\alpha=0$, where $A$ is an $r$ by
$n$ matrix with integral entries. More precisely, we want to solve
the following system of equations:
\begin{align*}
a_{1,1}\alpha_1 +a_{1,2} \alpha_2+\cdots +a_{1,n} \alpha_n &=0 \\
a_{2,1}\alpha_1 +a_{2,2} \alpha_2+\cdots +a_{2,n} \alpha_n &=0 \\
\cdots \cdots \qquad &=0\\
a_{r,1}\alpha_1 +a_{r,2} \alpha_2+\cdots +a_{r,n} \alpha_n &=0.
\end{align*}

Let $C_i$ be the $i$th column vector of $A$. Then the above system
is the same as
$$C_1 \alpha_1+C_2\alpha_2+\cdots +C_n\alpha_n=0.$$

Now let $E$ and $\bar{E}$ be the sets of all such solutions in
$\NN^n$ and $\PP^n$ respectively. It is natural to study the
generating functions of $E$ and $\bar{E}$:
\begin{align}
E(\mb{x})&=E(x_1,\dots ,x_n) =\sum_{\alpha\in E}
\mb{x}^\alpha,\quad
    \bar{E}(\mb{x})= \bar{E}(x_1,\dots ,x_n) =\sum_{\alpha\in \bar{E}} \mb{x}^\alpha
\end{align}
where if $\alpha=(\alpha_1,\dots ,\alpha_n)$, then
$\mb{x}^\alpha:= x_1^{\alpha_1}\cdots x_n^{\alpha_n}.$

Using MacMahon's partition analysis, we can realize the the $r$
linear constraints by introducing $\lambda_1, \lambda_2,\dots
,\lambda_r$ and then taking the constant terms. We have
\begin{align}
E(x) &= \ct_\Lambda \sum_{\alpha\in \NN^n} \lambda_1^{a_{1,1}
\alpha_1+\cdots +a_{1,n}\alpha_n}\cdots
\lambda_r^{a_{r,1} \alpha_1+\cdots +a_{r,n}\alpha_n} x^\alpha \nonumber\\
&= \ct_\Lambda \prod_{i=1}^n \frac{1}{1-
\lambda_1^{a_{1,i}}\lambda_2^{a_{2,i}}\cdots \lambda_r^{a_{r,i}}
x_i}=\ct_\Lambda \prod_{i=1}^n \frac{1}{1-\Lambda^{C_i}x_i},
\label{e-mac-Ex}
\end{align}
where we are working in $\CC[\Lambda,\Lambda^{-1}][[\mb{x}]]$,
which can be embedded into $\CC\ll \Lambda,\mb{x}\gg$. Similarly:
\begin{align}
\label{e-mac-Ebx} \bar{E}(x)=\ct_\Lambda \prod_{i=1}^n
\frac{\Lambda^{C_i}x_i}{1- \Lambda^{C_i}x_i}.
\end{align}

The well-known reciprocity theorem  \citep[]{EC1} for homogeneous
linear diophantine equations is the following:
\begin{thm}[Reciprocity Theorem]\label{t-3-recip}
Let $E$ and $A$ be as above. If the rank of $A$ is $r$, and
$\bar{E}$ is nonempty, then as rational functions
\begin{align}\label{e-3-recip}
    E(\mb{x})= (-1)^{n-r} \bar{E}(\mb{x^{-1}}).
\end{align}
\end{thm}
Previous proofs of this theorem use simplex decompositions, but we
want to give a proof directly from \eqref{e-mac-Ex} and
\eqref{e-mac-Ebx}. We will use the Elliott reduction identity to
derive this result.

We shall see that all of the work is done algebraically. First,
let us see some facts. Exchanging column $i$ and $j$ corresponds
to exchanging $x_i$ and $x_j$. Row operations, which will not
change the solutions of $A\alpha=0$, are equivalent to multiplying
$A$ on the left by an invertible matrix. This fact can be obtained
by applying the residue theorem. (In fact, Corollary
\ref{c-3-monomial}.)

 We
define $\mathcal{E}(\mb{x})$ to be the crucial generating function
of $E(x)$:
\begin{align}
    \mathcal{E}(\mb{x})
    =
    \prod_{i=1}^n \frac{1}{1- \Lambda^{C_i}
    x_i}.
\end{align}
The crucial generating function $\bar{\mathcal{E}}(\mb{x})$ of
$\bar{E}(\mb{x})$ is defined similarly. Now apply the residue
theorem, in fact Corollary \ref{c-3-monomial}, by changing
variables in equation \eqref{e-mac-Ebx} by $\Lambda\to
\Lambda^{-1}$, i.e., $\lambda_i\to \lambda_i^{-1}$ for all $i$.
Then we have:
\begin{align*}
\bar{E}(\mb{x})=\ct_{\Lambda^{-1}} \prod_{i=1}^n
\frac{\Lambda^{-C_i}x_i}{1- \Lambda^{-C_i} x_i}=\ct_{\Lambda^{-1}}
\prod_{i=1}^n \frac{1}{\Lambda^{C_i} x^{-1}_i-1}.
\end{align*}
Now we are taking the constant term of an element in the ring
$\CC[ \Lambda^{-1},\Lambda][[\mb{x}]]$, which is the same as the
original ring. So we have:
\begin{align}
\label{e-mac-Ebx1} \bar{E}(\mb{x})=\ct_\Lambda \prod_{i=1}^n
\frac{1}{ \Lambda^{C_i} x^{-1}_i-1}=(-1)^n\ct_\Lambda
\mathcal{E}(\mb{x^{-1}}).
\end{align}
Note that in the denominator of the right side of
\eqref{e-mac-Ebx1}, $1$ is not the initial term.

Now if we replace $\mb{x}$ with $\mb{\xx}$, then we have shown
that the reciprocity theorem is the full rank case of the
following proposition:

\begin{prop}\label{p-3-recip}
Suppose that $\bar{E}$ is nonempty. Then
\begin{align}\label{e-3-recip1}
    \ct_\Lambda  \mathcal{E}(x) =(-1)^{\mathrm{rank}(A)}\ct_\Lambda \uprho \mathcal{E}(x),
    \end{align}
    where $\rho$ is the endomorphism defined by $\rho(x_i)=x_i^{-1}$ and
    $\rho(\lambda_i)=\lambda_i$.
\end{prop}
On the other hand, it is easy to deal with the case of $\rank(A)<
r$. So the reciprocity Theorem \ref{t-3-recip} is equivalent to
Proposition \ref{p-3-recip}.

Before we give the proof of this proposition, let us see the
simple case of
 $r=1$. In this case, $\mathcal{E}(\mb{x}) $ has the form:
$$ \mathcal{E}(\mb{x})=\prod_{i=1} ^n \frac{1}{1-\lambda^{a_i}x_i}.$$
The condition that $\bar{E} $ is nonempty is equivalent to saying
that some of $a_i$ have to be positive and some of $a_i$ have to
be negative. Thus when written in the normal form of a rational
function in $\lambda$, $\mathcal{E}(\mb{x})$ is proper and its
numerator divides $\lambda$. So Proposition \ref{p-3-recip}
follows from Corollary \ref{c-mac-rec}.

The general case does not seem to work along this line because of
two problems. One is how to use the conditions that $\bar{E}$ is
nonempty, and the other is how to connect to the rank of $A$. The
proof we are going to give uses induction and Elliott's reduction
identity.

Clearly if $a_{11},\dots ,a_{1,n}$ are all positive or are all
negative, then $\bar{E}$ is empty. So we can assume that
$a_{11}>0$ and $ a_{12}<0$. Applying Elliott's reduction identity
on $\lambda_1$, we get:
\begin{align*}
\mathcal{E}(\mb{x})&= \frac{1}{1-\Lambda^{C_1+C_2}x_1x_2}
\left(\frac{1}{1-\Lambda^{C_1}x_1}+\frac{1}{1-\Lambda^{C_2}x_2}-1
\right) \prod_{i\ge 3}\frac{1}{1-\Lambda^{C_i}x_i}
\end{align*}
Now expand $\mathcal{E}(\mb{x})$ according to the middle term, and
denote the resulting three summans by $\mathcal{E}_1$,
$\mathcal{E}_2$, and $\mathcal{E}_3$ respectively. We have
\begin{align}\label{e-3-E123}
\mathcal{E}(\mb{x})=\mathcal{E}_1(x_1,x_1x_2,x_3,\dots )+
\mathcal{E}_2(x_1x_2,x_2,x_3,\dots)-
\mathcal{E}_3(x_1x_2,x_3,\dots ).
\end{align}
Then these $\mathcal{E}_i$ are very similar to $\mathcal{E}$.
Correspondingly, they are associated to matrices, and hence
solution spaces that lie in $\NN^n$ and $\PP^n$. More precisely,
$\mathcal{E}_i$, $i=1,2,3$, are associated to $A_1=
(C_1,C_1+C_2,C_3,\dots,C_n)$, $A_2=(C_1+C_2,C_2,C_3,\dots,C_n)$,
and $A_3=(C_1+C_2,C_3,\dots,C_n)$ respectively. Thus $E_i,
E_i(\mb{x})$ and $\bar{E}_i,\bar{E}_i(\mb{x})$ are defined
correspondingly.

Now the matrix $A_1$ is obtained from $A$ by adding the second
column to the first; the matrix $A_2$ is obtained from $A$ by
adding the first column to the second. They are obtained from $A$
through a column operation. So the rank of $A_1$ and $A_2$ are
both equal to that of $A$. The rank of $A_3$ might not equal the
rank of $A$.

Applying $\ct_\Lambda$ and $(-1)^n\ct_\Lambda^{\rho}$ to
\ref{e-3-E123} respectively, we get our key induction equations.
\begin{align}
E(\mb{x})&=E_1(x_1,x_1x_2,x_3,\dots )+ E_2(x_1x_2,x_2,x_3,\dots)-
E_3(x_1x_2,x_3,\dots ),\label{e-3-Ex123}\\
\bar{E}(\mb{x})&=\bar{E}_1(x_1,x_1x_2,x_3,\dots )+
\bar{E}_2(x_1x_2,x_2,x_3,\dots)\nonumber \\
&\qquad \qquad \qquad \qquad \qquad \quad +
(-1)^{\rank(A)-\rank(A_3)}\bar{E}_3(x_1x_2,x_3,\dots
).\label{e-3-Ebx123}
\end{align}

Looking more closely at these $E_i$, we can see that up to
isomorphism, $E_1$, $E_2$, and $E_3$ are obtained from $E$ by
intersecting the half spaces $\alpha_1\ge \alpha_2$, $\alpha_1\le
\alpha_2$, and the hyperplane $\alpha_1=\alpha_2$ respectively.
For instance, $(\alpha_1,\alpha_2,\dots,)$ belongs to $E$ with
$\alpha_1\ge \alpha_2$ if and only if
$(\alpha_1-\alpha_2,\alpha_2,\dots)$ belongs to $E_1$. Thus
Elliott's reduction identity in fact corresponds to a signed
decomposition of $E$. Equation \eqref{e-3-Ex123} and
\eqref{e-3-Ebx123} could be explained directly from geometry.

We need two more lemmas to give our proof of Proposition
\ref{p-3-recip}. We shall see that the condition on $\bar{E}$
plays an important role.

If $\bar{E}$ is nonempty, then $\dim E=\dim \bar{E}=n-\rank(A)$.
Clearly, the dimension of the solution space of $A\alpha=0$ is
$n-\rank(A)$. Let $\gamma\in \bar{E}$, and let $\Upsilon_1,\dots
,\Upsilon_{n-\rank(A)}$ be a $\ZZ$-basis of the solution space in
$\ZZ^n$ with $\Upsilon_1=\gamma$. Then for sufficiently large $m$,
$m \gamma +\Upsilon_1,\dots ,m\gamma+\Upsilon_{n-\rank(A)}$ will
be a linearly independent set in $\bar{E}$.

\begin{lem}\label{l-3-empty-23}
Suppose that $\bar{E}$ is nonempty, and that $\bar{E}_i$ is
defined as above for $i=1,2,3$. Then any two of the $\bar{E}_i$
being nonempty implies that they are all nonempty.
\end{lem}
\begin{proof}
Suppose that $\bar{E}_1$ and $\bar{E}_2$ are nonempty. Then we
have elements $\beta$ and $\gamma$ in $\bar{E}$ such that
$\beta=(\beta_1,\beta_2,\dots)$ with $\beta_1>\beta_2$ and
$\gamma=(\gamma_1,\gamma_2,\dots)$ with $\gamma_1<\gamma_2$. Then
$(\gamma_2-\gamma_1)\beta+(\beta_1-\beta_2)\gamma$ is in $\bar{E}$
with the first two entries being equal. This means $\bar{E}_3$ is
nonempty.

Suppose that $\bar{E_1}$ and $\bar{E}_3$ are nonempty. Then we
have elements $\beta$ and $\delta$ in $\bar{E}$ such that
$\beta=(\beta_1,\beta_2,\dots)$ with $\beta_1>\beta_2$ and
$\delta=(\delta_1,\delta_2,\dots)$ with $\delta_1=\delta_2$. Then
for sufficiently large $m$, $m\delta -\beta$ is in $\bar{E}$ with
the first entry being smaller than the second. This means
$\bar{E}_2$ is nonempty.

The case that $\bar{E_2}$ and $\bar{E}_3$ are nonempty is similar
to the previous case.
\end{proof}

\begin{lem}\label{l-3-rankall}
If all of the $\bar{E_i}$ are nonempty, then $\rank (A_3)=\rank
(A)$.
\end{lem}
\begin{proof}
By hypothesis, it is clear that $E$ is not contained in the
hyperplane $\alpha_1=\alpha_2$. Thus the intersection of $E$ with
the hyperplane has dimension $\dim E -1$. So $\dim E_3$ is also
$\dim E-1$ and the rank of $A_3$ equals $n-1-\dim E_3=\rank (A)$.
\end{proof}

\begin{proof}[Proof of Proposition \ref{p-3-recip}]
The base case, when $A$ is the zero matrix, is trivial.

By exchanging rows, we can assume that not all of the entries in
the first row of $A$ are zero. Moreover, since the entries can not
be all positive or negative, we can assume the first entry is
positive and the second is negative by exchanging columns.

We use induction on $S_1(A)$, which is defined to be the sum of
the absolute values of all the entries in the first row. Now the
above argument applies, and it is easy to see that
$S_1(A_i)<S_1(A)$ for $i=1,2,3$. Applying Lemma
\ref{l-3-empty-23}, we can reduce the seven cases of $E_i$ being
nonempty or not into the following four cases:

Case 1: only $\bar{E}_1$ is nonempty. Let $\beta$ in $\bar{E}$ be
such that $\beta_1>\beta_2$. We claim that all $\alpha$ with
$A\alpha=0$ satisfy the condition $\alpha_1
>\alpha_2$, so that $E_2(x_1x_2,x_2,x_3,\dots)$ equals $E_3(x_1x_2,x_3,\dots)$, and hence by induction we have
\begin{align*}
E(\mb{x})&=E_1(x_1,x_1x_2,x_3, \dots)\\
&= (-1)^{\rank(n-A_1)} \bar{E}_1(\xx _1,\xx _1 \xx_2, \xx _3,\dots)\\
&=(-1)^{n-\rank(A)}\bar{E}(\mb{x^{-1}}).
\end{align*}

If the claim does not hold, then  $\alpha_1\le \alpha_2$. But for
sufficiently large $m$, $m\beta -\alpha$ will produce an element
in $\bar{E}_2$  or $\bar{E}_3$, a contradiction.

Case 2: only $\bar{E}_2$ is nonempty. This is similar to case 1.

Case 3: only $\bar{E}_3$ is nonempty. This means that $E$ is
contained in the hyperplane $\alpha_1=\alpha_2$. Thus
$$E_1(x_1,x_1x_2,x_3,
\dots)=E_2(x_1x_2,x_2,x_3, \dots)=E_3(x_1x_2,x_3,\dots),$$ and we
have
$$\rank(A_3)=n-1-\dim(E_3)=n-\dim(E)-1=\rank(A)-1.$$
So
\begin{align*}
E(\mb{x})&=E_3(x_1x_2,x_3,\dots)\\
&=(-1)^{n-1-\rank(A_3)}\bar{E_3}(\xx _1\xx _2,\xx _3,\dots)\\
&=(-1)^{n-\rank(A)}\bar{E}(\mb{\xx}).
\end{align*}

Case 4: all of $\bar{E}_i$ are nonempty. By induction, we see that
$$E_i(\mb{x})=(-1)^{n-\rank(A_i)} \bar{E}_i(\mb{x^{-1}})$$
for $i=1,2$, and that
$$E_3(x_2,x_3,\dots)=(-1)^{n-1-\rank{A_3}}\bar{E}(\xx
_2,\xx_3,\dots).$$ From Lemma \ref{l-3-rankall},
$\rank(A_3)=\rank(A)$. Thus together with our key induction
equations \eqref{e-3-Ex123} and \eqref{e-3-Ebx123}, we get
\begin{align*}E(\mb{x})=&E_1(x_1,x_1x_2,x_3,\dots)
+ E_2(x_1x_2,x_2,x_3, \dots)-E_3(x_1x_2,x_3,\dots)\\
=&(-1)^{n-\rank(A)}\left(\bar{E}_1(\xx _1, \xx _1\xx _2,\xx
_3,\dots)\right.\\
&\qquad \qquad\qquad\qquad\qquad\left. +\bar{E}_2(\xx _1\xx_2, \xx
_2,\xx _3,\dots)+\bar{E_3}(\xx _1\xx
_2,\xx _3,\dots)\right)\\
=&(-1)^{n-\rank{A}}\bar{E}(\mb{x}).
\end{align*}
\end{proof}

\renewcommand{\theequation}{\thesection.\arabic{equation}}
\vfill\eject \setcounter{chapter}{3}
\chapter{Applications to Lattice Path Enumeration\label{s-lattice}}

In this section, we will use two methods to work on some lattice
path enumeration problems. One method is to use the bridge lemma,
which was used in \citep{bous}. The other method is to use the
factorization lemma, which was first discovered by \citep{ira},
and later rediscovered by \citep{bous}.

\section{Basic Concepts and the Bridge Lemma}

A path $\sigma$ in $\ZZ^2$ is a finite sequence of lattice points
$(a_0,b_0),\ldots , (a_n,b_n)$ in $\ZZ^2$, in which we call
$(a_0,b_0)$ the starting point, $(a_n,b_n)$ the ending point,
$(a_{i}-a_{i-1},b_i-b_{i-1})$ the steps of $\sigma$, and $n$ the
length of $\sigma$.

In what follows, the starting point of a path is always $(0,0)$
unless specified otherwise. The theory for other starting points
is similar.

Given two paths $\sigma_1$ and $\sigma_2$, we define their product
$\sigma_1 \sigma_2$ to be the path whose steps are those of
$\sigma_1$ followed by those of $\sigma_2$. If $\pi=\sigma_1
\sigma_2$, then we call $\sigma_1$ a {\em head} of $\pi$, and
$\sigma_2$ a {\em tail} of $\pi$.

For compatibility with the theory we developed on $\ZZ^2$, we
still use the reverse lexicographic order. Let $S$ be a
well-ordered subset of $\ZZ^2$. We are interested with paths  all
of whose steps lie in $S$.
 Denote by
$S^*$ the set of all such paths. Then any $\sigma\in S^*$ can be
uniquely factored as $\sigma=s_1s_2\cdots s_n$ for some $n\ge 0$,
and $s_i\in S$ for all $i$. The $s_i$'s are called the unit steps
of $\sigma$. Note that the empty path $\epsilon$ belongs to $S^*$.

The weight of a step $(a,b)\in S$ is defined to be
$\Gamma((a,b))=x^ay^bt$, and the weight of a path $\sigma
=s_1\cdots s_n$ is defined to be $\Gamma(\sigma)=
\Gamma(s_1)\cdots \Gamma(s_n)$. It is easy to see that for any two
paths $\sigma_1$ and $\sigma_2$, we have
$\Gamma(\sigma_1\sigma_2)=\Gamma(\sigma_1)\Gamma( \sigma_2)$. If
$P$ is a subset of $S^*$, then  we define
$$\Gamma(P)=\sum_{\sigma\in P} \Gamma({\sigma}).$$

In the special case that $P$ is the whole set $S^*$, we have
$$\Gamma(S^*)=\sum_{n\ge 0} (\Gamma(S))^n=(1-\Gamma(S))^{-1},$$
since each term in $(\Gamma(S))^n$ corresponds to a path of $n$
steps. This equation is interpreted as an identity in the field of
iterated Laurent series $\CC\ll x,y,t\gg$. In fact, we can do all
of the computations inside the ring $\CC\ll x,y\gg[[t]]$. Because
$\Gamma(S)$ is always the product of an element in $\CC\ll x,y\gg$
with $t$, it has positive order and $(1-\Gamma(S))^{-1}$ is in
$\CC\ll x,y\gg[[t]]$.

We can also write
\begin{align}\label{e-lattice-p}
\Gamma(P)=\sum_{n\ge 0} \sum_{i,j\in \ZZ} a_{i,j}(n) x^iy^j t^n,
\end{align}
where $a_{i,j}(n)$ is the number of paths in $P$ of length $n$
that end at $(i,j)$. The $a_{i,j}(n)$ in $\Gamma(P)$ is always
finite since this is true when $P=S^*$. We also call $\Gamma(P)$
the generating function of $P$ with respect to the ending points
and the lengths.

The requirement of $S$ being a well-ordered subset of $\ZZ^2$ is
necessary. For otherwise the $a_{i,j}(n)$ in $\Gamma(P)$ might be
infinite. For example, let $S=\{\, (m,0): m\in \ZZ\, \}$, and
consider the number of paths in $S^*$ that end at $(0,0)$ and have
length $2$.

Since $S$ is uniquely determined by $\Gamma(S)$, sometime we write
$\Gamma(S)$ instead of $S$. When $S$ is finite, $\Gamma(S)$ is a
Laurent polynomial, and $\Gamma(S^*)$ is in $
\CC[x,y,x^{-1},y^{-1}][[t]]$. Much work has been done in this
case.

Some operators on $\CC\ll x,y \gg[[t]]$ have simple combinatorial
interpretations. Let $P$ be a subset of $S^*$ with generating
function given by \eqref{e-lattice-p}.
\begin{enumerate}
\item The generating function for those paths in $P$
that end on the line $y=0$ is given by $\CT_y \Gamma(P)$.

\item The generating function for those paths in $P$
that end above the line $y=-1$ is given by $\pt_y \Gamma(P)$.

\item The generating function for those paths in $P$ that
end below the line $y=0$ is given by $\nt_y \Gamma(P)$.

\end{enumerate}
Similar properties hold for $x$. The diagonal generating function,
or the generating function for those paths in $P$ that end on the
line $y=x$, is given by
\begin{align}\label{e-lattice-diag}
\diag_{x=y} \Gamma(P) =\sum_{n\ge 0} \sum_{i\in \ZZ} a_{i,i}(n)
y^i t^n.
\end{align}
This can be written in terms of $\ct$. If we write
$F(x,y,t)=\Gamma(P)$, then
$$\diag_{x=y} \Gamma(P)= \ct_x F(x,y/x,t).$$
Note that $F(x/y,y,t)$ is not in the ring $\CC\ll x,y \gg [[t]]$,
unless $F(x,y,t)\in \CC((x,y))[[t]]$.

In the computations, we will use Theorem \ref{t-lagrange1}. Recall
that a positive root is a root with positive order.

\begin{exa}
If $S=\{\, (1,0),(-1,0),(0,1),(0,-1)\,\}$, then
$\Gamma(S)=t(x+y+x^{-1}+\yy)$. The generating function of paths in
$S^*$ that end on the line $y=0$ (or $x$-axis) can be computed as
follows.
\begin{align*}
\ct_y \Gamma(S^*) &= \ct_y \frac{1}{1-t(x+y+\xx+\yy)}\\
&=\ct_y \frac{y}{y-t(x+\xx)y -ty^2-t}\\
&=\frac1{1-t(x+\xx)-2tY},
\end{align*}
where $Y=Y(x,t)$ is the unique positive root  of $y-t(x+\xx)y
-ty^2-t$. This $Y$ can be found by quadratic formula. We have
$$Y= \frac{1-t\left (x+{x}^{-1}\right )-\sqrt{\left(1-t\left (x+{x}^{-1}\right )
\right)^2-4t^2}}{2t} .$$ After simplifying, the desired generating
function can be written as
$$\ct_y \Gamma(S^*)=\left[\left(1-t\left (x+{x}^{-1}\right )
\right)^2-4t^2\right]^{-1/2}.$$

Similarly, we can obtain the generating function for those paths
in $S^*$ that stay above the line $y=-1$. The generating function
for those paths in $S^*$ that stay below the line $y=0$ is
similar.
$$\pt_y \Gamma(S^*)=\frac1{1-t(x+\xx)-2tY} \frac{y}{y-Y} ,$$
\end{exa}

\vspace{3mm} The bridge lemma is a basic tool for studying
$S$-paths that never touch a certain area.

Let $S^R$ be a nonempty subset of $S^*$ with the condition that
$P\in S^R$ implies that all the heads of $P$ are also in $S^R$. We
call $R$ a {\em restriction} and say that $S^{R}$ is the set of
$S$-paths satisfying this restriction. Note that the empty path
belongs to $S^R$ since it is a head of any path.

The most common restrictions are avoiding some points. In
particular, $S^*$ is the set of $S$-paths avoiding nothing, or
with no restriction.

Let $S^{Rb}$ be the set of all $S$-paths that are not in $S^R$,
but for which deleting the last step results in a path in $S^R$.
The last step of a path in $S^{Rb}$ is like a bridge. Without it
we get a path in $S^R$, but with it we get a path not in $S^R$.

\begin{lem}[Bridge Lemma]\label{l-bridge}
Let $R$ be a restriction. Then we have
\begin{equation}
\Gamma(S^R)= \frac{1-\Gamma (S^{Rb})}{1-\Gamma (S)}.
\end{equation}
\end{lem}
\begin{proof}
For all $n\ge 1$, if we add another $S$-step to a path in $S^R$ of
length $n-1$, then we will get either a path in $S^R$ or a path in
$S^{Rb}$. On the other hand, by deleting the last step of a path
of length $n$ in $S^R$, we will get a path of length $n-1$ in
$S^R$. This bijection gives us the equation
$$1+\Gamma(S) \Gamma(S^R) =\Gamma(S^R)+\Gamma(S^{Rb}).$$
The lemma then follows by solving for $\Gamma(S^R)$.
\end{proof}

The bridge lemma is useful when $\Gamma(S^{Rb})$ has a certain
kind of nice form. This is best illustrated by examples.

We start by considering Catalan numbers, which are the most
frequently used numbers in combinatorics other than binomial
coefficients.
\begin{exa}
Let $S=\{\, (1,0),(0,1)\, \}$. In this case, if the ending point
of a path is $(i,j)$, then the path has length $i+j$. Hence we can
omit the parameter $t$, which keeps track of the length.

Let $R$ be the restriction that the path starts at $(0,0)$, and
never go above the line $y=x$. Then the number of paths in $S^R$
that end at $(n,n)$ is the Catalan number $C_n$.
\end{exa}

From the restriction, we see that the bridge paths must end on the
line $y=x+1$. Hence $\Gamma(S^{Rb})$ can be written as
$x^{-1}B(xy)$, where $B(z)$ is a power series in $z$ with $B(0)=0$
and $B'(0)=1$. Denote by $p(x,y)$ the generating function
$\Gamma(S^R)$. Then by the bridge lemma, we have
$$p(x,y)=\frac{1-x^{-1}B(xy)}{1-x-y}.$$

Since any path in $S^R$ never goes above the line $y=x$, its end
point cannot be $(n-1,n)$. Therefore $[x^{n-1}y^{n}] p(x,y)=0$ for
all $n\ge 0$. This is the same as $\diag\;  xp(x,y)=0$. Now we can
solve for $B$ as follows.
\begin{align*}
\diag \ xp(x,y)& =\CT_x xp(x,y/x) \\
&=\CT_x \frac{x(1-\xx B(y))}{1-x-y/x} \\
&=\CT_x \frac{x}{x-x^2-y} \cdot (x-B(y)) \\
&=\frac1{1-2X} (X-B(y)),
\end{align*}
where $X=\frac{1-\sqrt{1-4y}}{2}=yC(y)$ is the root of $x$ in
$x-x^2-y$ that belongs to $y\CC[[y]]$. Hence we have $B(y)=yC(y)$
and
$$p(x,y)=\frac{1-yC(xy)}{1-x-y}.$$
By setting $y=t$ and $x=t$, we will get the generating function
for paths in $S^R$ that are weighted by their lengths, but we do
not care where they end. We have
\begin{align*}
p(t,t)&=\left(1-t\frac{1-\sqrt{1-4t^2} }{2t^2}\right) \frac{1}{1-2t} \\
&=\frac{-1}{2t} +\frac{\sqrt{(1+2t)(1-2t)}}{2t(1-2t)} \\
&=\frac{-1}{2t}+\frac{1+2t}{2t\sqrt{1-4t^2}} \\
&=\sum_{n\ge 1} \frac{1}{2} \binom{2n}{n} t^{2n-1} +\sum_{n\ge 0}
\binom{2n}{n} t^{2n}.
\end{align*}
Hence the number of paths that start at $(0,0)$, with length $2n$,
and never go above the line $y=x$ is $\binom{2n}{n}$; the number
of paths that start at $(0,0)$, with length $2n-1$, and never go
above the line $y=x$ is $\binom{2n}{n}/2$.

\def\cxyt{\CC\ll x,y\gg[[t]]}

\vspace{3mm} Since the computation of the diagonal is always
converted into the computation of constant terms, it is convenient
to use Dyck paths, which are paths with northeast or southeast
steps, that never go below the horizontal axis. This time $S$ is
$\{\, (1,-1),(1,1)\,\}$, and we use $t$ for the horizontal
coordinate, and $y$ for vertical coordinate. Note that $x$ is
redundant since it also records the number of steps. The {\em
height} of a Dyck path $D$ is the $y$ coordinate of the highest
points in $D$.

\begin{exa}
Let $S=\{\, (1,-1),(1,1)\,\}$, and let $R$ be the restriction that
a path never goes below the line $y=0$ and never touches the line
$y=m$, where $m$ is a positive integer. Denote by $H_m(y,t)$ the
generating function $\Gamma(S^R)$. Then $H_m(y,t)$ is the
generating function for Dyck paths of height at most $m-1$. This
problem can be solved in many ways, but the method we are going to
give here has advantages. We use the bridge lemma and boundary
conditions to solve for it. The working field for this problem is
$\CC\ll y,t\gg$.
\end{exa}

It is easy to see that the paths in $S^{Rb}$ either end on the
line $y=-1$ or end on the line $y=m$, and paths in $S^{R}$ never
touch these two lines. Thus denote by $y^{-1}B(t)$ the generating
function for paths in $S^{Rb}$ that end on the bottom line $y=-1$,
and denote by $y^m T(t)$  the generating function for paths in
$S^{Rb}$ that end on the top line $y=m$. Then the bridge lemma and
the boundary conditions give us:
\begin{align*}
H_m(y,t) &=\frac{1-y^{-1}B(t)-y^mT(t)}{1-\yy t-yt} \\
\CT_y y H_m(y,t) & =0 \\
\CT_y y^{-m} H_m(y,t) &=0.
\end{align*}
To solve for $B$, $T$, and $H_m$, we have
\begin{align*}
\CT_y yH_m(y,t)&=\ct_y \frac{y\left(1-\yy B(t)-y^{m}T(t)\right)}{1-\yy t- yt} \\
&=\ct_y \frac{y}{y-ty^2-t} (y-B(t)-y^{m+1}T(t)) \\
&= \frac{1}{1-2tY}\left(Y-B(t)-Y^{m+1}T(t)\right),
\end{align*}
where $Y=\frac{1-\sqrt{1-4t^2}}{2t}=tC(t^2)$ is the root of $y$ in
$y-ty^2-t$ that belongs to $t\CC[[t]]$. So
\begin{equation}
\label{e-levelb} tC(t^2)-B(t)-(tC(t^2))^{m+1}T(t) =0
\end{equation}
Similarly,  we have
\begin{align*}
\CT_y y^{-m}H_m(y,t)&=\CT_{\yy} y^m H_m(\yy,t) \\
&=\frac{y^m(1-y B(t)-y^{-m}D(t)}{1-\yy t- yt} \\
&=\frac{y}{y-ty^2-t} (y^m-y^{m+1}B(t)-T(t)) \\
&= \frac{1}{1-2tY}(Y^m-Y^{m+1}B(t)-T(t)),
\end{align*}
where $Y=tC(t^2)$ is the unique root of $y$ in $y-ty^2-t$ that
belongs to $\CC[[t]]$. So
\begin{equation}
\label{e-levelt} (tC(t^2))^m-(tC(t^2))^{m+1}B(t)-T(t) =0.
\end{equation}
Now we can solve for $B(t)$ and $T(t)$ from \eqref{e-levelb} and
\eqref{e-levelt}. This gives us
$$ B(t)=tC(t^2) \frac{1-(tC(t^2))^{2m}}{1-(tC(t^2))^{2m+2}}, \qquad T(t)=\frac{(tC(t^2))^m-(tC(t^2))^{m+2}}
{1-(tC(t^2))^{2m+2}}.$$


\section{Applications to Walks on the Slit Plane}

Denote by $\mathcal{H}$ the half line $\{\, (-k,0); k\in \NN\,\}$.
Given  a well-ordered subset $\sy$ of $\ZZ^2$, {\em walks on the
slit plane} are paths that start at $(0,0)$ with steps in $\sy$
and never hit the half line $\mathcal{H}$ after the starting
point.

The problem of counting walks on the slit plane was first solved
by \citep{bous}. Much work has been done since then. See
\citep{bouso,bousslitc} The basic tools for solving this kind of
problem are the bridge lemma and the unique factorization lemma.
In the next section, we shall see that using the concept of
``Gessel pair", we can solve it directly by the unique
factorization lemma. We will work with walks on the slit plane by
using the bridge lemma.

The set of all walks on the slit plane is equal to $\sy^R$, where
$R$ is the restriction that a path never hits the half line
$\mathcal{H}$ after the starting point. Therefore, paths in
$\sy^{Rb}$ must end at some $(-k,0)$ for some $k\ge 0$, and
$\Gamma(\sy^{Rb})$ can be written as $B(\xx,t)$, which contains
only negative powers in $x$ except $1$. Denote by $S(x,y;t)$ the
generating function $\Gamma(\sy^{R})$. Then
$$S(x,y;t)=\sum_{n\ge 0}\sum_{i\in \ZZ}\sum_{j\in \ZZ} a_{i,j}(n) x^iy^j t^n,$$
where $a_{i,j}(n)$ is the number of $n$-step walks on the slit
plane that end at $(i,j)$. This is an element of $\cxyt$. When
$\sy$ is finite, $S(x,y;t)$ is a formal power series in $t$ with
coefficients in $\CC[x,\xx,y,\yy]$.

Applying the bridge lemma, we get the functional equation
\begin{align}\label{e-4-1}
S(x,y;t) &=\frac{1-B(\xx,t)}{1-\Gamma(\sy)},
\end{align}

Let $S_0(x,t)$ be the generating function of walks on the slit
plane that end on the line $y=0$. Then $S_0(x,t)$ contains only
positive powers in $x$ except $S_0(x,0)=1$. The boundary condition
is given by
\begin{align}
\label{e-4-2} \CT_y S(x,y;t)=S_0(x,t).
\end{align}

\citet{bouso} defined {\em bilateral walks} to be paths in $\sy^*$
that end on the $x$-axis. Let $S_x(x;t)$ be the generating
function of bilateral walks. Then we have
$$S_x(x;t)=\ct_y \frac{1}{1-\Gamma(\sy)}.$$

One important result for slit plane walks is the following
theorem, which was obtained in \citep{bouso} for the case of $\sy$
being a finite set. This result says that the $S_0(x,t)$,
$B(\xx,t)$, and $S(x,y;t)$ can be theoretically computed. In
practice, computing them is not a easy task. Only special cases
have been thoroughly studied.

\begin{thm}\label{t-4-slitg}
Let $\sy$ be a well-ordered set in $\ZZ^2$. Using notation as
above, we have:
\begin{align}
S_0(x,t) &=\left(S_x(x,t)\right)_+, \\
\frac{1}{1-B(\xx,t)} &= \left(S_x(x,t)\right)_0\left(S_x(x,t)\right)_-,\\
S(x,y;t)&=
\frac{1}{(1-\Gamma(\sy))\left(S_x(x,t)\right)_0\left(S_x(x,t)\right)_-}.\label{e-4-sxyt}
\end{align}
\end{thm}
\begin{proof}
Substituting \eqref{e-4-1} into \eqref{e-4-2}, we get
$$S_0(x,t)=\ct_y \frac{1-B(\xx,t)}{1-\Gamma(\sy)}=(1-B(\xx,t))\ct_y \frac{1}{1-\Gamma(\sy)},$$
which can be written as
\begin{equation} \label{e-slit-fac}
S_0(x,t) \frac{1}{1-B(\xx,t)} =S_x(x;t).
\end{equation}

Now we can check that $S_x(x;0)=1$, $S_0(x,0)=1$, and $
1-B(\xx,0)=1$. Recall that except for $1$, $S_0(x,t)$ contains
only positive powers in $x$, $(1-B(\xx,t))^{-1}$ contains only
negative powers in $x$. Thus the unique factorization lemma
applies, and the theorem follows.
\end{proof}

\begin{rem}
A combinatorial interpretation of equation \eqref{e-slit-fac} can
be given by using the cycle lemma \citep{bouso}. We will give
another interpretation in the next section.
\end{rem}

From the proof of the theorem, we see that $\log S_0(x,t)=\pt_x
\log S_x(x;t)$. Now if $\log S_0(x,t)$ has the form
$b(t)x^p+\text{higher degree terms}$, then so does $S_0(x,t)-1$.
This gives us the following result \citep[Proposition 4]{bouso}.
\begin{prop}
Let $p$ be the smallest positive integer such that there is a walk
on the slit plane that ends at $(p,0)$. Then the generating
function for such walks is given by
$$S_{p,0}(t)=[x^p] \log S_x(x;t).$$
\end{prop}
\bous\ shows in addition that $S_{k,0}$ is $D$-finite for every
$k$. We give an explicit example as follows, and we will discuss
this further in the next section.

\vspace{3mm} Example. We shall give an example where
$\Gamma(\sy)=t/(xy(1-x)(1-y))$. In this example, $S(x,y;t)$ does
not seem to be algebraic.

First we need to compute the generating function $S_x(x;t)$ of
bilateral walks.
\begin{align*}
S_x(x;t)=\ct_y \frac{1}{1-t/(xy(1-x)(1-y))}=\ct_y
\frac{y(1-y)}{y(1-y)-t/(x(1-x))}.
\end{align*}
Solving the denominator for $y$, we get a unique positive root
$$Y=Y(x;t)= \frac{1-\sqrt{1-4t/(x(1-x))}}{2}.$$
Applying Theorem \ref{t-lagrange1}, we get
$$S_x(x,t)=\frac{1-Y}{1-2Y}=\frac{1+\left(1-4t/(x(1-x))\right)^{-1/2}}{2}.$$
It is easy to expand this into a series. We have
$$S_x(x,t)=\frac12+\sum_{n\ge 0}\sum_{k\in \ZZ} \frac12 \binom{2n}{n}\binom{2n+k-1}{n+k-1}x^kt^n.$$
In particular, the constant term of $S_x(x,t)$ in $x$ is:
$$\ct_x S_x(x,t) = 1+\sum_{n\ge 1} \binom{2n}{n-1}^2 t^n.$$
 Though $S_x(x,t)$ is a
simple algebraic series, $S_0(x,t)=\left(S_x(x,t)\right)_+$ does
not seem to be algebraic.

We can get a formula for $\log S_0(x,t)$ by computing the series
expansion of $\log S_x(x,t)$, and then collecting all the terms
containing positive powers in $x$.

Now $\log S_x(x,t)$ is a power series in $t$ with constant term
$0$. Its series expansion can be obtained by finding the series
expansion of its derivative in $t$, and then integrating. We have
\begin{align*}
\frac{\partial}{\partial t} \log
S_x(x,t)&=\frac{\partial}{\partial t}\log
\frac{1+\left(1-4t/(x(1-x))\right)^{-1/2}}{2}\\
&= \frac1{2t}
\left(\frac{1}{1-4t/(x(1-x))}-\frac{1}{\sqrt{1-4t/(x(1-x))}}
\right),
\end{align*}
where the final formula is obtained after some simplification and
rationalization. Now it is easy to see the following:
$$\log S_0(x,t)=\sum_{n\ge 1}\sum_{k\ge 1} \frac1{2n} t^n \left(4^n-\binom{2n}{n}\right)
\binom{2n-1+k}{n-1}x^k. $$

In particular, \begin{align*} S_0(x,t)& =[x]\log S_0(x,t)
=\sum_{n\ge 1}
\frac{1}{2n}\left(4^n-\binom{2n}{n}\right) \binom{2n}{n-1}t^n \\
&=t+10{t}^{2}+110{t}^{3}+1302{t}^{4}+16212{t}^{5}+209352{t}^{6
}+\cdots.
\end{align*}

\vspace{3mm} Let $\sy$ be a set of steps. Then we say that $\sy$
satisfies the \emph{small height variation} condition if
$\text{for all } (i,j) \in \sy, |j|\le 1.$ In this case, it is
clear that we have
$$\Gamma (\sy )= t\sum_{(i,j)\in \sy}x^iy^j = tA_{-1}(x)\yy +tA_0(x)+tA_1(x) y.$$

Slit plane walks satisfying this condition have been thoroughly
studied. We have the following result, which is a slight variation
of \citep[Theorem 17]{bouso}.
\begin{thm}
Let $\sy$ be a well-ordered set of steps with small height
variations. Let
\begin{align}
\Delta(x;t)=(1-tA_0(x))^2-4t^2A_1(x)A_{-1}(x).
\end{align}
Then the generating function for bilateral walks is
$$S_x(x;t)= \frac{1}{\sqrt{\Delta(x;t)}}.$$
If $\Delta_-(\xx;t)\Delta_0(t)\Delta_+(x;t)$ is the third
decomposition of $\Delta$ with respect to $x$, then the generating
function for walks on the slit plan with steps in $\sy$ is
$$S(x,y;t)=\frac{\sqrt{\Delta_0(t)\Delta_-(\xx;t)}}{1-\Gamma(\sy)}. $$
\end{thm}
\begin{rem}
If $\Delta(x;t)$ is rational, then $\Delta_-, \Delta_0,$ and $
\Delta_+$ are algebraic, and hence $S(x,y;t)$ is algebraic. This
is always true when $\sy$ is a finite set, as has been discussed
in \citep{bouso}.
\end{rem}

\begin{proof}
First let us compute the generating function for bilateral walks.
We have
\begin{align*}
S_x(x;t)&=\ct_y \frac{1}{1-t(A_{-1}(x)\yy +A_0(x)+A_1(x) y)} \\
&= \ct_y \frac{y}{y-t(A_{-1}(x)+A_0(x)y+A_1(x)y^2)}
\end{align*}
By Theorem \ref{t-lagrange1}, if we let $Y=Y(x)$ be the unique
positive root for $y$ in the above denominator, i.e.,
$$ Y= \frac{1-tA_0(x)-\sqrt{\Delta(x;t)}}{2tA_1(x)}, $$
then
$$S_x(x;t)= \frac{1}{1-tA_0(x)-2tA_1(x)Y(x)}= \frac{1}{\sqrt{\Delta(x;t)}}.$$

The formula for $S(x,y;t)$ is obtained by applying the formula
\eqref{e-4-sxyt}.
\end{proof}

In \citep{bouso}, three examples were computed explicitly. They
are:
\begin{enumerate}
    \item The
example of the ordinary lattice, with \\
$\Gamma(\sy)=t(y+x+\xx+\yy)=t(x+y)(1+\xx \yy)$.
    \item The example of the
diagonal lattice, with \\
$\Gamma(\sy)= t(xy+\xx y+x\yy+\xx \yy)=t(x+\xx)(y+\yy)$.
    \item The example of the triangular lattice, with\\
$\Gamma(\sy)=t(y+xy+x+\xx+\yy+\xx\yy)$
\end{enumerate}
All of the above three examples are symmetric in $x$ and $y$. We
give another example as follows.

\begin{exa}
Let $\sy $ be given by $A_1(x)=A_{-1}(x)=\xx$ and
$A_0=1/(x(1-x))$.
\end{exa}
First we compute $\Delta$. \begin{align*} \Delta(x;t)&= \left(
1-{\frac {t}{x
\left( 1-x \right) }} \right) ^{2}-4{\frac {{ t}^{2}}{{x}^{2}}}\\
&={\frac {{x}^{4}-2\,{x}^{3}+ \left( 1-4\,{t}^{2}+2\,t \right)
{x}^{2}+
 \left( -2\,t+8\,{t}^{2} \right) x-3\,{t}^{2}}{{x}^{2} \left( x-1
 \right) ^{2}}}
.
\end{align*}

The four roots are given by
$$X_1(t)=\frac{1-\sqrt {1+4{t}^{2} } }{2}-t,\qquad X_2(t)=\frac{1-\sqrt {1-8t+4{t}^{2}}}{2}+t,$$
$$X_3(t)=\frac{1+\sqrt {1+4{t}^{2}}}{2}-t,\qquad X_4(t)= \frac{1+\sqrt {1-8t+4{t}^{2}}}{2}+t,
$$
where $X_1$ and $X_2$ have positive order, and $X_3$ and $X_4$
have zero order. Therefore
$$\frac{1}{\sqrt{\Delta(x;t)}}=\frac{1}{\sqrt{(1-X_1/x)(1-X_2/x)}}\cdot L(t) \cdot \frac{1-x}{\sqrt{(1-x/X_3)(1-x/X_4)}} $$
is the third decomposition of $S_x(x,t)$, where
$L(t)=\sqrt{X_3(t)X_4(t)}$ can be obtained by equating
coefficients of $x^2$ in $\Delta(x;t)$. Thus

\begin{align*}
S_0(x,t)&=\frac{1-x}{\sqrt{(1-x/X_3)(1-x/X_4)}},\\
S(x,y,t)&=\frac{\sqrt{(1-X_1/x)(1-X_2/x)/X_3(t)/X_4(t)}}{1-t(y/x+1/(x(1-x))+1/(xy))}.
\end{align*}

\section{Unique Factorization Lemma and Gessel Pairs\label{s-gesselpair}}
We now introduce the combinatorial interpretation to the
factorization lemma in terms of lattice paths in the plane. This
idea was first introduced in \citep{ira}. We modify this idea to
fit in a more general setting.

A {\em monoid} is a set $M$, equipped with a multiplication which
is associative, and having a unit element $1$.

For the set of paths, the multiplication of two paths is just the
product of two paths as we defined earlier. Thus the empty path is
the unit.

Let $H$ be a set of paths with steps in $S$ that start at $(0,0)$.
If $H$ is closed under multiplication of paths and contains the
empty path, then $H$ is a monoid. We call a nonempty path
$\sigma\in H$ a {\em prime} if it cannot be factored into two
nonempty paths in $H$. We say that $H$ is a {\em free monoid} if
any element in $H$ can be uniquely factored into products of
primes in $H$.

If $H$ is a free monoid, then for any $\sigma \in H$ with its
factorization into primes as $\sigma=h_1h_2\cdots h_m$, we say
that $h_1h_2\cdots h_i$ is an $H$ {\em head} of $\sigma$ for
$i=0,1,\ldots ,m$. If we let $P$ be the set of primes in $H$, then
$\Gamma(H)=1/(1-\Gamma(P))$.

For example, $S^*$ is a free monoid, whose primes are all the
elements in $S$.

The set of all paths in $S^*$ that end on the $x$-axis is a free
monoid, whose primes are those paths that return to the $x$-axis
only at the end point.

The set of all paths in $S^*$ that end at $(k,0)$ for some $k\ge
0$ is a free monoid. The primes are those paths that only return
the nonnegative half of the $x$-axis at the end point.

Let $\rho$ be a map from $H$ to $\ZZ$. We say that $\rho$ is a
{\em homomorphism} from $H$ to $\ZZ$ if $\rho(\epsilon)=0$ and for
all $\sigma_1,\sigma_2\in H$, $\rho(\sigma_1 \sigma_2)
=\rho(\sigma_1)+\rho(\sigma_2)$. The $\rho$ value of a path
$\sigma$ is $\rho(\sigma)$.

If $H$ is a free monoid, then any map from $H$ to $\ZZ$ defined on
the primes of $H$ induces a homomorphism. If in addition, $H$ is a
subset of $S^*$, then the natural map to the end point of a path
is a homomorphism from $H$ to $\ZZ^2$. Therefore, any homomorphism
from $\ZZ^2$ to $\ZZ$ induces a homomorphism from $H$ to $\ZZ$
through that natural map. The following two homomorphisms are
useful. Define $\rho_x(\sigma)$ to be  the $x$ coordinate of the
ending point of $\sigma$, then $\rho_x$ is clearly a homomorphism.
Similarly we can define $\rho_y$.

If $H$ is a free monoid, and $\rho$ is a
 homomorphism from $H$ to $\ZZ$, then we call $(H,\rho)$
a Gessel pair. For a Gessel pair $(H,\rho)$, we define:

A {\em minus-path} is either the empty path or a path whose $\rho$
value is negative and less than the $\rho$ values of all the other
$H$ heads.

A {\em zero-path} is a path with $\rho$ value $0$ and all of whose
$H$ heads have nonnegative $\rho$ values.

A {\em plus-path} is a path all of whose $H$ heads (except
$\epsilon$) have positive $\rho$ values.

For  a Gessel pair $(H,\rho)$, we denote by $H_ -$, $H_0$, and
$H_+$ respectively to be the sets of minus-, zero-, and plus-paths
in $H$. Note that the empty path, but no other path, belongs to
all three classes. The path $h_1h_2\cdots h_n$, where $h_i\in H$,
is a minus-path if and only if $h_nh_{n-1}\cdots h_1$ is a
plus-path; thus the theories of minus- and plus-paths are
identical.

\begin{lem}\label{l-pathfactor}
Let $(H,\rho)$ be a Gessel pair, and let $\pi$ be a path in $H$.
Then $\pi$ has a unique factorization $\pi_- \pi_0 \pi_+$, where
$\pi_-$ is a minus-path, $\pi_0$ is a zero-path, and $\pi_+$ is a
plus-path.
\end{lem}
\begin{proof}
Let $a$ be the smallest among all the $\rho$ values of the $H$
heads of $\pi$. Let $\pi_-$ be the shortest $H$ head of $\pi$
whose $\rho$ value equals $a$. Then if $\pi=\pi_- \sigma$, let
$\pi_-\pi_0$ be the longest $H$ head of $\pi$ whose $\rho$ value
equals $a$, and let $\pi_+$ be the rest of $\sigma$. It is easy to
see that this factorization satisfies the required conditions.

To see that it is unique, let $\tau_-\tau_0\tau_+$ be another
factorization of $\pi$. By definition, any $H$ head of $\tau_0
\tau_+$ has a nonnegative $\rho$ value. So the minimum $\rho$
value among all of the $H$ heads of $\pi$ is achieved in $\pi_-$.
By definition, it equals $\rho(\tau_-)$ and is unique in $\tau_-$.
Therefore, $\rho(\tau_-)=a$ and $\tau_-=\pi_-$ by the selection of
$\pi_-$. The reasons for $\pi_0=\tau_0$ and $\pi_+=\tau_+$ are
similar.

\end{proof}

\begin{prop}\label{p-pathfactor}
If $(H,\rho)$ is a Gessel pair, then $H_-$, $H_0$, and $H_+$ are
all free monoids. The map from $H$ to $H_-\times H_0\times H_+$
defined by $\pi\to(\pi_-,\pi_0,\pi_+)$ is a bijection.
\end{prop}
\begin{proof}
By Lemma \ref{l-pathfactor}, the map defined by
$\pi\to(\pi_-,\pi_0,\pi_+)$ is clearly a bijection. Now we show
that $H_-$, $H_0$, and $H_+$ are all free monoids.

It is easy to see that they are monoids. We only show that $H_-$
is free. The other parts are similar. Let $P$ be the subset of
$H_-$ such that $\sigma \in P$ if and only if
 $\rho(\sigma)$ is negative and every other
$H$ head of $\sigma$ has nonnegative $\rho$ value. We claim that
$P$ is the set of primes in $H_-$.

Clearly any $\sigma\in P$ cannot be factored as the product of two
nontrivial elements in $H_-$. Now let $\pi\in H_-$. In order to
factor $\pi$ into factors in $P$, we find the shortest $H$ head of
$\pi$ that has negative $\rho$ value, and denote it by $\sigma_1$.
Then $\pi$ is factored as $\pi=\sigma_1\pi'$ for some $\pi'$ in
$H$. From the definition of minus-path, $\rho(\sigma_1)$ is either
less than $\rho(\pi)$, in which case $\pi'$ is clearly in $ H_-$,
or $\rho(\sigma_1)=\rho(\pi)$, in which case $\pi'$ has to be the
unit and $\pi=\sigma_1$ is in $P$. So we can inductively obtain a
factorization of $\pi$ into elements in $P$.

The uniqueness of this factorization is clear.
\end{proof}

In a Gessel pair $(H,\rho)$, the weight of an element $\pi\in H$
is defined to be $\Gamma(\pi)z^{\rho(\pi)}$, where $z$ is a new
variable. When $H$ is also a subset of $S^*$ and we are
considering the Gessel pair $(H,\rho_x)$, the power in  $z$ is
always the same as the power in $x$ for any $\pi$ in $H$. So we
can replace $z$ by $1$ and let $x$ play the same role as $z$.
Since the factorization in $H$ is with respect to $\rho$, the
factorization of generating function is with respect to $z$.

\begin{thm}\label{t-pathfactor}
For any Gessel pair $(H,\rho)$, we have
$\Gamma(H_-)=[\Gamma(H)]_-$, $\Gamma(H_0)=[\Gamma(H)]_0,$ and
$\Gamma(H_+)=[\Gamma(H)]_+$.
\end{thm}
\begin{proof}
From Proposition \ref{p-pathfactor}, it follows that
$\Gamma(H)=\Gamma(H_-)\Gamma(H_0)\Gamma(H_+)$. Clearly except $1$,
which is the weight of the empty path, $\Gamma(H_-)$ contains only
negative powers in $z$, $\Gamma(H_0)$ is independent of $z$, and
$\Gamma(H_+)$ contains only positive power in $z$. The theorem
then follows from the unique Factorization Lemma with respect to
$z$.
\end{proof}

\citet{ira} gives many interesting examples about lattice paths on
the plane. We introduce the most classical example as the
following:
\begin{exa}
Let $S$ be $\{\, (1,r),(1,-1)\,\}$ with $r\ge 1$, and $H=S^*$.
Consider the Gessel pair $(H,\rho_y)$.
\end{exa}
Note that in this case the length of a path equals the $x$
coordinate of its end point. Replacing $x$ by $1$ will not lose
any information.

Clearly we have
$$\Gamma(H)=\Gamma(S^*)=\frac{1}{1-t(y^r+1/y)}.$$

We see that $H_+$ is the set of paths in $S^*$ that never go below
level $1$ after the starting point. The set $H_0$ contains all
paths in $S^*$ that end on level $0$ and never go below level $0$.
When $r=1$, this becomes Dyck paths.

To compute $\Gamma(H_0):=F(t)$, we let $Y(t)$ be the unique
positive root of $y-t(1+y^{r+1})$. By Theorem \ref{t-1-3rd0},
$F(t)=Y(t)/t)$. Now it is easy to see that
$F(t)=1+t^{r+1}F(t)^{r+1}$. So $F(t)$ equals the generating
function of complete $r+1$-ary trees.

\begin{exa}
Let $S$ be $\{\, (1,1), (1,-1) \,\}$, and let $H=S^*$. Let $\rho$
be determined by $\rho(1,1)=r$ and $\rho(1,-1)=-1$.
\end{exa}
It is easy to see that this example is isomorphic to the previous
one.

\begin{exa}\label{ex-4-hphl}
In general if $H=S^*$, then $(H,\rho_y)$ is a Gessel pair.
\end{exa}
We see that $H_+$ is the set of paths in $S^*$ that never go below
the line $y=1$ after the starting point.

If we let $J=H_+$, then $J$ is also a free monoid. The primes of
$J$ are paths that start at $(0,0)$, end at some positive level
$d$, and never hit level $d-1$ or lower.

The set $H_0$ contains all paths in $S^*$ that end on the line
$y=0$, and never go below the line $y=0$. In other words, $H_0$
contains all paths in $S^*$ that stays in the upper half plane and
end on the $x$-axis.

If we let $J=H_0$, then $(J,\rho_x)$ is a Gessel pair. The set
$J_+$ contains all paths in $J$ that avoiding the half line
$\mathcal{H}$ after the starting point. This is the same as walks
on the half plane avoiding the half line in \citep{bouso}.

The set $J_0$ contains all paths in $J $ that ending at $(0,0)$
and never touch the half line $\mathcal{H}$ except $(0,0)$.

\begin{exa}\label{ex-4-slit}
For any $S$, let $H$ be the set of paths that end on the $x$-axis.
Then $(H,\rho_x)$ is a Gessel pair.
\end{exa}
The set $H_+$ contains all paths that end on the $x$ axis and
never hit the half line $\mathcal{H}=\{\, (-k,0)\mid k\ge 0\,\}$
after the starting point. This is exactly the walks on the slit
plane that end on the $x$-axis.

The set $H_0$ contains all paths that end at $(0,0)$, and never
touch $(-k,0)$ for $k=1,2,\ldots$. This was call the set of loops
in \citep{bouso}.

As we proposed, we shall give a combinatorial explanation of
equation \eqref{e-slit-fac}. The set $H_-H_0$ is a free monoid. It
contain all paths that end at $(-k,0)$ for some $k\ge 0$. Its
primes are all paths that hit $(-k,0)$ only once at its end point.
These primes are exactly the bridge paths. So we have
$$\Gamma(H_-H_0)=\frac{1}{1-B(\xx;t)}, \quad \text{ and } \Gamma(H_+)=S_0(x,t).$$
Equation \eqref{e-slit-fac} then follows.

\begin{exa}\label{ex-4-hphl1}
For any $S$, let $H$ be the set of paths that end on the $x$-axis
and never go below the line $y=-d$ for some given $d>0$. Then it
is easy to check that $(H,\rho_x)$ is a Gessel pair.
\end{exa}

The set $H_+$ contains all paths that end on the $x$-axis, and
never hit the half line $\mathcal{H}$ after the starting point,
and never go below the line $y=-d$.

The set $H_0$ can be similarly described.

\begin{exa}\label{ex-4hphl2}
For any $S$, let $H$ be the set of paths that end on the $x$-axis
and never go below the line $y=-d$ and never go above the line
$y=f+1$ for some given positive integers $d$ and $e$. Then it is
easy to see that $(H,\rho_x)$ is a Gessel pair.
\end{exa}
This example is similar to the previous one.

\section{Explicit Examples\label{s-4-exa}}

We will do some explicit examples, several of which were proposed
in \citep{bouso}. Our task is to find a formula for $\log
\Gamma(H_+)$ for an algebraic $\Gamma(H)$ as previously described.

The idea is as follows. Let $P(x,y,t)$ be a polynomial and let
$Y(x;t)$ be the unique positive root of $y-tP(x,y,t)$ for $y$. The
problem will be reduced to finding the third decomposition of
$Q(x,Y(t),t)$ with respect to $x$ for some rational $Q$. We are
especially interested in $[x^p]Q_+(x,Y(t),t)$ for some positive
integer $p$, which is $D$-finite by the argument in \citep{bouso}.
This generating function can be obtained if we can get a nice form
of $\frac{\partial}{\partial t}\log Q(x,Y(t),t)$. Our approach to
finding such a nice form is to do all the computation implicitly.
It is best illustrated by examples.

\begin{exa}
Let $\sy$ be the set $\{ \, (1,0),(-1,0),(0,2)(0,-1)\, \}$, or
equivalently $\Gamma(\sy)=t(x+\xx+y^2+\yy)$. \citep{bouso}
proposed the problem of solving walks on the slit plane in this
model, or even replace the $2$ by a general positive integer $q$.
\end{exa}
Our method works for general $q$, but so far we have found a
reasonable formula only for $q=2$. We have:

\begin{prop}
The number of walks on the slit plane, with steps in \\
$\{ \, (1,0),(-1,0),(0,2),(0,-1)\, \}$, of length $N$, and ending
at $(1,0)$ equals
\begin{multline}\label{e-4-exay2s10}
a_{1,0}(N)= \binom{N}{\frac{N-1}2}+\sum_{n=1}^{\lfloor N/3\rfloor
}\frac{3^{3n-1}}{n2^{2n}}\binom{N-1}{3n-1}\binom{N-3n}{\frac{N-3n}{2}}+\\
\sum_{n,m,k}\frac{3^{3m+2}}{nN2^{2m+2}}
\binom{n}{k,2k+1,n-3k-1}\binom{N-n}{3m+2}
\binom{N-3m-3k-3}{\frac{N-3m-3k-4}{2}},
\end{multline}
where $\binom{A}{B+1/2}$ is interpreted as $0$ for all integers
$A,B$, and the second sum ranges over all $n,m,k$ such that $1\le
n\le N$, $0\le m\le \frac{N-n-2}{3} $, and $0\le k\le
\frac{n-1}{3}$.
\end{prop}

\begin{proof}
We proceed by computing $S_x(x;t)$. Let $b=x+\xx$. Then
$\Gamma(\sy)=t(b+y^2+\yy)$. We have
\begin{align*}
S_x(x;t)&=\ct_y \frac{y}{y-t(y^3+by+1)}=\frac{1}{1-tb-3tY^2},
\end{align*}
where $Y=Y(t)=Y(b,t)=Y(x,t)$ is the unique positive root of the
denominator for $y$. More precisely, $Y$ is the unique power
series in $t$ with constant term $0$ that satisfies
\begin{align}
\label{e-4-exay2} Y(t)-t(Y(t)^3+bY(t)+1)=0.
\end{align}
Using the Lagrange inversion formula we get
\begin{align}\label{4-exay2Y}
Y(t)=\sum_{n\ge 1} \sum_{k=0}^{\left\lfloor
\frac{n-1}{3}\right\rfloor} \binom{n}{k,2k+1,n-3k-1} b^{n-3k-1}
t^n.
\end{align}

We can compute $\log S_x(x;t)$ explicitly in order to obtain $\log
S_0(x,t)$. We have
\begin{align}
\frac{\partial }{\partial t}\log S_x(x;t)&=\frac{b+3 \left( Y
\left( t \right) \right) ^{2}+6tY \left( t \right) {
\frac{\partial }{\partial t}}Y \left( t
\right) }{1-tb-3tY(t)^2}\nonumber \\
&=\frac{b-t{b}^{2}+3\, \left( Y \left( t \right)  \right) ^{2}-3\,
\left( Y
 \left( t \right)  \right) ^{4}t+6\,tY \left( t \right)
}{(1-tb-3tY(t)^2)^2}\label{e-4-exay2-lns}
\end{align}
where
$$ \frac{\partial }{\partial t}Y \left( t \right) = \frac{1+bY(t)+Y(t)^3}{1-tb-3tY(t)^2}$$ is
determined implicitly by equation \eqref{e-4-exay2}.

Since $Y(t)$ satisfying \eqref{e-4-exay2}, we can rewrite
\eqref{e-4-exay2-lns} as $C_0+C_1Y(t)+C_2Y(t)^2$, where $C_i$ are
rational functions in $b$ and $t$. This can be done by Maple, and
we get
\begin{align}\label{e-4-exay2-lnsf}
\frac{\partial }{\partial t}\log S_x(x;t)= \frac{\left(
4\,{b}^{3}+27 \right) {t}^{2}-8\,t{b}^{2}+4\,b }{4(1-bt)^3-27t^3}+
\frac{9tY(t)}{4(1-bt)^3-27t^3}.
\end{align}
The first term has a simple form:
$$ \int \frac{\left( 4\,{b}^{3}+27 \right)
{t}^{2}-8\,t{b}^{2}+4\,b }{4(1-bt)^3-27t^3}dt = \log
\left(4(1-bt)^3-27t^3\right)^{-1/3}+C,
$$
where $C$ is independent of $t$. After some manipulation, we get
\begin{align*} \int \frac{\left( 4\,{b}^{3}+27
\right) {t}^{2}-8\,t{b}^{2}+4\,b }{4(1-bt)^3-27t^3}dt
=\log\frac{1}{1-bt}+\sum_{N\ge 1} \sum_{n=1}^{\lfloor
N/3\rfloor}\frac{3^{3n-1}}{n2^{2n}}\binom{N-1}{3n-1}b^{N-3n}t^N.
\end{align*}

For the second term, we have
\begin{align*}
\frac{9t}{4(1-bt)^3-27t^3} &= \frac{9t}{4(1-bt)^3}
\frac{1}{1-27t^3/(4(1-bt)^3)}.
\end{align*}
After some manipulation, we get
\begin{align*}
\frac{9t}{4(1-bt)^3-27t^3} &= \sum_{m\ge 0}
\frac{3^{3m+2}}{2^{2m+2}} \sum_{r\ge 0 } \binom{3m+r+2}{3m+2} b^r
t^{3m+r+1}.
\end{align*}
Thus together with the expansion of $Y(t)$ given by
\eqref{4-exay2Y}, we obtain
\begin{multline*}
    \int \frac{9tY(t)}{4(1-bt)^3-27t^3}dt =\sum_{N\ge
1}\sum_{n=1}^N\sum_{m=0}^{\left\lfloor \frac{N-n-2}{3}
\right\rfloor }
\sum_{k=0}^{\left\lfloor \frac{n-1}{3} \right\rfloor}\\
    \frac{3^{3m+2}}{nN2^{2m+2}}
\binom{n}{k,2k+1,n-3k-1}\binom{N-n}{3m+2}b^{N-3m-3k-3} t^N.
\end{multline*}
Note that the power in $b$ is always nonnegative. It is easy to
separate the negative power and positive powers in $b^M=(x+\xx)^M$
for every nonnegative integer $M$. Thus we can obtain a formula
for $\log S_0(x,t)$. In particular, from the formulas $[x]
(x+\xx)^M=\binom{M}{\frac{M-1}{2}}$ and $S_{1,0}(t)=[x]\log
S_x(x;t)$, we get \eqref{e-4-exay2s10}.
\end{proof}

\begin{exa}
We consider walks on the half plane avoiding half line; more
precisely, walks that never touch the half line $\mathcal{H}$ and
never hit a point $(i,j) $ with $j<0$. This is a continuation of
Example \ref{ex-4-hphl}. We denote by $HS(x,y;t)$ the generating
function for such paths.
\end{exa}

It turns out that this case is simpler than the previous one. We
obtain the following result, which includes \citep[Proposition
25]{bouso} as a special case.
\begin{thm}
For any well-ordered set $\sy$, let $p$ be the smallest positive
number such that there is an  $\sy$-path end at $(p,0)$. Then the
number of walks on the half plane avoiding the half line that end
at $(p,0)$ and are of length $n$ is equal to $1/n$th  times the
number of $\sy$-paths that end at $(p,0)$ and are of length $n$.
\end{thm}

\begin{proof}
We use the notation of Example \ref{ex-4-hphl}. From the Gessel
pair $(\sy^*, \rho_y)$, we have
$\Gamma(H_0)=\left(\Gamma(\sy^*)\right)_0$ and
\begin{align*}
\log \Gamma(H_0) =\ct_y \log \Gamma(\sy^*))=\ct_y \log
\frac{1}{1-\Gamma(\sy)}.
\end{align*}

Now let $J=H_0$ and consider the Gessel pair $(J, \rho_x)$. Then
$$ \log \Gamma(J_0J_+)= \pt_x  \log \Gamma(J). $$
In particular, we have
\begin{align*}
[x^p] \Gamma(J_+)&= [x^p] \log \Gamma(J)=[x^p] \log
\Gamma(H_0)=[x^p] \ct_y \log \Gamma(\sy^*).
\end{align*}
Therefore,
$$[x^p t^n] \Gamma(J_+) = [x^py^0t^n] \frac{1}{n}\Gamma(\sy)^n. $$
This prove the theorem.
\end{proof}

\section{Proof of a Conjecture about Walks on the Slit Plane\label{ss-conj}}


Let $a_{i,j}(n)$ denote the number of walks in $n$ steps from
$(0,0)$ to $(i,j)$, with steps $(\pm 1,0)$ and $(0,\pm 1)$, never
touching  a point $(-k,0)$ with $k\ge 0$ after
 the starting point. These are called {\em walks on the slit plane}.

Let $\xx$ denote $x^{-1}$ and $\yy$ denote $y^{-1}$.
\citet[Theorem 1]{bous} showed that
\begin{align}
\label{e-th1} S(x,y;t)&=\sum_{n\ge 0}\sum_{i,j\in \ZZ}
a_{i,j}(n)x^iy^jt^n  \nonumber\\
&=\frac{(1-2t(1+\xx)+\sqrt{1-4t})^{1/ 2}
(1+2t(1-\xx)+\sqrt{1+4t})^{1/ 2}} {2(1-t(x+\xx+y+\yy))},
\end{align}
where $S(x,y;t)$ is the complete generating function for walks on
the slit plane.

The authors also conjectured a closed form for $a_{-i,i}(2n)$ for
$i\ge 1$. By reflecting in the $x$-axis, we see that
$a_{-i,i}(2n)=a_{-i,-i}(2n)$, the closed form of which is given as
\eqref{conj1} in the following theorem.
\begin{thm}\label{t-conj}
For $i\ge 1 $ and $n\ge i$, we have
\begin{align}
a_{-i,-i}(2n) & =\frac{i}{2n}{2i\choose i}{n+i\choose
2i}\frac{{4n\choose 2n}}{
{2n+2i\choose 2i}}, \label{conj1} \\
a_{i,i}(2n) &=
a_{-i,-i}+4^n\frac{i}{n}\binom{2i}{i}\binom{2n}{n-i}
  \label{conj2}.
\end{align}
\end{thm}

We will prove this theorem in the next section. Theorem
\ref{t-lagrange1} is a basic tool to prove the conjecture.

There are two key steps in proving the conjecture that might be
worth mentioning: one is using Theorem \ref{t-lagrange1}
to obtain the generating function \eqref{conj3} that involves $a_{i,i}(2n)$ for all integers 
$i$; the other is guessing the formula \eqref{conj2}.

\vspace{3mm} Let $$C(t)=\sum_{n\ge 0} C_n t^n=
\frac{1-\sqrt{1-4t}}{2t}$$ be the Catalan generating function, and
let $$u=tC(t)C(-t)=\frac{\sqrt{1+4t}-1}{\sqrt{1-4t}+1}.$$

Much of the computation here involves rational functions of $u$.
We shall use the following facts from \citep{bous}.
$$\CC(u)=\CC(t,\sqrt{1-4t},\sqrt{1+4t}),$$
$$C(t)=\frac{1+u^2}{1-u}, \quad C(-t)=\frac{1+u^2}{1+u},
\quad C(4t^2)=\frac{(1+u^2)^2}{(1-u^2)^2}.$$

We shall prove Theorem \ref{t-conj} by computing the diagonal
generating function $F(y;t)$. More precisely, let
$$F(y;t)=\sum_{n\ge 0}\sum_{i\in \ZZ} a_{i,i}(2n) y^i t^{2n}.$$

Since $S(x,y\xx;t)$ belongs to $\CC[x,y,\xx,\yy][[t]]$, it is easy
to check that
\begin{align}\label{e-dd1}
F(y;t) =\CT_x S(x,y\xx;t)= \CT_x S(\xx, xy;t).
\end{align}

\begin{lem}\label{l-conj1}
\begin{multline}\label{conj3}
F(y;t)=\\
\frac{\left[\frac{1+\sqrt{1-4t}}{2}-t\left(1+\frac{1-\sqrt{1-4t^2(1+y)^2/y}}
{2t(1+y)}\right)\right]^{\frac{1}{2}}
\left[\frac{1+\sqrt{1+4t}}{2}+t\left(1-\frac{1-\sqrt{1-4t^2(1+y)^2/y}}{2t(1+y)}\right)
\right]^{\frac12} }{\sqrt{1-4t^2(1+y)^2/y}}.
\end{multline}
\end{lem}
\begin{proof}
Using \eqref{e-dd1} and \eqref{e-th1}, we get
\begin{align*}
F(y;&t)\\
=& \CT_x S(\xx, xy;t) \\
=&\CT_x \frac{(1-2t(1+x)+
\sqrt{1-4t})^{1/2}(1+2t(1-x)+\sqrt{1+4t})^{1/2}}
{2(1-t(x+\xx+x y+\xx \yy))}\\
=&\CT_x
\frac{x(1-2t(1+x)+\sqrt{1-4t})^{1/2}(1+2t(1-x)+\sqrt{1+4t})^{1/2}}{2(x-t(x^2+1+x^2y+\yy))}
.
\end{align*}

Applying Theorem \ref{t-lagrange1} with $R=\CC[y,\yy]$, this
becomes
$$ \frac{1}{2\left(1-t (2X+2X y)\right)}
(1-2t(1+X)+\sqrt{1-4t})^{1/2}(1+2t(1-X)+\sqrt{1+4t})^{1/2},
$$
where $X=X(t)$ is the unique solution in $tR[[t]]$ such that
$X=t(X^2+1+\yy+X^2y)$. We can solve for  $X$ by the quadratic
formula:
$$X= \frac{1-\sqrt{1-4t^2(1+y)^2/y}}{2t(1+y)}.$$

Equation \eqref{conj3}  then follows.
\end{proof}

It is clear that for any $G(y;t)\in R[y,\yy][[t]]$, there is a
unique decomposition $G(y;t)=G_+(y;t)+G_0(t)+G_-(\yy;t)$, such
that $G_+(y;t),G_-(y;t)\in yR[y][[t]]$ and $G_0(t)\in R[[t]]$.

Our task now is to find this decomposition of $F(y;t)$.
 There is no general theory to do this.
For this particular $F(y;t)$, thanks to the work of \bous\ and
Schaeffer, we can guess the formulas for $F_+$ and $F_-$ and prove
them.

The variable $s$ defined by the following is useful:
\begin{equation}
s=tC(4t^2)=\frac{u}{1-u^2} \text{ and } t=\frac{s}{1+4s^2}
\label{e-of-z}.
\end{equation}
Note that $s$ is also $S_{0,1}(t)$, the generating function of
walks on the slit plane that end at $(0,1)$. See \citep[P.
11]{bous}.

\begin{lem}\label{l-conj2}
We have the decomposition
$$F(y;t)=F_+(y,t)+1+F_-(\yy,t),$$
where
\begin{align}
F_+(y,t) & =F_-(y,t)+\frac{1}{2}((1-4s^2y)^{-1/2}-1),\label{e-f+} \\
F_-(y,t) &=\frac{(1-u^2)s^2yC(s^2y)}{1+u^2C^2(s^2y)s^2y}
\frac{1}{\sqrt{1-4s^2y}}\label{conj4} .
\end{align}
\end{lem}
\begin{proof}
Let
\begin{multline*}
    T(y;t)=\frac{(1-u^2)s^2yC(s^2y)}{1+u^2C^2(s^2y)s^2y}
\frac{1}{\sqrt{1-4s^2y}}+\frac{1}{2}((1-4s^2y)^{-1/2}-1)+1+\\
     \frac{(1-u^2)s^2\yy
C(s^2\yy)}{1+u^2C^2(s^2\yy)s^2\yy} \frac{1}{\sqrt{1-4s^2\yy}}.
\end{multline*}

From Lemma \ref{l-conj1}, the expression of $F(y;t)$ is:
\begin{align*}
\frac{\left[\frac{1+\sqrt{1-4t}}{2}-t\left(1+\frac{1-\sqrt{1-4t^2(1+y)^2/y}}
{2t(1+y)}\right)\right]^{1/2}
\left[\frac{1+\sqrt{1+4t}}{2}+t\left(1-\frac{1-\sqrt{1-4t^2(1+y)^2/y}}{2t(1+y)}\right)
\right]^{1/2} }{\sqrt{1-4t^2(1+y)^2/y}}.
\end{align*}

Therefore, it suffices to show that $T(y;t)=F(y;t)$. Since it is
easy to see that $T(y;0)=F(y;0)=1$, the proof will be completed by
showing that $T^2(y;t)-F^2(y;t)=0$.

Using the variable $u$, we can get rid of the radicals
$\sqrt{1-4t}$ and $\sqrt{1+4t}$ by the following:
$$\sqrt{1-4t}=\frac{1-2u-u^2}{1+u^2}, \text{ and } \sqrt{1+4t}=\frac{1+2u-u^2}{1+u^2}.$$

The radicals left are $D=\sqrt{1-4s^2y}$, $E=\sqrt{1-4s^2\yy}$,
and $\sqrt{1-4t^2(1+y)^2/y}$, which is easily checked to be equal
to $DE$.

Rewriting $T^2-F^2$ in terms of $u,D,E$, we get a rational
function of $u,D,E$. For $i=1,2$ (the degrees in $D$ and $E$ are
both $4$), replacing $D^{2i}$ by $(1-4s^2y)^i$, $D^{2i+1}$ by
$(1-4s^2y)^iD$, $E^{2i}$ by $(1-4s^2\yy)^i$, and $E^{2i+1}$ by
$(1-4s^2\yy)^iE$, we find that the expression reduces to $0$.
\end{proof}

Now we need to show the following.
\begin{lem}\label{l-conj3}
\begin{align}\label{conj7}
F_-(y,t)&=\sum_{n\ge 0} \sum_{i\ge 1} b_{i}(2n) t^ny^i ,
\end{align}
where
\begin{align}
b_i(2n)&=\frac{i}{2n}{2i\choose i}{n+i\choose 2i}\frac{{4n\choose
2n}}{ {2n+2i\choose 2i}}.
\end{align}
\end{lem}
We will give two proofs of this lemma. The first one starts from a
formula in \citep{bous}. We include it here as an example of
computing the generating function by Theorem \ref{t-lagrange1}.
The second proof is self-contained, and is simpler.

Let
\begin{align}
\label{e-dd2} f(y,t)=\sum_{n\ge 1} \sum_{i\ge 1} b_i(2n) t^ny^i .
\end{align}

We need to show that $F_-(y,t)=f(y,t)$.
\begin{proof}[First Proof of Lemma \ref{l-conj3}.]
It was stated in \citep{bous} that
\begin{equation}\label{conj8}
\sum_{n\ge 0} b_i(2n)
t^n=\frac{(-1)^i}{(1-u^2)^{2i-1}}\sum_{k=i}^{2i-1} {2i-1\choose k}
(-1)^k u^{2k}.
\end{equation}

Let $s$ be as in \eqref{e-of-z}. Using the following fact
$${n\choose k}= \ct_\alpha \frac{1}{\alpha^k} (1+\alpha)^n,$$
we can compute $f(y,t)$ by Theorem \ref{t-lagrange1}:
\begin{align*}
f(y,t) &=\sum_{i\ge 1}
\frac{(-1)^i}{(1-u^2)^{2i-1}}\sum_{k=i}^{2i-1}
{2i-1\choose k} (-1)^k u^{2k} y^i \\
&=\sum_{i\ge 1} \frac{(1-u^2)(-1)^i}{(1-u^2)^{2i}}
\sum_{r=0}^{i-1}
{2i-1\choose i+r} (-1)^{i+r} u^{2i+2r} y^i , \mbox{ where } r=k-i\\
&= (1-u^2) \sum_{r\ge 0} (-1)^ru^{2r}\sum_{i\ge r+1}
{2i-1\choose i-1-r} \frac{u^{2i}} {(1-u^2)^{2i}}y^i\\
&= (1-u^2)\sum_{r\ge 0}(-u^2)^r\sum_{i\ge r+1} \CT_\alpha
(1+\alpha )^{2i-1}\left(\frac{1}{\alpha}\right)^{i-1-r} (s^2y)^i \\
&= \CT_\alpha (1-u^2)\sum_{r\ge 0} \frac{\alpha }{1+\alpha }
(-u^2)^r\alpha ^r \sum_{i\ge r+1} \frac{(1+\alpha)^{2i}}{\alpha^i} (s^2y)^i \\
&= \CT_\alpha \frac{\alpha}{1+\alpha}(1-u^2)\sum_{r\ge
0}(-u^2\alpha)^r \left(\frac{(1+\alpha)^2}{\alpha}
s^2y\right)^{r+1}
\frac{1}{1-\displaystyle\frac{(1+\alpha)^2}{\alpha}s^2y} \\
&= \CT_\alpha (1-u^2)(1+\alpha) s^2y
\frac{1}{1+u^2(1+\alpha)^2s^2y} \cdot
\frac{1}{1-\displaystyle\frac{(1+\alpha)^2}{\alpha}s^2y} .
\end{align*}

Now
$$(1-u^2)(1+\alpha) s^2y \displaystyle\frac{1}{1+u^2(1+\alpha)^2s^2y}
$$ is a power series in $t$ with coefficients in $\CC[y][\alpha]$,
and
\begin{align*}
\frac{1}{1-\frac{(1+\alpha)^2}{\alpha}s^2y} &=
\frac{\alpha}{\alpha-(1+\alpha)^2s^2y}.
\end{align*}

Solving the denominator for $\alpha$, we get two solutions:
$$\frac{1-2s^2y+\sqrt{1-4s^2y}}{2s^2y}\text{ and } \frac{1-2s^2y-\sqrt{1-4s^2y}}{2s^2y}.$$
Only the latter is a power series in $t$ with constant term $0$,
which can also be written as $A=C(s^2y)-1$.

Thus we can apply Theorem \ref{t-lagrange1} to get
\begin{align*}
f(y,t)=& \CT_\alpha \frac{\alpha}{\alpha-(1+\alpha)^2s^2y}
(1-u^2)(1+\alpha) s^2y \frac{1}{1+u^2(1+\alpha)^2s^2y}
 \\
=& (1-u)^2s^2y(1+A)\frac{1}{1+u^2(1+A)^2 s^2y}\frac{1}{1-2s^2y(A+1)} \\
=& (1-u^2)s^2yC(s^2y)\frac{1}{1+u^2C^2(s^2y)s^2y}
\frac{1}{\sqrt{1-4s^2y}},
\end{align*}
which completes the proof.
\end{proof}

The second proof derives a different form of $F_-(y;t)$.

\begin{proof}[Second Proof of Lemma \ref{l-conj3}.]
We begin with finding the generating function of $2nb_i(n)$, which
equals
 $t{\partial\over \partial t} f(y,t)$.

We claim that
\begin{align}\label{e-proof2}
\sum_{n\ge 0} {n+i\choose 2i}\frac{{4n\choose 2n}}{ {2n+2i\choose
2i}} t^{2n} =\frac{\sqrt{1+4s^2} s^{2i}}{1-4s^2},
\end{align}
where the relation between $t$ and $s$ is given in \eqref{e-of-z}.

It is easy to check that
$${n+i\choose 2i}\frac{{4n\choose 2n}}{
{2n+2i\choose 2i}} = \binom{2n-1/2}{n-i}4^{n-i}.$$

In the well-known formula
$$ \frac{C(x)^k}{\sqrt{1-4x}} =\sum_{n\ge 0} \binom{2n+k}{n} x^n,$$
by setting $x=4t^2$, and $k=2i-1/2$, we  get
$$\sum_{n\ge 0} {n+i\choose 2i}\frac{{4n\choose 2n}}{
{2n+2i\choose 2i}}
t^{2n}=t^{2i}\frac{C(4t^2)^{2i-1/2}}{\sqrt{1-16t^2}}.$$ Using
\eqref{e-of-z} to write the above in terms of $s$, we get
\eqref{e-proof2}.

Now we have
$$t{\partial \over \partial t} f(y;t) =
\sum_{i\ge 1}\sum_{n\ge 0} i \binom{2i}{i} {n+i\choose
2i}\frac{{4n\choose 2n}}{ {2n+2i\choose 2i}} t^{2n} y^i =
\frac{2s^2y}{(1-4s^2y)^{{3/ 2}}} \frac{\sqrt{1+4s^2}}{1-4s^2}.$$

Hence
\begin{align*}
f(y;t)=&\int \frac{2s^2y}{(1-4s^2y)^{{3/ 2}}}
\frac{\sqrt{1+4s^2}}{1-4s^2} \frac{dt}{t} \\
=& \int \frac{2s^2y}{(1-4s^2y)^{{3/ 2}}}
\frac{\sqrt{1+4s^2}}{1-4s^2} \frac{1-4s^2}{s(1+4s^2)} ds \\
=& \frac{y\sqrt{1+4s^2}}{2(1+y)\sqrt{1-4s^2y}} +\text{constant},
\end{align*}
where the constant is independent of $t$. By setting $t=0$, and
hence $s=0$, we get $f(y;0)=\frac{y}{2(1+y)}+\text{constant}$.

Recalling equation \eqref{e-dd2}, we see that $f(y;0)=0$. Thus the
constant equals $-\frac{y}{2(1+y)}$. This gives another form of
$f(y;t)$:
$$f(y;t)=
\frac{y\sqrt{1+4s^2}}{2(1+y)\sqrt{1-4s^2y}}-\frac{y}{2(1+y)}
=\frac{y(1+u^2)}{2(1+y)(1-u^2)\sqrt{1-4s^2y}}-\frac{y}{2(1+y)},$$
which is easily checked to be equal to $F_-(y;t)$ as given in
\eqref{conj4}.
\end{proof}

\begin{proof}[Proof of Theorem \ref{t-conj}]
We gave a formula for the generating function
$$F(y;t)=\sum_{n\ge 0} \sum_{i\in \ZZ} a_{i,i}(2n) y^it^{2n} $$
in Lemma \ref{l-conj1}. In Lemma \ref{l-conj2}, we showed that
$$F_-(y,t)=\sum_{n\ge 0} \sum_{i> 0} a_{-i,-i}(2n) y^it^{2n}$$
has a formula as given in \eqref{conj4}. The proof of
\eqref{conj1} is thus accomplished by Lemma \ref{l-conj3}.

For equation \eqref{conj2}, once we get the formula \eqref{e-f+},
it is an easy exercise to show that
$$\frac12 ((1-4s^2y)^{-1/2}-1)=\sum_{n\ge 1}\sum_{i\ge 1} 
4^n\frac{i}{n}\binom{2i}{i}\binom{2n}{n-i} y^i t^{2n}.$$

\end{proof}

\section{Walks on the Quarter Plane\label{s-quart}}


{\em Walks on the quarter plane} are walks that stays in the first
quadrant $x> 0, y> 0$. Note that in some literature, the quarter
plane refers to $x\ge 0, y\ge 0$. Walks on the quarter plane has
be studied by many authors. See, e.g., \citep{fayolle}.
\citet{bousquart} used a functional equation  approach to solve
the enumeration problems for walks in the quarter plane. See also
\citep{bousquartD,bousquartf}. Here we will use the same ideas,
but work by our theory.

Let $\sy$ be a finite subset of $\ZZ^2$, and let $R$ be the
restriction that starts at $(1,1)$ and stays in the first quadrant
$x>0,y>0$. Then the walks in the quarter plane problem is to study
the properties of those paths in $\sy^R$. One basic problem is to
give a formula for the generating function $\Gamma(\sy^R)$.

Denote by $Q(x,y;t)$ the generating function $\Gamma(\sy^R)$. Then
it can be written as
$$Q(x,y;t)=\sum_{n\ge 0}t^n\sum_{i,j>0}a_{i,j}(n) x^iy^j,$$
where $a_{i,j}(n)$ is the number of walks in the first quadrant
that start at $(1,1)$, end at $(i,j)$.

Using the bridge lemma, we can get a functional equation. But we
do not have a general theory to solve this kind of functional
equation. Up to now, we can only deal with some simple situation.

To make things simpler, we suppose that $\sy$ contains only
$(r,s)$ with $-1\le r,s \le 1$. (Even in this situation, some
problems are left unsolved.) Such $\sy$ is said to be having small
lengths. In this case, $\Gamma(\sy)$  can be written as:
$$ \Gamma(\sy) =t(A(x) \yy +B(x)+C(x) y),$$
where $xA(x)$, $xB(x)$, and $xC(x)$ are polynomials in $x$ of
degree at most $2$.

If $\sy$  has small lengths, then the bridge paths must end at
$(k,0)$ or $(0,k)$ for some $k\ge 0$. Now let $H(x,t)$ be the
generating function for bridge paths that end  at $(k,0)$ for some
$k>0$, let $V(y,t)$ be the generating function for bridge paths
that end at $(0,k)$ for some $k>0$, and let $O(t)$ be the
generating function for bridge paths that end at $(0,0)$. Then the
bridge lemma gives us the following functional equation:

\begin{align} \label{e-qfunequ}
Q(x,y;t)=\frac{xy-H(x,t)-V(y,t)-O(t)}{1-\Gamma(\sy)},
\end{align}
where the $xy$ in the numerator is the weight of the starting
point $(1,1)$.

The boundary conditions are $\ct_y Q(x,y;t)=0$ and $\ct_x
Q(x,y;t)=0$.

Now it is routine to apply equation \eqref{e-qfunequ} to these two
boundary conditions. From the first boundary condition, we get
\begin{align*}
\ct_y Q(x,y) &= \ct_y \frac{y}{y-t(A_{-1}(x)+A_0(x)y+A_{1}(x)y)} (xy-H(x)-V(y)-O(t))\\
&= \frac{1}{1-2Y} (xY-H(x)-V(Y)-O(t)),
\end{align*}
where $Y=Y(x)$ is the unique positive root for $y$, in the
denominator $y-t(A_{-1}(x)+A_0(x)y+A_{1}(x)y)$. This denominator
can also be written as $y\Gamma(\sy)$, and $Y$ can be found by the
quadratic formula.

Hence we get our first functional equation:
\begin{align}\label{e-qfuny}
xY-H(x)-V(Y)-O(t)=0.
\end{align}

Similarly, from the second boundary condition,
 we derive our second functional equation:
\begin{align}\label{e-qfunx}
Xy-H(X)-V(y)-O(t)=0,
\end{align}
where $X$ is the unique positive root for $x$, in
$x-x\Gamma(\sy)$.

From the composition law, both $X\circ Y$ and $Y\circ X$ are well
defined. One can check that $X\circ Y=x$  and $Y\circ X =y$. A
simple reason for this to be true is that both $X$ and $Y$ are
solved from $xy\Gamma(\sy)$.

Using the above fact, we deduce that equations \eqref{e-qfunx} and
\eqref{e-qfuny} are equivalent, because \eqref{e-qfunx} can be
obtained from \eqref{e-qfuny} by replacing $x$ with $X$, and
\eqref{e-qfuny} can be obtained from \eqref{e-qfunx} by replacing
$y$ with $Y$.

Now the problem is how to solve the functional equation
\eqref{e-qfunx} for $H,V,O$. This can be done in some simple
cases.

Case $1$: If additionally $\sy$ is symmetric in $y$, i.e.
$$\Gamma(\sy)=\left.\Gamma(\sy)  \right|_{y=\yy},$$
then we know how to solve \eqref{e-qfunx}. The case that $\sy$ is
symmetric in $x$ is similar.

It is clear that in this case $X(y)=X(\yy)$. Substituting $y$ by
$\yy$ in \eqref{e-qfunx}, we get
\begin{align}\label{e-qfunx-1}
X \yy-H(X)-V(\yy)-O(t)=0.
\end{align}

Taking the difference on both sides of equations \eqref{e-qfunx}
and \eqref{e-qfunx-1}, we get
\begin{align}\label{e-qfunx-f}
V(y)-V(\yy)= Xy-X\yy.
\end{align}
Since $V(y)\in yt\CC[y][[t]]$, $V(y)+0+(-V(\yy))$ is the first
decomposition of $Xy-X\yy$. Hence $V(y)$ equals the positive part
of $Xy-X\yy$.

Similarly we can solve for $H(x)$. Then $O(t)$ can be obtained
from equation \eqref{e-qfunx}.

\begin{exa}
If $\sy=\{\, (1,0),(-1,0),(0,1),(0,-1)\,\}$, then this is called
ordinary lattice paths.
\end{exa}
The corresponding generating function $Q(x,y)$ has the form
$$Q(x,y)=\frac{xy-H(x)-V(y)-O(t)}{1-t(x+y+\xx+\yy)}.$$
By symmetry in $x$ and $y$, $H=V$, and it is easy to see that
$O(t)=0$.

Solve for $y$ in $1-t(x+y+\xx+\yy)=0$, we get
$$Y= \frac{1-t\left (x+{x}^{-1}\right )-\sqrt{\left(1-t\left (x+{x}^{-1}\right )
\right)^2-4t^2}}{2t}.$$

So by extracting the positive part of $xY-\xx Y$, we will get
$V(y)$.

The above argument is in fact the algebraic version of the
well-known reflection principle.

Using the reflection principle, we can solve the case that $\sy$
is symmetric in $y$ and has small lengths in $x$. In other words,
this is to say that $(s,t)\in \sy$ implies that $|s|\le 1$ and
that $(-s,t)\in \sy$. Of course we require that $\sy$ be a
well-ordered subset of $\ZZ^2$.

Let $a(i,j,n)$ be the number of paths of length $n$ that start at
$(1,1)$, end at $(i,j)$, and stay inside the quarter plane. Let
$p(i,j,n)$ be the number of paths of length $n$ that start at
$(1,1)$, end at $(i,j)$, and stay above the line $y=0$. Then among
all paths of length $n$ from $(1,1)$ to $(i,j)$ that stays above
the line $y=0$, those paths that never touches the line $x=0$ are
counted by $a(i,j,n)$, and those paths that touch the line $x=0$
are counted by $P(i+2,j,n)$, since they are the same as the number
of paths of length $n$ that start at $(-1,1)$, end at $(i,j)$, and
stay above the line $y=0$ by the reflection principle.

Let $P(x,y;t)$ be the generating function of paths that start at
$(0,0)$ and stay above the line $y=-1$. Then the above argument
gives us the equation
$$ V(x;t)=\pt_x  x yP(x,y;t)-x^{-1} yP(x,y;t).$$

\newpage
$ $ \thispagestyle{empty}
\newpage


\clearpage
\newpage
\clearpage
\bibliographystyle{ametsoc}


\end{document}